UNIVERSITÄT
LEIPZIG

# Normal Form of Equivariant Maps and Singular Symplectic Reduction in Infinite Dimensions with Applications to Gauge Field Theory

Von der Fakultät für Physik und Geowissenschaften
der Universität Leipzig
genehmigte

## Dissertation

zur Erlangung des akademischen Grades
Doktor der Naturwissenschaften
(Dr. rer. nat.)

vorgelegt von

M. Sc. Tobias Diez
geboren am 15.08.1989 in Leipzig

1. Gutachter: Prof. Dr. Gerd Rudolph
2. Gutachter: Prof. Dr. Markus Pflaum

# Bibliographic details

---

Tobias Diez

Normal Form of Equivariant Maps and Singular Symplectic Reduction in Infinite Dimensions with Applications to Gauge Field Theory

---

Inspired by problems in gauge field theory, this thesis is concerned with various aspects of infinite-dimensional differential geometry.

In the first part, a local normal form theorem for smooth equivariant maps between tame Fréchet manifolds is established. Moreover, an elliptic version of this theorem is obtained. The proof these normal form results is inspired by the Lyapunov–Schmidt reduction for dynamical systems and by the Kuranishi method for moduli spaces, and uses a slice theorem for Fréchet manifolds as the main technical tool. As a consequence of this equivariant normal form theorem, the abstract moduli space obtained by factorizing a level set of the equivariant map with respect to the group action carries the structure of a Kuranishi space, i.e., such moduli spaces are locally modeled on the quotient by a compact group of the zero set of a smooth map.

In the second part of the thesis, the theory of singular symplectic reduction is developed in the infinite-dimensional Fréchet setting. By refining the above construction, a normal form for momentum maps similar to the classical Marle–Guillemin–Sternberg normal form is established. Analogous to the reasoning in finite dimensions, this normal form result is then used to show that the reduced phase space decomposes into smooth manifolds each carrying a natural symplectic structure.

Finally, the singular symplectic reduction scheme is further investigated in the situation where the original phase space is an infinite-dimensional cotangent bundle. The fibered structure of the cotangent bundle yields a refinement of the usual orbit-momentum type strata into so-called seams. Using a suitable normal form theorem, it is shown that these seams are manifolds. Taking the harmonic oscillator as an example, the influence of the singular seams on dynamics is illustrated.

The general results stated above are applied to various gauge theory models. The moduli spaces of anti-self-dual connections in four dimensions and of Yang–Mills connections in two dimensions is studied. Moreover, the stratified structure of the reduced phase space of the Yang–Mills–Higgs theory is investigated in a Hamiltonian formulation after a $(3+1)$-splitting.

## Bibliographische Beschreibung

---

Tobias Diez

Normal Form of Equivariant Maps and Singular Symplectic Reduction in Infinite Dimensions with Applications to Gauge Field Theory

---

Von Problemen der Eichfeldtheorie inspiriert, beschäftigt sich diese Arbeit mit verschiedenen Aspekten der unendlich-dimensionalen Differentialgeometrie.

Im ersten Teil wird ein Satz über die Normalform von glatten äquivarianten Abbildungen zwischen zahmen Fréchet-Mannigfaltigkeiten aufgestellt. Ergänzend wird eine elliptische Version dieses Theorems bewiesen. Der Beweis dieser Normalformsätze greift die Lyapunov–Schmidt-Reduktion für dynamische Systeme und die Kuranishi-Methode für Modulräume unter Verwendung eines Scheibensatzes für Fréchet-Mannigfaltigkeiten als technisches Hauptwerkzeug auf. Als Konsequenz dieses äquivarianten Normalformsatzes trägt der abstrakte Modulraum, der als Quotient einer Niveaufläche der äquivarianten Abbildung bezüglich der Gruppenwirkung entsteht, die Struktur eines Kuranishi-Raums.

Im zweiten Teil der Arbeit wird die Theorie der singulären symplektischen Reduktion im unendlich-dimensionalen Fréchet-Kontext entwickelt. Durch Verfeinern der obigen Konstruktion erhält man eine Normalform für Impulsabbildungen ähnlich der klassischen Marle–Guillemin–Sternberg-Normalform. Dieses Normalformergebnis wird verwendet um zu zeigen, dass sich der reduzierte Phasenraum in glatte symplektische Mannigfaltigkeiten zerlegen lässt. Anschließend wird das singuläre symplektische Reduktionsschema für den Fall, dass der ursprüngliche Phasenraum ein unendlich-dimensionales Kotangentialbündel ist, untersucht. Die kanonische Faserstruktur des Kotangentialbündels ergibt eine Verfeinerung der üblichen Orbit-Impuls-Strata in sogenannte Seams. Mit Hilfe eines geeigneten Normalformsatzes wird nachgewiesen, dass diese Seams Mannigfaltigkeiten sind. Den Einfluss der einzelnen Seams auf die Dynamik veranschaulicht das Beispiel des harmonischen Oszillators.

Diese allgemeinen Ergebnisse werden auf verschiedene Modelle der Eichtheorie angewendet. So werden die Modulräume von anti-selbstdualen Zusammenhängen in vier Dimensionen und von Yang–Mills-Zusammenhängen in zwei Dimensionen analysiert. Weiterführend wird die stratifizierte Struktur des reduzierten Phasenraums der Yang–Mills–Higgs-Theorie in der Hamiltonschen (3+1)-Formulierung untersucht.

# Contents

# Acknowledgments

In his acceptance speech of the Nobel Prize, Ernest Hemingway affirmed: "Writing, at its best, is a lonely life." And yet, I have been fortunate that this dissertation has been a collaborative and cooperative process rather than a lonely endeavor.

I am deeply grateful to my mentor, Gerd Rudolph, for finding the perfect middle course between letting me wander off to explore my own ideas and giving me intellectual guidance when needed. His enlightening comments and suggestions have been immeasurable. Every page of this thesis (at least the pedagogically good ones) bear the deep imprint of his invaluable teaching experience. I would like to thank him for his generosity with his time and wealth of knowledge. I want to extend this feeling of gratitude to Johannes Huebschmann and Tudor Ratiu for their constant support, encouragement and hospitality in Lille, Lausanne and Shanghai. Working with them has been a pleasure and I learned a lot during our stimulating discussions.

Moreover, I am very thankful for the support I have received from so many sides during my time as a PhD student. In particular, I am very much indebted to Matthias Schmidt for many thoughtful comments and suggestions at various stages of the project and for reading parts of the manuscript. I would also like to thank Bas Janssens, Cornelia Vizman and Karl-Hermann Neeb for our inspiring and helpful conversations. To Christian Blohmann, Marius Crainic and Peter Albers, I am grateful for their interest in my work and their hospitality. Many sincere thanks to the members of my special task force — Alex and Damaris — for helping me along in this final phase by reading the manuscript. I gratefully acknowledge the support I received from the Max Planck Institute for Mathematics in the Sciences (Leipzig) and from the German National Academic Foundation (DAAD).

My colleagues and friends at the Max Planck Institute have contributed to my growth as a scientist and as a person in many ways. In particular, I am very grateful to Claudia, Nicólò, Thilo, Tim and Vishal for creating such a stimulating environment. I would like to thank my close friends who have borne with me when I was stressed and shared my excitement when things were going well. Among them, I would particularly like to express my special thanks to Alex, Anne, Bele, Damaris, Erik, Fine, Iris, Lisa and Lovis for the infinitely many good memories.

Lastly, none of this would have been possible without the constant encouragement and unconditional love of my family. I am deeply indebted to my parents for their unwavering support for all of my endeavors. Love you both. Furthermore, I am immensely grateful to my sister for absorbing my monothematic babbling and for being there even when I was not aware of the fact that I needed help. Beatrice, you are the best sister a brother could ask for.

# 1 Introduction

Nonlinear partial differential equations play an essential role in theoretical physics as well as in pure mathematics. In geometry, such equations share the common theme of having solutions that are sensitive to the topology and geometry of the underlying manifold. Already in the eighties, it was recognized that the same applies in reverse, that is, one can use differential-geometric equations to deduce some astounding topological and geometrical properties of manifolds. These insights have been accompanied by a shift in perspective: instead of studying the solutions themselves, the geometric structure of the space of solutions has increasingly become the subject of interest and has moved into the focus of mathematicians and physicists alike. In particular, topological and inherently non-perturbative aspects of the physical system are often manifested in the geometry of the solution space.

There is another general feature: the equations arising in geometry and physics frequently have large symmetry groups. For example, the Yang–Mills equation is invariant under the group of gauge transformations and Einstein's equation is invariant under the group of diffeomorphisms. In this situation, the right object of interest is the moduli space of solutions obtained by taking the quotient of the space of solutions by the action of the symmetry group. Thus, for applications in both geometry and physics, a deep understanding of the local structure of this quotient is essential. Of particular importance is the formation of singularities due to objects having a non-trivial stabilizer under the symmetry group action. Although relying on similar techniques, the analysis of these fundamental features so far only happened on a case by case basis for each moduli space separately.

*In the first part of this thesis, we provide a general framework that gives a unified approach to these differential-geometric moduli spaces.* Specifically, we establish a convenient normal form for a large class of non-linear differential equations with symmetries. Furthermore, we show that the corresponding moduli space of solutions can be endowed with the structure of a Kuranishi space, which roughly speaking means that it can be locally identified with the quotient of the zero set of a smooth map by the linear action of a compact group. Our approach is inspired by the ideas that underlie the Lyapunov–Schmidt reduction for dynamical systems and the Kuranishi method for moduli spaces in differential geometry. The results are phrased in terms of equivariant maps between infinite-dimensional manifolds endowed with actions of infinite-dimensional Lie groups.

As another advantage, the developed methods make it possible to investigate additional features of the moduli space under consideration. For example, it is known that some differential equations, like the constraint equations in Yang–Mills theory or general relativity, are tightly connected to the symmetry of the model and can often be recast in terms of the corresponding Noether current. The appropriate setting to study these phenomena is provided by the Hamiltonian framework formulated in the language of infinite-dimensional symplectic manifolds.

As a blueprint for developing this approach, we use the well-established finite-dimensional theory for classical particle mechanics. In this setting, the phase space is described by a finite-dimensional symplectic manifold and the symmetry of the system is encoded in a symplectic action of the symmetry group. The corresponding conserved quantities are captured in the momentum map. Following the classical result of Marsden, Weinstein and Meyer, the momentum map can be used to eliminate a number of variables and thereby pass to the reduced phase space, in which the symmetries are divided out. This elimination process is called symplectic reduction, see [MW74; SL91]. In finite dimensions, the equivariant Darboux theorem yields a convenient normal form for the original Hamiltonian system and thereby provides a fundamental tool to analyze the local structure of the reduced phase space. This analysis shows that the reduced phase space is a symplectic manifold again if the action of the symmetry group is free. However, in the situation where the action is not free, the reduced space has conic singularities and it is thus a stratified space with each stratum being a symplectic manifold. Starting in the early nineties [ER90], various case studies in finite dimensions have shown that these singularities have an influence on the properties of the quantum theory, see [RS17, Chapter 8] for a detailed discussion. In particular, they may carry information about the spectrum of the Hamiltonian, see [HRS09]. These studies show that information about the singular structure of the reduced phase space are an indispensable prerequisite for a deeper understanding of the corresponding quantum theory.

*In the second part of this thesis, we develop the theory of singular symplectic reduction in an infinite-dimensional setting.* This allows us to study the aforementioned moduli spaces and their additional features, which are derived from the symplectic structure. When passing to the infinite-dimensional setting, a simple counterexample by Marsden [Mar72] shows that the Darboux theorem fails spectacularly for weakly symplectic Banach manifolds. Thus, in order to establish an infinite-dimensional version of the symplectic reduction theorem, we pursue a different strategy, which is in sharp contrast to the usual proof in finite dimensions. Our derivation of this reduction theorem is based on the observation that it is not essential to bring the symplectic structure into a normal form. Instead, our focus will be primarily on momentum maps; using the symplectic structure only as a secondary tool. By refining and adapting the techniques developed in the first part of the thesis, we construct a normal form for equivariant momentum maps in the spirit of the classical

Marle–Guillemin–Sternberg normal form. With the help of this normal form, we then prove a singular symplectic reduction theorem including the analysis of the stratification into symplectic manifolds. Finally, the reduction scheme is applied to the symmetry reduction of various gauge theory models.

The symmetry reduction for gauge field theory is part of an ambitious program to rigorously develop quantum field theories. While perturbative methods in quantum field theory yield a satisfactory description of high energy processes in particle physics, many low energy processes are dominated by non-perturbative effects and, so far, have eluded a rigorous theoretical explanation. Many attempts to develop a mathematically rigorous, non-perturbative quantum gauge theory are inspired by the general approach of constructive quantum field theory. In its classical version, this approach is based on the Euclidean path integral concept. We refer to [tH005; JW] for status reports on the case of Yang–Mills theory. A different option is provided by the Hamiltonian formalism, see [BDI74; KS75; KR02; KR05; GR17]. The fundamental problem when dealing with Yang–Mills theory is the elimination of the unphysical gauge degrees of freedom. Within the Hamiltonian framework, this is accomplished via symplectic reduction. In both the path integral and the Hamiltonian approach, one usually analyzes a corresponding lattice approximation of the model as an intermediate step. The issues of the path integral approach now reappear as problems related to taking the continuum limit of vanishing lattice spacing. Concerning the Hamiltonian point of view, one may think of another strategy: instead of using a lattice approximation, one can study the classical continuum gauge theory as an infinite-dimensional Hamiltonian system with a gauge symmetry, reduce the symmetry and, then, pass to the quantum world by extending methods from geometric quantization to infinite-dimensional systems. The results of this thesis contribute to this approach by developing a rigorous reduction scheme and thereby clarify the structure of the reduced phase space on the classical level. The final and admittedly hardest step for constructing the corresponding quantum field theory is to quantize the reduced system, which is clearly outside the scope of this thesis.

In field theory and global analysis, the maps under consideration are usually given by partial differential operators between spaces of sections and thus they give rise to smooth maps between appropriate Sobolev completions. On the other hand, the symmetry action often involves compositions of maps and thus fails to be differentiable as a map between spaces of sections of a given Sobolev class. For example, the group of diffeomorphisms of a fixed Sobolev regularity is a Banach manifold as well as a topological group but not a Lie group, because the group operation is not differentiable. When working with smooth sections these problems disappear and the group of smooth diffeomorphisms is a bona fide Lie group modeled on a Fréchet space. In order to include these important examples, we assume throughout the thesis that all infinite-dimensional manifolds are modeled on Fréchet or even more general locally convex spaces. The approach via Fréchet spaces has also the advantage that certain geometric arguments are simpler, because one does not have to deal with issues originating in the low regularity of the geometric objects under study. We note that the polyfold framework recently proposed by Hofer, Wysocki, and Zehnder [HWZ17b] offers another approach to these issues of differentiability. We comment further on the relation of the polyfold theory to our Fréchet approach in Remarks 3.2.2 (iii) and A.1.4 (i) below.

The main results of this dissertation are as follows.

Chapter 2. Normal Forms of Maps This chapter lays the foundation for the study of equivariant normal forms of maps. We begin by considering the linear setting and determine under which conditions a continuous linear map between locally convex spaces factorizes through a topological isomorphism. We say that an operator is regular if it admits such a factorization. As we will demonstrate, an operator is regular if and only if it possesses a generalized inverse. For Fredholm operators, the existence of a generalized inverse is tightly connected to the existence of a parametrix. This implies, in particular, that Fredholm operators are regular. As a preparation for the non-linear case, we extend the discussion of regularity to families of linear maps depending continuously on a parameter and to chain complexes. Most results concerning regular operators are well-established in the Banach setting but their extension to more general locally convex spaces represents original work.

Next, we discuss the local behavior of a smooth map between locally convex manifolds. Unifying the concepts of immersion, submersion and subimmersion in one framework, the notion of a normal form of a non-linear map is introduced. Using versions of the Inverse Function Theorem, we establish Theorems 2.2.6, 2.2.9, 2.2.10, 2.2.13 and 2.2.14 which show that a given map can be brought into such a normal form in various functional-analytic settings under suitable conditions. These normal form theorems provide a unified approach to the immersion theorem, the level set theorem and the constant rank theorem in the setting of locally convex manifolds and tame Fréchet manifolds.

Chapter 3. Moduli Spaces Based on this preparation, we introduce the concept of an equivariant normal form and provide suitable conditions which ensure that an equivariant map can be brought into such a normal form (see Theorem 3.1.6 and its variants). Besides the normal form results of Chapter 2, the main technical tool is a version of the slice theorem for Fréchet manifolds as proved in [DR18c]. Then, we investigate the local structure of the moduli space obtained by taking the quotient of a level set of the equivariant map by the group action. Under the assumption that the map can be brought into a normal form, we show that the corresponding moduli space can be endowed with the structure of a Kuranishi space, which roughly speaking means that it can be locally identified with the quotient of the zero set of a smooth map with respect to the linear action of a compact group. Moreover, we find additional conditions on the normal form which ensure that the moduli space is stratified by orbit types. Finally, to show the utility of this novel framework, we apply the general theory to the example of the moduli space of anti-self-dual Yang–Mills connections.

Chapter 4. Singular Symplectic Reduction In this chapter, we generalize the well-known results concerning singular symplectic reduction to our infinite-dimensional setting. We begin by discussing some algebraic results of linear symplectic geometry needed for the symmetry reduction scheme. As a fundamental tool to handle weakly symplectic forms, we introduce and study a class of topologies associated to the symplectic form.

Next, we discuss symplectic manifolds and momentum maps in our infinite-dimensional setting. Based on our joint work with T. Ratiu [DR] on actions of diffeomorphism groups, we introduce the notion of a group-valued momentum map which unifies several other notions of generalized momentum maps. The group structure of the target allows to encode discrete topological information; a feature that is especially relevant to the action of geometric automorphism groups, which are sensitive to the topology of the spaces they live on. Using the results of Chapter 3 about normal forms of equivariant maps as a foundation, we establish a refined normal form result for momentum maps in the spirit of the classical Marle–Guillemin–Sternberg normal form.

This important technical tool then serves as the basis for our infinite-dimensional version of the Singular Symplectic Reduction Theorem (see Theorem 4.3.5). As in finite dimensions, it states that the reduced phase space, obtained by factorizing a level set of the momentum map with respect to the symmetry group, decomposes into smooth manifolds each carrying a natural symplectic structure. If a suitable approximation property holds, then this decomposition of the reduced phase space is a stratification. Moreover, the dynamics of the system reduces to the quotient and restricts to a Hamiltonian dynamics on each symplectic stratum. To our knowledge, these results concerning the normal form of a momentum map as well as the structure of singular symplectic quotients are new even for (weakly) symplectic Banach manifolds.

Finally, we apply these general results to the example of symplectic reduction in the context of the Yang–Mills equation over a Riemannian surface. This is based on joint work with J. Huebschmann [DH18].

Chapter 5. Singular Cotangent Bundle Reduction In most applications in physics, the phase space is a cotangent bundle over the configuration space of the system. In finite dimensions, Perlmutter, Rodriguez-Olmos, and Sousa-Dias [PRS07] have shown that the fibered structure of the cotangent bundle yields a refinement of the usual orbit-momentum type strata into so-called seams. The principal seam is symplectomorphic to a cotangent bundle while the singular seams are coisotropic submanifolds of the corresponding symplectic stratum.

In the first part of this chapter, we extend these results concerning singular symplectic reduction of cotangent bundles, including the analysis of the secondary stratification into seams, to the case of infinite-dimensional Fréchet manifolds (see Theorem 5.4.8). The proof relies on the construction of a suitable normal form for momentum maps of lifted actions, for which we directly exploit the additional structure of the cotangent bundle. As a by-product, our theory provides a much simpler approach to singular cotangent bundle reduction in the finite-dimensional setting. The possible significance of the seams for the dynamics of the system is demonstrated by analyzing the finite-dimensional example of the harmonic oscillator.

In the second part of the chapter, we apply our general theory to the singular cotangent bundle reduction of Yang–Mills–Higgs theory. The Singular Symplectic Reduction Theorem implies that the reduced phase space of the theory is a stratified symplectic space. We study this stratification in more detail and find that an inclusion of the singular strata leads to a refinement of what is called the resolution of the Gauß constraint in the physics literature. Finally, we further analyze the secondary stratification in the concrete example of the Higgs sector of the Glashow–Weinberg–Salam model. In this context, we find that the configuration space has only two orbit types. The singular stratum is characterized by the remarkable physical property of the absence of $W$-bosons, i.e., the Z-boson is the only non-trivial intermediate vector boson on the singular stratum. We then discuss the stratification of the phase space. The secondary stratification turns out to be similar to that of the harmonic oscillator in the sense that there are only three strata: two cotangent bundles which are glued together by one seam. The non-generic cotangent bundle is the phase space of a sub-theory consisting of electrodynamics described by a photon, the theory of a massive vector boson described by the $Z$-boson and the theory of a self-interacting real scalar field described by the Higgs boson. The seam is characterized by the condition that the $W$-boson field vanishes but its conjugate momentum is non-trivial. In contrast, all intermediate vector bosons of the model are present on the generic cotangent bundle. Finally, we study the structure of the strata of the reduced phase space in terms of gauge invariant quantities for the theory on $S^3$. By implementing unitary and Coulomb gauge fixing in a geometric fashion using momentum maps, we show that the singular structure of the reduced phase space is encoded in a finite-dimensional U(1)-Lie group action. The results of this chapter are published in a slightly extended form in [DR18b].

Our results concerning equivariant normal forms and singular symplectic reduction in infinite dimensions open many interesting avenues for further research, some of which are described in Chapter 6.

In the appendix, we summarize without proofs the relevant background material. Appendix A outlines the calculus of infinite-dimensional manifolds with a primary focus on the Inverse Function Theorem as well as Lie group actions. In Appendix B, we give an overview of the theory of dual pairs, which forms the basis for our discussion of linear symplectic geometry in Section 4.1. For the convenience of the reader, a list of mathematical symbols[1] is included at the end of the thesis.

[1] In the electronic version of this thesis, many mathematical symbols in the text are supplied with a hyperlink pointing to the corresponding entry in the list of symbols.

# Normal Forms of Maps 2

## 2.1 Normal form of a linear map

In this section, we are concerned with the normal form of continuous linear maps between locally convex spaces. Recall that every $m \times n$ matrix $T$ with rank $r$ can be written in the form

$$T = P \begin{pmatrix} 0 & 0 \\ 0 & \mathbb{1}_{r\times r} \end{pmatrix} Q, \tag{2.1.1}$$

where $P$ and $Q$ are invertible $m \times m$ and $n \times n$ matrices, respectively. As we will see, a similar factorization is possible for continuous linear maps between locally convex spaces, which are relatively open and whose kernel and image are closed complemented subspaces. We call such operators regular. In the finite-dimensional context, it is well-known that the factorization (2.1.1) can be used to construct a so-called generalized inverse of $T$. Guided by this construction, we will show that every regular operator between locally convex spaces possesses a generalized inverse. With a view towards applications, we give a brief overview of the theory of Fredholm operators and of elliptic operators in the locally convex framework and, in particular, show that these operators are regular. Finally, as a preparation for the non-linear case, we extend the discussion of regularity to families of linear maps depending continuously on a parameter and to chain complexes. Most results of this section are well-established in the Banach setting, but their extension to more general locally convex spaces represents original work (if not otherwise indicated).

### 2.1.1 Regular operators and generalized inverses

Let $X$ and $Y$ be locally convex vector spaces and let $T\colon X \to Y$ be a continuous linear map. Assume that the kernel and the image of $T$ are closed complemented subspaces. Thus, there exist topological decompositions[1]

$$X = \operatorname{Ker} T \oplus \operatorname{Coim} T, \qquad Y = \operatorname{Coker} T \oplus \operatorname{Im} T, \tag{2.1.2}$$

[1] Let $X$ and $Y$ be closed subspaces of a locally convex space $Z$ such that $X \cap Y = \{0\}$. If the linear map $X \times Y \to Z$ given by $(x, y) \mapsto x + y$ is a topological isomorphism, then we say that $Z$ is the topological direct sum of $X$ and $Y$, and write $Z = X \oplus Y$. Depending on the context, we then write elements of $Z$ either as a pair $(x, y)$ or as a sum $x + y$ with $x \in X$ and $y \in Y$.

where $\operatorname{Coim} T$ and $\operatorname{Coker} T$ are closed subspaces[1] of $X$ and $Y$, respectively. Then, $T$ factorizes as

$$\begin{array}{ccc} X & \xrightarrow{\quad T \quad} & Y, \\ \searrow & & \nearrow \\ \operatorname{Coim} T & \xrightarrow{\hat{T}} & \operatorname{Im} T \end{array} \tag{2.1.3}$$

where $\hat{T}\colon \operatorname{Coim} T \longrightarrow \operatorname{Im} T$ is a continuous linear bijection. The induced operator[2] $\hat{T}$ is called the *core* of $T$. Recall that a continuous linear map $T\colon X \longrightarrow Y$ is called relatively open if it is open as a map $T\colon X \longrightarrow \operatorname{Im} T$. In other words, the image under $T$ of every open subset of $X$ is open in the relative topology of $\operatorname{Im} T$. Thus, a continuous bijection is relatively open if and only if it is a topological isomorphism. In particular, $T$ is relatively open if and only if its core $\hat{T}$ is a topological isomorphism.

Definition 2.1.1 A continuous linear map $T\colon X \longrightarrow Y$ between locally convex spaces is called *regular* if $T$ is relatively open and $\operatorname{Ker} T$ as well as $\operatorname{Im} T$ are closed complemented subspaces of $X$ and $Y$, respectively. ◊

Thus, the core $\hat{T}\colon \operatorname{Coim} T \longrightarrow \operatorname{Im} T$ of a regular operator $T\colon X \longrightarrow Y$ is a topological isomorphism. Moreover, the above discussion shows that every regular operator $T$ can be written in a form similar to (2.1.1):

$$T = P \begin{pmatrix} 0 & 0 \\ 0 & \hat{T} \end{pmatrix} Q, \tag{2.1.4}$$

where $Q\colon X \longrightarrow \operatorname{Ker} T \oplus \operatorname{Coim} T$ and $P\colon \operatorname{Coker} T \oplus \operatorname{Im} T \longrightarrow Y$ are the natural isomorphisms determined by the decompositions (2.1.2). For this reason, we also say that a regular operator can be brought into a normal form.

Lemma 2.1.2 *Let $T\colon X \longrightarrow Y$ be a continuous linear map between locally convex spaces $X$ and $Y$.*

*(i) Assume that $X$ and $Y$ are Fréchet spaces. Then, $T$ is regular if and only if* $\operatorname{Ker} T$ *as well as* $\operatorname{Im} T$ *are closed complemented subspaces.*

*(ii) If* $\operatorname{Im} T$ *is finite-dimensional, then $T$ is regular. This is the case, in particular, when $X$ or $Y$ are finite-dimensional.* ◊

---

[1] In general, the coimage and cokernel are defined as $\operatorname{Coim} T = X/\operatorname{Ker} T$ and $\operatorname{Coker} T = Y/\operatorname{Im} T$, respectively. There exists, of course, no canonical realization of these quotient spaces as subspaces of $X$ and $Y$. Nonetheless, the choice of complements $A$ and $B$ of $\operatorname{Ker} T$ and $\operatorname{Im} T$, respectively, leads to the identifications $\operatorname{Coim} T \simeq A$ and $\operatorname{Coker} T \simeq B$. It is in this sense and with a slight abuse of notation that we view $\operatorname{Coim} T$ and $\operatorname{Coker} T$ as subspaces of $X$ and $Y$, respectively. The reader should keep in mind that these subspaces are not canonically associated to $T$.

[2] We use the word "operator" interchangeably with "continuous linear map".

*Proof.* The first claim follows directly from the open mapping theorem [Tre67, Theorem 17.1], which states that an operator between Fréchet spaces is relatively open if and only if its image is closed.

Suppose now that $\operatorname{Im} T$ is finite-dimensional. Then, by [Köt83, Proposition 20.5.5], $\operatorname{Im} T$ is closed and has a topological complement. Moreover, $\operatorname{Ker} T$ has finite codimension in $X$ and hence is topologically complemented according to [Köt83, Proposition 15.8.2]. The core of $T$ is a continuous linear bijection between finite-dimensional spaces and hence is a topological isomorphism. Thus, $T$ is regular. □

The core $\hat{T}$ of a regular operator $T$ has a continuous inverse. Hence, it is natural to expect that $T$ itself can be inverted in a certain sense. This is in fact the case, as we will see below, and the right notion turns out to be that of a generalized inverse.

Definition 2.1.3 Let $T\colon X \longrightarrow Y$ be a continuous linear map between locally convex spaces. A continuous linear map $S\colon Y \longrightarrow X$ satisfying

$$T \circ S \circ T = T \tag{2.1.5}$$

is called a *generalized inverse* of $T$. If, additionally,

$$S \circ T \circ S = S \tag{2.1.6}$$

holds, then $S$ is referred to as a *reflexive generalized inverse*. ◊

In finite dimensions, generalized inverses were first studied by Rao [Rao62] and have since become a valuable tool to find solutions to singular systems of linear equations. In general, a (reflexive) generalized inverse is not unique unless further conditions are imposed. For example, in the context of Hilbert spaces, it is convenient to require a certain hermiticity condition, which leads to the notion of a pseudoinverse independently introduced by Moore [Moo20] and Penrose [Pen55].

The following relationship between regularity and generalized inverses can be found (in parts) in [Har87, Theorem 3.8.2] for maps between normed spaces, but the proof can be generalized to maps between arbitrary locally convex spaces without much effort.

Proposition 2.1.4 *For a continuous linear map $T\colon X \longrightarrow Y$ between locally convex spaces the following are equivalent:*

*(i) $T$ is regular.*

*(ii) There exist topological decompositions $X = \operatorname{Ker} T \oplus \operatorname{Coim} T$ and $Y = \operatorname{Coker} T \oplus \operatorname{Im} T$ such that the core $\hat{T}\colon \operatorname{Coim} T \longrightarrow \operatorname{Im} T$ of $T$ is a topological isomorphism.*

*(iii) There exists a generalized inverse of $T$.*

*(iv) There exists a reflexive generalized inverse of $T$.* ◊

*Proof.* The equivalence of statements (i) and (ii) has been discussed above. We will establish the following implications: (ii) → (iv) → (iii) → (i).

First, assume that we have topological decompositions $X = \operatorname{Ker} T \oplus \operatorname{Coim} T$ and $Y = \operatorname{Coker} T \oplus \operatorname{Im} T$ such that the induced map $\hat{T}\colon \operatorname{Coim} T \to \operatorname{Im} T$ is a topological isomorphism. The continuous linear map $S\colon Y \to X$ defined by $S = \hat{T}^{-1} \circ \operatorname{pr}_{\operatorname{Im} T}$ satisfies

$$T \circ S = \hat{T} \circ \hat{T}^{-1} \circ \operatorname{pr}_{\operatorname{Im} T} = \operatorname{pr}_{\operatorname{Im} T}, \tag{2.1.7a}$$

$$S \circ T = \hat{T}^{-1} \circ \hat{T} \circ \operatorname{pr}_{\operatorname{Coim} T} = \operatorname{pr}_{\operatorname{Coim} T}, \tag{2.1.7b}$$

where $\operatorname{pr}_A$ denotes the projection onto the subspace $A$. Hence, by composing with $T$ and $S$, respectively, we obtain

$$T \circ S \circ T = \operatorname{pr}_{\operatorname{Im} T} \circ T = T, \tag{2.1.8a}$$

$$S \circ T \circ S = \operatorname{pr}_{\operatorname{Coim} T} \circ S = S. \tag{2.1.8b}$$

Thus, $S$ is a reflexive generalized inverse.

The existence of a reflexive generalized inverse clearly implies the existence of a generalized inverse.

Now, let $S$ be generalized inverse of $T$. The identity $T \circ S \circ T = T$ implies that $T \circ S$ and $S \circ T$ are projections. Obviously, $\operatorname{Ker} T \subseteq \operatorname{Ker}(S \circ T)$. We claim that both kernels actually coincide. Indeed, for $x \in \operatorname{Ker}(S \circ T)$ we have $Tx = T \circ (S \circ T)x = 0$, because $S$ is a generalized inverse. Similarly, $Tx = (T \circ S)(Tx)$ for $x \in X$ implies $\operatorname{Im}(T \circ S) = \operatorname{Im} T$. Hence, $\operatorname{Ker} T$ and $\operatorname{Im} T$ are images of continuous projections and thus they are closed and topologically complemented according to [Köt83, Proposition 15.8.1]. As above, denote the complements by $\operatorname{Coim} T$ and $\operatorname{Coker} T$, respectively. It remains to show that $\hat{T}\colon \operatorname{Coim} T \to \operatorname{Im} T$ has a continuous inverse. Let $\hat{S} \coloneqq \operatorname{pr}_{\operatorname{Coim} T} \circ S_{\restriction \operatorname{Im} T}\colon \operatorname{Im} T \to \operatorname{Coim} T$. Using $S \circ T = \operatorname{pr}_{\operatorname{Coim} T}$ and $T \circ S = \operatorname{pr}_{\operatorname{Im} T}$, it is straightforward to see that $\hat{S}$ is inverse to $\hat{T}$. This establishes the last remaining implication in the above chain and so completes the proof. □

We note in passing that the equivalence of (iii) and (iv) in the previous proposition can be established more directly. Indeed, given a generalized inverse $S$ of $T$, the operator $S \circ T \circ S$ is a reflexive generalized inverse of $T$.

Example 2.1.5 (Exterior differential) Let $M$ be a compact finite-dimensional manifold without boundary. Endow the space $\Omega^k(M)$ of differential $k$-forms on $M$ with its natural Fréchet topology. By Hodge theory, every Riemannian metric on $M$ yields the topological decompositions

$$\Omega^k(M) = \operatorname{Ker} \mathrm{d} \oplus \operatorname{Im} \mathrm{d}^*, \qquad \Omega^{k+1}(M) = \operatorname{Ker} \mathrm{d}^* \oplus \operatorname{Im} \mathrm{d}, \tag{2.1.9}$$

where $\mathrm{d}^*$ denotes the codifferential. Thus, the exterior differential $\mathrm{d}\colon \Omega^k(M) \longrightarrow \Omega^{k+1}(M)$ is a regular operator with $\operatorname{Coim}\mathrm{d} = \operatorname{Im}\mathrm{d}^*$ and $\operatorname{Coker}\mathrm{d} = \operatorname{Ker}\mathrm{d}^*$. Moreover, a straightforward calculation shows that the operator

$$\mathrm{d}^* \circ \triangle^{-1} \circ \mathrm{pr}_{\operatorname{Im}\mathrm{d}}\colon \Omega^{k+1}(M) \longrightarrow \Omega^k(M) \tag{2.1.10}$$

is a reflexive generalized inverse of d. Here, the Laplace operator $\triangle$ is viewed as a topological automorphism of $\operatorname{Im}\mathrm{d}$. ◊

Originally, generalized inverses were introduced to solve singular matrix equations. In the same spirit, generalized inverses in the locally convex setting can be used to solve singular boundary value problems. In this context, a generalized inverse is often called a *generalized Green's operator*. We limit ourselves here to providing a simple example that illustrates the main idea and refer the reader to [BG03; BS04] for further details.

Example 2.1.6 Following [KRR11, Section 4], let us consider the following boundary value problem:

$$u'' = f, \tag{2.1.11a}$$

$$u'(1) = 0 = u'(0) \tag{2.1.11b}$$

where $u, f \in C^\infty([0,1])$. It is clear that this problem is solvable only if $f$ satisfies the constraint $\int_0^1 f(y)\,\mathrm{d}y = 0$. Moreover, the solution is unique only up to a translation by a constant. Thus, to have a unique solution, the boundary value problem (2.1.11) needs to be supplemented by further conditions, say $u(1) = 0$ and $\int_0^1 f(y)\,\mathrm{d}y = 0$. In order to make the connection with the operator-theoretic setting, let $X \subseteq C^\infty([0,1])$ be the closed subspace of functions $u$ satisfying the boundary constraints $u'(1) = 0 = u'(0)$, set $Y = C^\infty([0,1])$ and let $T\colon X \longrightarrow Y$ be the continuous operator sending $u$ to $u''$. Then, $u$ is a smooth solution of the boundary value problem (2.1.11) if and only if $Tu = f$. As just noted, $T$ is neither injective nor surjective. Define the operator $S\colon Y \longrightarrow X$ by

$$S(f)(x) = x\int_0^x f(y)\,\mathrm{d}y - \int_0^x y f(y)\,\mathrm{d}y - \frac{1}{2}(x^2+1)\int_0^1 f(y)\,\mathrm{d}y + \int_0^1 y f(y)\,\mathrm{d}y \tag{2.1.12}$$

for $x \in [0,1]$. A straightforward calculation shows that $S$ is a reflexive generalized inverse of $T$. Note that we have

$$S \circ T(u)(x) = u(x) - u(1), \tag{2.1.13a}$$

$$T \circ S(f)(x) = f(x) - \int_0^1 f(y)\,\mathrm{d}y. \tag{2.1.13b}$$

Hence, the additional conditions on $u$ and $f$, which ensured that the supplemented boundary value problem had a unique solution, are encoded using the generalized inverse in terms of the projection maps $S \circ T$ and $T \circ S$. ♢

Proposition 2.1.4 gives a characterization of regular operators which is internal in the sense that it focuses on subspaces and the factorization of the operator. Following ideas of [SZ07], we now give an external characterization. The first result we establish towards this goal is the following construction of a generalized inverse as a part of a proper inverse of an extended operator, cf. [SZ07, Proposition 2.1].

Lemma 2.1.7 *Let $T\colon X \longrightarrow Y$, $T^+\colon Z^+ \longrightarrow Y$, $T^-\colon X \longrightarrow Z^-$ and $T^{+-}\colon Z^+ \longrightarrow Z^-$ be continuous linear maps between locally convex spaces such that the operator*

$$\begin{pmatrix} T & T^+ \\ T^- & T^{+-} \end{pmatrix} : X \oplus Z^+ \longrightarrow Y \oplus Z^- \tag{2.1.14}$$

*is invertible with a continuous inverse of the form*

$$\begin{pmatrix} S & S^- \\ S^+ & 0 \end{pmatrix}, \tag{2.1.15}$$

*where $S\colon Y \longrightarrow X$, $S^-\colon Z^- \longrightarrow X$ and $S^+\colon Y \longrightarrow Z^+$ are continuous linear maps. Then, $S$ is a generalized inverse of $T$ and the following identities hold:*

$$S^- \circ T^- = \mathrm{pr}_{\operatorname{Ker} T}, \qquad T^+ \circ S^+ = \mathrm{pr}_{\operatorname{Coker} T}. \tag{2.1.16}$$

*Moreover, $S$ is a reflexive generalized inverse of $T$ if and only if one of the following equivalent conditions is satisfied:*

*(i)* $\mathrm{pr}_{\operatorname{Ker} T} \circ S = 0$,

*(ii)* $S\restriction_{\operatorname{Coker} T} = 0$,

*(iii)* $S^- \circ T^{+-} \circ S^+ = 0$. ♢

*Proof.* We clearly have

$$\begin{pmatrix} S & S^- \\ S^+ & 0 \end{pmatrix} \circ \begin{pmatrix} T & T^+ \\ T^- & T^{+-} \end{pmatrix} = \begin{pmatrix} S \circ T + S^- \circ T^- & S \circ T^+ + S^- \circ T^{+-} \\ S^+ \circ T & S^+ \circ T^+ \end{pmatrix}, \tag{2.1.17a}$$

$$\begin{pmatrix} T & T^+ \\ T^- & T^{+-} \end{pmatrix} \circ \begin{pmatrix} S & S^- \\ S^+ & 0 \end{pmatrix} = \begin{pmatrix} T \circ S + T^+ \circ S^+ & T \circ S^- \\ T^- \circ S + T^{+-} \circ S^+ & T^- \circ S^- \end{pmatrix}. \tag{2.1.17b}$$

Since $\left(\begin{smallmatrix} T & T^+ \\ T^- & T^{+-} \end{smallmatrix}\right)$ and $\left(\begin{smallmatrix} S & S^- \\ S^+ & 0 \end{smallmatrix}\right)$ are inverse, these equations define various relations between the involved operators. For example, we read off that $T \circ S^- = 0$ and hence

$$T = T \circ S \circ T + T \circ S^- \circ T^- = T \circ S \circ T. \tag{2.1.18}$$

Thus, $S$ is a generalized inverse of $T$. Moreover, $S^+ \circ T^+ = \mathrm{id}_{Z^+}$ and $T^- \circ S^- = \mathrm{id}_{Z^-}$ imply that $T^+ \circ S^+$ and $S^- \circ T^-$ are projections. The relations (2.1.16) follow from

$$\mathrm{pr}_{\mathrm{Im}\,T} = T \circ S = \mathrm{id}_Y - T^+ \circ S^+, \tag{2.1.19a}$$

$$\mathrm{pr}_{\mathrm{Coim}\,T} = S \circ T = \mathrm{id}_X - S^- \circ T^-. \tag{2.1.19b}$$

Finally, we find

$$S \circ T \circ S + \mathrm{pr}_{\mathrm{Ker}\,T} \circ S = S \tag{2.1.20}$$

and

$$\mathrm{pr}_{\mathrm{Ker}\,T} \circ S = S^- \circ T^- \circ S = -S^- \circ T^{+-} \circ S^+ = S \circ T^+ \circ S^+ = S \circ \mathrm{pr}_{\mathrm{Coker}\,T}. \tag{2.1.21}$$

These identities imply that $S$ is a reflexive generalized inverse if and only if one, and thus all, of the above conditions are fulfilled. □

Lemma 2.1.7 suggests that the existence of a generalized inverse is tightly connected to the invertibility of an extended operator. In fact, we have the following characterization.

Proposition 2.1.8 *Let $T\colon X \longrightarrow Y$ be a continuous linear map between locally convex spaces. A continuous linear map $S\colon Y \longrightarrow X$ is a generalized inverse of $T$ if and only if there exist locally convex spaces $Z^\pm$ and continuous linear maps $T^+\colon Z^+ \longrightarrow Y$, $T^-\colon X \longrightarrow Z^-$ and $T^{+-}\colon Z^+ \longrightarrow Z^-$ such that*

$$\begin{pmatrix} T & T^+ \\ T^- & T^{+-} \end{pmatrix}^{-1} = \begin{pmatrix} S & S^- \\ S^+ & 0 \end{pmatrix}, \tag{2.1.22}$$

*where $S^-\colon Z^- \longrightarrow X$ and $S^+\colon Y \longrightarrow Z^+$ are continuous linear maps. Moreover, $S$ is a reflexive generalized inverse of $T$ if and only if we can choose $T^{+-} = 0$.* ◇

*Proof.* Suppose that $S\colon Y \longrightarrow X$ is a generalized inverse of $T$. Then, by Proposition 2.1.4, we have topological decompositions $X = \mathrm{Ker}\,T \oplus \mathrm{Coim}\,T$ and $Y = \mathrm{Coker}\,T \oplus \mathrm{Im}\,T$. Set $Z^+ = \mathrm{Coker}\,T$ and $Z^- = \mathrm{Ker}\,T$. Let $T^+\colon Z^+ \longrightarrow Y$, $S^-\colon Z^- \longrightarrow X$ and $T^-\colon X \longrightarrow Z^-$, $S^+\colon Y \longrightarrow Z^+$ be the natural injections and projections, respectively. Set $T^{+-} = -\mathrm{pr}_{\mathrm{Ker}\,T} \circ S\restriction_{\mathrm{Coker}\,T}$. A straightforward calculation using $T \circ S = \mathrm{pr}_{\mathrm{Im}\,T}$, $S \circ T = \mathrm{pr}_{\mathrm{Coim}\,T}$ and (2.1.17) shows that (2.1.22) holds. If, in addition, $S$ is a reflexive generalized inverse, then $S = S \circ T \circ S = S \circ \mathrm{pr}_{\mathrm{Im}\,T}$ implies $S \circ \mathrm{pr}_{\mathrm{Coker}\,T} = 0$ and so $T^{+-} = 0$. The converse direction has already been established in Lemma 2.1.7. □

An important class of examples of regular operators is given by Fredholm operators. Fredholm operators are usually studied as maps between Banach spaces (or Hilbert spaces) but most results extend to the locally convex setting. Here, we only give a brief account of the general theory of Fredholm operators between locally convex spaces and refer the reader to [Sch56; Sch59] for more details.

Definition 2.1.9 A continuous linear map $T\colon X \longrightarrow Y$ between locally convex spaces is called a *Fredholm operator* if $T$ is relatively open, the kernel of $T$ is a finite-dimensional subspace of $X$, and the image of $T$ is a finite-codimensional closed subspace of $Y$. The *index* $\operatorname{ind} T$ of a Fredholm operator $T$ is defined by

$$\operatorname{ind} T = \dim \operatorname{Ker} T - \dim \operatorname{Coker} T. \tag{2.1.23}$$

◊

Since finite-dimensional subspaces and finite-codimensional closed subspaces of a locally convex space are always topologically complemented according to [Köt83, Propositions 15.8.2 and 20.5.5], every Fredholm operator is regular.

As in the Banach setting, Fredholm operators are invertible modulo compact operators. A continuous linear map $K\colon X \longrightarrow Y$ between locally convex spaces is called *compact* if there exists a neighborhood $U$ of $0$ in $X$ such that the closure of $K(U)$ is a compact subset of $Y$. We refer to [Edw65, Chapter 9] for the theory of compact operators in the locally convex setting. According to [Sch56, Satz 12], a continuous linear map $T\colon X \longrightarrow Y$ between locally convex spaces is Fredholm if and only if there exists a continuous linear operator $S\colon Y \longrightarrow X$ such that

$$S \circ T = \mathrm{id}_X - K_1, \tag{2.1.24a}$$
$$T \circ S = \mathrm{id}_Y - K_2 \tag{2.1.24b}$$

holds for compact operators $K_1\colon X \longrightarrow X$ and $K_2\colon Y \longrightarrow Y$. The operator $S$ is called a *parametrix* of $T$. Every generalized inverse $S$ of a Fredholm operator $T$ is a parametrix, because we have

$$\mathrm{id}_X - S \circ T = \operatorname{Ker} T, \tag{2.1.25a}$$
$$\mathrm{id}_Y - T \circ S = \operatorname{Coim} T, \tag{2.1.25b}$$

and finite-rank operators are compact as a consequence of the Bolzano–Weierstrass Theorem.

### 2.1.2 Uniform regularity of operator families

As a preparation for the non-linear case, we extend the discussion of regularity to families of linear maps depending continuously on a parameter. Consider the following setup. Let $X$ and $Y$ be locally convex spaces, and let $P$ be an open neighborhood of $0$ in some locally convex space. A continuous map $T\colon P \times X \longrightarrow Y$ is called a *continuous family of linear maps* if, for all $p \in P$, the induced map $T_p \equiv T(p, \cdot)\colon X \longrightarrow Y$ is linear.

Recall that an operator $T\colon X \longrightarrow Y$ is regular if $\operatorname{Ker} T$ and $\operatorname{Im} T$ are closed complemented subspaces and the core $\hat{T} = \mathrm{pr}_{\mathrm{Im}\,T} \circ T_{\restriction \mathrm{Coim}\,T}$ is a topological isomorphism. This formulation suggests the following generalization of the notion of regularity to operator families.

Definition 2.1.10 A continuous family $T\colon P \times X \to Y$ of linear maps is called *uniformly regular* (at 0) if there exist topological decompositions

$$X = \operatorname{Ker} T_0 \oplus \operatorname{Coim} T_0, \qquad Y = \operatorname{Coker} T_0 \oplus \operatorname{Im} T_0, \tag{2.1.26}$$

and, for every $p \in P$, the map $\tilde{T}_p = \mathrm{pr}_{\mathrm{Im}\,T_0} \circ (T_p)_{\restriction \mathrm{Coim}\,T_0}\colon \operatorname{Coim} T_0 \to \operatorname{Im} T_0$ is a topological isomorphism such that the inverses form a continuous family $P \times \operatorname{Im} T_0 \to \operatorname{Coim} T_0$. ◇

Note that, in particular, $T_0$ is a regular operator. If the space of invertible maps is open in the space of all continuous linear maps, then, for every continuous family $T\colon P \times X \to Y$ with $T_0$ being regular, one can shrink $P$ to pass to a uniformly regular family. However, this openness property fails when one leaves the Banach realm as the following example demonstrates.

Example 2.1.11 Let $X = \mathrm{C}^\infty([0,1], \mathbb{R})$ and consider the continuous family $L\colon \mathbb{R} \times X \to X$ of linear differential operators defined by

$$L_r(g)(x) = g - rxg' \tag{2.1.27}$$

for $r \in \mathbb{R}$, $g \in X$ and $x \in [0,1]$. Clearly, $L_0$ is the identity operator on $X$. On the other hand, for every $n \in \mathbb{N}$, the operator $L_{1/n}$ annihilates the function $g_n(x) = x^n$. Hence, $L_r$ fails to be injective for arbitrarily small $r$ and so $L$ is not uniformly regular. Note that if we try to formulate this problem in Banach spaces and view $L_r$ as an operator from $\mathrm{C}^1([0,1])$ to $\mathrm{C}^0([0,1])$, then $L_0$ is only the inclusion of a dense subspace. ◇

It is clear from the definition that the property of uniform regularity of a continuous family is invariant under reparametrization. For ease of reference, let us record this.

Lemma 2.1.12 *Let $T\colon P \times X \to Y$ be a continuous uniformly regular family of linear maps and let $Q$ be an open neighborhood of 0 in some locally convex space. For every continuous map $L\colon Q \to P$ with $L(0) = 0$, the induced family*

$$T \circ (L \times \mathrm{id}_X)\colon Q \times X \to Y \tag{2.1.28}$$

*is uniformly regular (at 0).* ◇

Let $T\colon P \times X \to Y$ be a continuous family of relatively open linear maps. If $\operatorname{Ker} T_p = \operatorname{Ker} T_0$ and $\operatorname{Im} T_p = \operatorname{Im} T_0$ hold for all $p \in P$, then $\tilde{T}_p$ is clearly a topological isomorphism. Hence, in this case, the family $T$ is uniformly regular. In the converse direction, uniform regularity implies a semi-continuity property of the kernel and the image.

Proposition 2.1.13 *Let $T\colon P \times X \longrightarrow Y$ be a continuous family of linear maps. If $T$ is uniformly regular, then the following holds:*

*(i) The kernel of $T$ is upper semi-continuous at $0$ in the sense that $\operatorname{Ker} T_p \subseteq \operatorname{Ker} T_0$ for all $p \in P$.*

*(ii) The image of $T$ is lower semi-continuous at $0$ in the sense that $\operatorname{Im} T_p \supseteq \operatorname{Im} T_0$ for all $p \in P$.* ◊

*Proof.* The inclusions $\operatorname{Ker} T_p \subseteq \operatorname{Ker} T_0$ and $\operatorname{Im} T_p \supseteq \operatorname{Im} T_0$ need to be valid, because otherwise $\tilde{T}_p = \mathrm{pr}_{\operatorname{Im} T_0} \circ (T_p)\restriction_{\operatorname{Coim} T_0}$ cannot be an isomorphism from $\operatorname{Coim} T_0$ to $\operatorname{Im} T_0$. □

Similar semi-continuity properties are well-known to hold for families of Fredholm operators between Banach spaces, see e.g. [Hör07, Corollary 19.1.6]. Before we come to the uniform regularity of Fredholm operators, it is convenient to establish a characterization of uniform regularity analogous to Proposition 2.1.8. We need the following basic result concerning the invertibility of block matrices in terms of the Schur complement (whose proof is an easy exercise left to the reader, cf. [SZ07, Lemma 3.1]).

Lemma 2.1.14 *Let $A_{11}\colon X_1 \longrightarrow Y_1$, $A_{12}\colon X_2 \longrightarrow Y_1$, $A_{21}\colon X_1 \longrightarrow Y_2$ and $A_{22}\colon X_2 \longrightarrow Y_2$ be continuous linear maps between locally convex spaces such that*

$$\begin{pmatrix} A_{11} & A_{12} \\ A_{21} & A_{22} \end{pmatrix}^{-1} = \begin{pmatrix} B_{11} & B_{12} \\ B_{21} & B_{22} \end{pmatrix} \tag{2.1.29}$$

*for continuous linear maps $B_{ij}$ for $i, j = 1, 2$. If $B_{22}$ is a topological isomorphism, then so is $A_{11}$ and the inverse is given by*

$$A_{11}^{-1} = B_{11} - B_{12} B_{22}^{-1} B_{21}. \tag{2.1.30}$$

◊

With this preparation at hand, we can give the following characterization of uniformly regular operator families analogous to Proposition 2.1.8.

Proposition 2.1.15 *Let $T\colon P \times X \longrightarrow Y$ be a continuous family of linear maps. Then, the following are equivalent:*

*(i) $T$ is uniformly regular.*

*(ii) There exist locally convex spaces $Z^\pm$, continuous linear maps $T^+\colon Z^+ \longrightarrow Y$ and $T^-\colon X \longrightarrow Z^-$, and continuous families of linear maps $S\colon P \times Y \longrightarrow X$, $S^-\colon P \times Z^- \longrightarrow X$, $S^+\colon P \times Y \longrightarrow Z^+$ and $S^{-+}\colon P \times Z^- \longrightarrow Z^+$ with $S_0^{-+} = 0$ such that*

$$\begin{pmatrix} T_p & T^+ \\ T^- & 0 \end{pmatrix}^{-1} = \begin{pmatrix} S_p & S_p^- \\ S_p^+ & S_p^{-+} \end{pmatrix}, \tag{2.1.31}$$

*holds for all $p \in P$ and such that the operators*

$$\Gamma_p \equiv \begin{pmatrix} \mathrm{pr}_{\mathrm{Ker}\,T_0} \circ (S_p)_{\restriction \mathrm{Coker}\,T_0} & \mathrm{pr}_{\mathrm{Ker}\,T_0} \circ S_p^- \\ (S_p^+)_{\restriction \mathrm{Coker}\,T_0} & S_p^{-+} \end{pmatrix} : \mathrm{Coker}\,T_0 \oplus Z^- \to \mathrm{Ker}\,T_0 \oplus Z^+ \tag{2.1.32}$$

*are invertible for all $p \in P$ and their inverses form a continuous family.* ◇

*Proof.* First, suppose that $T$ is a uniformly regular family of linear maps. Then, by definition, we have topological decompositions $X = \mathrm{Ker}\,T_0 \oplus \mathrm{Coim}\,T_0$, $Y = \mathrm{Coker}\,T_0 \oplus \mathrm{Im}\,T_0$ and, for every $p \in P$, the map

$$\tilde{T}_p = \mathrm{pr}_{\mathrm{Im}\,T_0} \circ (T_p)_{\restriction \mathrm{Coim}\,T_0} \colon \mathrm{Coim}\,T_0 \to \mathrm{Im}\,T_0 \tag{2.1.33}$$

is a topological isomorphism. Set $Z^+ = \mathrm{Coker}\,T_0$ and $Z^- = \mathrm{Ker}\,T_0$, and let $T^+\colon Z^+ \to Y$ and $T^-\colon X \to Z^-$ be the canonical inclusion and projection, respectively. Moreover, set $S_p := \tilde{T}_p^{-1} \circ \mathrm{pr}_{\mathrm{Im}\,T_0} \colon Y \to X$ and define $S_p^{-+}\colon Z^- \to Z^+$ by

$$S_p^{-+} = \mathrm{pr}_{\mathrm{Coker}\,T_0} \circ (T_p \circ \tilde{T}_p^{-1} \circ \mathrm{pr}_{\mathrm{Im}\,T_0} - \mathrm{id}_Y) \circ T_p \circ \mathrm{pr}_{\mathrm{Ker}\,T_0}\,. \tag{2.1.34}$$

Finally, define $S_p^\pm$ by

$$S_p^+ = \mathrm{pr}_{\mathrm{Coker}\,T_0} \circ (\mathrm{id}_Y - T_p \circ S_p) \colon Y \to Z^+, \tag{2.1.35a}$$

$$S_p^- = (\mathrm{id}_X - S_p \circ T_p)_{\restriction \mathrm{Ker}\,T_0} \colon Z^- \to X. \tag{2.1.35b}$$

Since, by definition, $\tilde{T}^{-1}$ is a continuous family $P \times \mathrm{Im}\,T_0 \to \mathrm{Coim}\,T_0$, the families $S, S^\pm, S^{-+}$ defined above are continuous. Furthermore, a direct calculation yields

$$T_p \circ S_p = \mathrm{pr}_{\mathrm{Im}\,T_0} + \mathrm{pr}_{\mathrm{Coker}\,T_0} \circ T_p \circ \tilde{T}_p^{-1} \circ \mathrm{pr}_{\mathrm{Im}\,T_0}, \tag{2.1.36a}$$

$$S_p \circ T_p = \mathrm{pr}_{\mathrm{Coim}\,T_0} + \tilde{T}_p^{-1} \circ \mathrm{pr}_{\mathrm{Im}\,T_0} \circ (T_p)_{\restriction \mathrm{Ker}\,T_0}. \tag{2.1.36b}$$

Hence, we obtain

$$T_p \circ S_p \circ T_p = T_p + S_p^{-+}. \tag{2.1.37}$$

Using the identities (2.1.36) and (2.1.37), it is straightforward to check that

$$\begin{pmatrix} T_p & T^+ \\ T^- & 0 \end{pmatrix}^{-1} = \begin{pmatrix} S_p & S_p^- \\ S_p^+ & S_p^{-+} \end{pmatrix} \tag{2.1.38}$$

holds for every $p \in P$. Moreover, for $\Gamma_p$, we find here

$$\Gamma_p = \begin{pmatrix} \mathrm{pr}_{\mathrm{Ker}\,T_0} \circ (S_p)_{\restriction \mathrm{Coker}\,T_0} & \mathrm{pr}_{\mathrm{Ker}\,T_0} \circ S_p^- \\ (S_p^+)_{\restriction \mathrm{Coker}\,T_0} & S_p^{-+} \end{pmatrix} = \begin{pmatrix} 0 & \mathrm{id}_{\mathrm{Ker}\,T_0} \\ \mathrm{id}_{\mathrm{Coker}\,T_0} & S_p^{-+} \end{pmatrix}, \tag{2.1.39}$$

which is clearly invertible with continuous inverse

$$\Gamma_p^{-1} = \begin{pmatrix} -S_p^{-+} & \mathrm{id}_{\operatorname{Coker} T_0} \\ \mathrm{id}_{\operatorname{Ker} T_0} & 0 \end{pmatrix}. \tag{2.1.40}$$

Conversely, let $T^{\pm}, S, S^{\pm}$ and $S^{-+}$ satisfying the assumptions of the second point above. Since $S_0^{-+} = 0$, Proposition 2.1.8 implies that $S_0$ is a reflexive generalized inverse of $T_0$. Hence, by Proposition 2.1.4, there exist topological decompositions $X = \operatorname{Ker} T_0 \oplus \operatorname{Coim} T_0$ and $Y = \operatorname{Coker} T_0 \oplus \operatorname{Im} T_0$. It remains to show that $\tilde{T}_p = \mathrm{pr}_{\operatorname{Im} T_0} \circ (T_p)\restriction_{\operatorname{Coim} T_0}$ is a topological isomorphism and that the inverses form a continuous family. For this purpose, we write all operators in block form with respect to the decompositions $X = \operatorname{Coim} T_0 \oplus \operatorname{Ker} T_0$ and $Y = \operatorname{Im} T_0 \oplus \operatorname{Coker} T_0$ (note the non-standard order of the summands). Using this convention, the identity (2.1.31) becomes:

$$\begin{aligned} &\begin{pmatrix} \tilde{T}_p & \mathrm{pr}_{\operatorname{Im} T_0} \circ (T_p)\restriction_{\operatorname{Ker} T_0} & \mathrm{pr}_{\operatorname{Im} T_0} \circ T^+ \\ \mathrm{pr}_{\operatorname{Coker} T_0} \circ (T_p)\restriction_{\operatorname{Coim} T_0} & \mathrm{pr}_{\operatorname{Coker} T_0} \circ (T_p)\restriction_{\operatorname{Ker} T_0} & \mathrm{pr}_{\operatorname{Coker} T_0} \circ T^+ \\ (T^-)\restriction_{\operatorname{Coim} T_0} & (T^-)\restriction_{\operatorname{Ker} T_0} & 0 \end{pmatrix}^{-1} = \\ &\begin{pmatrix} \mathrm{pr}_{\operatorname{Coim} T_0} \circ (S_p)\restriction_{\operatorname{Im} T_0} & \mathrm{pr}_{\operatorname{Coim} T_0} \circ (S_p)\restriction_{\operatorname{Coker} T_0} & \mathrm{pr}_{\operatorname{Coim} T_0} \circ S_p^- \\ \mathrm{pr}_{\operatorname{Ker} T_0} \circ (S_p)\restriction_{\operatorname{Im} T_0} & \mathrm{pr}_{\operatorname{Ker} T_0} \circ (S_p)\restriction_{\operatorname{Coker} T_0} & \mathrm{pr}_{\operatorname{Ker} T_0} \circ S_p^- \\ (S_p^+)\restriction_{\operatorname{Im} T_0} & (S_p^+)\restriction_{\operatorname{Coker} T_0} & S_p^{-+} \end{pmatrix}, \end{aligned} \tag{2.1.41}$$

which should be read as an operator from $\operatorname{Im} T_0 \oplus \operatorname{Coker} T_0 \oplus Z^-$ to $\operatorname{Coim} T_0 \oplus \operatorname{Ker} T_0 \oplus Z^+$. Note that the lower right corner of the right-hand side coincides with the operator $\Gamma_p$. Since $\Gamma_p$ is invertible, Lemma 2.1.14 shows that $\tilde{T}_p$ is invertible, too. Moreover, the inverse is given by

$$\begin{aligned} \tilde{T}_p^{-1} &= \mathrm{pr}_{\operatorname{Coim} T_0}(S_p)\restriction_{\operatorname{Im} T_0} \\ &\quad - \begin{pmatrix} \mathrm{pr}_{\operatorname{Coim} T_0} \circ (S_p)\restriction_{\operatorname{Coker} T_0} \\ \mathrm{pr}_{\operatorname{Coim} T_0} \circ S_p^- \end{pmatrix} \Gamma_p^{-1} \begin{pmatrix} \mathrm{pr}_{\operatorname{Ker} T_0} \circ (S_p)\restriction_{\operatorname{Im} T_0} \\ (S_p^+)\restriction_{\operatorname{Im} T_0} \end{pmatrix} \end{aligned} \tag{2.1.42}$$

and thus forms a continuous family $P \times \operatorname{Im} T_0 \longrightarrow \operatorname{Coim} T_0$. Hence, $T$ is uniformly regular. □

Remark 2.1.16 Note that the relations (2.1.36) contain additional terms in contrast to the simpler form (2.1.7). These additional terms imply that, for $p \neq 0$, $S_p$ is not a generalized inverse of $T_p$ but that instead (2.1.37) hold. One can use this identity as a starting point to formalize the notion of a generalized inverse of an operator family. Since we do not need this in the remainder, we do not explore it further here. ◊

If $T_0$ is a Fredholm operator, then Proposition 2.1.15 takes a slightly simpler form.

Corollary 2.1.17 *Let $T\colon P \times X \longrightarrow Y$ be a continuous family of linear maps such that $T_0$ is a Fredholm operator. Then, $T$ is uniformly regular if and only if there* *exist finite-dimensional spaces $Z^{\pm}$ and continuous linear maps $T^+\colon Z^+ \longrightarrow Y$ and $T^-\colon X \longrightarrow Z^-$ such that, after possibly shrinking $P$,*

$$\begin{pmatrix} T_p & T^+ \\ T^- & 0 \end{pmatrix}^{-1} = \begin{pmatrix} S_p & S_p^- \\ S_p^+ & S_p^{-+} \end{pmatrix} \tag{2.1.43}$$

*holds for all $p \in P$, where $S\colon P \times Y \longrightarrow X$, $S^-\colon P \times Z^- \longrightarrow X$, $S^+\colon P \times Y \longrightarrow Z^+$ and $S^{-+}\colon P \times Z^- \longrightarrow Z^+$ are continuous families of linear maps with $S_0^{-+} = 0$.* ◊

*Proof.* If $T$ is uniformly regular, then the proof of Proposition 2.1.15 shows that one can choose $Z^+ = \operatorname{Coker} T_0$ and $Z^- = \operatorname{Ker} T_0$. Both spaces are finite-dimensional, because $T_0$ is a Fredholm operator. This establishes one direction.

Conversely, let $T^{\pm}, S, S^{\pm}, S^{-+}$ be given as stated above. By Proposition 2.1.15, it suffices to show that the operator

$$\Gamma_p = \begin{pmatrix} \mathrm{pr}_{\operatorname{Ker} T_0} \circ (S_p)\restriction_{\operatorname{Coker} T_0} & \mathrm{pr}_{\operatorname{Ker} T_0} \circ S_p^- \\ (S_p^+)\restriction_{\operatorname{Coker} T_0} & S_p^{-+} \end{pmatrix} : \operatorname{Coker} T_0 \oplus Z^- \longrightarrow \operatorname{Ker} T_0 \oplus Z^+ \tag{2.1.44}$$

is invertible and that the inverses form a continuous family. Since $S_0^{-+} = 0$, Lemma 2.1.7 implies that $S_0$ is a reflexive generalized inverse of $T_0$ and, moreover, that $S_0^- \circ T^- = \mathrm{pr}_{\operatorname{Ker} T_0}$ and $T^+ \circ S_0^+ = \mathrm{pr}_{\operatorname{Coker} T_0}$ hold. A straightforward calculation using these identities shows that we have

$$\Gamma_0 = \begin{pmatrix} 0 & \mathrm{pr}_{\operatorname{Ker} T_0} \circ S_0^- \\ (S_0^+)\restriction_{\operatorname{Coker} T_0} & 0 \end{pmatrix}, \tag{2.1.45}$$

and

$$\Gamma_0^{-1} = \begin{pmatrix} 0 & \mathrm{pr}_{\operatorname{Coker} T_0} \circ T^+ \\ (T^-)\restriction_{\operatorname{Ker} T_0} & 0 \end{pmatrix}. \tag{2.1.46}$$

Since, for every $p \in P$, $\Gamma_p$ is an operator between finite-dimensional space and $\Gamma_0$ is invertible, we can shrink $P$ in such a way that $\Gamma_p$ is invertible for all $p \in P$ and that the inverses form a continuous family. Thus, Proposition 2.1.15 implies that $T$ is uniformly regular. □

Proposition 2.1.18 *Let $T\colon P \times X \longrightarrow Y$ be a continuous family of linear maps such that $T_0$ is a Fredholm operator. If $T$ is uniformly regular, then $T_p$ is a Fredholm operator and* $\operatorname{ind} T_p = \operatorname{ind} T_0$ *for all $p \in P$.* ◊

*Proof.* By Corollary 2.1.17, there exist finite-dimensional spaces $Z^{\pm}$ and continuous linear maps $T^+\colon Z^+ \longrightarrow Y$ and $T^-\colon X \longrightarrow Z^-$ such that

$$\begin{pmatrix} T_p & T^+ \\ T^- & 0 \end{pmatrix} \tag{2.1.47}$$

is invertible for all $p \in P$. This is only possible if $\operatorname{Ker} T_p$ and $\operatorname{Coker} T_p$ are finite-dimensional so that $T_p$ has to be Fredholm. Moreover, using the invariance of the index under finite-rank perturbations, we have

$$0 = \operatorname{ind} \begin{pmatrix} T_p & T^+ \\ T^- & 0 \end{pmatrix} = \operatorname{ind} T_p + \dim Z^+ - \dim Z^- = \operatorname{ind} T_p - \operatorname{ind} T_0, \qquad (2.1.48)$$

which establishes the formula for the index. □

Remark 2.1.19 (Uniform regularity in the tame Fréchet category) In the previous two sections, we considered the general setting of continuous linear maps between locally convex spaces. It is clear that similar results hold in the tame Fréchet category if the word "tame" is inserted in the right places (see Appendix A.1 for a brief overview of the main concepts of tame Fréchet spaces). Since uniform regularity in the Fréchet setting will play a major role later, let us spell out some of the details. Let $X$ and $Y$ be tame Fréchet spaces and let $T\colon P \times X \to Y$ be a tame smooth family of linear maps. Then, $T$ is called *uniformly tame regular* if there exist tame decompositions[1] $X = \operatorname{Ker} T_0 \oplus \operatorname{Coim} T_0$ and $Y = \operatorname{Coker} T_0 \oplus \operatorname{Im} T_0$, and, for every $p \in P$, the map $\tilde{T}_p = \operatorname{pr}_{\operatorname{Im} T_0} \circ (T_p)\restriction_{\operatorname{Coim} T_0}\colon \operatorname{Coim} T_0 \to \operatorname{Im} T_0$ is a tame isomorphism such that the inverses form a tame smooth family $P \times \operatorname{Im} T_0 \to \operatorname{Coim} T_0$. ◊

### 2.1.3 Families of elliptic operators

We now discuss an important class of families of Fredholm operators which are uniformly regular.

For this purpose, let $E \to M$ and $F \to M$ be finite-dimensional vector bundles over a compact manifold $M$ without boundary. Endow the spaces $\mathcal{E}$ and $\mathcal{F}$ of smooth sections of $E$ and $F$, respectively, with the compact-open $C^\infty$-topology. With respect to this topology, these section spaces are tame Fréchet spaces, see [Ham82, Theorem II.2.3.1]. A linear partial differential operator $L\colon \mathcal{E} \to \mathcal{F}$ of degree $r$ assigns to every section $\phi$ of $E$ a section $L(\phi)$ of $F$ in such a manner that $L(\phi)$ depends only on the derivatives of $\phi$ of degree $r$ or less. Hence, $L$ factors through the jet bundle $J^r E$ as follows:

$$L\colon \mathcal{E} \xrightarrow{\ j^r\ } \Gamma^\infty(J^r E) \xrightarrow{\ f_*\ } \mathcal{F}, \qquad (2.1.49)$$

where $j^r$ denotes the $r$-th jet prolongation and $f_*$ is the push-forward by some vertical morphism of vector bundles $f\colon J^r E \to F$. We refer to $f$ as the coefficients of $L$. Conversely, every vertical vector bundle morphism $f\colon J^r E \to F$ induces a differential operator $L_f : \mathcal{E} \to \mathcal{F}$ of degree $r$ and with coefficients $f$. To each differential operator $L_f$ of degree $r$, we associate the principal symbol $\sigma_f$ by, roughly speaking, replacing each partial derivative in the highest order term by a variable in the cotangent bundle of $M$. The principal symbol $\sigma_f$ is a homogeneous polynomial of degree $r$ on $\mathrm{T}^*M$ with values in the bundle $\mathrm{L}(E, F)$, whose fiber over $m \in M$ consists of linear maps $E_m \to F_m$. Equivalently, we may view $\sigma_f$ as a vector bundle map $\overset{\star}{\tau}{}^*E \to \overset{\star}{\tau}{}^*F$ over the cotangent bundle $\overset{\star}{\tau} : \mathrm{T}^*M \to M$. A differential operator $L_f$ with coefficients $f$ is called *elliptic* if its symbol is invertible; that is, for each nonzero $p \in \mathrm{T}^*M$, the bundle map $\sigma_f(p, \ldots, p) \in \mathrm{L}(E, F)$ is invertible. It is a standard result in elliptic theory that every elliptic differential operator is a Fredholm operator between appropriate Sobolev spaces. Exactly the same arguments also show that an elliptic operator is Fredholm as an operator between spaces of *smooth* sections. Indeed, every elliptic operator $L_f$ can be inverted up to smoothing operators (see, e.g., [Wel07, Theorem IV.4.4]). That is, there exists a pseudo-differential operator $S : \mathcal{F} \to \mathcal{E}$ such that

$$S \circ L_f = \mathrm{id}_{\mathcal{E}} - K_1, \tag{2.1.50a}$$

$$L_f \circ S = \mathrm{id}_{\mathcal{F}} - K_2 \tag{2.1.50b}$$

holds for smoothing operators $K_1 : \mathcal{E} \to \mathcal{E}$ and $K_2 : \mathcal{F} \to \mathcal{F}$. On compact manifolds, every smoothing operator is a compact operator [Wel07, Proposition IV.4.5] and thus $S$ is a parametrix. In particular, $L_f$ is a Fredholm operator.

[1] Let $Z$ be a tame Fréchet space and assume that $Z$ is the topological direct sum of closed subspaces $X$ and $Y$. We say that the sum $Z = X \oplus Y$ is tame if the map $X \times Y \to Z$ given by $(x, y) \mapsto x + y$ is a tame isomorphism.

The parameterization of a partial differential operator by their coefficients yields a tame smooth family

$$L : \Gamma^\infty\big(\mathrm{L}(\mathrm{J}^r E, F)\big) \times \mathcal{E} \to \mathcal{F}, \quad (f, \phi) \mapsto L_f(\phi) \tag{2.1.51}$$

of linear operators. Let $f_0 \in \Gamma^\infty\big(\mathrm{L}(\mathrm{J}^r E, F)\big)$ be such that $L_{f_0}$ is an elliptic differential operator. By [Ham82, Theorem II.3.3.3], there exist an open neighborhood $\mathcal{U}$ of $f_0$ in $\Gamma^\infty\big(\mathrm{L}(\mathrm{J}^r E, F)\big)$, finite-dimensional vector spaces $Z^\pm$ and continuous linear maps $L^+ : Z^+ \to Y$ and $L^- : X \to Z^-$ such that

$$\begin{pmatrix} L_f & L^+ \\ L^- & 0 \end{pmatrix} : \mathcal{E} \times Z^+ \to \mathcal{F} \times Z^- \tag{2.1.52}$$

is invertible for all $f \in \mathcal{U}$. Moreover, the inverses form a tame smooth family $\mathcal{U} \times \mathcal{F} \times Z^- \to \mathcal{E} \times Z^+$ of linear operators. Hence, by Corollary 2.1.17, $L\restriction_{\mathcal{U}}$ is uniformly tame regular at $f_0$. We have thus proved the following.

Theorem 2.1.20 *The tame smooth family $L$ defined in* (2.1.51) *is uniformly tame regular in a neighborhood $\mathcal{U}$ of every $f_0 \in \Gamma^\infty\big(\mathrm{L}(\mathrm{J}^r E, F)\big)$ for which $L_{f_0}$ is an elliptic differential operator.* ◇

Example 2.1.21 (Covariant exterior differential) Let $P \longrightarrow M$ be a finite-dimensional principal $G$-bundle over a compact manifold $M$. The space $\mathcal{C}(P)$ of connections on $P$ is an affine tame Fréchet space. Let $\underline{E}$ be a finite-dimensional vector space that is endowed with a representation of $G$ and let $E = P \times_G \underline{E}$ be the associated vector bundle. For a connection $A \in \mathcal{C}(P)$, denote the covariant exterior differential on differential $k$-forms with values in $E$ by $\mathrm{d}_A : \Omega^k(M, E) \longrightarrow \Omega^{k+1}(M, E)$. The same type of arguments as in Example 2.1.5 show that on 0-forms $\mathrm{d}_A : \Omega^0(M, E) \longrightarrow \Omega^1(M, E)$ is a regular operator. In particular, we have topological decompositions[1]

$$\Omega^0(M, E) = \operatorname{Ker} \mathrm{d}_A \oplus \operatorname{Im} \mathrm{d}_A^*, \qquad \Omega^1(M, E) = \operatorname{Ker} \mathrm{d}_A^* \oplus \operatorname{Im} \mathrm{d}_A, \tag{2.1.53}$$

where $\mathrm{d}_A^*$ denotes the codifferential, cf. [RS17, Theorem 6.1.9]. Now, consider the continuous family

$$T : \mathcal{C}(P) \times \Omega^0(M, E) \to \Omega^1(M, E), \qquad (A, \alpha) \mapsto \mathrm{d}_A \alpha \tag{2.1.54}$$

of linear operators. We claim that $T$ is uniformly regular at every $A_0 \in \mathcal{C}(P)$. To see this, let us introduce the so-called *Faddeev–Popov operator* (cf. [RS17, eq. (8.4.8)])

$$\triangle_{AA_0} = \mathrm{d}_{A_0}^* \circ \mathrm{d}_A : \Omega^0(M, E) \to \Omega^0(M, E) \tag{2.1.55}$$

for $A \in \mathcal{C}(P)$. Note that, for $A = A_0$, the operator $\triangle_{A_0 A_0}$ coincides with the covariant Laplacian $\mathrm{d}_{A_0}^* \mathrm{d}_{A_0}$ on 0-forms and thus is an elliptic operator. Clearly, $\operatorname{Ker} \triangle_{A_0 A_0} = \operatorname{Ker} \mathrm{d}_{A_0}$ and $\operatorname{Im} \triangle_{A_0 A_0} = \operatorname{Im} \mathrm{d}_{A_0}^*$. By Theorem 2.1.20, the family $\triangle_{AA_0}$ is uniformly tame regular at $A = A_0$. Thus, there exists an open neighborhood $\mathcal{U}$ of $A_0$ in $\mathcal{C}(P)$ such that

$$L_{AA_0} = \mathrm{pr}_{\operatorname{Im} \mathrm{d}_{A_0}^*} \circ \left(\triangle_{AA_0}\right)_{\restriction \operatorname{Im} \mathrm{d}_{A_0}^*} = \left(\mathrm{d}_{A_0}^* \circ \mathrm{d}_A\right)_{\restriction \operatorname{Im} \mathrm{d}_{A_0}^*} \tag{2.1.56}$$

is an automorphism of $\operatorname{Im} \mathrm{d}_{A_0}^* \subseteq \Omega^0(M, E)$ for all $A \in \mathcal{U}$. Now a straightforward calculation shows that

$$\left(L_{AA_0}^{-1} \circ \mathrm{d}_{A_0}^*\right)_{\restriction \operatorname{Im} \mathrm{d}_{A_0}} : \Omega^1(M, E) \supseteq \operatorname{Im} \mathrm{d}_{A_0} \to \operatorname{Im} \mathrm{d}_{A_0}^* \subseteq \Omega^0(M, E) \tag{2.1.57}$$

is a continuous inverse of $\tilde{T}_A = \mathrm{pr}_{\operatorname{Im} \mathrm{d}_{A_0}} \circ (\mathrm{d}_A)_{\restriction \operatorname{Im} \mathrm{d}_{A_0}^*}$ for all $A \in \mathcal{U}$. Hence, $T$ is uniformly regular at $A_0$. ◊

### 2.1.4 Elliptic complexes

In this section, the notion of regularity and uniform regularity will be extended to linear chain complexes. The main focus lies on elliptic complexes. The presentation is inspired to a large extend by [AB67], where elliptic complexes are studied in the language of Fréchet spaces. The results concerning uniform regularity of elliptic complexes are original work.

[1] Similar decompositions do not hold for higher $k$-forms, because $\mathrm{d}_A^2$ does not need to vanish.

Let $X_0, X_1, \dots, X_N$ be a sequence of locally convex vector spaces and let $T_i \colon X_i \to X_{i+1}$ be a sequence of continuous linear maps. We adopt the convention that the sequence is extended to all $i \in \mathbb{Z}$ by setting $X_i = \{0\}$ and $T_i = 0$ for $i < 0$ or $i > N$. We call the pair $(X_i, T_i)$ a *chain* and shall write it also in the following form:

$$\cdots \longrightarrow X_i \xrightarrow{T_i} X_{i+1} \longrightarrow \cdots . \tag{2.1.58}$$

A chain $(X_i, T_i)$ is called a *chain complex* if $T_{i+1} \circ T_i = 0$ for all $i \in \mathbb{Z}$.

Remark 2.1.22 Every continuous linear map $T \colon X \to Y$ can be viewed as a chain complex

$$0 \longrightarrow X \xrightarrow{T} Y \longrightarrow 0. \tag{2.1.59}$$

The reader may find it instructive to check that all concepts we are about to introduce for complexes boil down to the corresponding notions for the linear map $T$. ◊

Recall that a continuous linear map $T \colon X \to Y$ is regular if its kernel and image are closed complemented subspaces and $T$ factors trough a topological isomorphism. The following notion provides a natural generalization to chain complexes.

Definition 2.1.23 A chain complex $(X_i, T_i)$ of locally convex spaces is called *regular* if the following holds for every $i \in \mathbb{Z}$:

(i) The subspace $\operatorname{Im} T_{i-1}$ is closed in $X_i$ and there exist closed subspaces $H_i$ and $\operatorname{Coim} T_i$ of $X_i$ such that

$$X_i = \operatorname{Im} T_{i-1} \oplus \operatorname{Coim} T_i \oplus H_i \tag{2.1.60}$$

is a topological isomorphism.

(ii) The map $T_i$ factors as

$$\begin{array}{ccc} X_i & \xrightarrow{\quad\quad T_i \quad\quad} & X_{i+1}\,, \\ \searrow & & \nearrow \\ \operatorname{Coim} T_i & \xrightarrow{\hat{T}_i} & \operatorname{Im} T_i \end{array} \tag{2.1.61}$$

where $\hat{T}_i \colon \operatorname{Coim} T_i \to \operatorname{Im} T_i$ is a topological isomorphism.

If additionally, for every $i \in \mathbb{Z}$, $X_i$ is a tame Fréchet space, $T_i$ is a tame linear map, the decomposition (2.1.60) of $X_i$ is tame and $\hat{T}_i$ is a tame isomorphism, then we say that the chain complex $(X_i, T_i)$ is *tame regular*. ◊

By definition, for a regular chain complex $(X_i, T_i)$, we have

$$\operatorname{Ker} T_i = \operatorname{Im} T_{i-1} \oplus H_i\,, \tag{2.1.62}$$

which justifies the notion Coim $T_i$ for the subspace complementary to Im $T_{i-1} \oplus H_i$. Moreover, the subspaces $H_i$ are clearly identified with the *homology groups*, that is,

$$H_i \simeq \operatorname{Ker} T_i / \operatorname{Im} T_{i-1}. \tag{2.1.63}$$

A regular chain complex $(X_i, T_i)$ is called a *Fredholm* complex if all homology groups $H_i$ are finite-dimensional. In this case, the Euler characteristic of $(X_i, T_i)$ is defined as the alternating sum of its Betti numbers:

$$\chi(X_i, T_i) = \sum_{i \in \mathbb{Z}} (-1)^i \dim H_i. \tag{2.1.64}$$

Recall that a linear map is a Fredholm operator if and only if there exists a continuous parametrix. Generalizing this characterization, a chain complex $(X_i, T_i)$ is a regular Fredholm chain complex if and only if there exists a sequence of continuous linear maps $S_i\colon X_{i+1} \to X_i$ such that

$$T_{i-1} \circ S_{i-1} + S_i \circ T_i = \mathrm{id}_{X_i} - K_i \tag{2.1.65}$$

holds for compact operators $K_i\colon X_i \to X_i$, cf. [AB67, Proposition 6.5][1]. The sequence $(S_i)$ is called a *parametrix* of the chain complex $(X_i, T_i)$.

Similarly to the case of Fredholm maps, elliptic complexes constitute an important class of examples of Fredholm complexes. Consider the following setup. Let $E_0, E_1, \ldots, E_N$ be a sequence of finite-dimensional vector bundles over a compact manifold $M$ and let $\mathcal{E}_i$ be the Fréchet space of smooth sections of $E_i$. Let $L_i\colon \mathcal{E}_i \to \mathcal{E}_{i+1}$ be a sequence of differential operators satisfying $L_{i+1} \circ L_i = 0$. The chain complex $(\mathcal{E}_i, L_i)$ is called *elliptic* if the sequence of principal symbols

$$\cdots \longrightarrow \overset{\star}{\tau}{}^* E_i \xrightarrow{\sigma(L_i)} \overset{\star}{\tau}{}^* E_{i+1} \longrightarrow \cdots \tag{2.1.66}$$

is exact outside of the zero section of the cotangent bundle $\overset{\star}{\tau}\colon \mathrm{T}^*M \to M$. By [AB67, Proposition 6.1], every elliptic chain complex possesses a continuous parametrix $(S_i)$. The construction of $(S_i)$ in [AB67] crucially involves the parametrix of an elliptic operator. Thus, [Ham82, Theorem II.3.3.3] implies that the operators $S_i$ are tame. Hence, we obtain the following.

Proposition 2.1.24 *Every elliptic chain complex is a tame regular Fredholm chain complex.* ◇

---

[1] [AB67, Proposition 6.5] is formulated for elliptic complexes but only the parametrix is used in the proof.

As a preparation for the non-linear case, we will now study chains depending on parameters. The following example shows that a deformation of a chain complex is in general not a complex.

Example 2.1.25 We continue with the setting and notation of Example 2.1.21. Thus, let $P \to M$ be a principal $G$-bundle and let $E = P \times_G \underline{E}$ be an associated vector bundle. Every connection $A \in \mathcal{C}(P)$ induces via the covariant exterior differential $\mathrm{d}_A$ a chain

$$\cdots \longrightarrow \Omega^i(M, E) \xrightarrow{\mathrm{d}_A} \Omega^{i+1}(M, E) \longrightarrow \cdots . \tag{2.1.67}$$

Let $F_A$ be the curvature of $A$ and let $\wedge$ denote the wedge product relative to the Lie algebra action $\mathfrak{g} \times \underline{E} \to \underline{E}$. Since we have $\mathrm{d}_A^2 \alpha = F_A \wedge \alpha$ for every $\alpha \in \Omega^i(M, E)$, the chain $(\Omega^i(M, E), \mathrm{d}_A)$ is a complex if $A$ is flat (or, more generally, if the curvature of $A$ takes values in the kernel of the action $\mathfrak{g} \to \mathrm{L}(\underline{E}, \underline{E})$). Of course, not every connection $A$ in a neighborhood of a flat connection $A_0$ is flat. Accordingly, the deformation $(\Omega^i(M, E), \mathrm{d}_A)$ of the complex $(\Omega^i(M, E), \mathrm{d}_{A_0})$ is in general not a complex. ◇

Thus, let us consider the following general setup. Let $P$ be an open neighborhood of 0 in some locally convex space, let $X_i$ be a sequence of locally convex vector spaces and let $T_i \colon P \times X_i \to X_{i+1}$ be a sequence of continuous families of linear maps such that $(X_i, T_{i,0})$ is a complex. We say that $(P, X_i, T_i)$ is a *continuous family of chains*.

Definition 2.1.26 A continuous family of chains $(P, X_i, T_i)$ is called *uniformly regular* (at 0) if $(X_i, T_{i,0})$ is a regular chain complex with

$$X_i = \operatorname{Im} T_{i-1,0} \oplus \operatorname{Coim} T_{i,0} \oplus H_i \tag{2.1.68}$$

and if, for every $i \in \mathbb{Z}$ and $p \in P$, the map

$$\tilde{T}_{i,p} = \mathrm{pr}_{\operatorname{Im} T_{i,0}} \circ (T_{i,p})\restriction_{\operatorname{Coim} T_{i,0}} \colon \operatorname{Coim} T_{i,0} \to \operatorname{Im} T_{i,0} \tag{2.1.69}$$

is a topological isomorphism such that the inverses form a continuous family $P \times \operatorname{Im} T_0 \to \operatorname{Coim} T_0$.

If $X$ and $Y$ are tame Fréchet spaces, the complex $(X_i, T_{i,0})$ is tame regular and $\tilde{T}_{i,p}$ is a tame isomorphism such that the inverses form a tame smooth family, then $(P, X_i, T_i)$ is called *uniformly tame regular*. ◇

The notion of a uniformly regular chain is a direct generalization of Definition 2.1.10. For the applications we have in mind, the following equivalent characterization turns out to be more convenient.

Proposition 2.1.27 *A continuous family of chains $(P, X_i, T_i)$ is uniformly regular if and only if, for every $i \in \mathbb{Z}$, the subspace* $\operatorname{Im} T_{i-1,0}$ *of $X_i$ is closed and topologically* *complemented, say* $X_i = \operatorname{Im} T_{i-1,0} \oplus \operatorname{Coker} T_{i-1,0}$*, and the continuous family* $p \mapsto (T_{i,p})\restriction_{\operatorname{Coker} T_{i-1,0}}$ *of linear maps is uniformly regular.* ◊

*Proof.* The claim is a simple consequence of the observation that the image of $(T_{i,0})\restriction_{\operatorname{Coker} T_{i-1,0}}$ coincides with the image of $T_{i,0}$ and that

$$H_i \simeq \operatorname{Ker}(T_{i,0})\restriction_{\operatorname{Coker} T_{i-1,0}} \tag{2.1.70}$$

holds, because $T_{i,0}$ is a complex. □

Roughly speaking, a family of chains $(P, X_i, T_i)$ is uniformly regular if each family $T_i$ of linear maps is uniformly regular after factoring-out the image of the direct predecessor $T_{i-1,0}$.

Let us now turn to deformations of elliptic complexes. Let $E_0, E_1, \ldots, E_N$ be a sequence of finite-dimensional vector bundles over a compact manifold $M$ and let $\mathcal{E}_i$ be the tame Fréchet space of smooth sections of $E_i$. Moreover, let $P$ be an open neighborhood of 0 in some tame Fréchet space and let $L_i \colon P \times \mathcal{E}_i \to \mathcal{E}_{i+1}$ be a sequence of differential operators parametrized by points of $P$. We assume that, for every $i \in \mathbb{Z}$, the parameterization factors through the space of coefficients as follows

$$P \times \mathcal{E}_i \xrightarrow{\hat{L}_i \times \mathrm{id}_{\mathcal{E}_i}} \Gamma^\infty\big(\mathrm{L}(\mathrm{J}^{r_i} E_i, E_{i+1})\big) \times \mathcal{E}_i \longrightarrow \mathcal{E}_{i+1}, \tag{2.1.71}$$

where $\hat{L}_i \colon P \mapsto \Gamma^\infty\big(\mathrm{L}(\mathrm{J}^{r_i} E_i, E_{i+1})\big)$ is a tame smooth map and the second map was defined in (2.1.51). For simplicity, let us assume that the degree $r_i$ of the differential operator $L_{i,p} \colon \mathcal{E}_i \to \mathcal{E}_{i+1}$ is the same for all $p \in P$ and $i \in \mathbb{Z}$. We will refer to this setting by saying that $(P, \mathcal{E}_i, L_i)$ is a *tame family of chains of differential operators*.

As a generalization of Theorem 2.1.20, we have the following result concerning deformations of elliptic complexes.

Theorem 2.1.28 *Let* $(P, \mathcal{E}_i, L_i)$ *be a tame family of chains of differential operators. If* $(\mathcal{E}_i, L_{i,0})$ *is an elliptic complex, then* $(P, \mathcal{E}_i, L_i)$ *is uniformly tame regular (after possibly shrinking* $P$*).* ◊

*Proof.* Our proof follows closely the proof of [AB67, Proposition 6.1], where a parametrix of an elliptic complex is constructed by reducing the problem to the construction of a parametrix of an elliptic operator. Similarly, we will reduce the question of the uniform tame regularity of the chain to the uniform regularity of a family of differential operators, for which we can employ Theorem 2.1.20.

For this purpose, fix a Riemannian metric on $M$ and a fiber Riemannian metric on every vector bundle $E_i$. These data define a natural $\mathrm{L}^2$-inner product on $\mathcal{E}_i$. By partial integration, we see that the adjoints $L^*_{i,p} \colon \mathcal{E}_{i+1} \to \mathcal{E}_i$ of $L_{i,p}$ with respect these inner products yield a tame family of chains of differential operators. For every $i \in \mathbb{Z}$, define the tame family $\triangle_i \colon P \times \mathcal{E}_i \longrightarrow \mathcal{E}_i$ by

$$\triangle_{i,p} = L^*_{i,0} \circ L_{i,p} + L_{i-1,p} \circ L^*_{i-1,0}. \tag{2.1.72}$$

Clearly, $\triangle_i$ is a family of differential operator of order $2r$. Moreover, $\triangle_{i,0}$ is an elliptic operator, because $(\mathcal{E}_i, L_{i,0})$ is an elliptic complex by assumption. Thus, Theorem 2.1.20 implies that the family $\triangle_i$ is uniformly tame regular. In particular, $\triangle_{i,0}$ is regular and self-adjoint so that we get the following topological decomposition

$$\mathcal{E}_i = \operatorname{Ker} \triangle_{i,0} \oplus \operatorname{Im} \triangle_{i,0}\,. \tag{2.1.73}$$

Moreover, $\tilde{\triangle}_{i,p} = \mathrm{pr}_{\operatorname{Im}\triangle_{i,0}} \circ (\triangle_{i,p})\restriction_{\operatorname{Im}\triangle_{i,0}}$ is a tame automorphism of $\operatorname{Im}\triangle_{i,0}$ for every $p \in P$ (after possibly shrinking $P$) in such a way that the inverses form a tame smooth family. Using the decomposition (2.1.73) and a standard argument from linear algebra, we conclude that the images of $L_{i-1,0}$ and $L^*_{i,0}$ are closed and fit into the topological decomposition

$$\mathcal{E}_i = \operatorname{Im} L_{i-1,0} \oplus \operatorname{Im} L^*_{i,0} \oplus H_i\,, \tag{2.1.74}$$

where $H_i \equiv \operatorname{Ker}\triangle_{i,0} = \operatorname{Ker} L_{i,0} \cap \operatorname{Ker} L^*_{i-1,0}$. It remains to show that, for every $i \in \mathbb{Z}$ and $p \in P$, the operator

$$\tilde{L}_{i,p} = \mathrm{pr}_{\operatorname{Im} L_{i,0}} \circ (L_{i,p})\restriction_{\operatorname{Im} L^*_{i,0}} \colon \operatorname{Im} L^*_{i,0} \longrightarrow \operatorname{Im} L_{i,0} \tag{2.1.75}$$

has a tame inverse. For this purpose, consider the tame smooth family $G_i$ defined by

$$G_{i,p} = L^*_{i,0} \circ (\tilde{\triangle}^{-1}_{i+1,p})\restriction_{\operatorname{Im} L_{i,0}} \colon \operatorname{Im} L_{i,0} \longrightarrow \operatorname{Im} L^*_{i,0}\,. \tag{2.1.76}$$

Using the decomposition (2.1.74), we obtain

$$\begin{aligned}
\tilde{L}_{i,p} \circ G_{i,p} &= \mathrm{pr}_{\operatorname{Im} L_{i,0}} \circ L_{i,p} \circ L^*_{i,0} \circ (\triangle^{-1}_{i+1,p})\restriction_{\operatorname{Im} L_{i,0}} \\
&= \mathrm{pr}_{\operatorname{Im} L_{i,0}} \circ \big(L_{i,p} \circ L^*_{i,0} + L^*_{i+1,p} \circ L_{i+1,p}\big) \circ (\tilde{\triangle}^{-1}_{i+1,p})\restriction_{\operatorname{Im} L_{i,0}} \\
&= \mathrm{pr}_{\operatorname{Im} L_{i,0}} \circ \triangle^{-1}_{i+1,p} \circ (\tilde{\triangle}^{-1}_{i+1,p})\restriction_{\operatorname{Im} L_{i,0}} \\
&= \mathrm{id}_{\operatorname{Im} L_{i,0}}.
\end{aligned} \tag{2.1.77}$$

Similarly, the commutation identity

$$\tilde{\triangle}_{i,p} \circ (L^*_{i,0})\restriction_{\operatorname{Im}\triangle_{i,0}} = L^*_{i,0} \circ L_{i,p} \circ L^*_{i,0} = L^*_{i,0} \circ \tilde{\triangle}_{i+1,p} \tag{2.1.78}$$

implies

$$\begin{aligned} G_{i,p} \circ \tilde{L}_{i,p} &= L^*_{i,0} \circ \tilde{\triangle}^{-1}_{i+1,p} \circ \mathrm{pr}_{\mathrm{Im}\, L_{i,0}} \circ (L_{i,p})\restriction_{\mathrm{Im}\, L^*_{i,0}} \\ &= \tilde{\triangle}^{-1}_{i,p} \circ L^*_{i,0} \circ (L_{i,p})\restriction_{\mathrm{Im}\, L^*_{i,0}} \\ &= \tilde{\triangle}^{-1}_{i,p} \circ (\tilde{\triangle}_{i,p})\restriction_{\mathrm{Im}\, L^*_{i,0}} \\ &= \mathrm{id}_{\mathrm{Im}\, L^*_{i,0}} . \end{aligned} \tag{2.1.79}$$

Thus, $G_i$ is a tame smooth family of inverses of $\tilde{L}_i$. This completes the proof that $(P, \mathcal{E}_i, L_i)$ is uniformly tame regular. □

Example 2.1.29 Recall the setting of Example 2.1.21 and consider the chain

$$\cdots \longrightarrow \Omega^i(M, E) \xrightarrow{\mathrm{d}_A} \Omega^{i+1}(M, E) \longrightarrow \cdots \tag{2.1.80}$$

induced by the covariant derivative of the connection $A \in \mathcal{C}(P)$. As noted in Example 2.1.25, this chain is an elliptic complex if $A$ is a flat connection. Thus, Theorem 2.1.28 entails that the family of chains $\big(\mathcal{C}(P), \Omega^k(M, E), \mathrm{d}_A\big)$ is uniformly tame regular in a neighborhood of every flat connection $A_0 \in \mathcal{C}(P)$. In this context, the operator defined in (2.1.72) takes the form

$$\triangle_{A_0 A} = \mathrm{d}^*_{A_0} \mathrm{d}_A + \mathrm{d}^*_A \mathrm{d}_{A_0} : \Omega^k(M, E) \to \Omega^k(M, E) \tag{2.1.81}$$

for $A \in \mathcal{C}(P)$. Note that $\triangle_{A_0 A}$ is a natural extension of the Faddeev–Popov operator to forms of higher degree. An operator similar to $\triangle_{A_0 A}$ played a central role in [DH18, p. 405] for the study of the curvature map $F\colon \mathcal{C}(P) \to \Omega^2(M, \mathrm{Ad} P)$ near a flat connection. We will take up this example again in Section 4.4. ◊

## 2.2 Normal form of a non-linear map

In this section, we study the local behavior of a smooth map $f\colon M \to N$ between manifolds. In particular, we introduce the concept of a normal form and find suitable conditions that ensure that $f$ can be brought into such a normal form.

In the linear setting, we have seen that every (regular) operator factorizes through a linear isomorphism. Thus, one expects that every non-linear map can be represented locally by a linear isomorphism up to some higher order error term.

Definition 2.2.1 An *abstract normal form* consists of a tuple $(X, Y, \hat{f}, f_{\mathrm{sing}})$, where

(i) $X$ and $Y$ are locally convex vector spaces with topological decompositions[1] $X = \mathrm{Ker} \oplus \mathrm{Coim}$ and $Y = \mathrm{Coker} \oplus \mathrm{Im}$,

[1] In these decompositions, Ker, Coim, etc. denote abstract spaces. Below, we will identify them

(ii) $\hat{f}\colon$ Coim $\to$ Im is a linear topological isomorphism,

(iii) $f_{\text{sing}}\colon X \supseteq U \to$ Coker is a smooth map defined on an open neighborhood $U$ of 0 in $X$ such that $f_{\text{sing}}(0, x_2) = 0$ holds for all $x_2 \in U \cap$ Coim and such that the derivative $\mathrm{T}_{(0,0)} f_{\text{sing}}\colon X \to$ Coker of $f_{\text{sing}}$ at $(0, 0)$ vanishes. ◊

Given an abstract normal form $(X, Y, \hat{f}, f_{\text{sing}})$, let

$$f_{\text{NF}} = \hat{f} + f_{\text{sing}}\colon X \supseteq U \to Y. \tag{2.2.1}$$

Note that the 0-level set of $f_{\text{NF}}$ is given by

$$f_{\text{NF}}^{-1}(0) = \{(x_1, 0) \in X : f_{\text{sing}}(x_1, 0) = 0\}. \tag{2.2.2}$$

Since $\mathrm{T}_{(0,0)} f_{\text{sing}} = 0$, the level set $f_{\text{NF}}^{-1}(0)$ is in general not a smooth manifold Its singular structure is completely determined by $f_{\text{sing}}$. For this reason, we refer to $f_{\text{sing}}$ as the *singular part* of $f_{\text{NF}}$.

The predominant theme of this section is the systematic reduction of a smooth map between manifolds to an abstract normal form by a convenient choice of coordinates. This idea is formalized in the following.

Definition 2.2.2 We say that a smooth map $f\colon M \to N$ between manifolds can *be brought into the normal form* $(X, Y, \hat{f}, f_{\text{sing}})$ at the point $m \in M$ if there exist charts[1] $\kappa\colon M \supseteq U' \to U \subseteq X$ at $m$ and $\rho\colon N \supseteq V' \to V \subseteq Y$ at $f(m)$ such that $f(U') \subseteq V'$ and

$$\rho \circ f_{\restriction U'} = f_{\text{NF}} \circ \kappa \tag{2.2.3}$$

hold. ◊

Assume that the smooth map $f\colon M \to N$ can be brought into a normal form $(X, Y, \hat{f}, f_{\text{sing}})$ in a neighborhood of $m$ using diffeomorphisms $\kappa\colon U' \to U$ and $\rho\colon V' \to V$. Since $\kappa$ is a diffeomorphism, $\mathrm{T}_m\kappa\colon \mathrm{T}_m M \to X$ is a topological isomorphism. Similarly, $\mathrm{T}_{f(m)}\rho$ identifies $\mathrm{T}_{f(m)}N$ with $Y$. Under these identifications, the abstract spaces Ker and Im in the decomposition of $X$ and $Y$ coincide with the kernel and the image of $\mathrm{T}_m f$, respectively, because $\mathrm{T}_0 f_{\text{NF}} = \hat{f}$ holds and $\hat{f}\colon$ Coim $\to$ Im is an isomorphism.

In certain cases, a normal form amounts to a linearization of the map under consideration. We say that a smooth map $f\colon M \to N$ is a *submersion* at $m \in M$ if it is equivalent to a linear projection in a neighborhood of $m$. Similarly, $f$ is called an *injection* at $m$ if it is equivalent to a linear injection in a neighborhood of $m$. More generally, $f$ is a *subimmersion* at $m$ if it is equivalent to a linear map in a neighborhood of $m$.

---

with the kernel, coimage, etc. of the tangent map $\mathrm{T}_m f$, respectively.

[1] Throughout this work, we follow the convention that a chart $\kappa\colon M \supseteq U' \to U \subseteq X$ at a point $m \in M$ satisfies $\kappa(m) = 0$. Moreover, in this notation, $U'$ and $U$ are understood to be open neighborhoods of $m$ in $M$ and of 0 in $X$, respectively.

Proposition 2.2.3 *Let $f\colon M \longrightarrow N$ be a smooth map. Assume that $f$ can be brought into a normal form in a neighborhood $U'$ of $m \in M$. Then, the following holds:*

*(i) (Submersion) $f$ is a submersion at $m$ if and only if $\mathrm{T}_m f$ is surjective.*

*(ii) (Immersion) $f$ is an immersion at $m$ if and only if $\mathrm{T}_m f$ is injective.*

*(iii) (Constant rank) $f$ is a subimmersion at $m$ if, for all $p \in U'$, $\mathrm{T}_p f$ is a finite-rank operator*[1] *satisfying* $\operatorname{rk} \mathrm{T}_p f = \operatorname{rk} \mathrm{T}_m f$. ◊

*Proof.* Let $(X, Y, \hat{f}, f_{\mathrm{sing}})$ be a normal form of $f$ at $m$. Since the claims are of local nature, it is sufficient to consider the normal form $f_{\mathrm{NF}} = \hat{f} + f_{\mathrm{sing}}\colon X \supseteq U \longrightarrow Y$ of $f$. For simplicity, we continue writing $f$ for $f_{\mathrm{NF}}$. If $\mathrm{T}_0 f$ is surjective, then Coker is trivial and hence $f_{\mathrm{sing}} = 0$. Similarly, if $\mathrm{T}_0 f$ is injective, then $\mathrm{Coim} = X$ and thus $f_{\mathrm{sing}}(0, x_2) = 0$ for all $x_2 \in U \cap \mathrm{Coim}$ implies $f_{\mathrm{sing}} = 0$. The converse direction is clear.

Suppose now that $\mathrm{T}_x f$ is a finite-rank operator and that $\operatorname{rk} \mathrm{T}_x f = \operatorname{rk} \mathrm{T}_0 f$ holds for all $x \in U$. We have

$$\mathrm{T}_x f(v_1, v_2) = \big(\mathrm{T}_x f_{\mathrm{sing}}(v_1, v_2), \hat{f}(v_2)\big) \in \mathrm{Coker} \oplus \mathrm{Im} = Y \tag{2.2.4}$$

for all $x \in U$, $v_1 \in \mathrm{Ker}$ and $v_2 \in \mathrm{Coim}$. Since $\operatorname{rk} \mathrm{T}_x f = \operatorname{rk} \mathrm{T}_0 f = \dim \mathrm{Im}$, we conclude that $\mathrm{T}_{(x_1,x_2)} f_{\mathrm{sing}}(v_1, 0) = 0$ holds for all $(x_1, x_2) \in U$ and $v_1 \in \mathrm{Ker}$. Hence, $f_{\mathrm{sing}}(x_1, x_2)$ does not depend on $x_1$ and so

$$f_{\mathrm{sing}}(x_1, x_2) = f_{\mathrm{sing}}(0, x_2) = 0. \tag{2.2.5}$$

In other words, $f$ is locally represented by the linear map $\hat{f}$ and thus it is a subimmersion. □

Remark 2.2.4 It is straightforward to verify that the pull-back of a submanifold along a submersion is a submanifold, see [Glö15, Theorem C]. In particular, the level set $f^{-1}(\mu)$ of a smooth map $f\colon M \longrightarrow N$ corresponding to the value $\mu \in N$ is a submanifold of $M$ if $f$ is a submersion at every $m \in f^{-1}(\mu)$. Moreover, if $f\colon M \longrightarrow N$ is an immersion and a topological embedding, then $f(M)$ is a submanifold of $N$ and $f$ yields a diffeomorphism of $M$ onto $f(M)$, see [Glö15, Lemma I.13]. ◊

The upshot of Proposition 2.2.3 is that, for any normal form theorem, one obtains a corresponding version of the submersion, regular value, immersion and constant rank theorem. *A normal form theorem thus unifies these fundamental theorems under one umbrella.*

[1] The statement generalizes to maps whose derivatives have a constant but not necessarily finite-dimensional image, cf. [MRA02, Theorem 2.5.15] for a Banach version.

REMARK 2.2.5 (Finite-dimensional reduction and Sard–Smale theorem) Let $f: M \to N$ be a smooth map which can be brought into the normal form $(X, Y, \hat{f}, f_{\mathrm{sing}})$ at $m \in M$. Denote $\mu = f(m)$. In terms of the normal form, the level set $f^{-1}(\mu)$ in a neighborhood of $m$ is given by (cf. (2.2.2))

$$\{(x_1, 0) \in X : \ x_1 \in \mathrm{Ker}, f_{\mathrm{sing}}(x_1, 0) = 0\}. \tag{2.2.6}$$

Hence, all singularities of $f^{-1}(\mu)$ that might arise are encoded in the map

$$f_{\mathrm{sing}}(\cdot, 0) \colon \ \mathrm{Ker} \to \mathrm{Coker}\,. \tag{2.2.7}$$

If $\mathrm{T}_m f$ is a Fredholm operator, then Ker and Coker are finite-dimensional spaces. Thus, in this case, the study of the singular structure of $f^{-1}(\mu)$ is reduced to a question in finite dimensions. Exploiting this reduction to finite dimensions, we may mimic the usual proof of the Sard–Smale theorem [Sma65] to obtain a version of this theorem for Fredholm maps between *locally convex* manifolds. We leave the details to the reader. ◊

The aim of the remainder of the section is to find suitable conditions on the derivative of $f: M \to N$ which ensure that $f$ can be brought into a normal form.

### 2.2.1 Banach version

To illustrate the main idea and to somewhat reduce the functional analytic complexity, we first restrict attention to the Banach setting.

THEOREM 2.2.6 (Normal form — Banach version) *Let $f: M \to N$ be a smooth map between Banach manifolds and let $m \in M$. If $\mathrm{T}_m f : \mathrm{T}_m M \to \mathrm{T}_{f(m)} N$ is a regular operator, then $f$ can be brought into a normal form around $m$. In particular, every smooth map between finite-dimensional manifolds can be brought into a normal form around every point.* ◊

*Proof.* Since the claim is of local nature, we can use charts $\tilde{\kappa}: M \supseteq U' \to U \subseteq X$ at $m$ and $\tilde{\rho}: N \supseteq V' \to V \subseteq Y$ at $f(m)$ to replace $f$ by its local representative $f: X \supseteq U \to Y$. To economize on notation, we abbreviate $T \equiv \mathrm{T}_0 f: X \to Y$. Since $T$ is regular by assumption, there exist topological decompositions $X = \mathrm{Ker}\, T \oplus \mathrm{Coim}\, T$ and $Y = \mathrm{Coker}\, T \oplus \mathrm{Im}\, T$. Moreover, the core $\hat{T}: \ \mathrm{Coim}\, T \to \mathrm{Im}\, T$ of $T$ is a topological isomorphism.

Define the smooth map $\psi: X \supseteq U \to X$ by

$$\psi(x_1, x_2) = \big(x_1, \hat{T}^{-1} \circ \mathrm{pr}_{\mathrm{Im}\, T} \circ f(x_1, x_2)\big), \tag{2.2.8}$$

with $x_1 \in \mathrm{Ker}\, T$ and $x_2 \in \mathrm{Coim}\, T$. Note that $\psi(0) = 0$. Since $\mathrm{T}_0 \psi = \mathrm{id}_X$, it follows from the Inverse Function Theorem A.1.1 that we can shrink $U$ in such a way that $\psi(U)$ is an open neighborhood of 0 in $X$ and $\psi: U \to \psi(U)$ is a diffeomorphism.

By possibly shrinking $V$, we may assume that $V \cap \operatorname{Im} T \subseteq \hat{T} \circ \psi(U \cap \operatorname{Coim} T)$. Define the smooth map $\phi\colon Y \supseteq V \to Y$ by

$$\phi(y_1, y_2) = \big(y_1 + \mathrm{pr}_{\operatorname{Coker} T} \circ f \circ \psi^{-1}(0, \hat{T}^{-1} y_2), y_2\big) \tag{2.2.9}$$

with $y_1 \in \operatorname{Coker} T$ and $y_2 \in \operatorname{Im} T$. A moment's reflection convinces us that $\phi(0) = 0$ and $\mathrm{T}_0\phi = \mathrm{id}_Y$ hold. Thus, the Inverse Function Theorem A.1.1 implies that we can shrink $V$ such that $\phi(V)$ is an open neighborhood of 0 in $Y$ and $\phi\colon V \to \phi(V)$ is a diffeomorphism. By possibly shrinking $U$, we may assume $f(U) \subseteq V$ and $f(U) \subseteq \phi(V)$.

Set $\hat{f} = \hat{T}\colon \operatorname{Coim} T \to \operatorname{Im} T$ and define the smooth map $f_{\mathrm{sing}}\colon X \supseteq \psi(U) \to$ Coker by

$$f_{\mathrm{sing}} \circ \psi(x_1, x_2) = \mathrm{pr}_{\operatorname{Coker} T} \left(f(x_1, x_2) - f \circ \psi^{-1}\big(0, \hat{T}^{-1} \circ \mathrm{pr}_{\operatorname{Im} T} \circ f(x_1, x_2)\big)\right). \tag{2.2.10}$$

A straightforward computation shows that the following diagram commutes:

$$\begin{array}{ccc} X \supseteq U & \xrightarrow{\ f\ } & \phi(V) \subseteq Y \\ {\scriptstyle\psi}\big\downarrow & & \big\uparrow{\scriptstyle\phi} \\ X \supseteq \psi(U) & \xrightarrow[\hat{f} + f_{\mathrm{sing}}]{} & V \subseteq Y. \end{array} \tag{2.2.11}$$

Thus, in the charts $\kappa = \psi \circ \tilde{\kappa}$ and $\rho = \phi^{-1} \circ \tilde{\rho}$, the map $f$ coincides with $f_{\mathrm{NF}} = \hat{f} + f_{\mathrm{sing}}$. Finally, let us verify the asserted properties of the singular part $f_{\mathrm{sing}}$. We clearly have $\mathrm{T}_0 f_{\mathrm{sing}} = 0$. Moreover, for all $x_2 \in U \cap \operatorname{Coim} T$, we get

$$\begin{aligned} f_{\mathrm{sing}} \circ \psi(0, x_2) &= \mathrm{pr}_{\operatorname{Coker} T} \big(f(0, x_2) - f \circ \psi^{-1}(0, \hat{T}^{-1} \circ \mathrm{pr}_{\operatorname{Im} T} \circ f(0, x_2))\big) \\ &= \mathrm{pr}_{\operatorname{Coker} T} \big(f(0, x_2) - f \circ \psi^{-1} \circ \psi(0, x_2)\big) \\ &= 0. \end{aligned} \tag{2.2.12}$$

From $\psi(0, x_2) \in \operatorname{Coim} T$ it follows that $f(0, x_2') = 0$ holds for all $x_2' \in \psi(U) \cap \operatorname{Coim} T$. □

Remarks 2.2.7

(i) A weaker version of this normal form theorem can be found in [MD92, Theorem 5.1.8; MRA02, Theorem 2.5.14]. There, the chart on $N$ is not modified and hence the additional property $f_{\mathrm{sing}}(0, \cdot) = 0$ of the singular part is not deduced. Note that this property was crucial in the proof of Proposition 2.2.3 (ii) to show that a smooth map with injective differential is an immersion.

(ii) The first part of the proof of Theorem 2.2.6 is inspired by the Lyapunov–

Schmidt reduction procedure. To establish the link, let us give a brief outline of this procedure, mostly ignoring the peculiarities of the infinite-dimensional setting, see e.g. [Cha05, Section 1.3] for a textbook treatment. Given Banach spaces $X$ and $Y$, and a smooth map $f\colon X \supseteq U \to Y$ defined on an open neighborhood $U$ of 0 in $X$ with $f(0) = 0$, we are interested in solutions of the non-linear equation

$$f(x) = 0 \tag{2.2.13}$$

near the solution $x = 0$. The Lyapunov–Schmidt scheme consists of the following steps:

(i) Split $X$ and $Y$ into direct sums $X = \operatorname{Ker} T \oplus \operatorname{Coim} T$ and $Y = \operatorname{Coker} T \oplus \operatorname{Im} T$, where $T = \mathrm{T}_0 f\colon X \to Y$ as above. The equation (2.2.13) is then equivalent to the system

$$\begin{aligned} \mathrm{pr}_{\operatorname{Coker} T} \circ f(x_1, x_2) &= 0, \\ \mathrm{pr}_{\operatorname{Im} T} \circ f(x_1, x_2) &= 0, \end{aligned} \tag{2.2.14}$$

with $x_1 \in \operatorname{Ker} T$ and $x_2 \in \operatorname{Coim} T$.

(ii) The Implicit Function Theorem shows that, after possibly shrinking $U$, the second equation in (2.2.14) has for each $x_1 \in U \cap \operatorname{Ker} T$ a unique solution $x_2(x_1) \in U \cap \operatorname{Coim} T$.

(iii) Substituting this solution of the second equation in the first equation of (2.2.14) yields the reduced equation

$$\mathrm{pr}_{\operatorname{Coker} T} \circ f(x_1, x_2(x_1)) = 0. \tag{2.2.15}$$

for the unknown $x_1 \in U \cap \operatorname{Ker} T$.

In this way, the non-linear equation (2.2.13) is reduced to the non-linear equation (2.2.15), which often happens to be a set of finitely many equations for a finite number of unknowns. When comparing this reduction scheme with our construction of the chart deformation $\psi$ in the proof of Theorem 2.2.6, the only conceptual difference is our usage of the Inverse Function Theorem in place of the Implicit Function Theorem employed in the Lyapunov–Schmidt procedure. Both methods rely fundamentally on the fact that the map $\mathrm{pr}_{\operatorname{Im} T} \circ f\colon X \supseteq U \to \operatorname{Im} T$ has a surjective derivative at 0. In our language, the reduced equation (2.2.15) takes the form

$$f_{\mathrm{sing}}(x_1, 0) = 0 \tag{2.2.16}$$

for $x_1 \in U \cap \operatorname{Ker} T$.

Similar ideas are also used in the study of deformations of geometric objects, see e.g. [Kur65] concerning deformations of complex structures and [Don83, Section II.2; Tau82, Section 6] in the gauge theoretic setting; see also [DK97, Section 4.2.5]. In this context, the counterpart of the local diffeomorphism $\psi$ is usually referred to as the Kuranishi map.

(iii) In the setting of Theorem 2.2.6, additionally assume that $f$ is equivariant with respect to actions of a compact Lie group $G$ on $M$ and $N$. If there exist charts $\tilde{\kappa}\colon M \supseteq U' \to U \subseteq X$ at $m$ and $\tilde{\rho}\colon N \supseteq V' \to V \subseteq Y$ at $f(m)$ which are $G$-equivariant with respect to *linear* actions of $G$ on $X$ and $Y$, respectively, then the chart deformations $\psi$ and $\phi$ introduced in the proof of Theorem 2.2.6 are $G$-equivariant. To see this, consider the topological decompositions $X = \operatorname{Ker} T \oplus \operatorname{Coim} T$ and $Y = \operatorname{Coker} T \oplus \operatorname{Im} T$, where $T = \mathrm{T}_0 f$, as above. Since $G$ is compact and $T$ is $G$-equivariant, the complements $\operatorname{Coim} T$ and $\operatorname{Coker} T$ can be chosen to be $G$-invariant according to Lemma A.2.5. Therefore, the map $\psi$ defined in (2.2.8) is a composition of $G$-equivariant maps and hence is $G$-equivariant. The map $\phi$ defined in (2.2.9) and the singular part $f_{\text{sing}}$ are $G$-equivariant, too, for similar reasons. We will refer to this situation by saying that $f$ can be *brought into a normal form in a G-equivariant way.* ◊

By Proposition 2.2.3 (ii), we get a corresponding normal form theorem for immersions. Note that in our construction of the normal form the chart on the domain is always deformed, which is in contrast to the classical immersion theorem (e.g., [MRA02, Theorem 2.5.12]). Hence, the constructed submanifold charts differ from the usual ones. This difference is illustrated by the following simple example.

Example 2.2.8 Consider the embedding $f\colon \mathbb{R} \to \mathbb{R}^2$, $\vartheta \mapsto (\cos\vartheta, \sin\vartheta)$ of the circle in $\mathbb{R}^2$. According to the construction in the proof of Theorem 2.2.6, the normal form of $f$ near $\vartheta = 0$ is given as follows. We have $X = \mathbb{R} = \operatorname{Coim}$ and

$$Y = \mathbb{R}^2 = \underbrace{\mathbb{R}\begin{pmatrix}1\\0\end{pmatrix}}_{\text{Coker}} \oplus \underbrace{\mathbb{R}\begin{pmatrix}0\\1\end{pmatrix}}_{\text{Im}}. \tag{2.2.17}$$

Moreover, the chart deformations defined in (2.2.8) and (2.2.9) here take the form

$$\psi(\vartheta) = \sin\vartheta \tag{2.2.18}$$

and

$$\phi(y_1, y_2) = \left(y_1 + \sqrt{1 - y_2^2}, y_2\right). \tag{2.2.19}$$

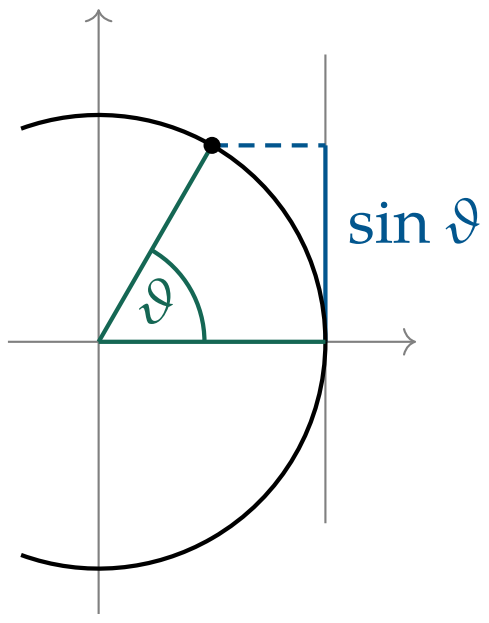

Figure 2.1: Illustration of the submanifold charts $\phi$ (blue) and $\Phi$ (green).

A direct inspection shows that the following diagram commutes:

$$\begin{array}{ccc} \vartheta & \overset{f}{\longmapsto} & (\cos\vartheta, \sin\vartheta) \\ {\scriptstyle\psi}\downarrow & & \uparrow{\scriptstyle\phi} \\ \sin\vartheta & \overset{\hat{f}}{\longmapsto} & (0, \sin\vartheta), \end{array} \tag{2.2.20}$$

where $\hat{f}$ is the inclusion of $\mathbb{R}$ in $\mathbb{R}^2$ as the second component. In contrast, the chart on $\mathbb{R}^2$ constructed in the usual immersion theorem [MRA02, Theorem 2.5.12] is given by

$$\Phi(y_1, y_2) = f(y_2) + (y_1, 0) = (y_1 + \cos y_2, \sin y_2) \tag{2.2.21}$$

and fits into the following commutative diagram:

$$\begin{array}{ccc} \vartheta & \overset{f}{\longmapsto} & (\cos\vartheta, \sin\vartheta) \\ {\scriptstyle\mathrm{id}}\downarrow & & \uparrow{\scriptstyle\Phi} \\ \vartheta & \overset{\hat{f}}{\longmapsto} & (0, \vartheta). \end{array} \tag{2.2.22}$$

Hence, the normal form $\hat{f}$ coincides in both constructions but the submanifold charts constructed are different. In our approach, the point on the circle is characterized by its second component, while in the usual construction the point is determined by the angle of rotation (see Figure 2.1). ◊

In the following, we give generalizations of the Banach normal form theorem to different analytic settings. Beyond the Banach context, the classical Banach Inverse Function Theorem used in the proof of Theorem 2.2.6 has to be replaced by a different version. We will use Glöckner's Inverse Function Theorem for maps between Banach spaces with parameters in a locally convex space and the Nash–Moser Theorem in the tame Fréchet setting. Extrapolating from these cases, it will become clear that analogous normal form theorems can be obtained based on other Inverse Function Theorems. The reader may consult

Appendix A.1 for an overview of different versions of the Inverse Function Theorem and related background material.

### 2.2.2 Banach target or domain

Let us start with the following generalization of Theorem 2.2.6 to more general domains.

Theorem 2.2.9 (Normal form — Banach target) *Let $f: M \to N$ be a smooth map between manifolds, where $N$ is a Banach manifold. If, for some $m \in M$, the differential $\mathrm{T}_m f: \mathrm{T}_m M \to \mathrm{T}_{f(m)} N$ is a regular operator, then $f$ can be brought into a normal form at $m$. In particular, every smooth map $f: M \to N$ with finite-dimensional target $N$ can be brought into a normal form at every point.* ◊

*Proof.* The proof of Theorem 2.2.6 carries over word by word except for the part where the Inverse Function Theorem was used to show that the map $\psi$ defined in (2.2.8) is a local diffeomorphism. The idea is to replace the classical Inverse Function Theorem A.1.1 by the more advanced Theorem A.1.2 due to [Glö06]. For this purpose, we continue to work in the setting and notation of the proof of Theorem 2.2.6. Define the smooth map $\bar{\psi}: X \supseteq U \to \operatorname{Coim} T$ by

$$\bar{\psi}(x_1, x_2) = \hat{T}^{-1} \circ \mathrm{pr}_{\operatorname{Im} T} \circ f(x_1, x_2). \tag{2.2.23}$$

The partial derivative of $\bar{\psi}$ at 0 with respect to the second component is given by $\mathrm{T}^2_0 \bar{\psi} = \mathrm{id}_{\operatorname{Coim} T}$. Since $\hat{T}: \operatorname{Coim} T \to \operatorname{Im} T$ is an isomorphism, $\operatorname{Coim} T$ is a Banach space. Hence, considering $x_1 \in \operatorname{Ker}$ as a parameter, the Inverse Function Theorem A.1.2 shows that the map

$$\psi(x_1, x_2) = \big(x_1, \bar{\psi}(x_1, x_2)\big) \tag{2.2.24}$$

is a local diffeomorphism. The rest of the proof of Theorem 2.2.6 goes through without modification and we thus conclude that $f$ can be brought into a normal form.

For the second part of the claim, recall from Lemma 2.1.2 (ii) that every linear continuous map with finite-dimensional target is a regular operator. □

In combination with Proposition 2.2.3, we recover the submersion theorem and the constant rank theorem [Glö15, Theorem A and F] for maps with values in a finite-dimensional manifold as a special case of Theorem 2.2.9.

We also have the following version of the normal form theorem where not the target but the domain is a Banach manifold.

Theorem 2.2.10 (Normal form — Banach domain) *Let $f: M \to N$ be a smooth map between manifolds. Assume that $M$ is a Banach manifold. If, for some $m \in M$, the differential $\mathrm{T}_m f: \mathrm{T}_m M \to \mathrm{T}_{f(m)} N$ is a regular operator, then $f$ can be brought* *into a normal form around $m$. In particular, every smooth map $f: M \longrightarrow N$ with $M$ being finite-dimensional can be brought into a normal form around every point.* ◇

*Proof.* The proof proceeds in complete analogy to the one of Theorem 2.2.9 with the modification that Glöckner's Inverse Function Theorem A.1.2 has to be applied on the domain and not on the target. Details are left to the reader. □

In conjunction with Proposition 2.2.3 (ii), we recover the immersion theorem [Glö15, Theorem H] for maps from a Banach manifold into a locally convex manifold.

### 2.2.3 Nash–Moser version

We now establish a normal form theorem in the tame Fréchet category. For the convenience of the reader, the main notions of the tame category are briefly summarized in Appendix A.1. One of the major complications of the Nash–Moser Inverse Function Theorem is that the derivative has to be invertible at a point and in its neighborhood. This additional condition is due to the fact that the subset of invertible operators is no longer open in the space of all linear continuous operators. In Section 2.1.2, we have introduced the notion of uniform regularity to address a similar problem related to the lack of openness. The following results hints at uniform regularity playing a major role in the context of normal forms (in fact, it was the main inspiration for this notion).

Proposition 2.2.11 *For every normal form $(X, Y, \hat{f}, f_{\mathrm{sing}})$, the family $\mathrm{T}_x f_{\mathrm{NF}}: X \longrightarrow Y$ of linear maps parametrized by $x \in U$ is uniformly regular at $0$.* ◇

*Proof.* First note that $\mathrm{T}_0 f_{\mathrm{NF}} = \hat{f}$, because $\mathrm{T}_0 f_{\mathrm{sing}} = 0$. Thus, $\mathrm{T}_0 f_{\mathrm{NF}}$ is a regular operator. With respect to the decompositions $X = \mathrm{Ker} \oplus \mathrm{Coim}$ and $Y = \mathrm{Coker} \oplus \mathrm{Im}$, the derivative of $f_{\mathrm{NF}}$ is given in block form as

$$\mathrm{T}_x f_{\mathrm{NF}} = \begin{pmatrix} (\mathrm{T}_x f_{\mathrm{sing}})\restriction_{\mathrm{Ker}} & (\mathrm{T}_x f_{\mathrm{sing}})\restriction_{\mathrm{Coim}} \\ 0 & \hat{f} \end{pmatrix}. \tag{2.2.25}$$

Thus, the map $\mathrm{pr}_{\mathrm{Im}} \circ (\mathrm{T}_x f_{\mathrm{NF}})\restriction_{\mathrm{Coim}}$ coincides with the isomorphism $\hat{f}$, which confirms that $\mathrm{T}_x f_{\mathrm{NF}}$ is uniformly regular at $0$. □

In order to generalize the notion of uniform regularity to the setting of manifolds, consider a morphism $T: E \longrightarrow F$ of vector bundles over a manifold $M$. We will denote the induced operator on the fibers by $T_m: E_m \longrightarrow F_m$ for every $m \in M$. Suppose that $E$ and $F$ are trivialized over an open subset $U \subseteq M$, that is, we are given vector bundle isomorphism $E\restriction_U \longrightarrow U \times \underline{E}$ and $F\restriction_U \longrightarrow U \times \underline{F}$, where $\underline{E}$ and $\underline{F}$ are locally convex spaces. With respect to these trivializations, we view $T_m$ for $m \in U$ as an operator $\underline{E} \longrightarrow \underline{F}$ and accordingly identify $T$ with a family $T\restriction_U: U \times \underline{E} \longrightarrow \underline{F}$ of linear maps.

Definition 2.2.12 A morphism $T\colon E \to F$ of vector bundles over $M$ is called *uniformly regular* at $m \in M$ if there exist local trivializations of $E$ and $F$ on an open neighborhood $U$ of $m$ such that the induced family $T_{\restriction U}\colon U \times \underline{E} \to \underline{F}$ of linear maps is uniformly regular at 0 in the sense of Definition 2.1.10. Similarly, a tame morphism $T\colon E \to F$ between tame Fréchet vector bundles is called *uniformly tame regular* at $m$, if there exist tame local trivializations such that the induced family $T_{\restriction U}\colon U \times \underline{E} \to \underline{F}$ is uniformly tame regular at 0. ◊

Phrased in this language, Proposition 2.2.11 entails that the derivative of $f\colon M \to N$, viewed as a vector bundle map $\mathrm{T}f\colon \mathrm{T}M \to f^*\mathrm{T}N$ over $M$, is uniformly regular at $m \in M$ if $f$ can be brought into a normal form at $m$. The following theorem shows that, in the tame category, uniform regularity of the derivative is also a sufficient condition for the existence of a normal form.

Theorem 2.2.13 (Normal form — Tame Fréchet) *Let $f\colon M \to N$ be a tame smooth map between tame Fréchet manifolds and let $m \in M$. If $\mathrm{T}f\colon \mathrm{T}M \to f^*\mathrm{T}N$ is uniformly tame regular at $m$, then $f$ can be brought into a normal form at $m$.* ◊

*Proof.* The proof follows the same line of arguments as in the proof of Theorem 2.2.6 except that we will use the Nash–Moser Theorem A.1.3 to show that the chart deformations (2.2.8) and (2.2.9) are local diffeomorphisms. Continuing in the notation of the proof of Theorem 2.2.6, abbreviate $T_x \equiv \mathrm{T}_x f\colon X \to Y$ for every $x \in U$. The assumption of uniform tame regularity of $\mathrm{T}f$ implies that the family $T\colon U \times X \to Y$ is uniform tame regular at 0.

The derivative at $x \in U$ of the map $\psi$ defined in (2.2.8) evaluates to

$$\mathrm{T}_x\psi = \begin{pmatrix} \mathrm{id}_{\operatorname{Ker} T_0} & 0 \\ \hat{T}_0^{-1} \circ \mathrm{pr}_{\operatorname{Im} T_0} \circ (T_x)_{\restriction \operatorname{Ker} T_0} & \hat{T}_0^{-1} \circ \tilde{T}_x \end{pmatrix}, \tag{2.2.26}$$

where $\tilde{T}_x = \mathrm{pr}_{\operatorname{Im} T_0} \circ (T_x)_{\restriction \operatorname{Coim} T_0}\colon \operatorname{Coim} T_0 \to \operatorname{Im} T_0$. Since $T$ is uniformly tame regular, $\tilde{T}_x$ is invertible for all $x \in U$ and the inverses form a tame smooth family. Hence, $\mathrm{T}_x\psi$ has a tame smooth family of inverses given by

$$\left(\mathrm{T}_x\psi\right)^{-1} = \begin{pmatrix} \mathrm{id}_{\operatorname{Ker} T_0} & 0 \\ -\tilde{T}_x^{-1} \circ \mathrm{pr}_{\operatorname{Im} T_0} \circ (T_x)_{\restriction \operatorname{Ker} T_0} & \tilde{T}_x^{-1} \circ \hat{T}_0 \end{pmatrix}. \tag{2.2.27}$$

Thus, we can apply the Nash–Moser Theorem A.1.3 to conclude that $\psi$ is a local diffeomorphism at 0. For the map $\phi$ defined in (2.2.9) things are somewhat simpler. Indeed, using (2.2.27), the derivative of $\phi$ can be written in block form as

$$\mathrm{T}_{y_1,y_2}\phi = \begin{pmatrix} \mathrm{id}_{\operatorname{Coker} T_0} & \mathrm{pr}_{\operatorname{Coker} T_0} \circ T_x \circ \tilde{T}_x^{-1} \\ 0 & \mathrm{id}_{\operatorname{Im} T_0} \end{pmatrix}, \tag{2.2.28}$$

where $x = \psi^{-1}\big(0, \hat{T}^{-1}(y_2)\big)$. A direct calculation verifies that $\mathrm{T}_{y_1,y_2}\phi$ is invertible

with inverse given by

$$\left(\mathrm{T}_{y_1,y_2}\phi\right)^{-1} = \begin{pmatrix} \mathrm{id}_{\operatorname{Coker} T_0} & -\mathrm{pr}_{\operatorname{Coker} T_0} \circ T_x \circ \tilde{T}_x^{-1} \\ 0 & \mathrm{id}_{\operatorname{Im} T_0} \end{pmatrix}. \tag{2.2.29}$$

Hence, the inverses $\left(\mathrm{T}_{y_1,y_2}\phi\right)^{-1}$ parametrized by $(y_1, y_2) \in V$ form a tame smooth family. Therefore, the Nash–Moser Theorem A.1.3 implies that $\phi$ is a local diffeomorphism. The remainder of the proof of Theorem 2.2.6 goes through without modification. □

### 2.2.4 Elliptic version

In problems coming from geometry and physics, one is usually interested in differential operators between spaces of geometric objects. Let $E \to M$ and $F \to M$ be finite-dimensional fiber bundles over the compact manifold $M$. Denote the space of smooth sections of $E$ and $F$ by $\mathcal{E}$ and $\mathcal{F}$, respectively. By [Ham82, Theorem II.2.3.1], the spaces $\mathcal{E}$ and $\mathcal{F}$ are tame Fréchet manifolds. Following [Pal68, Definition 15.3], we call a map $L\colon \mathcal{E} \to \mathcal{F}$ a *non-linear differential operator* of degree $r$ if it can be factored as the composition

$$L\colon \mathcal{E} \xrightarrow{\mathrm{j}^r} \Gamma^\infty(\mathrm{J}^r E) \xrightarrow{f_*} \mathcal{F}, \tag{2.2.30}$$

where $\mathrm{j}^r$ denotes the $r$-th jet prolongation and $f\colon \mathrm{J}^r E \to F$ is a vertical morphism of fiber bundles. As for linear differential operators, we write $L_f$ for the differential operator induced by the bundle morphism $f\colon \mathrm{J}^r E \to F$. Every non-linear differential operator $L\colon \mathcal{E} \to \mathcal{F}$ is a tame smooth map according to [Ham82, Corollary II.2.2.7].

Let us describe the non-linear differential operator $L_f$ in local coordinates. For this purpose, choose local coordinates $(x_1, \ldots, x_n)$ on an open subset $U$ of $M$ and $(y_1, \ldots, y_m)$ along the fibers of $E$ over $U$. We then have natural coordinates $y_j^\alpha$ along the fibers of $\mathrm{J}^r E$ over $U$, with $\alpha$ being a multi-index, defined by

$$y_j^\alpha(\mathrm{j}^r \phi) = \frac{\partial^{|\alpha|}\phi_j}{\partial^\alpha x_i} \equiv \partial^\alpha \phi_j, \tag{2.2.31}$$

where $\phi_j$ denotes the local section $\phi$ expressed in the coordinates $(x_i, y_j)$. Moreover, choose fiber coordinates $(v_1, \ldots, v_q)$ in $F$ over $U$. With respect to these coordinates, the bundle morphism $f\colon \mathrm{J}^r E \to F$ is written as

$$v_k \circ f = f_k(x_i, y_j, y_j^\alpha), \tag{2.2.32}$$

where $1 \leq k \leq q$. Then, the local expression of the non-linear differential

operator $L_f$ is

$$\big(L_f(\phi)\big)_k(x_i) = f_k\big(x_i, \phi_j(x_i), \partial^\alpha \phi_j(x_i)\big). \tag{2.2.33}$$

As one would expect, the linearization of a non-linear differential operator is a linear differential operator. To see this, recall that $\mathrm{T}_\phi \mathcal{E} \simeq \Gamma^\infty(\phi^* \mathrm{V}E)$ for all $\phi \in \mathcal{E}$, where $\mathrm{V}E$ denotes the subbundle of $\mathrm{T}E$ consisting of vectors tangent to the fibers. By [Pal68, Theorem 17.1], there exist a natural isomorphism of $(\mathrm{j}^r\phi)^*\mathrm{V}(\mathrm{J}^r E)$ and $\mathrm{J}^r(\phi^*\mathrm{V}E)$ such that the following diagram commutes

$$\begin{array}{ccc} \Gamma^\infty(\phi^*\mathrm{V}E) & \xrightarrow{\mathrm{T}_\phi \mathrm{j}^r} & \Gamma^\infty((\mathrm{j}^r\phi)^*\mathrm{V}(\mathrm{J}^r E)) \\ \big\downarrow \mathrm{id} & & \big\downarrow \simeq \\ \Gamma^\infty(\phi^*\mathrm{V}E) & \xrightarrow{\mathrm{j}^r} & \Gamma^\infty(\mathrm{J}^r(\phi^*\mathrm{V}E)). \end{array} \tag{2.2.34}$$

The chain rule implies that the linearization $\mathrm{T}_\phi L_f$ of the non-linear differential operator $L_f$ thus factors as the composition

$$\begin{array}{rccc} \mathrm{T}_\phi L_f \colon \Gamma^\infty(\phi^*\mathrm{V}E) & \xrightarrow{\mathrm{j}^r} & \Gamma^\infty(\mathrm{J}^r(\phi^*\mathrm{V}E)) & \\ & & \big\downarrow \simeq & \\ & & \Gamma^\infty((\mathrm{j}^r\phi)^*\mathrm{V}(\mathrm{J}^r E)) & \xrightarrow{(\mathrm{V}f)_*} \Gamma^\infty((L_f(\phi))^*\mathrm{V}F), \end{array} \tag{2.2.35}$$

where $\mathrm{V}f \colon \mathrm{V}(\mathrm{J}^r E) \to f^*\mathrm{V}F$ denotes the vertical derivative of $f$. Hence, $\mathrm{T}_\phi L_f$ is a linear differential operator with coefficients $\mathrm{V}f$. In terms of the local coordinates introduced above, we obtain using the chain rule the following expression for the linearized operator:

$$\begin{aligned} \big(\mathrm{T}_\phi L_f(\sigma)\big)_k(x_i) &= \sum_j \frac{\partial f_k}{\partial y_j}\big(x_i, \phi_j(x_i), \partial^\alpha \phi_j(x_i)\big) \cdot \sigma_j(x_i) \\ &\quad + \sum_{j,\alpha} \frac{\partial f_k}{\partial y_j^\alpha}\big(x_i, \phi_j(x_i), \partial^\alpha \phi_j(x_i)\big) \cdot \frac{\partial^{|\alpha|}\sigma_j}{\partial^\alpha x_j}(x_i) \end{aligned} \tag{2.2.36}$$

for $\sigma \in \Gamma^\infty(\phi^*\mathrm{V}E)$, see [Pal68, Theorem 17.6].

We say that $L_f$ is *elliptic* if the linear differential operator $\mathrm{T}_\phi L_f$ is elliptic for all $\phi \in \mathcal{E}$. In Theorem 2.1.20, we have seen that a family of elliptic operators is uniformly tame regular. Hence, ellipticity provides an important class of examples for which the Normal Form Theorem 2.2.13 holds.

Theorem 2.2.14 (Normal form — Elliptic version) *Let $E \to M$ and $F \to M$ be finite-dimensional fiber bundles over the compact manifold $M$. Every non-linear elliptic differential operator $L_f \colon \mathcal{E} \to \mathcal{F}$ can be brought into a normal form at every $\phi \in \mathcal{E}$.* ◇

*Proof.* According to Theorem 2.2.13, we have to show that the bundle map $\mathrm{T}L_f\colon \mathrm{T}\mathcal{E} \to (L_f)^*\mathrm{T}\mathcal{F}$ is uniformly tame regular at $\phi \in \mathcal{E}$. Using tubular neighborhoods of $\operatorname{Im}\phi \subseteq E$ and $\operatorname{Im} L_f(\phi) \subseteq F$, it suffices to consider the case where $E$ and $F$ are vector bundles. In this linear setting, the vertical tangent bundle of $\mathrm{J}^r E$ is identified with $\mathrm{V}(\mathrm{J}^r E) \simeq \mathrm{J}^r E \times_M \mathrm{J}^r E$. Accordingly, the vertical derivative of $f$ is identified with a map $\mathrm{V}f\colon \mathrm{J}^r E \times_M \mathrm{J}^r E \to F$ and we have to show that the tame smooth family

$$T_\varphi\colon \mathcal{E} \to \mathcal{F}, \qquad \sigma \mapsto \mathrm{V}f(\mathrm{j}^r\varphi, \mathrm{j}^r\sigma) \tag{2.2.37}$$

of linear maps parametrized by $\varphi \in \mathcal{E}$ is uniformly tame regular at $\varphi = \phi$. For this purpose, note that the associated map $T\colon \mathcal{E} \times \mathcal{E} \to \mathcal{F}$ factorizes as the composition of the tame smooth map $\mathcal{E} \ni \varphi \mapsto (\mathrm{V}f)_*(\mathrm{j}^r\varphi, \cdot) \in \Gamma^\infty(\mathrm{L}(\mathrm{J}^r E, F))$ and the family of differential operators

$$\Gamma^\infty(\mathrm{L}(\mathrm{J}^r E, F)) \times \mathcal{E} \to \mathcal{F}, \qquad (\Lambda, \sigma) \mapsto \Lambda(\mathrm{j}^r\sigma). \tag{2.2.38}$$

By assumption, the differential operator with coefficients $\Lambda = (\mathrm{V}f)_*(\mathrm{j}^r\phi, \cdot)$ is elliptic. Thus, Lemma 2.1.12 and Theorem 2.1.20 imply that the family $T_\varphi$ is uniformly tame regular at $\varphi = \phi$. □

In applications, one sometimes encounters geometric spaces which are not realized as spaces of sections of fiber bundles, but whose linearization is still modeled on spaces of sections of vector bundles. A prime example is the Fréchet Lie group of diffeomorphisms of a compact manifold whose Lie algebra is the space of vector fields. In order to include such examples, we follow [Sub84, p. 57] and introduce the following subclass of Fréchet manifolds.

Definition 2.2.15 A tame Fréchet manifold $M$ is said to be *geometric* if it is locally modeled on the space of smooth sections of some vector bundle over a compact manifold.

A tame smooth map $f\colon M \to N$ is called *geometric* if for every point $m \in M$ there exist vector bundles $E$ and $F$ over the same compact manifold, and local trivializations $(\mathrm{T}M)_{\restriction U} \simeq U \times \Gamma^\infty(E)$ and $(f^*\mathrm{T}N)_{\restriction U} \simeq U \times \Gamma^\infty(F)$ in a neighborhood $U$ of $m$ such that in this trivialization the derivative $\mathrm{T}f\colon \mathrm{T}M \to f^*\mathrm{T}N$ factorizes as the composition of a tame smooth map $U \to \Gamma^\infty(\mathrm{L}(\mathrm{J}^r E, F))$ and the family of differential operators

$$\Gamma^\infty(\mathrm{L}(\mathrm{J}^r E, F)) \times \Gamma^\infty(E) \to \Gamma^\infty(F), \qquad (\Lambda, \sigma) \mapsto \Lambda(\mathrm{j}^r\sigma). \tag{2.2.39}$$

The notions of a geometric vector bundle and a geometric vector bundle morphism are defined in a similar way. ◊

Roughly speaking, a tame smooth map $f\colon M \to N$ is geometric if its linearization at a point $m \in M$ is a linear differential operator whose coefficients depend tamely on $m$. As we have seen in the proof of Theorem 2.2.14, every non-linear differential operator is a geometric map.

A slight reformulation of the proof of Theorem 2.2.14 gives the following slightly more general normal form theorem.

Theorem 2.2.16 (Normal form — Elliptic) *Let $f : M \to N$ be a geometric map between geometric Fréchet manifolds. If $\mathrm{T}_m f : \mathrm{T}_m M \to \mathrm{T}_{f(m)} N$ is an elliptic differential operator for some $m \in M$, then $f$ can be brought into a normal form at $m$.* ◊

# Moduli Spaces 3

In abstract terms, a moduli space is a space whose points parametrize isomorphism classes of geometric objects. Usually, one is mainly interested in a subclass of geometric objects satisfying an additional condition, which is often phrased in the form of a partial differential equation. This interest in moduli spaces stems from the fact that a better understanding of the geometric structure of the moduli space yields a deeper insight into the geometry of the objects and their deformations. Moreover, moduli spaces of solutions of geometric partial differential equations often reflect the intrinsic features of the underlying base manifold and thus can be used as a tool to extract topological invariants, which leads to striking geometric applications.

In the following, we will consider moduli spaces fitting into the following general setup. Let $f\colon M \to N$ be an equivariant map between $G$-manifolds. For every $\mu \in N$, let $M_\mu \equiv f^{-1}(\mu)$ be the $\mu$-level set of $f$ and set

$$\check{M}_\mu \equiv f^{-1}(\mu)/G_\mu, \tag{3.0.1}$$

where $G_\mu$ is the stabilizer subgroup of $\mu$ under the $G$-action on $N$. In the context of a moduli problem, $M$ is the space of geometric objects, the equation $f(m) = \mu$ describes the additional properties of these objects one is interested in, and the $G$-action implements the notion of equivalence. Thus, $\check{M}_\mu$ can be viewed as an abstract moduli space. In order to highlight the comprehensive and flexible nature of this general setting, let us briefly outline how some well-known moduli spaces fit into it.

Projective space: Let $f\colon \mathbb{C}^{n+1} \to \mathbb{R}$ be given by $f(z) = |z|^2$. Clearly, $f$ is equivariant with respect to the natural action of $\mathrm{U}(1)$ by rotation on $\mathbb{C}^{n+1}$ and the trivial action on $\mathbb{R}$. We recover the complex projective space $\mathbb{CP}^n$ as the associated moduli space

$$f^{-1}(1)/\mathrm{U}(1) = S^{2n+1}/\mathrm{U}(1), \tag{3.0.2}$$

whose points parametrize complex lines trough the origin in $\mathbb{C}^{n+1}$.

Gauge theory: Consider a finite-dimensional principal $G$-bundle $P$ over the compact manifold $\Sigma$. Let $f\colon \mathcal{C}(P) \to \Omega^2(\Sigma, \mathrm{Ad}P)$ be the map which assigns to a connection $A \in \mathcal{C}(P)$ on $P$ its curvature $F_A$. The map $f$ is equivariant with respect to the natural actions of the group $\mathcal{G}au(P)$ of gauge transformations of $P$. The associated moduli space

$$f^{-1}(0)/\mathcal{G}au(P) \tag{3.0.3}$$

is the moduli space of flat connections on $P$. This space, and the related moduli space of (central) Yang–Mills connections, has been extensively studied by Atiyah and Bott [AB83]. We will take up the discussion of this example in Section 4.4.

A similar setup yields the moduli space of anti-self-dual Yang–Mills connections, which will be the subject of Section 3.4. Moreover, the space of solutions of the Cauchy problem for the Yang–Mills–Higgs equation in (3 + 1)-dimensions will be studied in detail in Section 5.6.

Pseudo-holomorphic curves: Let $(M, \omega)$ be a finite-dimensional symplectic manifold endowed with a compatible almost complex structure $J$. For every compact Riemann surface $(\Sigma, j)$ with complex structure $j$, consider the Cauchy–Riemann operator

$$\bar{\partial}_{j,J} u = \frac{1}{2}(\mathrm{T}u + J \circ \mathrm{T}u \circ j), \tag{3.0.4}$$

where $u \colon \Sigma \to M$ is a smooth map. In order to realize $\bar{\partial}_{j,J}$ as a map between infinite-dimensional manifolds, let us introduce the Fréchet vector bundle $\mathcal{E} \to C^\infty(\Sigma, M)$ whose fiber over $u \in C^\infty(\Sigma, M)$ is the space $\Omega^{0,1}(\Sigma, u^*\mathrm{T}M)$ of smooth $(j, J)$-antilinear 1-forms on $\Sigma$ with values in $u^*\mathrm{T}M$. The Cauchy–Riemann operator yields a smooth section $f$ of $\mathcal{E}$ by setting $f(u) = \bar{\partial}_{j,J} u$. Note that $f$ is equivariant with respect to the natural reparametrization actions of the group $\mathcal{Diff}_j(\Sigma)$ of diffeomorphisms of $\Sigma$ preserving $j$. The associated moduli space

$$f^{-1}(0) / \mathcal{Diff}_j(\Sigma) \tag{3.0.5}$$

is the moduli space of pseudo-holomorphic curves[1]. We refer to [MS99] for further information.

These examples show that the moduli space $\check{M}_\mu = f^{-1}(\mu)/G_\mu$ may have a complicated geometry. In the simplest case, when $\mu$ is a regular value of $f$ and $G_\mu$ acts freely, the space $\check{M}_\mu$ is a smooth manifold. This is what happens, for example, for the projective space. However, in other examples, $f$ is not a submersion and the $G_\mu$-action is not free. Thus, the moduli space has singularities in these cases.

In this chapter, we will develop a framework to study the singular geometry of general moduli spaces of the form $\check{M}_\mu = f^{-1}(\mu)/G_\mu$. For this purpose, we first introduce the concept of an equivariant normal form and give suitable conditions which ensure that an equivariant map can be brought into such a normal form. In Sections 3.2 and 3.3, we will investigate the local structure of the moduli space $\check{M}_\mu = f^{-1}(\mu)/G_\mu$ under the assumption that $f$ can be brought into a normal form. In this case, $\check{M}_\mu$ can be locally identified with the quotient of the zero set of a smooth map by the linear action of a compact group, i.e. it has the structure of a Kuranishi space. We also find additional conditions on the normal form which ensure that $\check{M}_\mu$ is stratified by orbit types. Finally, in Section 3.4, we apply the general theory to the example of the moduli space of anti-self-dual Yang–Mills connections.

[1] Depending on the context, one usually restricts attention to simple or stable maps $u$ in a given homology class.

## 3.1 Normal form of an equivariant map

In this section, we study the local properties of a smooth equivariant map. Consider the following setup. Let $G$ be a locally convex Lie group and let $M$ and $N$ be locally convex $G$-manifolds. We refer the reader to Appendix A.2 for relevant background material and for our conventions regarding Lie group actions in infinite dimensions. Let $f\colon M \to N$ be a smooth equivariant map. Choose a point $m \in M$ and denote its image under $f$ by $\mu = f(m) \in N$. We assume throughout this section that the stabilizer subgroup $G_\mu$ of $\mu$ under the $G$-action on $N$ is a Lie subgroup of $G$ and that the induced $G_\mu$-action on $M$ is proper.

Let us start by describing the general idea of how to construct an equivariant normal form of $f$. Recall that a slice provides a normal form for the $G$-action in a neighborhood of a given orbit (see Appendix A.2 for background information concerning slices and tubular neighborhoods). Thus, it naturally comes to mind as a tool to study equivariant maps. An initial idea would be to split-off the $G$-action with the help of slices $S$ and $R$ for the $G$-action at the points $m \in M$ and $\mu \in N$, respectively. Every point $p \in M$ in a slice neighborhood of $m$ is of the form $p = g \cdot s$ with $g \in G$ and $s \in S$. Moreover, by $G$-equivariance, we have $f(p) = g \cdot f(s)$. If the slices satisfy $f(S) \subseteq R$, then the part of $f$ that does not come from the group action is captured in the map $f^R_S = f\restriction_S \colon S \to R$ between the slices (see Figure 3.1). For the reduced map $f^R_S$, we can use the normal form results of Section 2.2 to arrive at an equivariant normal form for $f$. In the finite-dimensional setting, a similar strategy has been used in [PW17, Proposition 2.6] to establish an equivariant submersion theorem for maps that are equivariant with respect to the action of a compact group.

However, the assumptions concerning the slices $S$ and $R$ are too restrictive for the applications we are interested in. To illustrate the shortcomings, for the moment, we restrict attention to the case when $f$ is an equivariant momentum map $J\colon M \to \mathfrak{g}^*$ for a symplectic $G$-action on a finite-dimensional symplectic manifold $(M, \omega)$. For non-compact groups $G$, the coadjoint action is never proper and thus we cannot expect it to admit a slice. Indeed, [GLS96, Point 5 in Section 2.3.1] gives a examples where the infinitesimal orbit $\mathfrak{g} \,.\, \mu$ has not even an $G_\mu$-invariant complement. Therefore, if we want to include symplectic quotients by non-compact groups, we cannot assume that the $G$-action on the target $N$ has a slice $R$.

As a solution, one may try to reduce this approach to its well-functioning elements and only use a slice for the $G$-action on $M$ while leaving the action on $N$ untouched. In the context of symplectic reduction, this idea is at the heart of the proof of the Marle–Guillemin–Sternberg normal form of an equivariant momentum map. In this case, the slice is constructed using an equivariant version of Darboux's theorem. However, a close inspection of the proof of the singular symplectic reduction theorem reveals that the Marle–Guillemin–Sternberg form is actually *not* the most convenient normal form for the study of the geometry of momentum map level sets. In the end, the problem arises from using a slice for the $G$-action but taking the quotient with respect to the subgroup $G_\mu$. This asymmetry was counterbalanced in the proof of the reduction theorem [OR03, Proposition 8.1.2] by deforming the momentum map $J$ using a local diffeomorphism of $\mathfrak{g}^*$. Instead of changing coordinates on the target, we will take a different approach and use a slice on $M$ for the action of $G_\mu$ instead of $G$. This has also the advantage that $G_\mu$ is often considerably smaller than $G$ and, thus, a slice for the subgroup action is easier to construct.

Finally, the symplectic setting is special in that the action of $G$ on $\mathfrak{g}^*$ is linear. When moving beyond the momentum map example, the $G$-action on the target $N$ is usually non-linear. Moreover, we cannot use a slice to control it as we have explained above. On the other hand, the subgroup $G_m$ is compact as the stabilizer of a proper action and its action on $N$ leaves $\mu$ invariant due to $G$-equivariance of $f$. Hence, we can at least hope to linearize the induced action of $G_m$ around the fixed point $\mu$. In summary, we will work with a $G_\mu$-slice $S$ on $M$ and with a $G_m$-slice on $N$, and then bring the restriction of $f$ to the slice $S$ in a normal form using the results of Section 2.2.

Remark 3.1.1 The main reason why slices are helpful for the study of local properties of equivariant maps is the basic observation that an equivariant map is locally completely determined by its restriction to the slice. Indeed, the restriction of a $G$-equivariant map $f\colon M \to N$ to a $G$-slice $S$ at the point $m \in M$ is a $G_m$-equivariant map $f_{\restriction S}: S \to N$. Conversely, every $G_m$-equivariant map $g: S \to N$ extends uniquely to $G$-equivariant map $g_{\mathrm{ext}}: G \cdot S \to N$ defined on the open neighborhood $G \cdot S$ of the orbit $G \cdot m$. Restriction and extension are obviously compatible in the sense that $(f_{\restriction S})_{\mathrm{ext}} = f_{\restriction G \cdot S}$ holds. ◊

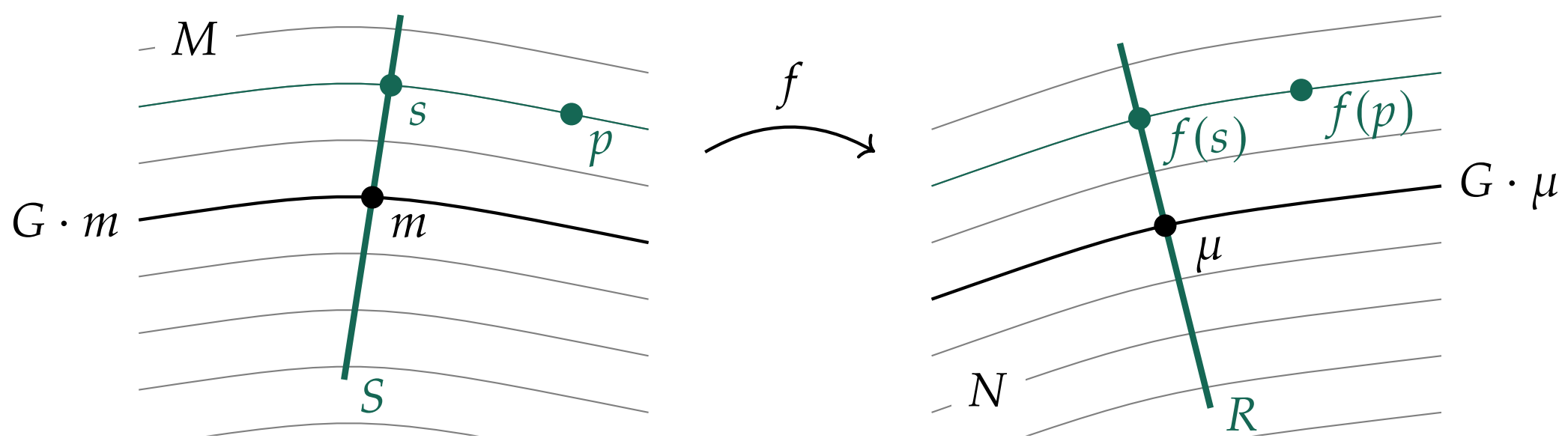

Figure 3.1: Illustration of the idea to capture the part of $f$ that does not come from the group action in a map $f_S^R\colon S \to R$ between the slices $S$ and $R$. If the slices satisfy $f(S) \subseteq R$, then every $p \in M$ in a slice neighborhood of $m$ can be written as $p = g \cdot s$ with $g \in G$ and $s \in S$ such that $f(p) = g \cdot f(s)$ with $f(s) \in R$.

### 3.1.1 General normal form theorem

Using the strategy outlined above, we will show that a wide class of equivariant maps can be brought into the following equivariant normal form.

Definition 3.1.2 An *abstract equivariant normal form* is a tuple $(H, X, Y, \hat{f}, f_{\mathrm{sing}})$ consisting of:

(i) a compact Lie group $H$,

(ii) locally convex vector spaces $X$ and $Y$ which are endowed with linear $H$-actions and which admit topological $H$-invariant decompositions

$$X = \mathrm{Ker} \oplus \mathrm{Coim} \qquad \text{and} \qquad Y = \mathrm{Coker} \oplus \mathrm{Im}, \tag{3.1.1}$$

(iii) an $H$-equivariant linear topological isomorphism $\hat{f}:\ \mathrm{Coim} \to \mathrm{Im}$,

(iv) a smooth $H$-equivariant map $f_{\mathrm{sing}}: X \supseteq U \mapsto \mathrm{Coker}$ defined in a $H$-invariant open neighborhood $U$ of $0$ such that $f_{\mathrm{sing}}(0, x_2) = 0$ for all $x_2 \in U \cap \mathrm{Coim}$ and $\mathrm{T}_0 f_{\mathrm{sing}} = 0$. ◊

For an equivariant normal form $(H, X, Y, \hat{f}, f_{\mathrm{sing}})$ and a Lie group $G_\mu$ with $H \subseteq G_\mu$, define the smooth map $f_{\mathrm{NF}}: G_\mu \times_H U \to G_\mu \times_H Y$ by

$$f_{\mathrm{NF}}\big([g, x_1, x_2]\big) = \big[g, \hat{f}(x_2) + f_{\mathrm{sing}}(x_1, x_2)\big] \tag{3.1.2}$$

for $g \in G_\mu$, $x_1 \in U \cap \mathrm{Ker}$ and $x_2 \in U \cap \mathrm{Coim}$.

Similar to the non-equivariant case studied in Section 2.2, the main theme is to reduce an equivariant map to an appropriate abstract equivariant normal form. In this context, a slice for the $G_\mu$-action on $M$ plays a fundamental role. Recall from Appendix A.2, that a $G_\mu$-slice at $m \in M$ is a $G_m$-invariant submanifold $S$ of $M$ which is transverse to the orbit $G_\mu \cdot m$ and which possesses a few further properties, see Definition A.2.2. In particular, the natural fibration $G_\mu \to G_\mu / G_m$ is a locally trivial principal bundle and it admits a local section $\chi: G_\mu / G_m \supseteq U \to G_\mu$ defined on an open neighborhood $U$ of the identity coset $[e]$ in such a way that the map

$$\chi^S: U \times S \to M, \qquad ([g], s) \mapsto \chi([g]) \cdot s \tag{3.1.3}$$

is a diffeomorphism onto an open neighborhood $V \subseteq M$ of $m$, i.e. the slice yields convenient product coordinates in a neighborhood of $m$. For every

$G_\mu$-slice $S$ at $m$, the tube map

$$\chi^{\mathrm{T}}: G_\mu \times_{G_m} S \to M, \qquad [g,s] \mapsto g \cdot s \tag{3.1.4}$$

is a $G_\mu$-equivariant diffeomorphism onto an open, $G_\mu$-invariant neighborhood of $G_\mu \cdot m$ in $M$, see Proposition A.2.3.

Definition 3.1.3 Let $M$ be a $G$-manifold and let $f: M \to N$ be a smooth $G$-equivariant map. Choose $\mu \in N$ and assume that the stabilizer $G_\mu$ is a Lie subgroup[1] of $G$. We say that $f$ can be *brought into the equivariant normal form* $(H, X, Y, \hat{f}, f_{\mathrm{sing}})$ at $m \in f^{-1}(\mu)$ if $H = G_m$ and there exist

(i) a linear slice $S$ at $m$ for the $G_\mu$-action,

(ii) a $G_m$-equivariant diffeomorphism $\iota_S: X \supseteq U \to S \subseteq M$ and

(iii) a $G_m$-equivariant chart $\rho: N \supseteq V' \to V \subseteq Y$ at $\mu$ with $f(S) \subseteq V'$

such that the following diagram commutes:

$$\begin{array}{ccc} M & \xrightarrow{\quad f \quad} & N \\ {\scriptstyle \chi^{\mathrm{T}} \circ (\mathrm{id}_{G_\mu} \times \iota_S)} \uparrow & & \uparrow {\scriptstyle \rho^{\mathrm{T}}} \\ G_\mu \times_{G_m} X \supseteq G_\mu \times_{G_m} U & \xrightarrow{\; f_{\mathrm{NF}} \;} & G_\mu \times_{G_m} V \subseteq G_\mu \times_{G_m} Y, \end{array} \tag{3.1.5}$$

where $\chi^{\mathrm{T}}: G_\mu \times_{G_m} S \to M$ is the tube diffeomorphism associated to $S$ and $\rho^{\mathrm{T}}: G_\mu \times_{G_m} V \to N$ is defined by $\rho^{\mathrm{T}}([g, y]) = g \cdot \rho^{-1}(y)$. ◇

We usually suppress the slice isomorphism $\iota_S$ and then view $f_{\mathrm{NF}}$ as a map defined on $G_\mu \times_{G_m} S$.

For the notion of an equivariant normal form, the linear action of the stabilizer $G_m$ on $X$ and $Y$ plays an important role. In finite dimensions, a classical result by Bochner [Boc45, Theorem 1] entails that every action of a compact group can be linearized near a fixed point. Here, under a *linearization* of the $G$-action on the manifold $M$ at a fixed point $m$ we understand a $G$-invariant open neighborhood $U$ of $m$ in $M$, a chart $\kappa: M \supseteq U \to X$ at $m$ and a linear $G$-action on $X$ such that $\kappa$ is $G$-equivariant. We have the following generalization of Bochner's theorem to the action of a compact Lie group on a Banach manifold.

Proposition 3.1.4 (Bochner's Linearization Theorem) *Every action of a compact Lie group $G$ on a Banach manifold $M$ can be linearized at a fixed point.* ◇

[1] It is not known, even for Banach Lie group actions, whether the stabilizer is always a Lie subgroup, see [Nee06, Problem IX.3.b].

*Proof.* We closely follow the proof in the finite-dimensional setting, see e.g. [DK99, Theorem 2.2.1].

Let $m \in M$ be a fixed point of the $G$-action on $M$. Choose an arbitrary chart $\tilde{\kappa}\colon M \supseteq \tilde{U} \to \tilde{V} \subseteq X$ at $m$. Since the action is continuous and since $m$ is a fixed point, for each $g \in G$ there exist open neighborhoods $W_g$ of $g$ in $G$ and $U_g$ of $m$ in $M$ such that $W_g \cdot U_g \subseteq \tilde{U}$. As $G$ is compact it can be covered by finitely many such open sets $W_{g_i}$. Then, $U := G \cdot \left(\bigcap_{g_i} U_{g_i}\right)$ is an open $G$-invariant neighborhood of $m$ contained in $\tilde{U}$.

Since $m$ is a fixed point, the linearization $\mathrm{T}_m \Upsilon_g\colon \mathrm{T}_m M \to \mathrm{T}_m M$ of the action $\Upsilon_g\colon M \to M$ of $g \in G$ endows $\mathrm{T}_m M$ with a linear $G$-action, which we carry over to $X$ using the isomorphism $\mathrm{T}_m \tilde{\kappa}\colon \mathrm{T}_m M \to X$. Define the map $\kappa\colon U \to X$ as the average

$$\kappa(m') := \int\limits_G a^{-1} \cdot \tilde{\kappa}(a \cdot m') \, \mathrm{d}a \tag{3.1.6}$$

for $m' \in U$, where $\mathrm{d}a$ is the normalized Haar measure on $G$. The integral exist, because $G$ is compact and $X$ is complete. A standard calculation shows that $\kappa$ intertwines the $G$-action on $U$ and the linear action on $X$. Moreover, for all $v \in \mathrm{T}_m M$ we have

$$\mathrm{T}_m \kappa(v) = \int\limits_G a^{-1} \cdot \mathrm{T}_{a \cdot m} \tilde{\kappa}(a \,.\, v) \, \mathrm{d}a = \mathrm{T}_m \tilde{\kappa}(v) \tag{3.1.7}$$

since $m$ is a fixed point. Thus, by the Inverse Function Theorem, $\kappa$ is a local diffeomorphism. □

Remark 3.1.5 The proof of Proposition 3.1.4 generalizes to arbitrary locally convex manifolds modeled on Mackey complete spaces except for the last argument, which uses the Inverse Function Theorem. We refrain from stating the tame Fréchet version, because the application of the Nash–Moser inverse function needs the existence of a chart $\tilde{\kappa}$ at $m$ with the equivariance property

$$\mathrm{T}_{a \cdot m'} \tilde{\kappa}(a \,.\, X) = a \cdot \mathrm{T}_{m'} \tilde{\kappa}(X) \tag{3.1.8}$$

not only at the point $m' = m$ but for $m'$ ranging over a neighborhood of $m$. Finding such a chart is almost as hard as directly guessing a chart in which the action is linear. Despite the fact that these problems concerning the Inverse Function Theorem prevent a generalization of Bochner's Theorem beyond Banach manifolds, it is often possible to find a linearization at a fixed point by direct means, e.g. by constructing the required chart by hand. ◊

After this little excursus about linearization of actions at fixed points, we come back to the main theme of an equivariant normal form.

Theorem 3.1.6 (Equivariant normal form — General) *Let $f\colon M \to N$ be an equivariant map between Fréchet G-manifolds. Let $m \in M$ and $\mu = f(m)$. Assume that the following conditions hold:*

*(i) The stabilizer subgroup $G_\mu$ of $\mu$ is a Lie subgroup of $G$.*

*(ii) The induced $G_\mu$-action on $M$ is proper and admits a slice $S$ at $m$.*

*(iii) The induced $G_m$-action on $N$ can be linearized at $\mu$.*

*(iv) The restriction $f_S \equiv f_{\restriction S}\colon S \to N$ of $f$ to $S$ can be brought into a normal form at $m$ in a $G_m$-equivariant way in the sense of Definition 2.2.1 and Remark 2.2.7 (iii).*

*Then, $f$ can be brought into an equivariant normal form at $m$.* ◇

*Proof.* The proof will proceed in a number of steps.

Step 1: Slice on $M$ for the $G_\mu$-action Let $S$ be a slice at $m$ for the induced $G_\mu$-action on $M$. Recall from Proposition A.2.3 that the tube map $\chi^{\mathrm{T}}\colon G_\mu \times_{G_m} S \to M$ defined by $\chi^{\mathrm{T}}([g, s]) = g \cdot s$ is a $G_\mu$-equivariant diffeomorphism onto an open neighborhood of $G_\mu \cdot m$ in $M$. By $G$-equivariance of $f$, the following diagram commutes:

$$\begin{array}{ccc} M & \xrightarrow{\quad f \quad} & N \\ {\scriptstyle \chi^{\mathrm{T}}}\uparrow & & \uparrow{\scriptstyle \mathrm{id}} \\ G_\mu \times_{G_m} S & \xrightarrow{\quad f_\chi \quad} & N, \end{array} \tag{3.1.9}$$

where the smooth map $f_\chi$ is defined by $f_\chi([g, s]) = g \cdot f(s)$ for $g \in G_\mu$ and $s \in S$. Thus, the map $f$ decomposes into the $G_\mu$-action and the restriction $f_S = f_{\restriction S}\colon S \to N$ of $f$ to the slice. Note that $f_S$ is $G_m$-equivariant.

Step 2: Linearization of the $G_m$-action on $N$ at $\mu$ By assumption, the $G_m$-action on $N$ can be linearized at $\mu$. That is, there exist a $G_m$-invariant open neighborhood $V'$ of $\mu$ in $N$, a chart $\rho\colon N \supseteq V' \to Y$ and a linear $G_m$-action on $Y$ such that $\rho$ is $G_m$-equivariant. Set $V = \rho(V') \subseteq X$ and define $\rho^{\mathrm{T}}\colon G_\mu \times_{G_m} V \to N$ by $\rho^{\mathrm{T}}([g, y]) = g \cdot \rho^{-1}(y)$. By possibly shrinking $S$, we may assume that $f(S) \subseteq V'$. Due to $G_m$-equivariance of $f_S$ and $\rho$, we can define $f_{\chi\rho}\colon G_\mu \times_{G_m} S \to G_\mu \times_{G_m} V$ by $f_{\chi\rho}([g, s]) = [g, \rho \circ f_S(s)]$ for $g \in G_\mu$ and $s \in S$. Then, the following diagram commutes:

$$\begin{array}{ccc} G_\mu \times_{G_m} S & \xrightarrow{\quad f_\chi \quad} & N \\ {\scriptstyle \mathrm{id}}\downarrow & & \uparrow{\scriptstyle \rho^{\mathrm{T}}} \\ G_\mu \times_{G_m} S & \xrightarrow{\quad f_{\chi\rho} \quad} & G_\mu \times_{G_m} V. \end{array} \tag{3.1.10}$$

Abbreviate $f_{S\rho} = \rho \circ f_S \colon S \to Y$ so that $f_{\chi\rho}([g,s]) = [g, f_{S\rho}(s)]$. At this point, we have completely split-off the $G_\mu$-action and reduced the problem to constructing a normal form of the map $f_{S\rho}$. Note that $f_{S\rho}$ is equivariant with respect to *linear* actions of the compact Lie group $G_m$, which is a considerable simplification of the general situation we started with.

Step 3: Bringing $f_{S\rho}$ into a normal form By Theorem A.2.4, we may regard $S$ as a $G_m$-equivariant open neighborhood of 0 in some Fréchet space $X$. By possibly shrinking $S$, we may assume that $f_{S\rho}$ takes values in an open $G_m$-invariant neighborhood $W$ of 0 in $Y$. Since, by assumption, $f_S$ can be brought into a normal form at $m$, there exist topological decompositions $X = \mathrm{Ker} \oplus \mathrm{Coim}$ and $Y = \mathrm{Coker} \oplus \mathrm{Im}$, and local diffeomorphisms $\psi \colon X \supseteq S \to X$ and $\phi \colon Y \supseteq W \to Y$ such that the following diagram commutes:

$$\begin{array}{ccc} X \supseteq S & \xrightarrow{f_{S\rho}} & \phi(W) \subseteq Y \\ \Big\downarrow \psi & & \phi \Big\uparrow \\ X \supseteq \psi(S) & \xrightarrow{\hat{f} + f_{\mathrm{sing}}} & W \subseteq Y, \end{array} \tag{3.1.11}$$

where $\hat{f} \colon \mathrm{Coim} \to \mathrm{Im}$ is a linear isomorphism and $f_{\mathrm{sing}} \colon S \to \mathrm{Coker}$ is a smooth map satisfying $\mathrm{T}_0 f_{\mathrm{sing}} = 0$ and $f_{\mathrm{sing}}(0, x_2) = 0$ for all $x_2 \in S \cap \mathrm{Coim}$. Without loss of generality, we may assume $\phi(W) \subseteq V$. Recall that $f_{S\rho}$ is $G_m$-equivariant with respect to the linear $G_m$-action on $Y$ introduced above. Thus, as explained in Remark 2.2.7 (iii), we may assume that the above decompositions of $X$ and $Y$ are $G_m$-invariant[1], and that $\psi$, $\phi$ and $f_{\mathrm{sing}}$ are $G_m$-equivariant maps. These equivariance properties allow us to define the map $f_{\mathrm{NF}} \colon G_\mu \times_{G_m} \psi(S) \to G_\mu \times_{G_m} W$ by

$$f_{\mathrm{NF}}([g, s']) = \big[g, \hat{f} \circ \mathrm{pr}_{\mathrm{Coim}}(s') + f_{\mathrm{sing}}(s')\big] \tag{3.1.12}$$

for $g \in G_\mu$ and $s' \in \psi(S)$. Since $\psi$ is a $G_m$-equivariant diffeomorphism, $\psi(S)$ is a new linear slice for the $G_\mu$-action at $m$.

[1] This is the only place in the proof where we need completeness of $X$ and $Y$.

Putting everything together Combining these three steps yields the following commutative diagram:

$$
\begin{array}{ccc}
M & \xrightarrow{\quad f \quad} & N \\
{\scriptstyle \chi^{\mathrm{T}}}\uparrow & & \uparrow{\scriptstyle \mathrm{id}} \\
G_\mu \times_{G_m} S & \xrightarrow{\quad f_\chi \quad} & N \\
{\scriptstyle \mathrm{id}}\downarrow & & \uparrow{\scriptstyle \rho^{\mathrm{T}}} \\
G_\mu \times_{G_m} S & \xrightarrow{f_{\chi\rho}} & G_\mu \times_{G_m} \phi(W) \\
{\scriptstyle \mathrm{id}\times\psi}\downarrow & & \uparrow{\scriptstyle \mathrm{id}\times\phi} \\
G_\mu \times_{G_m} \psi(S) & \xrightarrow{f_{\mathrm{NF}}} & G_\mu \times_{G_m} W.
\end{array}
\quad
\begin{array}{l}
\left.\vphantom{\begin{array}{c}a\\b\\c\end{array}}\right\} \text{1. Slice for } G_\mu \\
\left.\vphantom{\begin{array}{c}a\\b\\c\end{array}}\right\} \text{2. Slice for } G_m \\
\left.\vphantom{\begin{array}{c}a\\b\\c\end{array}}\right\} \text{3. Normal form of } f_{S\rho}
\end{array}
\tag{3.1.13}
$$

This finishes the proof of Theorem 3.1.6. □

Remarks 3.1.7

(i) In the setting of Theorem 3.1.6, assume additionally that $G_\mu$ is a split[1] Lie subgroup of $G$. In this case, we get control over the behavior of $f_S$ with respect to the residual action of $G/G_\mu$, at least on the infinitesimal level. Indeed, since $G_\mu$ is split, there exists a closed subspace $\mathfrak{q} \subseteq \mathfrak{g}$ such that

$$\mathfrak{g} = \mathfrak{g}_\mu \oplus \mathfrak{q} \tag{3.1.14}$$

is a topological decomposition. Intersecting the topological decomposition

$$\mathrm{T}_m M = \mathfrak{g}_\mu \,.\, m \oplus \mathrm{T}_m S \tag{3.1.15}$$

that comes from the slice property (SL3) of Definition A.2.2 with $\mathfrak{g} \,.\, m = \mathfrak{g}_\mu \,.\, m \oplus \mathfrak{q} \,.\, m$ yields $\mathfrak{q} \,.\, m \subseteq \mathrm{T}_m S$. Since $f$ is $G$-equivariant, the restriction of $\mathrm{T}_m f_S$ to $\mathfrak{q} \,.\, m$ takes the form

$$\mathrm{T}_m f_S \colon \mathfrak{q} \,.\, m \to \mathrm{T}_\mu N, \quad A \,.\, m \mapsto A \,.\, \mu. \tag{3.1.16}$$

By definition, $\mathfrak{q}$ is the complement of the stabilizer Lie algebra $\mathfrak{g}_\mu$ and hence $\mathrm{T}_m f_S$ yields an isomorphism of $\mathfrak{q} \,.\, m$ onto the orbit $\mathfrak{g} \,.\, \mu = \mathfrak{q} \,.\, \mu$. This isomorphism is, of course, part of the isomorphism $\hat{f} \colon \operatorname{Coim} \to \operatorname{Im}$ of the normal form.

(ii) Recall from (2.2.8) the local diffeomorphism $\psi \colon X \supseteq S \to X$ which brings $f_{S\rho}$ into a normal form. When $S$ is viewed as a submanifold of $M$, then

[1] A Lie subgroup $H \subseteq G$ is said to be split if its Lie algebra $\mathfrak{h}$ is topologically complemented in $\mathfrak{g}$.

$\Psi = \psi^{-1}$ satisfies $\Psi(0,0) = m$,

$$\mathrm{T}_{\Psi(x_1,x_2)} M = \mathfrak{g} \,.\, \Psi(x_1, x_2) + \operatorname{Im} \mathrm{T}_{(x_1,x_2)} \Psi \tag{3.1.17}$$

for all $(x_1, x_2) \in \psi(S)$, and

$$\mathrm{T}_m f \circ \mathrm{T}_{(0,0)} \Psi = \hat{f} \circ \mathrm{pr}_{\operatorname{Coim} \mathrm{T}_m f_S}. \tag{3.1.18}$$

In the context when $f$ is the momentum map for a symplectic $G$-action, a map with these properties is called a slice map in [Cho+03, Definition 2.1]. Although moving in a similar circle of ideas, the construction in [Cho+03, Proposition 2.2] of such a slice map makes no use of slices. Thus, our approach has the advantage of bringing the $G$-action into a normal form simultaneously.

(iii) Recall that in the finite-dimensional setting, for a proper action, a slice always exists [Pal61]. In the infinite-dimensional case, this may no longer be true and additional assumptions have to be made (see [Sub86; DR18c] for general slice theorems in infinite dimensions and [Ebi70; ACM89; CMM91] for constructions of slices in concrete examples). We refer the reader to Appendix A.2 for more details. Having in mind the application to Yang–Mills theory, it is especially important to note that, in particular, the action of the group of gauge transformations on the space of connections admits a slice. ◇

In the remainder of this section, we present variations of the Normal Form Theorem 3.1.6 which have assumptions that are often easier to verify in applications. Similar to the discussion in Section 2.2, we start with the simplest case and then work through various levels of functional-analytic settings, finishing with the tame Nash–Moser category.

### 3.1.2 Normal form theorem with finite-dimensional target or domain

Let us start by considering the finite-dimensional setting.

Theorem 3.1.8 (Equivariant normal form — finite-dimensional version) *Let $G$ be a finite-dimensional Lie group and let $f\colon M \to N$ be a smooth $G$-equivariant map between finite-dimensional $G$-manifolds. If the $G$-action on $M$ is proper, then $f$ can be brought into an equivariant normal form at every point.* ◇

*Proof.* Let $m \in M$ and $\mu = f(m)$. Since the stabilizer $G_\mu$ of $\mu$ is a closed subgroup of a finite-dimensional Lie group, it is a Lie subgroup. The induced $G_\mu$-action on $M$ is proper and thus the Slice Theorem of Palais [Pal61] implies that there exists a slice $S$ for the $G_\mu$-action at $m$. Properness of the action also implies that $G_m$ is compact and thus Bochner's Linearization Theorem implies that the $G_m$-action on $N$ can be linearized around the fixed point $\mu$. Let $f_S = f\restriction_S : S \to N$. By Theorem 2.2.6 and Remark 2.2.7 (iii), $f_S$ can be brought into a normal form in a $G_m$-equivariant way. Hence, all assumptions of Theorem 3.1.6 are verified and so $f$ can be brought into an equivariant normal form. □

This result for the finite-dimensional case can be generalized such that only one of the manifolds involved needs to be finite-dimensional.

Theorem 3.1.9 (Equivariant normal form — finite-dimensional target) *Let $G$ be a Lie group and let $f : M \to N$ be a smooth $G$-equivariant map between $G$-manifolds, where $N$ is finite-dimensional. Let $m \in M$ and $\mu = f(m)$. If the stabilizer subgroup $G_\mu$ of $\mu$ is a Lie subgroup of $G$ such that the induced $G_\mu$-action on $M$ is proper and admits a slice $S$ at $m$, then $f$ can be brought into an equivariant normal form at $m$.* ◇

*Proof.* The proof is similar to the one of Theorem 3.1.8 except for the fact that Theorem 2.2.9 has to be used instead of Theorem 2.2.6. □

Theorem 3.1.10 (Equivariant normal form — finite-dimensional domain) *Let $G$ be a finite-dimensional Lie group and let $f : M \to N$ be a smooth $G$-equivariant map between $G$-manifolds, where $M$ is finite-dimensional. Let $m \in M$ and $\mu = f(m)$. If the $G$-action on $M$ is proper and the induced $G_m$-action on $N$ can be linearized at $\mu$, then $f$ can be brought into an equivariant normal form at $m$.* ◇

*Proof.* The proof is analogous to the one of Theorem 3.1.8, only replacing Theorem 2.2.6 by Theorem 2.2.10. □

### 3.1.3 Tame Fréchet and elliptic normal form theorem

Let us now come to a version of Theorem 3.1.6 which lives in the tame Fréchet category.

Theorem 3.1.11 (Equivariant normal form — tame Fréchet version) *Let $G$ be a tame Fréchet Lie group and let $f : M \to N$ be an equivariant map between tame Fréchet $G$-manifolds. Let $m \in M$ and $\mu = f(m)$. Assume that the following conditions hold:*

*(i) The stabilizer subgroup $G_\mu$ of $\mu$ is a tame Fréchet Lie subgroup of $G$.*

*(ii) The induced $G_\mu$-action on $M$ is proper and admits a tame slice $S$ at $m$.*

*(iii) The induced $G_m$-action on $N$ can be linearized at $\mu$.*

*(iv) The chain*

$$0 \longrightarrow \mathfrak{g}_\mu \longrightarrow \mathrm{T}_s M \xrightarrow{\mathrm{T}_s f} \mathrm{T}_{f(s)} N \longrightarrow 0 \tag{3.1.19}$$

*of linear maps parametrized by $s \in S$ is uniformly tame regular at $m$. Here, the first map is the Lie algebra action given by $\xi \mapsto \xi \,.\, s$ for $\xi \in \mathfrak{g}_\mu$.*

*Then, $f$ can be brought into an equivariant normal form at $m$.* ◊

*Proof.* Let $S$ be a tame slice at $m$ modeled on the tame Fréchet space $X$. According to Theorem 3.1.6, we have to show that the map $f_S\colon S \to N$ can be brought into a normal form in a $G_m$-equivariant way. For this purpose, consider slice coordinates as in (SL3) of Definition A.2.2. That is, let $\chi\colon U \to G_\mu$ be a local section of $G_\mu \to G_\mu / G_m$ defined on an open neighborhood $U$ of the identity coset $[e]$ such that the map $\chi^S\colon U \times S \to M$ defined by $\chi^S([g], s) = \chi([g]) \cdot s$ is a local diffeomorphism. Clearly, $\mathrm{T}_{[e]}\chi\colon \mathfrak{g}_\mu / \mathfrak{g}_m \to \mathfrak{g}_\mu$ yields a continuous splitting of the exact sequence

$$0 \longrightarrow \mathfrak{g}_m \longrightarrow \mathfrak{g}_\mu \longrightarrow \mathfrak{g}_\mu / \mathfrak{g}_m \longrightarrow 0 \tag{3.1.20}$$

and thus induces a topological isomorphism $\mathfrak{g}_\mu \simeq \mathfrak{g}_\mu / \mathfrak{g}_m \times \mathfrak{g}_m$. With respect to this decomposition, we write elements $\xi \in \mathfrak{g}_\mu$ as pairs $([\xi], \xi_{\mathfrak{g}_m})$ with $[\xi] \in \mathfrak{g}_\mu / \mathfrak{g}_m$ and $\xi_{\mathfrak{g}_m} \in \mathfrak{g}_m$. Let $\rho\colon N \supseteq V' \to V \subseteq Y$ be a local chart at $\mu$ which linearizes the $G_m$-action on $N$. As above, let $f_{S\rho}\colon X \supseteq S \to Y$ denote the chart representation of $f_S\colon S \to N$. With respect to the local trivialization induced by $\chi^S$ and $\rho$, the chain (3.1.19) parametrized by $s \in S$ takes the form

$$0 \longrightarrow \mathfrak{g}_\mu \xrightarrow{\Gamma_s} \mathfrak{g}_\mu / \mathfrak{g}_m \times X \xrightarrow{\Xi_s} \times Y \longrightarrow 0 \tag{3.1.21}$$

where $\Gamma\colon S \times \mathfrak{g}_\mu \to \mathfrak{g}_\mu / \mathfrak{g}_m \times X$ and $\Xi\colon S \times \mathfrak{g}_\mu / \mathfrak{g}_m \times X \to Y$ are tame smooth families of linear maps defined by

$$\Gamma_s(\xi) = ([\xi], \xi_{\mathfrak{g}_m} \,.\, s) \tag{3.1.22}$$

for $s \in S$ and $\xi \in \mathfrak{g}_\mu$, and

$$\Xi_s([\xi], v) = \mathrm{T}_{f(s)}\rho\big((\mathrm{T}_{[e]}\chi([\xi]))\,.\,f(s)\big) + \mathrm{T}_s f_{S\rho}(v) \tag{3.1.23}$$

for $s \in S$, $[\xi] \in \mathfrak{g}_\mu / \mathfrak{g}_m$ and $v \in X$. Since the chain (3.1.19) is uniformly tame regular at $m$, we can assume that the chain (3.1.21) is uniformly tame regular at $m$. Note that $\operatorname{Im} \Gamma_m = \mathfrak{g}_\mu / \mathfrak{g}_m \times \{0\}$. Thus, Proposition 2.1.27 implies that the family

$$(\Xi_s)_{\restriction \{[0]\} \times X} = \mathrm{T}_s f_{S\rho}\colon X \mapsto Y \tag{3.1.24}$$

of linear maps parametrized by $s \in S$ is uniformly tame regular at $m$. Now, Theorem 2.2.13 shows that $f_{S\rho}$ can be brought into a normal form. The deforming diffeomorphisms can be chosen to be $G_m$-equivariant according to Remark 2.2.7 (iii). □

For the following elliptic version of the Equivariant Normal Form Theorem, the reader might want to recall the notion of a geometric Fréchet manifold from Definition 2.2.15.

Theorem 3.1.12 (Equivariant normal form — elliptic) *Let $G$ be a tame Fréchet Lie group and let $f\colon M \to N$ be an equivariant map between tame Fréchet $G$-manifolds. Let $m \in M$ and $\mu = f(m)$. Assume that the following conditions hold:*

*(i) The stabilizer subgroup $G_\mu$ of $\mu$ is a geometric tame Fréchet Lie subgroup of $G$.*

*(ii) The induced $G_\mu$-action on $M$ is proper and admits a geometric slice $S$ at $m$.*

*(iii) The induced $G_m$-action on $N$ can be linearized at $\mu$.*

*(iv) The chain*

$$0 \longrightarrow \mathfrak{g}_\mu \longrightarrow \mathrm{T}_s M \xrightarrow{\mathrm{T}_s f} \mathrm{T}_{f(s)} N \longrightarrow 0 \tag{3.1.25}$$

*is a chain of geometric linear maps parametrized by $s \in S$, which is an elliptic complex at $m$.*

*Then, $f$ can be brought into an equivariant normal form at $m$.* ◊

*Proof.* The claim follows directly as a special case of Theorem 3.1.11, because, according to Theorem 2.1.28, a tame family of chains of differential operators is uniformly tame regular in a neighborhood of a point at which the chain is an elliptic complex. □

## 3.2 Kuranishi structures

Kuranishi spaces were introduced by Fukaya and Ono [FO99, Section 1.5] in their study of the geometry of moduli spaces of pseudo-holomorphic curves. The notion of a Kuranishi structure builds on ideas of Kuranishi [Kur65] for the moduli space of complex structures and of Donaldson [Don83, Section II.2] and Taubes [Tau82, Section 6] for moduli problems in gauge theory, see also [DK97, Section 4.2.5]. As we will see, Kuranishi structures play an important role in our general setting, too.

Definition 3.2.1 Let $X$ be a topological space. A *Kuranishi chart* at a point $x \in X$ is a tuple $(V, E, F, H, s, \kappa)$ consisting of the following data:

(i) an open neighborhood $V$ of $0$ in a locally convex vector space $E$,

(ii) a locally convex vector space $F$,

(iii) a compact Lie group $H$ acting linearly and continuously on $V$ and $F$,

(iv) a smooth $H$-equivariant map $s\colon V \to F$,

(v) a homeomorphism $\kappa$ from $s^{-1}(0)/H$ to a neighborhood of $x$ in $X$.

The bundle $F \times V \to V$ is called the *obstruction bundle* and $s$ is referred to as the *obstruction map*. If $E$ and $F$ are finite-dimensional vector spaces, then the Kuranishi chart is said to be *finite-dimensional*. ◊

Informally speaking, a space admitting Kuranishi charts is locally modeled on the quotient of the zero set of a smooth map by a compact group. As such, Kuranishi charts extend the notion of an orbifold and of a manifold chart. Indeed, if the map $s$ vanishes and $H$ is a finite group acting faithfully on $E$, then the Kuranishi chart reduces to an (infinite-dimensional) orbifold chart, i.e., $X$ is locally modeled on the quotient of $E$ by a finite group action. If, in addition, the $H$-action on $E$ is trivial, then we obtain an ordinary manifold chart on $X$. For a finite-dimensional Kuranishi chart $(V, E, F, H, s, \kappa)$ at $x \in X$, the number

$$\dim E - \dim F - \dim H \tag{3.2.1}$$

is called the *virtual dimension* of $X$ (at $x$). This is motivated by the fact that $s^{-1}(0)/H$ forms a manifold of this dimension at regular points of $s$ for which the $H$-action is free.

Remarks 3.2.2

(i) Our notion of a finite-dimensional Kuranishi chart is a slight generalization of the one proposed by Oh et al. [Oh+09, Definition A.1.1]. There, $H$ is assumed to be a finite group (acting effectively on $V$). Finiteness of $H$ is a natural assumption in the context of moduli spaces of pseudo-holomorphic curves, but it is too restrictive in our more general setting. Moreover, we do not require $X$ to be compact nor to be endowed with a metric.

Usually, only finite-dimensional Kuranishi charts are discussed in the literature. However, for general moduli spaces, one cannot expect the Kuranishi chart to be finite-dimensional. As we will see below in Remark 3.2.4, this amounts to requiring that a certain complex is Fredholm.

(ii) In order to define a Kuranishi structure on $X$, in a similar spirit to the smooth structure of a manifold, one needs to introduce coordinate transition maps, which explain how to glue together different Kuranishi charts. A variety of definitions of suitable chart transitions are proposed in the literature, each with their own advantages and functorial properties. We refer the reader to [Joy14, Appendix A] for a recent review on this matter.

(iii) Recently, Hofer, Wysocki, and Zehnder introduced the polyfold framework as a different approach to deal with the analytic and geometric issues occurring in the study of moduli spaces in symplectic field theory, see [HWZ17a; HWZ17b] and references therein. We refer to [Yan14] for an extensive discussion of the relation of Kuranishi structures and polyfold theory. Moreover, in Remark A.1.4 (i), we compare the theory of tame Fréchet spaces to the so-called scale calculus that underlies the polyfold framework. ◊

Theorem 3.2.3 *Let $f\colon M \to N$ be an equivariant map between (locally convex) G-manifolds. Let $\mu \in N$ be such that $M_\mu \equiv f^{-1}(\mu)$ is not empty and that the stabilizer subgroup $G_\mu$ of $\mu$ is a Lie subgroup of $G$. If $f$ can be brought into an equivariant normal form at every point of $M_\mu$, then there exists a Kuranishi chart on $\check{M}_\mu \equiv f^{-1}(\mu)/G_\mu$ at every point.* ◊

*Proof.* Let $m \in M_\mu$ and let $(S, \rho, \hat{f}, f_{\mathrm{sing}})$ be an equivariant normal form of $f$ at $m$ in the sense of Definition 3.1.3. Then, there exists maps $\chi$ and $\rho$ such that the following diagram commutes:

$$\begin{array}{ccc} M & \xrightarrow{\quad f \quad} & N \\ \chi^{\mathrm{T}}\uparrow & & \uparrow\rho^{\mathrm{T}} \\ G_\mu \times_{G_m} X \supseteq G_\mu \times_{G_m} S & \xrightarrow{f_{\mathrm{NF}}} & G_\mu \times_{G_m} V \subseteq G_\mu \times_{G_m} Y, \end{array} \tag{3.2.2}$$

where $f_{\mathrm{NF}}$ is defined by $f_{\mathrm{NF}}([g, x_1, x_2]) = [g, \hat{f}(x_2) + f_{\mathrm{sing}}(x_1, x_2)]$ for $g \in G_\mu$, $x_1 \in S \cap \operatorname{Ker}$ and $x_2 \in S \cap \operatorname{Coim}$. Recall that the map $\rho^{\mathrm{T}}\colon G_\mu \times_{G_m} V \to N$ is given by $\rho^{\mathrm{T}}([g, y]) = g \cdot \rho^{-1}(y)$, where $\rho\colon N \supseteq V' \to V \subseteq Y$ is a chart on $N$ satisfying $\rho(\mu) = 0$. Hence, $(\rho^{\mathrm{T}})^{-1}(\mu) = G_\mu \times_{G_m} \{0\}$. Using the commutative diagram (3.2.2), we obtain

$$\begin{aligned} (\chi^{\mathrm{T}})^{-1}(M_\mu) &= (\chi^{\mathrm{T}})^{-1}\big(f^{-1}(\mu)\big) \\ &= f_{\mathrm{NF}}^{-1}\big((\rho^{\mathrm{T}})^{-1}(\mu)\big) \\ &= f_{\mathrm{NF}}^{-1}\big(G_\mu \times_{G_m} \{0\}\big) \\ &= G_\mu \times_{G_m} \{(x_1, 0) \in S : \; f_{\mathrm{sing}}(x_1, 0) = 0\}. \end{aligned} \tag{3.2.3}$$

By $G_\mu$-equivariance, the tube diffeomorphism $\chi^{\mathrm{T}}$ thus induces a local homeomorphism of $\check{M}_\mu = f^{-1}(\mu)/G_\mu$ with

$$\{(x_1, 0) \in S : \; f_{\mathrm{sing}}(x_1, 0) = 0\}/G_m. \tag{3.2.4}$$

This local homeomorphism provides a Kuranishi chart on $\check{M}_\mu$ with $V = S \cap \operatorname{Ker}$, $F = \operatorname{Coker}$, $H = G_m$ and $s = f_{\mathrm{sing}}(\cdot, 0)\colon V \to F$ in the notation of Definition 3.2.1. □

Remark 3.2.4 In the setting of Theorem 3.2.3, consider the chain complex

$$0 \longrightarrow \mathfrak{g}_\mu \xrightarrow{\mathrm{T}_e \Upsilon_m} \mathrm{T}_m M \xrightarrow{\mathrm{T}_m f} \mathrm{T}_\mu N \longrightarrow 0, \tag{3.2.5}$$

where $\Upsilon_m : G_\mu \longrightarrow M$ is the orbit map of the $G_\mu$-action at $m \in M$. The homology groups of this complex are identified with the spaces

$$\mathfrak{g}_m, \qquad \operatorname{Ker} \mathrm{T}_m f / \mathfrak{g}_\mu \,.\, m \simeq \operatorname{Ker}, \qquad \mathrm{T}_\mu N / \operatorname{Im} \mathrm{T}_m f \simeq \operatorname{Coker}, \tag{3.2.6}$$

respectively. Here, Ker and Coker refer to the spaces occurring in the equivariant normal form as in Definition 3.1.3. The Kuranishi charts on $\check{M}_\mu$ constructed in the proof of Theorem 3.2.3 had $H = G_m$, $E = \operatorname{Ker}$ and $F = \operatorname{Coker}$. Thus, these Kuranishi charts are finite-dimensional if and only if the complex (3.2.5) is Fredholm. In this case, the Euler characteristic (2.1.64) of the complex (3.2.5),

$$\dim \mathfrak{g}_m - \dim \operatorname{Ker} + \dim \operatorname{Coker}, \tag{3.2.7}$$

is (minus) the virtual dimension of $X$. ◇

Since every equivariant map between finite-dimensional $G$-manifolds can be brought into an equivariant normal form according to Theorem 3.1.8, we obtain the following corollary of Theorem 3.2.3.

Corollary 3.2.5 *Let $f : M \longrightarrow N$ be an equivariant map between finite-dimensional $G$-manifolds and let $\mu \in N$. If the $G_\mu$-action on $M$ is proper, then $\check{M}_\mu \equiv f^{-1}(\mu)/G_\mu$ admits a Kuranishi chart at every point.* ◇

## 3.3 Orbit type stratification

As we have seen in the previous section, the abstract moduli space $\check{M}_\mu = f^{-1}(\mu)/G_\mu$ admits Kuranishi charts if $f$ can be brought into an equivariant normal form. Under additional conditions on the equivariant normal form, $\check{M}_\mu$ carries even more structure.

Proposition 3.3.1 *Let $f : M \longrightarrow N$ be a smooth equivariant map between (locally convex) $G$-manifolds and let $\mu \in N$. Assume that the stabilizer subgroup $G_\mu$ is a Lie subgroup of $G$ and that the induced action of $G_\mu$ on $M$ is proper. Moreover, assume that $f$ can be brought into an equivariant normal form $(S, \rho, \hat{f}, f_{\mathrm{sing}})$ at every $m \in f^{-1}(\mu)$ such that the following holds:*

*(i) (submanifold property) The set*

$$\left\{(x_1, 0) \in S_{(G_m)} : f_{\mathrm{sing}}(x_1, 0) = 0\right\} = f_{\mathrm{sing}}^{-1}(0) \cap \operatorname{Ker} \cap S_{(G_m)} \tag{3.3.1}$$

*is a submanifold of $S_{(G_m)}$.*

*(ii) (approximation property) For every orbit type $(K) \leq (G_m)$ of the $G_\mu$-action on $M$, the point $0$ lies in the closure of $f_{\mathrm{sing}}^{-1}(0) \cap \operatorname{Ker} \cap S_{(K)}$ in $S$.*

*Then, the partition of* $M_\mu = f^{-1}(\mu)$ *into the orbit type subsets* $M_\mu \cap M_{(H)}$ *is a stratification. Moreover, under these assumptions, the decomposition of* $\check{M}_\mu = f^{-1}(\mu)/G_\mu$ *into the sets* $M_\mu \cap M_{(H)}/G_\mu$ *is a stratification, too.* ◇

*Proof.* Let $m \in M_\mu$ and let $(S, \rho, \hat{f}, f_{\mathrm{sing}})$ be an equivariant normal form of $f$ at $m$ satisfying the submanifold and approximation property. In the proof of Theorem 3.2.3, we have seen that the $G_\mu$-equivariant tube diffeomorphism $\chi^{\mathrm{T}}$ locally identifies the subset $M_\mu$ with

$$G_\mu \times_{G_m} \{(x_1, 0) \in S : \; f_{\mathrm{sing}}(x_1, 0) = 0\}. \tag{3.3.2}$$

Accordingly, $M_\mu \cap S_{(G_m)}$ is identified with the set

$$\big\{(x_1, 0) \in S_{(G_m)} : \; f_{\mathrm{sing}}(x_1, 0) = 0\big\}. \tag{3.3.3}$$

Since the latter is a submanifold of $S_{(G_m)}$ by the submanifold property, we conclude that $M_\mu \cap S_{(G_m)}$ is a submanifold of $S_{(G_m)}$. Moreover, the approximation property entails that, for every $G_\mu$-orbit type $(K) \leq (G_m)$, the point $m$ lies in the closure of $M_\mu \cap S_{(K)}$. Thus, the claims follow from Proposition A.2.8. □

## 3.4 Application: Gauge orbit types and anti-self-dual connections

In this section, we are concerned with the local properties of the moduli space of anti-self-dual Yang–Mills connections. The local structure of this moduli space is well-understood, see e.g. [DK97]. We will show how these structure results can be rederived with relatively small effort using the general framework developed in the previous sections.

Before we come to the anti-self-dual equation, let us recall a few general results concerning the orbit type stratification for the action of the group of gauge transformations on the space of connections. For this purpose, consider the principal $G$-bundle $\pi\colon P \to M$ over the compact $n$-dimensional Riemannian manifold $M$ with $G$ being a compact Lie group. A connection in $P$ is a $G$-equivariant splitting of the tangent bundle $\mathrm{T}P = \mathrm{V}P \oplus \mathrm{H}P$ into the canonical vertical distribution $\mathrm{V}P$ and a horizontal distribution $\mathrm{H}P$. Recall that $\mathrm{V}P$ is spanned by the Killing vector fields $p \mapsto \xi \,.\, p$ for $\xi \in \mathfrak{g}$ and hence it is isomorphic to the adjoint bundle $\mathrm{Ad}P$. A connection $A$ in $P$ yields a splitting of the Atiyah sequence

$$0 \longrightarrow \mathrm{Ad}P \longrightarrow \mathrm{T}P/G \xrightarrow{\mathrm{T}\pi} \mathrm{T}M \longrightarrow 0 \tag{3.4.1}$$

and we realize the space $\mathcal{C}(P)$ of connections on $P$ as that of vector bundle sections $\mathrm{T}M \to \mathrm{T}P/G$. Thus, the points of $\mathcal{C}(P)$ are in bijective correspondence with the sections of an affine bundle over $M$ and in this way $\mathcal{C}(P)$ carries a natural tame Fréchet manifold structure modeled on the vector space $\Omega^1(M, \mathrm{Ad}P)$, cf. [Abb+86]. The covariant derivative with respect to $A$ is denoted by $\mathrm{d}_A$ and the curvature of $A$ is written as $F_A$.

The group $\mathcal{G}au(P)$ of local gauge transformations on $P$ is a tame Fréchet Lie group, because it is realized as the space of smooth sections of the group bundle $P \times_G G$, see [CM85] for details. The natural left action of $\mathcal{G}au(P)$ on $\mathcal{C}(P)$ is given by

$$A \mapsto \mathrm{Ad}_\lambda A + \lambda\, \mathrm{d}\lambda^{-1}, \tag{3.4.2}$$

for $\lambda \in \mathcal{G}au(P)$. This action is proper, see [RSV02b; Die13]. Moreover, it admits a slice $\mathcal{S}_{A_0}$ at every $A_0 \in \mathcal{C}$ given by the Coulomb gauge condition. That is[1],

$$\mathcal{S}_{A_0} := \{A \in \mathcal{U} : \, \mathrm{d}^*_{A_0}(A - A_0) = 0\}, \tag{3.4.3}$$

where $\mathcal{U}$ is a suitable open neighborhood of $A_0$ in $\mathcal{C}$. One uses the Nash–Moser Inverse Function Theorem to show that $\mathcal{S}_{A_0}$ is a slice. The details can be found in [ACM89; Die13]. In the Banach context, the orbit type decomposition of $\mathcal{C}(P)$ has been extensively studied in [KR86], see also [RSV02b]. The proof that the decomposition satisfies the frontier condition [KR86, Theorem 4.3.5] carries over to our Fréchet setting without major changes. As a consequence, the decomposition of $\mathcal{C}(P)$ into gauge orbit types is a stratification, see Proposition A.2.8.

Finally, let us comment on the classification of the gauge orbit types. For this purpose, choose a point $p_0 \in P$. Evaluation at $p_0$ of a gauge transformation $\lambda \in \mathcal{G}au(P)$, seen as a $G$-equivariant map $P \to G$, yields a Lie group homomorphism

$$\mathrm{ev}_{p_0} \colon \mathcal{G}au(P) \to G. \tag{3.4.4}$$

A gauge transformation $\lambda \in \mathcal{G}au(P)$ leaves the connection $A$ invariant if and only if $\lambda$ is constant on every $A$-horizontal curve, that is, if it is constant on the holonomy bundle $P_A$ of $A$ (based at $p_0$). By $G$-equivariance, $\lambda \in \mathcal{G}au_A(P)$ is thus completely determined by its value at some point in $P_A$. Hence, the evaluation map $\mathrm{ev}_{p_0}$ restricts to an isomorphism of Lie groups between the stabilizer $\mathcal{G}au_A(P)$ and the centralizer $\mathrm{C}_G(\mathrm{Hol}_A)$ of the holonomy group of $A$ based at $p_0$, cf. [RSV02b, Theorem 2.1]. We usually suppress the evaluation map in our notation and view $\mathcal{G}au_A(P)$ directly as a Lie subgroup of $G$. If the stabilizer $\mathcal{G}au_A(P)$ of $A$ is trivial, then $A$ is said to be irreducible. Recall that a subgroup which can be written as a centralizer is called a Howe subgroup. In particular,

$$H = \mathrm{C}_G\big(\mathcal{G}au_A(P)\big) = \mathrm{C}^2_G(\mathrm{Hol}_A) \tag{3.4.5}$$

[1] Here, as usual, $\mathrm{d}^*_A \alpha := (-1)^k * \mathrm{d}_A * \alpha$ for a $k$-form $\alpha$.

is a Howe subgroup of $G$. Consider the $H$-principal bundle

$$P_H := P_A \times_{\mathrm{Hol}_A} H \equiv P_A \cdot H \subseteq P \tag{3.4.6}$$

associated to the holonomy bundle $P_A$. The bundle $P_H$ consists of all $p \in P$ obeying $\lambda(p) = \lambda(p_0)$ for every $\lambda \in \mathcal{G}au_A(P)$. Conversely, the stabilizer subgroup can be recovered from $P_H$ as the subgroup

$$\mathcal{G}au_A(P) = \{\lambda \in \mathcal{G}au(P) : \lambda_{\restriction P_H} = \mathrm{const}\}. \tag{3.4.7}$$

A bundle reduction of $P$ to a Howe subgroup is called a Howe subbundle. A bundle reduction $Q$ of $P$ is called holonomy-induced if there exists a connected reduction $\tilde{Q} \subseteq P$ to a subgroup $\tilde{H}$ such that $Q = \tilde{Q} \cdot \mathrm{C}^2_G(\tilde{H})$. The following proposition (cf. [RSV02b, Theorem 3.3]) is basic for the classification procedure.

Proposition 3.4.1 *Let $P$ be a principal $G$-bundle over the compact manifold $M$ with* $\dim M \geq 2$. *Then, the map*

$$[\mathcal{G}au_A(P)] \mapsto \left[P_A \cdot \mathrm{C}^2_G(\mathrm{Hol}_A)\right] \tag{3.4.8}$$

*from the set of gauge orbit types to the set of isomorphism classes of holonomy-induced Howe subbundles of $P$ (factorized by the action of $G$) is a bijection.* ◊

In [RSV02b], this proposition is proved in the context of Sobolev spaces, but the result clearly holds true in the smooth category as well.

Remark 3.4.2 By [RSV02a, Theorem 6.2], every Howe subbundle of a principal $\mathrm{SU}(n)$-bundle is automatically holonomy induced. However, for other classical groups this is not always the case, see the counterexample after Theorem 6.2 in [RSV02a]. ◊

Thus, to enumerate the gauge orbit types for a given principal $G$-bundle $P \to M$, one has to work through the following program:

(i) Classify the Howe subgroups up to conjugacy.

(ii) Classify the Howe subbundles up to isomorphy.

(iii) Extract the Howe subbundles which are holonomy-induced.

(iv) Factorize by the principal action.

(v) Determine the natural partial order.

The classification of the orbit types for all classical groups has been accomplished in [RSV02a; RSV02b; RSV02c; HRS10; HRS11]. For the convenience of the reader, we recall the result for the case of $G = \mathrm{SU}(n)$, see also [RSV02c] for the discussion of the partial order.

Proposition 3.4.3 *Let $P$ be a principal $\mathrm{SU}(n)$-bundle over a compact manifold $M$ of dimension $2, 3$, or $4$. The gauge orbit types of the space of connections on $P$ are in one-to-one correspondence with symbols $[(I;\alpha,\xi)]$, where*

*(i)* $I = \big((k_1,\dots,k_r),(m_1,\dots,m_r)\big)$ *is a pair of sequences of positive integers obeying*

$$\sum_{i=1}^{r} k_i m_i = n\,, \tag{3.4.9}$$

*(ii)* $\alpha = (\alpha_1,\dots,\alpha_r)$ *is a sequence of elements* $\alpha_i \in \mathrm{H}^*(M,\mathbb{Z})$ *representing admissible values of Chern classes of* $\mathrm{U}(k_i)$*-bundles over* $M$*,*

*(iii)* $\xi \in \mathrm{H}^1(M,\mathbb{Z}_d)$*, where $d$ is the greatest common divisor of* $(m_1,\dots,m_r)$*.*

*The cohomology elements $\alpha_i$ and $\xi$ are subject to the relations*

$$\sum_{i=1}^{r} \frac{m_i}{d}\alpha_{i,1} = \beta_d(\xi), \quad \alpha_1^{m_1} \cup \ldots \cup \alpha_r^{m_r} = \mathrm{c}(P), \tag{3.4.10}$$

*where $\mathrm{c}(P)$ is the total Chern class of $P$ and* $\beta_d\colon \mathrm{H}^1(M,\mathbb{Z}_d) \to \mathrm{H}^2(M,\mathbb{Z})$ *is the Bockstein homomorphism associated with the short exact sequence of coefficient groups* $0 \to \mathbb{Z} \to \mathbb{Z} \to \mathbb{Z}_d \to 0$*. For any permutation $\sigma$ of* $\{1,\dots r\}$*, the symbols* $[(I;\alpha,\xi)]$ *and* $[(\sigma I;\sigma\alpha,\xi)]$ *have to be identified.* ◇

With these general results on the gauge orbit type stratification in mind, we can now investigate the solution space of the anti-self-dual equation. For this purpose, let $M$ be a oriented compact Riemannian manifold of dimension 4. Consider a principal $G$-bundle $P \to M$, where $G$ is a compact, semi-simple Lie group. On a 4-dimensional manifold, the Hodge star operator $*$ associated to the Riemannian metric satisfies $** = \mathrm{id}$ on 2-forms and thus determines a decomposition

$$\Omega^2(M) = \Omega^2_+(M) \oplus \Omega^2_-(M) \tag{3.4.11}$$

of the space of 2-forms into the $\pm1$-eigenspaces. A similar decomposition holds for vector-valued 2-forms and, in particular, the curvature $F_A$ of a connection $A \in \mathcal{C}(P)$ can be written as

$$F_A = F_A^+ + F_A^- \tag{3.4.12}$$

with $F_A^\pm \in \Omega^2_\pm(M,\mathrm{Ad}P)$. A connection $A$ with $F_A^+ = 0$ is called *anti-self-dual (ASD)*. The Bianchi identity implies that an ASD connection satisfies the Yang–Mills equation. The self-dual part of the curvature gives a smooth map

$$F^+\colon \mathcal{C}(P) \to \Omega^2_+(M,\mathrm{Ad}P), \tag{3.4.13}$$

which is equivariant with respect to the natural actions of the group $\mathcal{G}au(P)$ of gauge transformations. The *moduli space of anti-self-dual connections* is, by

definition, the space

$$\mathcal{A} = (F^+)^{-1}(0)/\mathcal{G}au(P) \tag{3.4.14}$$

of anti-self-dual connections on $P$ modulo gauge equivalence. This clearly fits into the general framework considered in the previous sections. Let us verify that the assumptions of Theorem 3.1.12 are met for the model under consideration:

(i) The stabilizer of $\mu = 0 \in \Omega^2_+(M, \mathrm{Ad}P)$ is the whole group $\mathcal{G}au(P)$, which is a geometric tame Fréchet Lie group with Lie algebra $\mathfrak{gau}(P) = \Gamma^\infty(\mathrm{Ad}P)$.

(ii) The natural action of $\mathcal{G}au(P)$ on $\mathcal{C}(P)$ is proper and admits a slice $\mathcal{S}$ at $A \in \mathcal{C}(P)$ as discussed above.

(iii) The action of $\mathcal{G}au(P)$ on $\Omega^2_+(M, \mathrm{Ad}P)$ is clearly linear.

(iv) Let $A$ be an ASD connection and let $\mathcal{S}$ be a slice at $A$. The chain (3.1.19) here takes the form

$$0 \longrightarrow \Omega^0(M, \mathrm{Ad}P) \xrightarrow{\mathrm{d}_B} \Omega^1(M, \mathrm{Ad}P) \xrightarrow{\mathrm{d}_B^+} \Omega^2_+(M, \mathrm{Ad}P) \longrightarrow 0, \tag{3.4.15}$$

where $B \in \mathcal{S}$ and $\mathrm{d}_B^+$ denotes the self-dual part of the covariant derivative $\mathrm{d}_B$. This chain is clearly a chain of linear differential operators tamely parametrized by the connection $B$. The ASD condition for $A$ asserts that $\mathrm{d}_A^+ \circ \mathrm{d}_A = 0$ and so, at $B = A$, the chain (3.4.15) is a complex, called the *Yang–Mills complex*. Ellipticity of the Yang–Mills complex is well-known and follows from a straightforward computation in linear algebra, see e.g. [RS17, Lemma 6.5.2]. Moreover, the Atiyah–Singer Index Theorem shows that its Euler characteristic is given by

$$-2\,\mathrm{p}_1(\mathrm{Ad}P) + \frac{1}{2}(\chi_M - \sigma_M) \cdot \dim G \tag{3.4.16}$$

where $\mathrm{p}_1(\mathrm{Ad}P)$ is the Pontryagin index of the adjoint bundle, $\chi_M$ is the Euler number of $M$ and $\sigma_M$ is the signature[1] of $M$, see [RS17, Lemma 6.5.5].

Hence, by Theorem 3.1.12, the map $F^+$ can be brought into an equivariant normal form at every ASD connection $A \in \mathcal{C}(P)$. Moreover, as a consequence

[1] The intersection form of $M$ is the symmetric, non-degenerate bilinear form given by

$$\mathrm{H}^2(M) \times \mathrm{H}^2(M) \to \mathbb{R}, \qquad ([\alpha], [\beta]) \mapsto \int_M \alpha \wedge \beta. \tag{3.4.17}$$

The signature $\sigma_M$ of $M$ is defined as the difference between the number of positive eigenvalues and the number of negative eigenvalues of the quadratic form corresponding to the intersection form.

of Theorem 3.2.3 and Remark 3.2.4, we obtain the following description of the local geometry of the moduli space $\mathcal{A}$ of anti-self-dual connections.

Theorem 3.4.4 *Let $P$ be a principal $G$-bundle with a compact, semi-simple structure group $G$ over a 4-dimensional compact Riemannian manifold $M$. Then, the moduli space $\mathcal{A}$ of anti-self-dual connections on $P$ admits a finite-dimensional Kuranishi chart at every point $[A] \in \mathcal{A}$. Moreover, the virtual dimension of $\mathcal{A}$ is given by*

$$2\,\mathrm{p}_1(\mathrm{Ad}P) - \frac{1}{2}(\chi_M - \sigma_M) \cdot \dim G. \tag{3.4.18}$$

◊

Let us describe the constructed Kuranishi charts on $\mathcal{A}$ in more detail. For this purpose, let $A \in \mathcal{C}(P)$ be an ASD connection. According to Remark 3.2.4, the linear spaces occurring in definition of a Kuranishi chart at $[A] \in \mathcal{A}$ are given by the cohomology groups

$$\begin{aligned} \mathcal{E} &= \mathrm{H}^{1,+}_A(M, \mathrm{Ad}P) \equiv \operatorname{Ker} \mathrm{d}^+_A / \operatorname{Im} \mathrm{d}_A, \\ \mathcal{F} &= \mathrm{H}^{2,+}_A(M, \mathrm{Ad}P) \equiv \Omega^2_+(M, \mathrm{Ad}P) / \operatorname{Im} \mathrm{d}^+_A. \end{aligned} \tag{3.4.19}$$

These spaces are finite-dimensional, because the Yang–Mills complex is elliptic. Moreover, they are endowed with a natural linear action of the compact stabilizer subgroup $\mathcal{G}au_A(P)$ of $A$. Thus, the obstruction map is a $\mathcal{G}au_A(P)$-equivariant map

$$s\colon \mathrm{H}^{1,+}_A(M, \mathrm{Ad}P) \supseteq \mathcal{V} \to \mathrm{H}^{2,+}_A(M, \mathrm{Ad}P), \tag{3.4.20}$$

where $\mathcal{V}$ is a $\mathcal{G}au_A(P)$-invariant, open neighborhood of $0$ in $\mathrm{H}^{1,+}_A(M, \mathrm{Ad}P)$. Finally, the moduli space $\mathcal{A}$ in a neighborhood of $[A]$ is modeled on the quotient $s^{-1}(0)/\mathcal{G}au_A(P)$. In this way, we recover the well-known result [DK97, Proposition 4.2.23] concerning the local structure of $\mathcal{A}$.

Further insights into the geometry of $\mathcal{A}$ require a more precise control of the obstruction map $s$, which is rather difficult to obtain in full generality. However, in concrete examples, one can often find conditions which ensure that $s$ vanishes. For example, if $M$ is self-dual and has positive scalar curvature, then it can be shown using the Weitzenböck formula that $\mathrm{H}^{2,+}_A(M, \mathrm{Ad}P)$ is trivial for every irreducible ASD connection $A$, see e.g. [RS17, Lemma 6.5.4]. Thus, in this case, the moduli space of irreducible anti-self-dual connections is a smooth manifold. This important result was originally obtained by Atiyah, Hitchin, and Singer [AHS78, Theorem 6.1].

For the remainder of this section, we specialize to $G = \mathrm{SU}(2)$. This setting was used by Donaldson [Don83] to gain astounding insights into the topology and geometry of 4-manifolds. For $G = \mathrm{SU}(2)$, we have $\mathrm{p}_1(\mathrm{Ad}P) = 4k$, where $k$ is the so-called instanton number. Let us restrict our attention to the case $k = 1$. Moreover, we assume that $M$ is simply connected and that its intersection form is positive definite. These assumptions imply $\chi_M - \sigma_M = 1$ so that the virtual dimension of $\mathcal{A}$ is $8 - 3 = 5$. For $G = \mathrm{SU}(2)$, a reducible connection $A$ has a stabilizer group conjugate to $\mathrm{U}(1)$ or to $\mathbb{Z}_2$. Connections with a discrete stabilizer subgroup are flat as a consequence of the Ambrose–Singer Theorem and thus only connections with a stabilizer subgroup conjugate to $\mathrm{U}(1)$ are of interest for the geometry of $\mathcal{A}$. By refining the classification of orbit types in Proposition 3.4.3, one can show that there are only finitely many gauge-equivalence classes of ASD connections which are reducible to $\mathrm{U}(1)$, see [RS17, Proposition 6.5.11]. The Kuranishi charts constructed above determine the structure of $\mathcal{A}$ in a neighborhood of a singular point. Let $A$ be an ASD connection that is reducible to $\mathrm{U}(1)$. In the present setting, a straightforward application of the representation theory of $\mathrm{U}(1)$ yields the following isomorphisms of $\mathcal{G}au_A(P)$-representation spaces

$$\mathrm{H}^{1,+}_A(M, \mathrm{Ad}P) \simeq \mathbb{C}^p, \qquad \mathrm{H}^{2,+}_A(M, \mathrm{Ad}P) \simeq \mathbb{C}^q \tag{3.4.21}$$

for some $p$ and $q$ satisfying $p + q = 3$, where $\mathrm{U}(1)$ acts in the usual way on $\mathbb{C}^p$ and $\mathbb{C}^q$, see [FU84, Proposition 4.9]. Thus, if $\mathrm{H}^{2,+}_A(M, \mathrm{Ad}P)$ is trivial, then, in a neighborhood of the singular point $[A]$, the moduli space $\mathcal{A}$ is identified with the cone $\mathbb{C}^3/\mathrm{U}(1)$ over $\mathbb{CP}^2$. We thus recover the following important result of Donaldson [Don83].

Theorem 3.4.5 *Let $P$ be a principal $\mathrm{SU}(2)$-bundle with instanton number $1$ over a $4$-dimensional, compact, simply connected, oriented Riemannian manifold $M$, whose intersection form is positive definite. Assume that $\mathrm{H}^{2,+}_A(M, \mathrm{Ad}P) = \{0\}$ for every ASD connection $A$ on $P$. Then, the moduli space $\mathcal{A}$ of anti-self-dual connections on $P$ is a smooth $5$-dimensional manifold except at finitely many points. At a singular point, $\mathcal{A}$ has the geometry of a cone over $\mathbb{CP}^2$.* ◊

In particular, every singular connection in $\mathcal{A}$ can be approximated by a sequence of irreducible ASD connections. This implies that $\mathcal{A}$ is stratified by orbit types. The case $\mathrm{H}^{2,+}_A(M, \mathrm{Ad}P) \neq \{0\}$ is more complicated and a perturbation of the metric on $M$ is required, see [FU84, Section 4].

Further examples are discussed in [DK97, Section 4.2.6]. We note that, for the above example and for the ones discussed in [DK97], the moduli space of anti-self-dual connections is stratified by orbit types. We do not know whether this holds true in general.

# Singular Symplectic Reduction 4

This chapter focuses on symmetry reduction of infinite-dimensional Hamiltonian dynamical systems. In the modern geometric formulation of Hamiltonian systems, the phase space of the system is modeled as a symplectic manifold $(M, \omega)$ and symmetries are represented by the action of a Lie group $G$ on $M$ preserving the symplectic structure. The conserved quantities corresponding to the $G$-symmetry are encoded in the momentum map $J: M \to \mathfrak{g}^*$, which takes values in the dual space of the Lie algebra of $G$. The Marsden-Weinstein or symplectically reduced phase space at $\mu \in \mathfrak{g}^*$ is, by definition, the space

$$\check{M}_\mu = J^{-1}(\mu)/G_\mu, \tag{4.0.1}$$

where $G_\mu$ denotes the stabilizer subgroup of $\mu$ under the coadjoint action of $G$ on $\mathfrak{g}^*$. It is clear that this setting fits into the general framework of abstract moduli spaces studied in the previous chapter. In finite dimensions, the structure of $\check{M}_\mu$ is well understood: If $\mu$ is a regular value and the $G$-action is free and proper, then $\check{M}_\mu$ is a manifold and the symplectic form $\omega$ on $M$ descends to a symplectic structure $\check{\omega}$ on $\check{M}_\mu$. Moreover, the dynamics associated to a $G$-invariant Hamiltonian projects down to a Hamiltonian flow on $\check{M}_\mu$. In general, without assuming these regularity conditions, the space $\check{M}_\mu$ has singularities and is stratified by symplectic manifolds. The regular case is discussed in [Mey73; MW74] and the singular case is studied in [AGJ90; SL91]. The aim of this chapter is to extend these results concerning singular symplectic reduction to the infinite-dimensional setting using the results of Chapter 3 about normal forms of equivariant maps. We begin with the linear case of symplectic vector spaces in Section 4.1, where we discuss under which conditions many well-known results of linear symplectic geometry generalize to weakly symplectic vector spaces. As a fundamental tool, we introduce and study a class of topologies associated to the symplectic form. Next, in Section 4.2, we discuss symplectic manifolds and momentum maps in our infinite-dimensional setting. Based on joint work with T. Ratiu [DR] on actions of diffeomorphism groups, we introduce the notion of a group-valued momentum map which unifies several other notions of generalized momentum maps. The group structure of the target allows to encode discrete topological information; a fact that is especially relevant for the action of geometric automorphism groups, which are sensitive to the topology of the spaces they live on. With the help of a refined normal form result for momentum maps, in Section 4.3, we prove a theorem on singular symplectic reduction in our infinite-dimensional setting. To our knowledge, the results about the normal form of a momentum map and about the structure of singular symplectic quotients are new even for symplectic Banach manifolds[1]. Finally, we discuss the example of symplectic reduction in the context of the Yang–Mills equation over a Riemannian surface, which is based on joint work with J. Huebschmann [DH18].

## 4.1 Symplectic functional analysis

In the Banach space setting, the theory of symplectic geometry splits into two branches depending on whether the map $\omega^\flat \colon \mathrm{T}M \to \mathrm{T}^* M$ induced by the symplectic form $\omega$ on $M$ is an isomorphism or merely an injection. The former forms are called strongly symplectic and the latter are referred to as weakly symplectic. The well-consolidated building of finite-dimensional symplectic geometry generalizes almost without changes to strongly symplectic structures, but it is confronted with serious problems if weakly symplectic forms are considered. For example, the Darboux theorem holds for strongly symplectic forms but fails for weak ones, see [Mar72].

In the case of a Fréchet manifold, a 2-form cannot be strongly symplectic, because the dual of a Fréchet space is never a Fréchet space (except when it is a Banach space). Thus, for us, weakly symplectic forms are the norm rather than the exception and we need to find a way to address the problems that originate from the failure of $\omega^\flat$ to be surjective. At the root of our approach lies the observation that the symplectic form induces a natural topology on the tangent spaces. For strongly symplectic forms, this topology is equivalent to the one induced from the manifold topology; while for weakly symplectic forms the manifold topology is finer. With respect to the symplectic topology, the symplectic form behaves as if it were strongly symplectic and thus many standard results of the finite-dimensional setting carry over to statements relative to this topology. Then, tools from functional analysis, especially the theory of dual pairs, can be used to compare the original and the symplectic topology and in this way allow to transfer these conclusions to results relative to the original topology. From a different angle, we translate algebraic questions into a functional analytic setting and then employ the well-oiled machine of dual pairs to solve these problems. We are not aware of any previous systematic application of the theory of dual pairs to symplectic vector spaces, though some elements can be found in [BZ18]. For the convenience of the reader we summarize the relevant material concerning dual pairs in Appendix B.

Definition 4.1.1 A *symplectic vector space* $(X, \omega)$ is a locally convex vector space $X$ endowed with a jointly continuous, antisymmetric bilinear form $\omega$ which is

[1] The theory of regular symplectic reduction for strongly symplectic Banach manifolds was already discussed in [MW74].

non-degenerate in the sense that the induced map

$$\omega^\flat : X \to X', \quad x \mapsto \omega(x, \cdot) \tag{4.1.1}$$

is injective, where $X'$ denotes the topological dual of $X$. If $\omega^\flat$ is a bijection, then $\omega$ is called *strongly symplectic*. ◊

In the literature, symplectic forms $\omega$ for which $\omega^\flat$ is only injective are often called weakly symplectic. Since, however, a genuine Fréchet space is never isomorphic to its topological dual, there are no strongly symplectic forms beyond the Banach setting and the weak case is the generic one for us.

Example 4.1.2 Consider the Fréchet space $X = \Omega^1(M)$ of differential one-forms on a closed two-dimensional surface $M$. The integration pairing yields a symplectic form $\omega$ on $\Omega^1(M)$ by setting

$$\omega(\alpha, \beta) = \int_M \alpha \wedge \beta \tag{4.1.2}$$

for $\alpha, \beta \in \Omega^1(M)$. Indeed, given a Riemannian metric $g$ on $M$ with associated Hodge star operator $*$, the value

$$\omega_A(\alpha, * \alpha) = \int_M \|\alpha\|_g^2 \tag{4.1.3}$$

is positive for every non-vanishing $\alpha \in \Omega^1(M)$ and thus $\omega$ is non-degenerate. Moreover, the symplectic structure $\omega$ is only weakly non-degenerate, because the image of $\omega^\flat$ consists of regular functionals and not of all distributional 1-forms.

This example will reappear throughout the next sections serving as an illustration of the abstract theory. The series consists of Examples 4.2.2, 4.2.6, 4.3.6, 4.2.9, 4.1.10 and 4.2.28. ◊

If we ignore for a moment that $X$ already carries a locally convex topology, we are left with a vector space endowed with a bilinear form $\omega$. This setting is well-studied in functional analysis, where such a pair $(X, \omega)$ is called a *dual pair*, see Appendix B for background information. The bilinear form $\omega$ singles out certain topologies $\tau$ on $X$ for which it is strong, i.e., every $\tau$-continuous functional on $X$ is of the form $\omega(x, \cdot)$ for some $x \in X$. In the general theory such topologies are called compatible with the dual pair $(X, \omega)$. In the following, we will use the notation $(X, \tau)'$ for the space of all $\tau$-continuous functionals on $X$ to emphasize the dependence on the topology $\tau$ of $X$.

Definition 4.1.3 Let $(X, \omega)$ be a symplectic vector space. A locally convex topology $\tau$ on $X$ is called *compatible* with $\omega$ or a *symplectic topology* if $(X, \tau)' = \operatorname{Im} \omega^\flat$. ◊

For a general dual pair, the so-called weak, strong, and Mackey topologies play an important role, see Example B.1.5. In the context of a symplectic vector space $(X, \omega)$, we put the prefix "symplectic" in front. For example, the *weak symplectic topology*, denoted by $\sigma_\omega$, is determined by the seminorms

$$\|x\|_A := \sup_{y \in A} |\omega(y, x)| \tag{4.1.4}$$

indexed by a finite subset $A \subseteq X$. The *symplectic Mackey topology* is defined in a similar way except for the fact that the seminorms are indexed by convex, circled and $\sigma_\omega$-compact subsets $A \subseteq X$. Both the weak symplectic and the symplectic Mackey topology are compatible with the symplectic form. Moreover, we have the following symplectic version of the Mackey–Arens Theorem B.1.6.

Proposition 4.1.4 *A locally convex topology on a symplectic vector space $(X, \omega)$ is symplectic if and only if it is lies between the weak symplectic and the symplectic Mackey topology.* ◇

In order to distinguish symplectic topologies from the original one in which $\omega$ is jointly continuous, we call the latter the *original topology*. Note that, for a weakly symplectic form, a symplectic topology is strictly coarser than the original topology, because the latter has more continuous functionals.

Proposition 4.1.5 *Let $(X, \omega)$ be a symplectic vector space. The original topology on $X$ is symplectic if and only if $X$ is a normed space and $\omega$ is strongly symplectic.* ◇

*Proof.* For a strongly symplectic form, the original topology is clearly symplectic. Conversely, suppose that the original topology is symplectic. Then, by definition, $\omega^\flat : X \to X'$ is a bijection and it remains to show that $X$ is a normed space. For this purpose, endow $X'$ with the topology $\tau_\omega$ by declaring $\omega^\flat$ to be a homeomorphism. Accordingly, the inverse $\omega^\sharp : X' \to X$ of $\omega^\flat$ is a continuous linear map with respect to the topology $\tau_\omega$ on $X'$ and the original topology on $X$. Note that the canonical evaluation map $X' \times X \to \mathbb{R}$ can be written as

$$(\alpha, x) \mapsto \alpha(x) = \omega(\omega^\sharp(\alpha), x). \tag{4.1.5}$$

Thus, it is jointly continuous with respect to $\tau_\omega$ on $X'$ and the original topology on $X$. However, the evaluation map is only jointly continuous (for some vector space topology on $X'$) if $X$ is normable according to [Mai63]. □

In other words, the difference between the original topology and a symplectic topology is a measure of how much $\omega^\flat : X \to X'$ fails to be surjective.

By definition, symplectic topologies are closely tied to the symplectic form and thus they are often able to detect symplectic phenomena, which are hard to describe in the original topology. As we will see now, problems involving subspaces of symplectic spaces can be conveniently dealt within the framework of symplectic topologies. A first hint of this interplay can be gathered from the fact that all symplectic topologies have the same closed linear subspaces according to Proposition B.1.7. This observation allows us to introduce the following notion.

Definition 4.1.6 Let $(X, \omega)$ be a symplectic vector space. A linear subspace $V \subseteq X$ is called *symplectically closed*, or simply $\omega$-closed, if it is closed with respect to some (and hence all) symplectic topology on $X$. In a similar vein, $V$ is said to be *symplectically dense* if it is dense relative to a symplectic topology. ◇

Let $(X, \omega)$ be a symplectic vector space and let $V \subseteq X$ be a linear subspace. The *symplectic orthogonal* of $V$ is defined by

$$V^{\omega} := \{x \in X : \ \omega(x, v) = 0 \text{ for all } v \in V\}. \tag{4.1.6}$$

Note that $V^{\omega}$ is the symplectic counterpart of the polar (or annihilator) of $V$ in the general theory of dual pairs, see (B.1.2). Hence, the results about polars apply, in particular, to symplectic orthogonals as well. The Bipolar Theorem B.1.8 yields the following description of symplectic double orthogonals and thereby generalizes previous results of Booß-Bavnbek and Zhu [BZ18, Lemma 1.4].

Proposition 4.1.7 *Let $(X, \omega)$ be a symplectic vector space. For every linear subspace $V \subseteq X$, the symplectic double orthogonal $V^{\omega\omega}$ coincides with the closure of $V$ with respect to a symplectic topology. In particular, $V$ is symplectically closed if and only if $V^{\omega\omega} = V$.* ◇

This results in a shift of perspective as we may complement the algebraic approach to symplectic double orthogonals by powerful tools from topology. Illustrating this shift in philosophy, the following result is almost trivial from a topological point of view but not so straightforward to prove in an algebraic fashion.

Lemma 4.1.8 *Every finite-dimensional subspace of a symplectic vector space is symplectically closed.* ◇

*Proof.* Every finite-dimensional vector space $V \subseteq X$ is closed for whatever vector space topology we put on $X$, see [Köt83, Proposition 15.5.2]. □

In the literature, it is often assumed or even "proven" that in weakly symplectic Banach spaces every linear subspace $V$ which is closed with respect to the original topology on $X$ satisfies $V^{\omega\omega} = V$, see for example [Kob87, Lemma 7.5.9] or [Bam99, Lemma 3.2]. Extending a counterexample of Booß-Bavnbek and Zhu [BZ18, Example 1.6] we can, however, show that every genuinely weakly symplectic space has at least one closed subspace that is not symplectically closed. The following result is analogous to the fact that the dual of a non-reflexive Banach space contains subspaces that are both norm-closed and weak*-dense, see [DeV78, Fact 4.1.6].

Proposition 4.1.9 *Let $(X, \omega)$ be a symplectic vector space. If $\omega$ is not strongly symplectic, then there exists a closed proper linear subspace that is symplectically dense. In particular, every closed linear subspace of $X$ is symplectically closed if and only if $X$ is a normed space and $\omega$ is strongly symplectic.* ◊

*Proof.* Endow $X'$ with the weak topology determined by the dual pair $(X', X)$, see Example B.1.5. Suppose that $\omega$ is not strongly non-degenerate. Then, there exists a non-zero continuous linear functional $\alpha \in X' \setminus \operatorname{Im} \omega^\flat$. The subspace $A := \mathbb{R} \cdot \alpha \subseteq X'$ is finite-dimensional and hence closed, see [Köt83, Proposition 15.5.2]. Consider the subspace $V := A^\circ \subseteq X$, where the polar is taken with respect to the dual pair $(X', X)$. By definition, $V$ coincides with the kernel of $\alpha$ and thus it is a proper closed subspace of $X$. The Bipolar Theorem B.1.8 with respect to the dual pair $(X', X)$ and closedness of $A$ imply $V^\circ = A^{\circ\circ} = A = \mathbb{R} \cdot \alpha$. Since $\alpha$ is not contained in the image of $\omega^\flat$, we have $(\omega^\flat)^{-1}(V^\circ) = \{0\}$. As in Example B.2.2, we can view $\omega^\flat \colon X \to X'$ as a weakly continuous map between the dual pairs $\omega(X, X)$ and $(X', X)$ with the identity on $X$ as the adjoint map. Hence, Proposition B.2.3 (i) implies

$$V^\omega = (\mathrm{id}(V))^\omega = (\omega^\flat)^{-1}(V^\circ) = \{0\}. \tag{4.1.7}$$

Thus, $V^{\omega\omega} = X$. In summary, every genuinely weakly symplectic space has at least one closed proper subspace whose symplectic closure is the whole space. In other words, if every closed subspace is symplectically closed then the symplectic form has to be strongly symplectic. The latter is only possible if $X$ is normable according to Proposition 4.1.5. Conversely, if $(X, \omega)$ is a strongly symplectic space, then the original topology on $X$ is symplectic and thus every closed subspace is also symplectically closed as a consequence of Proposition B.1.7. □

Example 4.1.10 Continuing Example 4.1.2. Every tangent vector $X_m \in \mathrm{T}_m M$ yields a continuous Dirac-like functional $\delta_{X_m} \colon \Omega^1(M) \to \mathbb{R}$ by evaluation of a 1-form on $X_m$. Note that $\delta_{X_m}$ is singular in the sense of distributions and thus does *not* lie in the image of $\omega^\flat$. Consider the closed subspace

$$V := (\operatorname{span} \delta_{X_m})^\circ = \{\alpha \in \Omega^1(M) : \ \alpha(X_m) = 0\}. \tag{4.1.8}$$

Since integration is not sensitive to the behavior at a single point, we find $V^\omega = \{0\}$ and thus $V^{\omega\omega} = \Omega^1(M)$. In summary, $V$ is a closed but symplectically dense subspace. ◊

In their original article on symplectic reduction, Marsden and Weinstein [MW74] considered symplectic forms $\omega$ on reflexive Banach spaces whose associated musical isomorphism $\omega^\flat$ has a closed image. This setting, a priori, lies between weakly and strongly symplectic Banach spaces. However, under these assumptions, they showed that every closed subspace is also symplectically closed, see [MW74, Lemma on p. 123]. Hence, Proposition 4.1.9 implies that such symplectic forms are automatically strongly non-degenerate. Let us record this observation.

Proposition 4.1.11 *Let $X$ be a reflexive Banach space endowed with a symplectic form $\omega$. If $\omega^\flat$ has closed image in $X'$, then $\omega$ is a strongly symplectic form.* ◊

The restriction $\omega_V$ of the symplectic structure $\omega$ to a closed subspace $V \subseteq X$ is in general degenerate, the kernel of $\omega_V^\flat : V \to V'$ being $V \cap V^\omega$.

Definition 4.1.12 A closed subspace $V$ of a symplectic vector space $(X, \omega)$ is called *symplectic* if $V \cap V^\omega = \{0\}$. ◊

Accordingly, the restriction $\omega_V$ of $\omega$ to a symplectic subspace $V$ yields a symplectic form on $V$. We emphasize that the notions "symplectically closed subspace" and "symplectic subspace" should not be confused.

Example 4.1.13 Let $(X, \omega)$ be a symplectic vector space. Assume that $\omega$ is not strongly symplectic and consider a closed, symplectically dense, proper subspace $V \subset X$ (which always exists according to Proposition 4.1.9). For every $x \in V^\omega$, the functional $\omega(x, \cdot)$ on $X$ is continuous with respect to a symplectic topology on $X$ but it also vanishes on the symplectically dense subspace $V$. Thus, it has to vanish on the whole space $X$, which implies $x = 0$. Hence, $V^\omega = \{0\}$ and so $V \cap V^\omega = \{0\}$. In other words, every symplectically dense subspace is a symplectic subspace. ◊

In finite dimensions, every symplectic subspace $V \subseteq X$ induces a direct sum decomposition $X = V \oplus V^\omega$. The previous example shows that this is no longer the case in infinite dimensions (for a proper symplectically dense subspace $V$, the subspace $V \oplus V^\omega = V$ is a proper subspace of $X$). The following phenomenon, peculiar for the infinite-dimensional setting, is related. Given a subspace $V \subseteq X$ and $v \in V$, the functional $\omega(v, \cdot)$ on $X$ is continuous with respect to a symplectic topology $\tau_\omega$ on $X$. Hence, its restriction to a functional on $V$ is continuous relative to the subspace topology, the latter being denoted by $\tau_\omega$ as well. This furnishes a linear map

$$\Gamma_V : V \mapsto (V, \tau_\omega)', \qquad v \mapsto \omega(v, \cdot). \tag{4.1.9}$$

Clearly, $\Gamma_V$ is injective if and only if $V$ is a symplectic subspace. In finite dimensions, a count of dimensions shows that $\Gamma_V$ is bijective if $V$ is a symplectic subspace. However, when passing to the infinite-dimensional setting, the above example of a symplectically dense subspace $V \subseteq X$ shows that $\Gamma_V$ is in general not surjective (for every non-zero $y \in X \setminus V$, the $\tau_\omega$-continuous functional $\omega(y, \cdot)$ on $V$ is not representable by some $v \in V$).

Proposition 4.1.14 *Let $(X, \omega)$ be a symplectic vector space. For a closed subspace $V \subseteq X$, the following are equivalent:*

*(i) $V$ is a symplectic subspace and the restriction of every symplectic topology on $X$ yields a symplectic topology on $(V, \omega_V)$.*

*(ii) There exists an algebraic direct sum decomposition $X = V \oplus V^\omega$.*

*(iii) The linear map $\Gamma_V$ defined in (4.1.9) is surjective.* ◊

*Proof.* The implication (i) → (iii) is clear. Now assume that $\Gamma_V$ is surjective. Then, we have

$$\begin{aligned} V \cap V^\omega &= \{x \in V : \Gamma_V(v)(x) = \omega(v, x) = 0 \text{ for all } v \in V\} \\ &= \{x \in V : \alpha(x) = 0 \text{ for all } \alpha \in (V, \tau_\omega)'\} \\ &= \{0\}, \end{aligned} \tag{4.1.10}$$

where the last equality is a consequence of the Hahn–Banach Theorem [Köt83, Proposition 20.1.2], which implies that $\tau_\omega$-continuous functionals on $V$ separate points of $V$; that is, for each non-zero $x \in V$ there exists $\alpha \in (V, \tau_\omega)'$ such that $\alpha(x) = 1$. Moreover, by surjectivity of $\Gamma_V$, for every $x \in X$, there exists $v \in V$ such that the functionals $\omega(x, \cdot)$ and $\omega(v, \cdot)$ coincide on $V$. Thus, $\omega(x - v, \cdot)$ vanishes on $V$, that is, $x - v \in V^\omega$. Hence, $X$ is the algebraic direct sum of $V$ and $V^\omega$, which establishes the implication (iii) → (ii). For the last implication (ii) → (i), suppose now that $X = V \oplus V^\omega$ is an algebraic direct sum. Then, $V$ is a symplectic subspace, because the sum is direct. Moreover, the Hahn–Banach Theorem [Köt83, Proposition 20.1.1] implies that the restriction map $(X, \tau_\omega)' \to (V, \tau_\omega)'$ is surjective. Hence, using the definition of the symplectic topology, every $\tau_\omega$-continuous functional on $V$ is obtained as the restriction to $V$ of $\omega(x, \cdot)$ for some $x \in X$. Write $x$ as $x = v + w$ with $v \in V$ and $w \in V^\omega$. Then, the restrictions to $V$ of $\omega(x, \cdot)$ and $\omega(v, \cdot)$ coincide. In other words, every $\tau_\omega$-continuous functional on $V$ is of the form $\omega(v, \cdot)$ for some $v \in V$, which completes the proof. □

Finally, we come to what can be considered the linear toy example of symplectic reduction. Later on, the non-linear case will be reduced to this simple setting. Assume that a compact Lie group $G$ acts continuously and linearly on a symplectic vector space $(X, \omega)$. We say that the symplectic form $\omega$ is preserved by the $G$-action if

$$\omega(g \cdot x, g \cdot y) = \omega(x, y) \tag{4.1.11}$$

holds for all $g \in G$ and $x, y \in X$.

Proposition 4.1.15 *Let $(X, \omega)$ be a symplectic vector space and let a compact Lie group $G$ act continuously and linearly on $X$ in such a way that the symplectic form $\omega$ is preserved. Then, the subspace $X_G$ of $G$-invariant elements is a symplectic subspace of $X$ and $X = (X_G) \oplus (X_G)^\omega$ is an algebraic direct sum.* ◊

*Proof.* Let $\tau_\omega$ be a symplectic topology on $X$. We will denote the induced subspace topology on $X_G$ by $\tau_\omega$ as well. According to Proposition 4.1.14, we have to show that the linear map

$$\Gamma_{X_G} \colon X_G \mapsto (X_G, \tau_\omega)', \qquad v \mapsto \omega(v, \cdot) \tag{4.1.12}$$

is surjective. For this purpose, let $\alpha \in (X_G, \tau_\omega)'$. The functional $\alpha$ can be extended to a $\tau_\omega$-continuous functional $\bar{\alpha}$ defined on the whole of $X$ according to the Hahn–Banach theorem [Köt83, Proposition 20.1.1]. By taking the average over the compact Lie group $G$, we may assume that the extension $\bar{\alpha}$ is $G$-invariant. Since $\bar{\alpha}$ is $\tau_\omega$-continuous, there exists $y \in X$ such that $\bar{\alpha} = \omega(y, .)$. As $\bar{\alpha}$ and $\omega$ are $G$-invariant, $y$ has to be an element of $X_G$. Clearly, $\Gamma(y) = \alpha$, which shows that $\Gamma$ is surjective. □

## 4.2 Symplectic manifolds and momentum maps

In this section, we will introduce the notion of a symplectic structure on an infinite-dimensional manifold. While some parts will be devoted to certain general properties of symplectic manifolds, our main focus lies on the momentum map geometry. Although this topic is well-studied in finite dimensions and is the subject of many textbooks, we are not aware of any previous systematic treatments of the infinite-dimensional case, especially beyond the Banach real. Some aspects of the theory of Hamiltonian dynamics on symplectic Banach manifolds can be found in [CM74; RM99; MRA02].

### 4.2.1 Group-valued momentum maps

Much of the content of this subsection represents joint work with T. Ratiu and will be published as part of [DR].

Definition 4.2.1 Let $M$ be a smooth manifold. A differential 2-form $\omega$ on $M$ is called a *symplectic form* if it closed and, for every $m \in M$, the induced bilinear form $\omega_m \colon \mathrm{T}_m M \times \mathrm{T}_m M \to \mathbb{R}$ is a symplectic structure on $\mathrm{T}_m M$ in the sense of Definition 4.1.1. ◊

We remind the reader that, for a symplectic form $\omega$ on $M$, the associated map

$$\omega^\flat_m \colon \mathrm{T}_m M \to (\mathrm{T}_m M)', \qquad v \mapsto \omega_m(v, \cdot) \tag{4.2.1}$$

is required to be injective for all $m \in M$. If $\omega^\flat_m$ is a topological isomorphism, then we say that $\omega$ is a strongly symplectic form. According to Proposition 4.1.5, this may be the case only when $M$ is a Banach manifold.

Example 4.2.2 Continuing in the series of Example 4.1.2. Let $P \to M$ be a $\mathrm{U}(1)$-bundle on a closed surface $M$. The space $\mathcal{C}(P)$ of connections on $P$ is an affine space modeled on the Fréchet space $\Omega^1(M)$. Generalizing Example 4.1.2, the 2-form $\omega$ on $\mathcal{C}(P)$ defined by the integration pairing

$$\omega_A(\alpha, \beta) = \int_M \alpha \wedge \beta \tag{4.2.2}$$

for $A \in \mathcal{C}(P)$ and $\alpha, \beta \in \Omega^1(M)$ is a symplectic form. Indeed, $\omega$ is closed, because the right-hand side of (4.2.2) is independent of $A$, and non-degeneracy follows from the same arguments as in Example 4.1.2. ◊

Let $(M, \omega)$ be a symplectic manifold and let $\Upsilon$ be an action of a Lie group $G$ on $M$ by symplectic diffeomorphism, i.e. $\Upsilon_g^* \omega = \omega$ for all $g \in G$. We will refer to this setting by saying that $(M, \omega)$ is a *symplectic G-manifold*. In finite dimensions, the conserved quantities corresponding to the $G$-symmetry of the system are encoded in the momentum map, which is a map from $M$ to the dual space of the Lie algebra $\mathfrak{g}$ of $G$. To make sense of the notion of a momentum map in our infinite-dimensional manifold, we first need to explain what we understand by the dual space of the Lie algebra. Guided by the theory of dual pairs, we say that a *dual space of* $\mathfrak{g}$ is a locally convex vector space $\mathfrak{h}$, which is in duality with $\mathfrak{g}$ through a given weakly non-degenerate jointly continuous bilinear map $\kappa \colon \mathfrak{h} \times \mathfrak{g} \to \mathbb{R}$. Using notation stemming from functional analysis, we often write the dual pair as $\kappa(\mathfrak{h}, \mathfrak{g})$. Intuitively, we think of $\mathfrak{h}$ as the dual vector space of $\mathfrak{g}$ and, for this reason, we often write $\mathfrak{g}^* := \mathfrak{h}$, even though $\mathfrak{g}^*$ is not necessarily the (topological) dual of $\mathfrak{g}$.

Definition 4.2.3 Let $(M, \omega)$ be a symplectic $G$-manifold and let $\kappa(\mathfrak{g}^*, \mathfrak{g})$ be a dual pair. A map $J \colon M \to \mathfrak{g}^*$ is called a *momentum map* if

$$\xi^* \lrcorner\, \omega + \kappa(\mathrm{d}J, \xi) = 0 \tag{4.2.3}$$

holds for all $\xi \in \mathfrak{g}$, where $\xi^*$ denotes the fundamental vector field induced by the action of the Lie algebra $\mathfrak{g}$ on $M$. ◊

For every $\xi \in \mathfrak{g}$, the map $J_\xi \colon M \to \mathbb{R}$ defined by $J_\xi(m) = \kappa(J(m), \xi)$ is called the $\xi$-component of $J$.

Example 4.2.4 (As a generalization of Lie algebra-valued momentum maps) A jointly continuous, non-degenerate, symmetric bilinear form $\kappa \colon \mathfrak{g} \times \mathfrak{g} \to \mathbb{R}$ identifies the dual $\mathfrak{g}^*$ with $\mathfrak{g}$. This leads to the concept of a Lie algebra-valued momentum map. Although this notion can be found in the literature since the mid-seventies, it was recently formalized by Neeb, Sahlmann, and Thiemann [NST14, Definition 4.3]: a Lie algebra-valued momentum map is a smooth map $J \colon M \to \mathfrak{g}$ such that, for all $\xi \in \mathfrak{g}$, the component functions $J_\xi \colon M \to \mathbb{R}$ defined by $J_\xi(m) = \kappa(J(m), \xi)$ for $m \in M$ satisfy

$$\xi^* \lrcorner\, \omega + \mathrm{d}J_\xi = 0. \tag{4.2.4}$$

It is immediately clear from the definition that such a Lie algebra-valued momentum map can be regarded as a $\mathfrak{g}$-valued momentum map with respect to the dual pair $\kappa(\mathfrak{g}, \mathfrak{g})$. ◊

Example 4.2.5 Let $(X, \omega)$ be a symplectic vector space endowed with a continuous linear action of the compact finite-dimensional Lie group $G$. Assume that the symplectic form is preserved by the action, which in the linear setting simply means that

$$\omega(g \cdot x, g \cdot y) = \omega(x, y) \tag{4.2.5}$$

for all $x, y \in X$ and $g \in G$. As a consequence, the Lie algebra action of $\mathfrak{g}$ is skew-symmetric with respect to $\omega$. Although $\mathfrak{g}$ is finite-dimensional and thus there is no ambiguity of its dual space, we continue to use the dual pair notation $\kappa(\mathfrak{g}^*, \mathfrak{g})$. A straightforward calculation shows that the quadratic map $J\colon X \to \mathfrak{g}^*$ defined by

$$\kappa(J(x), \xi) = \frac{1}{2}\omega(x, \xi \,.\, x), \tag{4.2.6}$$

for $\xi \in \mathfrak{g}$, is a momentum map for the $G$-action. ◊

Example 4.2.6 Continuing in the setting of Example 4.2.2, let $P \to M$ be a principal U(1)-bundle on the closed surface $M$. The group $\mathcal{G}au(P)$ of gauge transformations of $P$ is identified with the space $\mathrm{C}^\infty(M, \mathrm{U}(1))$ and thus is a Fréchet Lie group with Lie algebra $\mathfrak{gau}(P) = \mathrm{C}^\infty(M)$. The natural pairing

$$\kappa(\alpha, \phi) = \int_M \phi\, \alpha, \tag{4.2.7}$$

for $\alpha \in \Omega^2(M)$ and $\phi \in \mathfrak{gau}(P)$, identifies $\Omega^2(M)$ as the dual of $\mathfrak{gau}(P)$.

The gauge group $\mathcal{G}au(P)$ acts on $\mathcal{C}(P)$ via gauge transformations

$$\mathcal{G}au(P) \times \mathcal{C}(P) \to \mathcal{C}(P), \qquad (\lambda, A) \mapsto A - \delta^R \lambda, \tag{4.2.8}$$

where $\delta^R \lambda \in \Omega^1(M)$ denotes the right logarithmic derivative of $\lambda \in \mathrm{C}^\infty(M, \mathrm{U}(1))$ defined by $\delta^R \lambda(v) = \mathrm{T}_m \lambda(v) \,.\, \lambda^{-1}(m)$ for $v \in \mathrm{T}_m M$. The infinitesimal action of $\mathfrak{gau}(P)$ coincides with minus the exterior differential $\mathrm{d}\colon \mathrm{C}^\infty(M, \mathbb{R}) \to \Omega^1(M)$. The curvature map

$$\mathcal{J}\colon \mathcal{C}(P) \to \Omega^2(M), \qquad A \mapsto -F_A \tag{4.2.9}$$

is the momentum map for this action with respect to the pairing $\kappa$ introduced above. Indeed, we have $\mathrm{T}_A \mathcal{J} = -\mathrm{d}$ as $F_{A+\alpha} = F_A + \mathrm{d}\alpha$ for every $\alpha \in \Omega^1(M)$,

and the calculation

$$\begin{aligned} \omega(\psi \,.\, A, \alpha) &= \omega(-\,\mathrm{d}\psi, \alpha) \\ &= -\int_M \mathrm{d}\psi \wedge \alpha = \int_M \psi \wedge \mathrm{d}\alpha \\ &= \kappa(\mathrm{d}\alpha, \psi) = -\kappa(\mathrm{T}_A \mathcal{J}(\alpha), \psi) \end{aligned} \tag{4.2.10}$$

for $\psi \in \mathrm{C}^\infty(M)$ and $\alpha \in \Omega^1(M)$ verifies the momentum map relation (4.2.3). ◊

Even in finite dimensions, a momentum map for a symplectic action does not need to exist. This is the case if the 1-form $\xi^* \lrcorner\, \omega$ is only closed and not exact for some $\xi \in \mathfrak{g}$. Since $\xi^* \lrcorner\, \omega$ is closed, it gives rise to the period homomorphism

$$\mathrm{per}_\xi \colon \mathrm{H}^1(M, \mathbb{Z}) \to \mathbb{R}, \qquad [\gamma] \mapsto \int_\gamma (\xi^* \lrcorner\, \omega), \tag{4.2.11}$$

where $\gamma$ is a closed curve in $M$. One could argue that the topological data encoded in the period homomorphism is conserved by the action and such conservation laws should be encoded in the momentum map as well, cf. Proposition 4.2.11 below. However, the classical momentum map takes values in a continuous vector space and thus there is no space to store (discrete) topological information. In order to capture this additional data, we introduced in [DR] the concept of a group-valued momentum map. At the heart of this generalized notion of a momentum map lies the observation that if the periods of $\xi^* \lrcorner\, \omega$ are integral, i.e. $\mathrm{per}_\xi([\gamma]) \in \mathbb{Z}$ for all closed curves $\gamma$ in $M$, then there exists a smooth map $J_\xi \colon M \to \mathrm{U}(1)$ such that the left logarithmic derivative of $J_\xi$ equals $\xi^* \lrcorner\, \omega$. The U(1)-valued component function $J_\xi$ can be viewed as a generalized primitive of $\xi^* \lrcorner\, \omega$ in extension of the real-valued maps components of an ordinary momentum map. The additional topological information is encoded in the winding number of $J_\xi$.

In order to formalize these ideas, we say that two Lie groups $G$ and $H$ are *dual* to each other if there exists a non-degenerate bilinear form $\kappa \colon \mathfrak{h} \times \mathfrak{g} \to \mathbb{R}$ relative to which the associated Lie algebras are in duality. We use the notation $\kappa(H, G)$ in this case. As for Lie algebras, we often write $G^* := H$, intuitively thinking of $G^*$ as the dual group, as in the theory of Poisson Lie groups [LW90]. In this thesis, we assume, for simplicity, that $G^*$ is abelian and write its group operation as addition and use $0 \in G^*$ for the identity element. The reader is referred to [DR] for the general non-abelian case.

Definition 4.2.7 Let $(M, \omega)$ be symplectic $G$-manifold. A *group-valued momentum map* is a tuple $(J, \kappa)$, where $\kappa(G^*, G)$ is a dual pair of Lie groups with abelian $G^*$ and $J \colon M \to G^*$ is a smooth map satisfying

$$\xi^* \lrcorner\, \omega + \kappa(\delta J, \xi) = 0 \tag{4.2.12}$$

for all $\xi \in \mathfrak{g}$, where $\delta J \in \Omega^1(M, \mathfrak{g}^*)$ denotes the left logarithmic derivative of $J$ defined by $\delta J(v) = J(m)^{-1} . \mathrm{T}_m J(v)$ for $v \in \mathrm{T}_m M$. ◇

Example 4.2.8 Consider the the action of U(1) on the torus $T^2 = \mathrm{U}(1) \times \mathrm{U}(1)$ by multiplication in the first factor. This action is symplectic with respect to the natural volume form on $T^2$, but it does not have a classical momentum map. On the other hand, the projection onto the second factor yields a U(1)-valued momentum map $J \colon T^2 \to \mathrm{U}(1)$. ◇

Example 4.2.9 Continuing with the setting of Example 4.2.6, let $P \to M$ be a principal U(1)-bundle on the compact connected surface $M$. Choose a point $m_0 \in M$. The evaluation $\mathrm{ev}_{m_0} \colon \mathcal{G}au(P) \to \mathrm{U}(1)$ at $m_0$ is a morphism of Lie groups. The kernel of $\mathrm{ev}_{m_0}$ is called the group of pointed gauge transformations and is denoted by $\mathcal{G}au_{m_0}(P)$. This group is a normal, locally exponential Lie subgroup of $\mathcal{G}au(P)$ due to [Nee06, Proposition IV.3.4]. The Lie algebra $\mathfrak{gau}_{m_0}(P)$ of $\mathcal{G}au_{m_0}(P)$ consists of all $\phi \in \mathfrak{gau}(P)$ vanishing at $m_0$. The integration pairing[1] suggests $\mathfrak{gau}_{m_0}(P)^* = \mathrm{d}\Omega^1(M)$ as a natural choice for the dual of $\mathfrak{gau}_{m_0}(P)$. However, the momentum map $\mathcal{J}$ for the $\mathcal{G}au(P)$-action defined in (4.2.9) does not yield a $\mathfrak{gau}_{m_0}(P)^*$-valued momentum map for the $\mathcal{G}au_{m_0}(P)$-action, because the curvature $F_A$ of $A \in \mathcal{C}(P)$ is in general not exact. Instead, the curvature is a closed 2-form with integral periods, i.e. $F_A \in \Omega^2_{\mathrm{cl},\mathbb{Z}}(M)$. Thus, we get a well-defined map

$$\mathcal{J}_{m_0} \colon \mathcal{C}(P) \to \Omega^2_{\mathrm{cl},\mathbb{Z}}(M), \qquad A \mapsto -F_A \tag{4.2.13}$$

and it can easily be verified that $\mathcal{J}_{m_0}$ is a $\Omega^2_{\mathrm{cl},\mathbb{Z}}(M)$-valued momentum map for the action of $\mathcal{G}au_{m_0}(P)$. Here, we view the abelian Lie group $\Omega^2_{\mathrm{cl},\mathbb{Z}}(M)$ with group multiplication given by addition as a dual group of $\mathcal{G}au_{m_0}(P)$, because the exact sequence

$$0 \longrightarrow \mathrm{d}\Omega^1(M) \longrightarrow \Omega^2_{\mathrm{cl},\mathbb{Z}}(M) \longrightarrow \mathrm{H}^2(M, \mathbb{Z}) \longrightarrow 0 \tag{4.2.14}$$

identifies the Lie algebra of $\Omega^2_{\mathrm{cl},\mathbb{Z}}(M)$ with $\mathrm{d}\Omega^1(M) = \mathfrak{gau}_{m_0}(P)^*$. The $\Omega^2_{\mathrm{cl},\mathbb{Z}}(M)$-valued momentum map $\mathcal{J}_{m_0}$ remembers the topological type of $P$ in form of the Chern class of $P$:

$$\mathrm{c}(P) = \int_M F_A = -\int_M \mathcal{J}_{m_0}(A) \in \mathbb{Z} \tag{4.2.15}$$

for some $A \in \mathcal{C}(P)$. The group-valued momentum map $\mathcal{J}_{m_0}$ simplifies to a classical momentum map if and only if the bundle $P$ is trivial. In summary,

[1] The map $\phi \mapsto (\phi - \phi(m_0), \phi(m_0))$ gives a topological decomposition $\mathfrak{gau}(P) \simeq \mathfrak{gau}_{m_0}(P) \oplus \mathbb{R}$ which is dual to the Hodge isomorphism $\Omega^2(M) = \mathrm{d}\Omega^1(M) \oplus \mathrm{H}^2(M)$.

for non-trivial bundles $P$, the group-valued momentum map $\mathcal{J}_{m_0}$ carries the topological structure of $P$. ◊

The example of the action of the group of pointed gauge transformations shows that valuable topological information is contained in a group-valued momentum map. This is especially important for actions of diffeomorphism groups, which by their very nature are sensitive to topological properties of the manifold. As the full story goes beyond the scope of this thesis, we give only a brief overview and refer the reader to [DR] for further information.

Example 4.2.10 Let $(M, \mu)$ be a closed finite-dimensional manifold endowed with a volume form $\mu$ and let $(F, \omega)$ be a symplectic manifold. The space $C^\infty(M, F)$ of smooth maps from $M$ to $F$ carries the symplectic form

$$\Omega_\phi(X, Y) = \int_M \omega_{\phi(m)}(X(m), Y(m))\, \mu(m), \tag{4.2.16}$$

where $\phi \in C^\infty(M, F)$ and $X, Y \in T_\phi C^\infty(M, F)$, i.e., $X, Y\colon M \to TF$ satisfy $X(m), Y(m) \in T_{\phi(m)}F$ for all $m \in M$. The natural action by precomposition of the group $\mathcal{D}iff_\mu(M)$ of diffeomorphisms of $M$ preserving the volume form $\mu$ leaves $\Omega$ invariant. If $\omega$ is exact, say with primitive $\vartheta \in \Omega^1(F)$, then the associated momentum map for the $\mathcal{D}iff_\mu(M)$-action is given by

$$\mathcal{J}\colon C^\infty(M, F) \to \mathfrak{X}_\mu(M)^*, \qquad \phi \mapsto \left[\phi^*\vartheta\right], \tag{4.2.17}$$

where the space of volume-preserving vector fields $\mathfrak{X}_\mu(M)$ (the vector fields whose $\mu$-divergence vanishes) is identified with the space of closed $(n-1)$-forms so that $\mathfrak{X}_\mu(M)^* = \Omega^1(M)/d\Omega^0(M)$. More generally, Gay-Balmaz and Vizman [GV12] showed that a (non-equivariant) momentum map exists also for the case when the pull-back of $\omega$ by all maps $\phi \in C^\infty(M, F)$ is exact; for example, this happens when $H^2(M)$ is trivial. However, without topological conditions on $M$ and without exactness of $\omega$, a classical momentum map does not exist. In contrast, our generalized group-valued momentum map no longer takes values in $\mathfrak{X}_\mu(M)^*$, but instead, in the abelian group $\hat{H}^2(M, U(1))$ that parametrizes principal circle bundles with connections modulo gauge equivalence. If $(F, \omega)$ has a prequantum bundle $(L, \vartheta)$, then the group-valued momentum map

$$\mathcal{J}\colon C^\infty(M, F) \to \hat{H}^2(M, U(1)), \qquad \phi \mapsto \phi^*(L, \vartheta) \tag{4.2.18}$$

sends $\phi$ to the pull-back bundle with connection $\phi^*(L, \vartheta)$. We see that no (topological) restrictions have to be made for $M$ and only the integrability condition of the symplectic form $\omega$ is needed for the existence of the group-valued momentum map. In contrast to the classical momentum map, the $\hat{H}^2(M, U(1))$-valued momentum map contains topological information. Indeed, the Chern class of the bundle, as a class in $H^2(M, \mathbb{Z})$, is available from the generalized momentum map. In our simple example, this is just the integral refinement of the closed 2-form $\phi^* \omega$ on $M$.

This additional topological information encoded in the group-valued momentum map has been used in [DR] to refine and extend many well-known geometric constructions. For example, Marsden and Weinstein [MW83] construct Clebsch variables for ideal fluids using an infinite-dimensional symplectic system similar to the one discussed above. It turns out that every vector field represented in those Clebsch variables has vanishing helicity, i.e., such a fluid configuration has trivial topology and no links or knots. The more general framework of group-valued momentum maps allows to construct generalized Clebsch variables for vector fields with integral helicity.

When applied to the space of Lagrangian immersions, the group-valued momentum map recovers the Liouville class as the conserved topological datum. This observation allows to realize moduli spaces of Lagrangian immersions (and modifications thereof) as symplectic quotients.

A wide range of interesting examples with geometric significance are obtained when the target $F$ is a fiber bundle over $M$ with structure group $G$ and typical fiber a symplectic homogeneous space $G/H$. In this case, the space $\mathcal{F}$ of smooth sections of $F$ carries a symplectic form defined in a similar way as in (4.2.16). Note that sections of $F$ correspond to a reduction of the $G$-bundle to $H$. Of special interest is the case where the fiber is $\mathrm{Sp}(2n, \mathbb{R})/\mathrm{U}(n)$, so that points of the symplectic manifold $\mathcal{F}$ correspond to almost complex structures compatible with a given symplectic structure. In this case, the group-valued momentum map for the group of symplectomorphisms assigns to an almost complex structure the anti-canonical bundle. It was already observed by Fujiki [Fuj92] and Donaldson [Don97] that the Hermitian scalar curvature furnishes a classical momentum map for the action of the group of *Hamiltonian* symplectomorphisms. Of course, the Hermitian scalar curvature is the curvature of the anti-canonical bundle. Thus, the group-valued momentum map combines the geometric curvature form with the topological data of the anti-canonical bundle. As a consequence, for the case of a 2-dimensional base manifold, the Teichmüller moduli space with the symplectic Weil–Petersson form can be realized as a symplectic (orbit) reduced space. ◊

Moreover, the concept of a group-valued momentum map generalizes and unifies several other notions such as the circle-valued momentum [PR12, Definition 1] and the cylinder-valued momentum map [CDM88], see [DR] for details.

Despite its general nature, a group-valued momentum map still captures conserved quantities of the dynamical system, i.e., it has the *Noether property* (see [OR03, Definition 4.3.1]). In finite dimensions, every Hamiltonian $h$ on a symplectic manifold $(M, \omega)$ induces a Hamiltonian flow. This no longer holds in an infinite-dimensional context. For one thing, the map $\omega^\flat \colon \mathrm{T}M \to \mathrm{T}^*M$ induced by the symplectic form $\omega$ on $M$ is, in general, only injective and not surjective. Hence, a Hamiltonian vector field $X_h$ associated to a Hamiltonian $h\colon M \to \mathbb{R}$ by the relation

$$X_h \lrcorner\, \omega + \mathrm{d}h = 0 \tag{4.2.19}$$

may not exist[1]. Even if $X_h$ exists, it may not have a unique flow. The construction of a flow requires the solution of an ordinary differential equation on $M$, which a priori is not guaranteed to exist and to be unique in infinite dimensions. Nonetheless, for concrete examples, one can often show that a unique flow exists. For example, in the gauge theory context studied in Section 5.6 below, existence and uniqueness of the Hamiltonian flow is equivalent to the well-posedness of the Cauchy problem for the Yang–Mills-Higgs theory.

Proposition 4.2.11 (Noether's theorem) *Let $(M, \omega)$ be a symplectic $G$-manifold and $\kappa(G, G^*)$ a dual pair of Lie groups. Assume that the $G$-action on $M$ has a $G^*$-valued momentum map $J\colon M \to G^*$. Let $h \in C^\infty(M)$ be a smooth function for which the Hamiltonian vector field $X_h$ exists and has a unique local flow. If $h$ is $G$-invariant, then $J$ is constant along the integral curves of $X_h$.* ◊

*Proof.* Let $\xi \in \mathfrak{g}$ and $m \in M$. Using the defining equation for the momentum map, we have

$$\begin{aligned}\kappa\big((\delta J)_m(X_h), \xi\big) &= -\omega_m\big(\xi^*, X_h\big) = -(\mathrm{d}h)_m(\xi^*) \\ &= -\frac{\mathrm{d}}{\mathrm{d}\varepsilon}\bigg|_0 h(\exp(\varepsilon\xi)\cdot m) = 0\end{aligned} \tag{4.2.20}$$

by $G$-invariance of $h$. Since $\xi \in \mathfrak{g}$ is arbitrary and the pairing $\kappa$ is non-degenerate, we conclude $\delta J(X_h) = 0$. Hence, $J$ is constant along integral curves of $X_h$. □

Fix a dual pair $\kappa(G^*, G)$ of Lie groups. Let $G$ act on the symplectic manifold $(M, \omega)$ and assume that the action has a group-valued momentum map $J\colon M \to G^*$. A natural question to ask is in which sense $J$ is equivariant. Recall that the coadjoint action is defined with respect to the duality pairing $\kappa$ by

$$\kappa(\mathrm{Ad}^*_g \mu, \xi) = \kappa(\mu, \mathrm{Ad}_{g^{-1}} \xi), \tag{4.2.21}$$

for $g \in G$, $\xi \in \mathfrak{g}$ and $\mu \in \mathfrak{g}^*$. In infinite dimensions, this relation only ensures uniqueness of $\mathrm{Ad}^*$ but not its existence, because $\kappa$ is in general only weakly non-degenerate. In the following, we will assume that $\mathrm{Ad}^*$ exists.

Definition 4.2.12 A left action $\Upsilon\colon G \times G^* \to G^*$ of $G$ on $G^*$ is called a *coconjugation action* if it integrates the coadjoint action, that is,

$$\delta_\eta \Upsilon_g(\eta \,.\, \mu) = \mathrm{Ad}^*_g \mu \tag{4.2.22}$$

[1] For example, consider a symplectic vector space $(X, \omega)$. The Hamiltonian vector field associated to a continuous linear functional $h\colon X \to \mathbb{R}$ exists if and only if $h \in X'$ lies in the image of $\omega^\flat\colon X \to X'$.

holds for all $g \in G$, $\eta \in G^*$ and $\mu \in \mathfrak{g}^*$, where $\eta \,.\, \mu \in \mathrm{T}_\eta G^*$ is the derivative at the identity of the (left) translation by $\eta \in G^*$ on $G^*$ in the direction $\mu \in \mathfrak{g}^*$ and $\delta_\eta \Upsilon_g(\eta \,.\, \mu) \in \mathfrak{g}^*$ denotes the (left) logarithmic derivative at $\eta$ of the map $\Upsilon_g \colon G^* \to G^*$ in the direction $\eta \,.\, \mu \in \mathrm{T}_\eta G^*$. If, in addition, $\Upsilon_g(\zeta + \eta) = \Upsilon_g(\zeta) + \Upsilon_g(\eta)$ holds for all $\zeta, \eta \in G^*$, then we say that the coconjugation action is *standard*.

For a given coconjugation action, we say that the group-valued momentum map $J$ is *equivariant* if it is $G$-equivariant as a map $J \colon M \to G^*$. ◊

Example 4.2.13 The coadjoint representation is a standard coconjugation action of $G$ on $G^* := \mathfrak{g}^*$. Moreover, every 1-cocycles $c \colon G \to \mathfrak{g}^*$ defines a (non-standard) coconjugation action by

$$(g, \mu) \mapsto \mathrm{Ad}^*_g \mu + c(g). \tag{4.2.23}$$

Recall that such affine actions play an important role for ordinary non-equivariant momentum maps (see [OR03, Definiton 4.5.23]). ◊

Every standard coconjugation action $\Upsilon$ of $G$ on $G^*$ defines a Poisson Lie structure $\Lambda \colon G^* \times \mathfrak{g} \to \mathfrak{g}^*$ on $G^*$ by

$$\Lambda(\eta, \xi) = \eta^{-1} \,.\, \mathrm{T}_e \Upsilon_\eta(\xi), \tag{4.2.24}$$

where $\Upsilon_\eta \colon G \to G^*$ is the orbit map through $\eta \in G^*$. Moreover, non-standard coconjugation actions are related to affine Poisson structures. A group-valued momentum map is equivariant with respect to a given standard coconjugation action $\Upsilon$ if and only if it is a Poisson map with respect to the induced Poisson structure $\Lambda$. This is essentially a consequence of the infinitesimal equivariance property

$$\delta_m J(\xi \,.\, m) = J(m)^{-1} \,.\, \big(\xi \,.\, J(m)\big) = \Lambda(J(m), \xi), \tag{4.2.25}$$

see [DR] for details. In the case when $G^* = \mathfrak{g}^*$ endowed with the coadjoint representation as the standard coconjugation action, the Poisson Lie structure defined in (4.2.24) is simply given by

$$\Lambda(\eta, \xi) = \mathrm{ad}^*_\xi \eta \tag{4.2.26}$$

and the equivariance condition (4.2.25) takes the usual form

$$\mathrm{T}_m J(\xi \,.\, m) = \mathrm{ad}^*_\xi \big(J(m)\big). \tag{4.2.27}$$

### 4.2.2 Bifurcation lemma

Let $(M, \omega, \Upsilon)$ be a symplectic $G$-manifold with group-valued momentum map $J \colon M \to G^*$. Throughout this section we assume that $J$ is $G$-equivariant with respect to a given coconjugation action on $G^*$. The momentum map relation (4.2.12) can be written in the form

$$\kappa\big(\delta_m J(v), \xi\big) = \omega_m\big(v, \mathrm{T}_e \Upsilon_m(\xi)\big) \tag{4.2.28}$$

for all $m \in M$, $v \in \mathrm{T}_m M$ and $\xi \in \mathfrak{g}$. Phrased in the language of dual pairs, this identity shows that $\delta_m J$ is the adjoint of the infinitesimal action $\mathrm{T}_e \Upsilon_m$ with respect to the dual pairs $\omega(\mathrm{T}_m M, \mathrm{T}_m M)$ and $\kappa(\mathfrak{g}^*, \mathfrak{g})$, cf. Appendix B.2. In particular, $\mathrm{T}_e \Upsilon_m$ and $\delta_m J$ are continuous with respect to the weak topologies induced by the dual pairs $\omega$ and $\kappa$. We will represent this situation diagrammatically as follows:

$$\begin{array}{ccc} \mathrm{T}_m M & \xleftarrow{\mathrm{T}_e \Upsilon_m} & \mathfrak{g} \\ \times_{\omega_m} & & \times_\kappa \\ \mathrm{T}_m M & \xrightarrow[\delta_m J]{} & \mathfrak{g}^*. \end{array} \tag{4.2.29}$$

Taken with a grain of salt, the momentum map $J$ "integrates" this diagram to

$$\begin{array}{ccc} M & \xleftarrow{\Upsilon_m} & G \\ \times_{\omega} & & \times_\kappa \\ M & \xrightarrow[J]{} & G^*. \end{array} \tag{4.2.30}$$

In this sense, we may view a (group-valued) momentum map as the symplectic adjoint of the group action.

The theory of adjoints allows us to generalize the following important result concerning the infinitesimal momentum map geometry to infinite dimensions (see, e.g., [OR03, Proposition 4.5.12] for the finite-dimensional version). The polar with respect to the dual pair $\kappa$ will be denoted by $\perp$, cf. (B.1.1).

Lemma 4.2.14 (Bifurcation Lemma) *Let $(M, \omega)$ be a symplectic $G$-manifold with equivariant momentum map $J\colon M \to G^*$. Then, for every $m \in M$, the following holds:*

*(i) (weak version)*

$$\operatorname{Ker} \delta_m J = (\mathfrak{g} \cdot m)^\omega, \tag{4.2.31a}$$

$$\mathfrak{g}_m = (\operatorname{Im} \delta_m J)^\perp. \tag{4.2.31b}$$

*(ii) (strong version) If, additionally, $\mathfrak{g} \cdot m = (\mathfrak{g} \cdot m)^{\omega\omega}$ holds, then we have*

$$(\operatorname{Ker} \delta_m J)^\omega = \mathfrak{g} \cdot m \tag{4.2.32a}$$

*and*

$$\operatorname{Ker} \delta_m J \cap (\operatorname{Ker} \delta_m J)^\omega = \mathfrak{g}_\mu \cdot m, \tag{4.2.32b}$$

*where* $\mu = J(m) \in G^*$. *Moreover, if* $\operatorname{Im} \delta_m J = (\operatorname{Im} \delta_m J)^{\perp\perp}$, *then*

$$\operatorname{Im} \delta_m J = \mathfrak{g}_m^{\perp}. \tag{4.2.32c}$$

◇

*Proof.* The formulae (4.2.31) follow from the general relation $\operatorname{Ker} T = (\operatorname{Im} T^*)^{\perp}$ relating the image of a linear operator $T$ and the kernel of its adjoint $T^*$, see Proposition B.2.3, applied to $T = \delta_m J$ and $T = \mathrm{T}_e \Upsilon_m$, respectively. In a similar vein, the general identity $(\operatorname{Im} T)^{\perp\perp} = (\operatorname{Ker} T^*)^{\perp}$ of Proposition B.2.3 (iii) implies (4.2.32a) and (4.2.32c). Using (4.2.32a), we find

$$\operatorname{Ker} \delta_m J \cap (\operatorname{Ker} \delta_m J)^{\omega} = \operatorname{Ker} \delta_m J \cap \mathfrak{g} \,.\, m \,. \tag{4.2.33}$$

The equivariance property (4.2.25) of $J$ implies that $\xi \,.\, m \in \operatorname{Ker} \delta_m J$ if and only if $\xi \,.\, \mu = 0$, i.e., $\xi \in \mathfrak{g}_{\mu}$. Hence, $\operatorname{Ker} \delta_m J \cap \mathfrak{g} \,.\, m = \mathfrak{g}_{\mu} \,.\, m$. □

Recall from Proposition 4.1.7 that the condition $(\mathfrak{g} \,.\, m)^{\omega\omega} = \mathfrak{g} \,.\, m$ is equivalent to $\mathfrak{g} \,.\, m$ being symplectically closed. In particular, by Lemma 4.1.8, every infinitesimal orbit of the action of a finite-dimensional Lie group is symplectically closed.

The Bifurcation Lemma is a helpful tool to study the bifurcations (surprise, surprise!) of Hamiltonian flows in the neighborhood of points with non-trivial stabilizer. With view towards symplectic reduction, its main significance is the existence of a symplectic slice at the infinitesimal level.

Proposition 4.2.15 *Let $(M, \omega)$ be a symplectic $G$-manifold with equivariant momentum map $J \colon M \to G^*$. Let $m \in M$ and $\mu = J(m) \in G^*$. Assume that the stabilizer $G_{\mu}$ of $\mu$ is a Lie subgroup of $G$ and that $\mathfrak{g} \,.\, m$ is symplectically closed. Then, for every slice $S$ at $m$ for the induced $G_{\mu}$-action,*

$$E \coloneqq \mathrm{T}_m S \cap \operatorname{Ker} \delta_m J \tag{4.2.34}$$

*is a symplectic subspace of $(\mathrm{T}_m M, \omega_m)$ in the sense of Definition 4.1.12.* ◇

*Proof.* Using (4.2.31) and Proposition B.1.3 (iv), we obtain

$$\begin{aligned} E \cap E^{\omega} &= \mathrm{T}_m S \cap \operatorname{Ker} \delta_m J \cap E^{\omega} \\ &= \mathrm{T}_m S \cap (\mathfrak{g} \,.\, m)^{\omega} \cap E^{\omega} \\ &= \mathrm{T}_m S \cap (\mathfrak{g} \,.\, m + E)^{\omega}. \end{aligned} \tag{4.2.35}$$

Since $\operatorname{Ker} \delta_m J \subseteq \mathfrak{g} \,.\, m + E$, the anti-monotony of the symplectic orthogonal, see Proposition B.1.3 (ii), implies $E \cap E^{\omega} \subseteq (\operatorname{Ker} \delta_m J)^{\omega}$. Thus, using the Bifurcation Lemma 4.2.14 (ii), we obtain

$$E \cap E^{\omega} \subseteq \mathrm{T}_m S \cap \operatorname{Ker} \delta_m J \cap (\operatorname{Ker} \delta_m J)^{\omega} = \mathrm{T}_m S \cap \mathfrak{g}_{\mu} \,.\, m = \{0\}, \tag{4.2.36}$$

[1] The momentum map fails to be a submersion exactly at points with non-trivial stabilizer.

because $S$ is a slice for the $G_\mu$-action. Hence, $E$ is a symplectic subspace of $(\mathrm{T}_m M, \omega_m)$. □

By $G$-equivariance of $J$, we have a chain complex

$$0 \longrightarrow \mathfrak{g}_\mu \xrightarrow{\mathrm{T}_e \Upsilon_m} \mathrm{T}_m M \xrightarrow{\delta_m J} \mathfrak{g}^* \longrightarrow 0. \tag{4.2.37}$$

Since, for every $G_\mu$-slice $S$ at $m$, $\mathrm{T}_m S$ is a topological complement of $\mathfrak{g}_\mu \,.\, m$ in $\mathrm{T}_m M$, the subspace $E$ defined in (4.2.34) is identified with the middle homology group of this complex. If the assumptions of the Strong Bifurcation Lemma hold, then the other homology groups are given by $\mathfrak{g}_m$ and $\mathfrak{g}^*_m$ according to (4.2.32c). Thus, if (4.2.37) is a Fredholm complex, then its Euler characteristic is given by

$$\dim \mathfrak{g}_m - \dim E + \dim \mathfrak{g}^*_m = 2 \dim \mathfrak{g}_m - \dim E. \tag{4.2.38}$$

According to Proposition 4.2.15, this number is always even. Keeping in mind that the Euler characteristic of (4.2.37) coincides (up to a sign) with the virtual dimension of the symplectic quotient, this observation is a first indication that the latter carries a symplectic structure.

Remark 4.2.16 Let us rephrase the result of Proposition 4.2.15 in a more familiar form by connecting it to the classical Witt–Artin decomposition. In the setting of Proposition 4.2.15, assume additionally that $G_\mu$ is a split Lie subgroup of $G$. Then, there exists a topological complement $\mathfrak{q}$ of $\mathfrak{g}_\mu$ in $\mathfrak{g}$. Since $S$ is a slice at $m$ for the $G_\mu$-action, we have a topological decomposition

$$\mathrm{T}_m M = \mathfrak{g}_\mu \,.\, m \oplus \mathrm{T}_m S. \tag{4.2.39}$$

Assume that $E$ has a topological complement $F'$ in $\mathrm{T}_m S$. Since $\mathfrak{q} \,.\, m \cap E = \{0\}$, $\mathfrak{q} \,.\, m$ is a subspace of $F'$. Assume that $\mathfrak{q} \,.\, m$ has a topological complement $F$ in $F'$. Then we obtain the topological decomposition:

$$\mathrm{T}_m M = \overbrace{\mathfrak{q} \,.\, m \oplus \mathfrak{g}_\mu \,.\, m}^{\mathfrak{g}.m} \oplus E \oplus F. \tag{4.2.40}$$

where $\underbrace{\mathfrak{g}_\mu \,.\, m \oplus E}_{\operatorname{Ker} \delta_m J}$.

For a classical momentum map, this decomposition is often referred to as the Witt–Artin decomposition (see e.g. [RS12, Equation 10.2.11]). We emphasize that, in addition to the conditions of Proposition 4.2.15, the derivation of the Witt–Artin decomposition in infinite dimensions relied on several additional assumptions concerning the existence of topological complements.

Let us now pass to the corresponding decomposition of $\mathfrak{g}^*$. For this purpose, suppose that $\mathfrak{g}_m$ is topologically complemented in $\mathfrak{g}_\mu$ so that we get the topological decomposition

$$\mathfrak{g} = \mathfrak{q} \oplus \underbrace{\mathfrak{m} \oplus \mathfrak{g}_m}_{\mathfrak{g}_\mu} . \tag{4.2.41}$$

If this decomposition is weakly continuous with respect to the dual pair $\kappa$, then, by [Köt83, Section 20.5], we obtain a dual decomposition of $\mathfrak{g}^*$:

$$\mathfrak{g}^* = \mathfrak{q}^* \oplus \mathfrak{m}^* \oplus \mathfrak{g}^*_m, \tag{4.2.42}$$

where $\mathfrak{q}^*$ is the annihilator of $\mathfrak{g}_\mu$ in $\mathfrak{g}$ and $\mathfrak{m}^*$ is the annihilator of $\mathfrak{g}_m$ in $\mathfrak{g}^*_\mu$. Moreover, assume that $\operatorname{Im} \delta_m J = (\operatorname{Im} \delta_m J)^{\perp\perp}$. Then, by Lemma 4.2.14 (ii), we have $\operatorname{Im} \delta_m J = \mathfrak{g}_m^\perp = \mathfrak{q}^* \oplus \mathfrak{m}^*$. Hence, in summary, we obtain the decomposition

$$\mathfrak{g}^* = \underbrace{\mathfrak{q}^* \oplus \mathfrak{m}^*}_{\operatorname{Im} \delta_m J} \oplus \mathfrak{g}^*_m. \tag{4.2.43}$$

Accordingly, $\delta_m J$ yields, by restriction, a bijection between $\mathfrak{q} \,.\, m \oplus F$ and $\mathfrak{q}^* \oplus \mathfrak{m}^*$, cf. (4.2.40). Due to the equivariance property (4.2.25), the restriction of $\delta_m J$ to $\mathfrak{q} \,.\, m$ yields the map

$$I_\mu \colon \mathfrak{q} \,.\, m \mapsto \mathfrak{q}^* \oplus \mathfrak{m}^*, \qquad \xi \,.\, m \mapsto \mu^{-1} \,.\, (\xi \,.\, \mu) = \Lambda(\mu, \xi). \tag{4.2.44}$$

Since $\mathfrak{q}$ is a complement of $\mathfrak{g}_\mu$, the map $I_\mu$ is injective. In the case when $G^* = \mathfrak{g}^*$ endowed with the coadjoint action as the coconjugation action, $I_\mu$ is given by

$$I_\mu \colon \mathfrak{q} \,.\, m \mapsto \mathfrak{q}^* \oplus \mathfrak{m}^*, \qquad \xi \,.\, m \mapsto \operatorname{ad}^*_\xi \mu, \tag{4.2.45}$$

see (4.2.27). Then, $I_\mu$ takes values in $\mathfrak{q}^*$, because $\operatorname{ad}^*_\xi \mu \in \mathfrak{g}_\mu^\perp = \mathfrak{q}^*$ for every $\xi \in \mathfrak{g}$. If $\mathfrak{q}$ is finite-dimensional, then a dimension count shows that $I_\mu$ yields an isomorphism between $\mathfrak{q} \,.\, m$ and $\mathfrak{q}^*$. Moreover, in this case, $\delta_m J$ restricts to a bijection between $F$ and $\mathfrak{m}^*$. We do not know if this observation generalizes to infinite-dimensional $\mathfrak{q}$ and/or to other coconjugation actions. ◇

Example 4.2.17 For the action of the group of gauge transformations on the space of connections on a principal U(1)-bundle $P$ discussed in Example 4.2.6, we have found that the infinitesimal action $\mathrm{d} \colon \Omega^0(M) \to \Omega^1(M)$ has $\delta_\alpha \mathcal{J} = \mathrm{d} \colon \Omega^1(M) \to \Omega^2(M)$ as its adjoint. The $L^2$-orthogonal Hodge decompositions

$$\Omega^1(M) = \mathrm{d}\, C^\infty(M) \oplus \mathrm{d}^* \Omega^2(M) \oplus \mathrm{H}^1(M, \mathbb{R}), \tag{4.2.46}$$

$$\Omega^2(M) = \mathrm{d}\Omega^1(M) \oplus \mathrm{H}^2(M, \mathbb{R}) \tag{4.2.47}$$

show that the closedness conditions $\mathfrak{g} \,.\, m = (\mathfrak{g} \,.\, m)^{\omega\omega}$ and $\operatorname{Im} \delta_m J = (\operatorname{Im} \delta_m J)^{\perp\perp}$ hold for the case under consideration. The identities (4.2.32a) and (4.2.32c) of the Strong Bifurcation Lemma are equivalent to

$$(\mathrm{d}\,\mathrm{C}^\infty(M) \oplus \mathrm{H}^1(M, \mathbb{R}))^\omega = \mathrm{d}\,\mathrm{C}^\infty(M), \tag{4.2.48}$$

$$\mathrm{d}\Omega^1(M) = \mathrm{H}^0(M, \mathbb{R})^\perp. \tag{4.2.49}$$

These identities can be verified also by hand using a simple application of Hodge theory.

Let us now determine the symplectic subspace $\mathcal{E}$ of Proposition 4.2.15 for this example. The orbit through $A \in \mathcal{C}(P)$ of the $\mathcal{G}au(P)$-action is given by

$$\mathcal{G}au(P) \cdot A = A + \Omega^1_{\mathrm{cl},\mathbb{Z}}(M), \tag{4.2.50}$$

where $\Omega^1_{\mathrm{cl},\mathbb{Z}}(M)$ denotes the space of closed 1-forms with integral periods. According to the Hodge decomposition (4.2.46), a natural choice of a slice $\mathcal{S}$ at $A$ is given by $\mathcal{S} = A + \mathcal{U}$, where $\mathcal{U}$ is a sufficiently small neighborhood of 0 in $\mathrm{d}^*\Omega^2(M) \oplus \mathrm{H}^1(M, \mathbb{R})$. For the subspace defined in (4.2.34), we hence obtain

$$\mathcal{E} = \mathrm{T}_A\mathcal{S} \cap \operatorname{Ker} \delta_A \mathcal{J} = \big(\mathrm{d}^*\Omega^2(M) \oplus \mathrm{H}^1(M, \mathbb{R})\big) \cap \operatorname{Ker} \mathrm{d} \simeq \mathrm{H}^1(M, \mathbb{R}). \tag{4.2.51}$$

It is a simple exercise to verify that $\mathcal{E}$ is a symplectic subspace, in accordance with Proposition 4.2.15. Moreover, the restriction of the symplectic structure $\omega$ on $\mathrm{T}_A\big(\mathcal{C}(P)\big)$ to $\mathcal{E}$ coincides with the intersection form on $\mathrm{H}^1(M, \mathbb{R})$:

$$\big([\alpha], [\beta]\big) \mapsto \int_M \alpha \wedge \beta. \tag{4.2.52}$$

◊

### 4.2.3 Normal form of a momentum map

Starting with the work of Arms [Arm81] on the Yang–Mills equation and of Fischer, Marsden, and Moncrief [FMM80] on the Einstein equation, it became clear that the solution spaces of field theories have singularities at points with internal symmetry. Shortly thereafter, a similar relation between internal symmetries and singularities of the momentum map was established in [AMM81]. The occurrence of these singularities of momentum map level sets was later explained by Marle [Mar83; Mar85] and Guillemin and Sternberg [GS84], who proved normal form theorems for momentum maps. These considerations were in a finite-dimensional setting or had a formal nature in the sense that the problems of the infinite-dimensional context were largely ignored. The aim of this section is to use and refine the construction of Section 3.1 to establish a normal form result for equivariant momentum maps in infinite dimensions.

At this point, the reader might want to recall the definition of an equivariant normal form, cf. Definition 3.1.3. In a similar spirit, we define the concept of a normal form of an equivariant momentum map.

Definition 4.2.18 An abstract *Marle–Guillemin–Sternberg* (MGS) normal form is a tuple $(H, X, \mathfrak{g}^*, \hat{J}, J_{\mathrm{sing}}, \bar{\omega})$ consisting of:

(i) an abstract equivariant normal form $(H, X, \mathfrak{g}^*, \hat{J}, J_{\mathrm{sing}})$ in the sense of Definition 3.1.3 satisfying $\mathrm{Coker} = \mathfrak{h}^*$, where Coker refers to the space in the decomposition (3.1.1), $\mathfrak{g}$ is a Lie algebra with a dual pair $\kappa(\mathfrak{g}^*, \mathfrak{g})$ and the Lie algebra $\mathfrak{h}$ of $H$ is identified with a subalgebra of $\mathfrak{g}$,

(ii) a closed $H$-invariant 2-form[1] $\bar{\omega}$ on $U \cap \mathrm{Ker}$ such that $\bar{\omega}_0$ is a symplectic form on Ker, where $U$ and Ker are as in Definition 3.1.3.

Moreover, we require that the momentum map identity

$$\bar{\omega}_x\big(\xi \,.\, x, w\big) + \kappa\big(\mathrm{T}_x J_{\mathrm{sing}}(w), \xi\big) = 0 \tag{4.2.53}$$

holds for all $x \in U \cap \mathrm{Ker}$, $w \in \mathrm{Ker}$ and $\xi \in \mathfrak{h}$. An abstract MGS normal form is called *strong* if $\bar{\omega}_x = \bar{\omega}_0$ holds and $J_{\mathrm{sing}}$ satisfies

$$\kappa\big(J_{\mathrm{sing}}(x), \xi\big) = \frac{1}{2}\bar{\omega}_0(x, \xi \,.\, x) \tag{4.2.54}$$

for all $x \in U \cap \mathrm{Ker}$ and $\xi \in \mathfrak{h}$. ◊

The important additional property of an MGS normal form in comparison to an abstract equivariant normal form is the fact that the singular part $J_{\mathrm{sing}}$ takes values in the dual of the Lie algebra of $H$ and that it is tamed by the 2-form $\bar{\omega}$ using the momentum map identity (4.2.53). Clearly, (4.2.54) is the integrated version of (4.2.53). Note that we do not require $\bar{\omega}$ to be symplectic away from the origin. Nonetheless, this is automatically the case if Ker is finite-dimensional, because then the space of invertible operators is open in the space of all linear maps. Moreover, in this case, one can use the Darboux theorem to pass from an MGS normal form to a strong one in the following sense.

Proposition 4.2.19 *Let $(H, X, \mathfrak{g}^*, \hat{J}, J_{\mathrm{sing}}, \bar{\omega})$ be an abstract MGS normal form. If the subspace* Ker *occurring in the decomposition* (3.1.1) *is finite-dimensional, then there exists an $H$-equivariant local diffeomorphism $\psi$ of $X$ such that $(H, X, \mathfrak{g}^*, \hat{J}, J_{\mathrm{sing}} \circ \psi, \psi^* \bar{\omega})$ is a strong MGS normal form.* ◊

*Proof.* Since Ker is finite-dimensional and $\bar{\omega}_0$ is symplectic, we may shrink $U$ such that $\bar{\omega}_x$ is a symplectic form on Ker for all $x \in U \cap \mathrm{Ker}$. The classical Darboux theorem for finite-dimensional symplectic manifolds yields a local diffeomorphism $\psi$ of $U \cap \mathrm{Ker}$ satisfying $\psi(0) = 0$ and $\mathrm{T}_0 \psi = \mathrm{id}_{\mathrm{Ker}}$ such that

$$\psi^* \bar{\omega} = \bar{\omega}_0. \tag{4.2.55}$$

[1] That is, $\bar{\omega}_x$ is an antisymmetric $H$-invariant bilinear form on Ker for every $x \in U \cap \mathrm{Ker}$.

By extending $\psi$ to the whole of $X$ using the identity map on Coim, we may regard $\psi$ as a local diffeomorphism of $X$. As $H$ is a compact Lie group, we may furthermore choose $\psi$ to be $H$-equivariant. The momentum map relation behaves naturally with respect to equivariant symplectomorphisms and thus $J_{\text{sing}} \circ \psi$ is a momentum map for the linear $H$-action on Ker with respect to the constant symplectic form $\bar{\omega}_0$. Since, in finite dimensions, the momentum map is unique (up to a constant which is fixed by the condition $J_{\text{sing}} \circ \psi(0) = J_{\text{sing}}(0) = 0$), we have

$$\kappa\big(J_{\text{sing}} \circ \psi(x), \xi\big) = \frac{1}{2}\bar{\omega}_0\big(x, \xi \,.\, x\big) \tag{4.2.56}$$

for every $x \in \text{Ker}$ and $\xi \in \mathfrak{h}$. In summary, $J_{\text{sing}} \circ \psi$ is in the strong MGS normal form according to (4.2.54). □

For an MGS normal form $(H, X, \mathfrak{g}^*, \hat{J}, J_{\text{sing}})$ and a Lie group $G_\mu$ with $H \subseteq G_\mu$, define the smooth map $J_{\text{NF}} \colon G_\mu \times_H U \to G_\mu \times_H \mathfrak{g}^*$ by

$$J_{\text{NF}}([g, x_1, x_2]) = [g, \hat{J}(x_2) + J_{\text{sing}}(x_1, x_2)] \tag{4.2.57}$$

for $g \in G_\mu$, $x_1 \in U \cap \text{Ker}$ and $x_2 \in U \cap \text{Coim}$.

Definition 4.2.20 Let $(M, \omega)$ be a symplectic $G$-manifold with equivariant momentum map $J \colon M \to G^*$. Let $m \in M$. Assume that the stabilizer $G_\mu$ of $\mu = J(m) \in G^*$ is a Lie subgroup of $G$. We say that $J$ can be *brought into the MGS normal form* $(H, X, \mathfrak{g}^*, \hat{J}, J_{\text{sing}}, \bar{\omega})$ at $m$ if $H = G_m$ and if there exists a linear slice $S$ at $m$ for the $G_\mu$-action, a $G_m$-equivariant diffeomorphism $\iota_S \colon X \supseteq U \to S \subseteq M$ and a $G_m$-equivariant chart $\rho \colon G^* \supseteq V' \to V \subseteq \mathfrak{g}^*$ at $\mu$, which bring $J$ into the equivariant normal form $J_{\text{NF}}$ according to Definition 3.1.3 and for which the restriction of $\iota_S^* \omega$ to Ker coincides with $\bar{\omega}$. ◇

If $J$ can be brought into the MGS normal form $(G_m, X, \mathfrak{g}^*, \hat{J}, J_{\text{sing}}, \bar{\omega})$ at $m \in M$, then according to Definition 3.1.3 there exists a commutative diagram of the type

$$\begin{array}{ccc} M & \xrightarrow{\quad J \quad} & G^* \\ {\scriptstyle \chi^{\mathrm{T}} \circ (\mathrm{id}_{G_\mu} \times \iota_S)} \big\uparrow & & \big\uparrow {\scriptstyle \rho^{\mathrm{T}}} \\ G_\mu \times_{G_m} X \supseteq G_\mu \times_{G_m} U & \xrightarrow{J_{\text{NF}}} & G_\mu \times_{G_m} V \subseteq G_\mu \times_{G_m} \mathfrak{g}^*, \end{array} \tag{4.2.58}$$

where $J_{\text{NF}}$ was defined in (4.2.57).

Remark 4.2.21 Assume that the equivariant momentum map $J \colon M \to G^*$ can be brought into an MGS normal form $(H, X, \mathfrak{g}^*, \hat{J}, J_{\text{sing}}, \bar{\omega})$ at $m \in M$ using the maps $\chi^T \colon G_\mu \times_{G_m} S \to M$, $\iota_S \colon X \supseteq U \to S$ and $\rho \colon G^* \supseteq V' \to V \subseteq \mathfrak{g}^*$. Under the natural isomorphisms $\mathrm{T}_0\iota_S \colon X \to \mathrm{T}_m S$ and $\mathrm{T}_\mu \rho \colon \mathrm{T}_\mu G^* \to \mathfrak{g}^*$, the Witt–Artin decomposition (cf. Remark 4.2.16) yields the following identification of the spaces occurring in the MGS normal form:

$$\mathrm{T}_m M = \mathfrak{g}_\mu \,.\, m \oplus \overbrace{\underbrace{E}_{\text{Ker}} \oplus \underbrace{F \oplus \mathfrak{q} \,.\, m}_{\text{Coim}}}^{X} \tag{4.2.59}$$

and

$$\mathfrak{g}^* = \underbrace{\mathfrak{q}^* \oplus \mathfrak{m}^*}_{\text{Im}} \oplus \underbrace{\mathfrak{g}^*_m}_{\text{Coker}} \,. \tag{4.2.60}$$

Moreover, under these identifications, suppressing $\iota_S$, the diagram (4.2.58) takes the form

$$\begin{array}{ccc} M & \xrightarrow{\quad J \quad} & G^* \\ {\scriptstyle \chi^{\mathrm{T}}}\uparrow & & \uparrow{\scriptstyle \rho^{\mathrm{T}}} \\ G_\mu \times_{G_m} (E \oplus F \oplus \mathfrak{q} \,.\, m) & \xrightarrow{J_{\mathrm{NF}}} & G_\mu \times_{G_m} \mathfrak{g}^* \,, \end{array} \tag{4.2.61}$$

where $J_{\mathrm{NF}}$ is given by $J_{\mathrm{NF}}([g, e, f, \xi] \,.\, m) = [g, \hat{J}(f, \xi \,.\, m) + J_{\mathrm{sing}}(e, f, \xi)]$. Furthermore, the restriction of the symplectic form $\omega_m$ on $\mathrm{T}_m M$ to $E$ is identified with the symplectic form $\bar{\omega}_0$ on Ker. ◊

Remark 4.2.22 (Classical Marle–Guillemin–Sternberg normal form) The classical Marle–Guillemin–Sternberg normal form in finite dimensions (see, e.g., [OR03, Theorem 7.5.5]) has a slightly different form and can be summarized in the following commutative diagram:

$$\begin{array}{ccc} M & \xrightarrow{\quad J \quad} & \mathfrak{g}^* \\ \uparrow & & \uparrow{\scriptstyle \Lambda} \\ G \times_{G_m} (E \times \mathfrak{m}^*) & \xrightarrow{J_{\mathrm{MGS}}} & G \times_{G_m} \mathfrak{g}^* \,, \end{array} \tag{4.2.62}$$

where $J_{\mathrm{MGS}}$ and $\Lambda$ are defined by

$$J_{\mathrm{MGS}}([g, e, \eta]) = [g, \eta + J_E(e)] \tag{4.2.63}$$

and $\Lambda([g, \nu]) = \mathrm{Ad}^*_{g^{-1}}(\mu + \nu)$, respectively. Here, $J_E : E \longrightarrow \mathfrak{g}^*_m$ is the quadratic momentum map for the symplectic linear $G_m$-action on $E$, see Example 4.2.5.

The most important difference is that the classical MGS normal form splits off the $G$-action while the $G_\mu$-action lies in the focus of our normal form. With view towards symmetry reduction, our formulation has advantages over the traditional approach, because in the context of symplectic reduction the group $G_\mu$ and not $G$ is the major player. In the classical MGS normal form $\mathfrak{m}^*$ is identified with $F$ via the map $\mathfrak{m}^* \ni \nu \mapsto f \in F$ implicitly defined by

$$\kappa(\nu, \xi) = \omega_m(\xi \,.\, m, f) \tag{4.2.64}$$

for $\xi \in \mathfrak{m}$. As already noted in Remark 4.2.16, we do not know whether this isomorphism generalizes to infinite dimensions. Furthermore, in the classical MGS normal form the non-linear term $J_E$ only depends on points in Ker while our normal form is weaker in this regard, allowing the singular part $J_{\mathrm{sing}}$ to additionally depend on points in $\mathfrak{q}$ and $F$. Finally, the classical Marle–Guillemin–Sternberg normal form brings both the symplectic structure and the momentum map simultaneously into a normal form, while we are mostly concerned with the momentum map. In fact, we do not have any control over the symplectic form in the Coim-direction. ◇

The following result is of fundamental importance for the construction of an MGS normal form.

Lemma 4.2.23 *Let $(M, \omega)$ be symplectic $G$-manifold with equivariant momentum map $J\colon M \to G^*$. Let $m \in M$ be such that $\mathfrak{g} \,.\, m$ is symplectically closed. Assume that $J$ can be brought into the $G_\mu$-equivariant normal form $(G_m, X, \mathfrak{g}^*, \hat{J}, J_{\mathrm{sing}})$ at $m$ using the maps $\chi^{\mathrm{T}}\colon G_\mu \times_{G_m} S \to M$ and $\rho\colon G^* \supseteq V' \to \mathfrak{g}^*$ in the sense of Definition 3.1.3. If, for all $\nu \in V' \subseteq G^*$, the left derivative*

$$\mathrm{T}^{\mathrm{L}}_\nu \rho\colon \mathfrak{g}^* \to \mathfrak{g}^*, \qquad \eta \mapsto \mathrm{T}_\nu \rho(\nu \,.\, \eta) \tag{4.2.65}$$

*of $\rho$ at $\nu$ restricts to the identity on $\mathfrak{g}^*_m$, then $J$ can be brought into the MGS normal form $(G_m, X, \mathfrak{g}^*, \hat{J}, J_{\mathrm{sing}}, \bar{\omega})$, where $\bar{\omega} = (\iota_S^* \omega)\restriction_{\mathrm{Ker}}$.* ◇

*Proof.* First, note that the commutative diagram (4.2.58) yields by restriction the following commutative diagram:

$$\begin{array}{ccc} M & \xrightarrow{\;J\;} & G^* \\ {\scriptstyle \iota_S}\uparrow & & \uparrow{\scriptstyle \rho^{-1}} \\ U \cap \mathrm{Ker} & \xrightarrow{\;J_{\mathrm{sing}}\;} & \mathfrak{g}^*_m \,. \end{array} \tag{4.2.66}$$

Moreover, by (A.0.4), we have $\delta_{\rho(\nu)}(\rho^{-1}) \circ \mathrm{T}^{\mathrm{L}}_\nu \rho = \mathrm{id}_{\mathfrak{g}^*}$. Thus, $\delta_{\rho(\nu)}(\rho^{-1})\colon \mathfrak{g}^* \to \mathfrak{g}^*$ restricts to the identity on $\mathfrak{g}^*_m$ for all $\nu \in V'$. Hence, for every $x \in U \cap \mathrm{Ker}$, $w \in \mathrm{Ker}$ and $\xi \in \mathfrak{g}_m$, we obtain

$$\begin{aligned} \kappa\big(\mathrm{T}_x J_{\mathrm{sing}}(w), \xi\big) &= \kappa\big(\delta_{\rho\,\circ J \circ \iota_S(x)}(\rho^{-1}) \circ \mathrm{T}_x J_{\mathrm{sing}}(w), \xi\big) \\ &= \kappa\big(\delta_{\iota_S(x)} J \circ \mathrm{T}_x \iota_S(w), \xi\big) \\ &= \omega_{\iota_S(x)}\big(\mathrm{T}_x \iota_S(w), \xi \,.\, \iota_S(x)\big) \\ &= \omega_{\iota_S(x)}\big(\mathrm{T}_x \iota_S(w), \mathrm{T}_x \iota_S(\xi \,.\, x)\big) \\ &= (\iota_S^* \omega)_x\big(w, \xi \,.\, x\big), \end{aligned} \tag{4.2.67}$$

where we used the momentum map relation and $G_m$-equivariance of $\iota_S$. Thus, if we let $\bar{\omega} = (\iota_S^* \omega)_{\restriction \mathrm{Ker}}$, then (4.2.53) holds. Moreover, $\bar{\omega}_0$ coincides with the restriction of $\omega_m$ to the subspace $E = \mathrm{T}_m S \cap \mathrm{Ker}\, \delta_m J$ (under the isomorphism $\mathrm{T}_0 \iota_S \colon X \to \mathrm{T}_m S$). Hence, $\bar{\omega}_0$ is a symplectic form according to Proposition 4.2.15. In summary, $J$ is brought into the MGS normal form $(G_m, X, \mathfrak{g}^*, \hat{J}, J_{\mathrm{sing}}, \bar{\omega})$ using the maps $\chi^{\mathrm{T}}$ and $\rho$. □

The assumption in Lemma 4.2.23 concerning the chart $\rho$ is rather easy to satisfy. For example, in the case when $G^* = \mathfrak{g}^*$, we can simply take $\rho$ to be the translation by $\mu$. More generally, if $G^*$ has an exponential map $\exp_{G^*}$ which is a local diffeomorphism, then $\rho \coloneqq \exp_{G^*}^{-1}$ satisfies $\mathrm{T}^{\mathrm{L}}_\nu \rho = \mathrm{id}_{\mathfrak{g}^*}$, because $G^*$ is abelian. This covers the cases we are interested in and thus, for simplicity, we will assume in the sequel that $\mathrm{T}^{\mathrm{L}}_\nu \rho = \mathrm{id}_{\mathfrak{g}^*}$.

Remark 4.2.24 The proof of Lemma 4.2.23 shows that the momentum map relation (4.2.53) holds even for all $x \in U$ and not just $x \in U \cap \mathrm{Ker}$ (with $\bar{\omega}$ being a 2-form on $U$). If Ker is finite-dimensional, this observation can be combined with a slightly modified version of Proposition 4.2.19 to show that $J_{\mathrm{sing}}$ even satisfies

$$\kappa\big(J_{\mathrm{sing}}(x_1, x_2), \xi\big) = \frac{1}{2} \bar{\omega}_{(0, x_2)}(x_1, \xi \,.\, x_1) \tag{4.2.68}$$

for all $(x_1, x_2) \in U$ and $\xi \in \mathfrak{g}_m$. Here, we used the property $J_{\mathrm{sing}}(0, x_2) = 0$ coming from the equivariant normal form to fix the constant of integration. We see that $J_{\mathrm{sing}}(x_1, x_2)$ is quadratic in $x_1$ for all $x_2$. In comparison to the classical MGS normal form, the only difference is that the quadratic form $x_1 \mapsto \bar{\omega}_{(0, x_2)}(\xi \,.\, x_1, x_1)$ on Ker still depends on $x_2$. ◇

With the help of Lemma 4.2.23, we can upgrade an equivariant normal form of $J$ in the sense of Definition 3.1.3 to an MGS normal form.

Theorem 4.2.25 (Marle–Guillemin–Sternberg Normal Form) *Let $(M, \omega)$ be symplectic $G$-manifold with equivariant momentum map $J \colon M \to G^*$. Let $m \in M$ and $\mu = J(m) \in G^*$. Assume that there exists a chart $\rho \colon G^* \supseteq V' \to \mathfrak{g}^*$ at $\mu$ which linearizes the $G_m$-action on $G^*$ and which satisfies $\mathrm{T}^{\mathrm{L}}_\nu \rho = \mathrm{id}_{\mathfrak{g}^*}$ for every $\nu \in V'$. Moreover, assume that $\mathfrak{g} \,.\, m$ is symplectically closed and that the image of $\delta_m J$ is weakly closed. If $J$ satisfies the assumptions of Theorem 3.1.6 (or of one of the other equivariant normal form theorems of Section 3.1), then $J$ can be brought into an MGS normal form.* ◇

*Proof.* Theorem 3.1.6 (or its other variations in Section 3.1) shows that $J$ can be brought into an equivariant normal form in a neighborhood of $m$ by deforming the chart $\rho$ using the local diffeomorphism $\phi$ defined in (2.2.9). A simple inspection shows that the restriction of $\mathrm{T}_{(y_1, y_2)} \phi \colon Y \to Y$ to $\mathrm{Coker}\, T$ is the identity map for all $(y_1, y_2) \in V$ (continuing in the notation of the proof of Theorem 2.2.6). Since $\mathrm{Im}\, \delta_m J$ is a weakly closed subspace of $\mathfrak{g}^*$, the strong Bifurcation Lemma 4.2.14 (ii) shows that in our setting $\mathrm{Coker}\, T$ is identified with $\mathfrak{g}^*_m$. Hence, the chart $\rho' = \rho \circ \phi^{-1}$ satisfies the assumptions of Lemma 4.2.23, which implies that $J$ can be brought into an MGS normal form. □

We note that, in finite dimensions and for classical momentum maps, all assumptions of Theorem 4.2.25 are met automatically except for properness of the action.

Remark 4.2.26 The proof of Theorem 4.2.25 is constructive in the sense that the diffeomorphisms bringing $J$ into the MGS normal form have explicit and relatively simple expressions. In contrast, the usual proof of the classical MGS normal form theorem relies on the relative Darboux theorem, which makes it difficult[1] to determine the deforming diffeomorphism (both analytically as well as numerically). Thus, one would expect that our proof of the MGS normal form theorem can be used to design new numerical discretization algorithms which preserve the momentum map geometry. This will be explored in further work. ◇

Recall from Proposition 4.2.15 that $\operatorname{Ker} \simeq E$ can be identified with the middle homology group of the complex

$$0 \longrightarrow \mathfrak{g}_\mu \xrightarrow{\mathrm{T}_e \Upsilon_m} \mathrm{T}_m M \xrightarrow{\delta_m J} \mathfrak{g}^* \longrightarrow 0, \tag{4.2.69}$$

where $\Upsilon_m \colon G \to M$ denotes the orbit map as before. If this complex is Fredholm, then $\operatorname{Ker} \simeq E$ is finite-dimensional and every MGS normal form can be upgraded to a strong MGS normal form according to Proposition 4.2.19. This is, in particular, the case for elliptic actions so that Theorems 3.1.12 and 4.2.25 yield the following MGS normal form theorem.

Theorem 4.2.27 (MGS Normal Form — elliptic version) *Let $G$ be a tame Fréchet Lie group, let $G^*$ be a geometric dual Lie group and let $(M, \omega, \Upsilon)$ be a symplectic geometric tame Fréchet $G$-manifold with equivariant momentum map $J \colon M \to G^*$. Let $m \in M$ and $\mu = J(m)$. Assume that the following conditions hold:*

*(i) The stabilizer subgroup $G_\mu$ of $\mu$ is a geometric tame Fréchet Lie subgroup of $G$.*

*(ii) The induced $G_\mu$-action on $M$ is proper and admits a slice $S$ at $m$.*

*(iii) The induced $G_m$-action on $G^*$ can be linearized at $\mu$ using a chart $\rho \colon G^* \supseteq V' \to \mathfrak{g}^*$ satisfying $\mathrm{T}^{\mathrm{L}}_\nu \rho = \mathrm{id}_{\mathfrak{g}^*}$ for every $\nu \in V'$.*

*(iv) The chain*

$$0 \longrightarrow \mathfrak{g}_\mu \xrightarrow{\mathrm{T}_e \Upsilon_s} \mathrm{T}_s M \xrightarrow{\delta_s J} \mathfrak{g}^* \longrightarrow 0 \tag{4.2.70}$$

[1] Even in the simplest cases, it is usually impossible to integrate in closed-form the differential flow equation underlying Moser's method.

*of linear maps parametrized by $s \in S$ is a chain of differential operators, which constitute an elliptic complex at $m$.*

*(v) The subspace $\mathfrak{g} \, . \, m \subseteq \mathrm{T}_m M$ is symplectically closed and the image of $\delta_m J$ is weakly closed in $\mathfrak{g}^*$.*

*Then, $J$ can be brought into a strong MGS normal form at $m$.* ◊

Example 4.2.28 Let us use Theorem 4.2.27 to show that the momentum map for the action of the group of gauge transformations on the space of connections on a principal $\mathrm{U}(1)$-bundle $P \to M$ over a surface $M$ can be brought into an MGS normal form. By Example 4.2.6, the map $\mathcal{J} : \mathcal{C}(P) \to \Omega^2(M)$ assigning the curvature $F_A$ to a connection $A \in \mathcal{C}(P)$ is the momentum map for the $\mathcal{G}au(P)$-action on $\mathcal{C}(P)$. Let $\mu \in \Omega^2(M)$. Since $\mathcal{G}au(P)$ acts trivially on $\Omega^2(M)$, the stabilizer group of $\mu$ coincides with the whole group $\mathcal{G}au(P)$, which clearly is a geometric tame Fréchet Lie group. As discussed in Example 4.2.17, the $\mathcal{G}au(P)$-action is proper and admits a slice $\mathcal{S}$ at every $A \in \mathcal{C}(P)$ of the form $\mathcal{S} = A + \mathcal{U}$, where $\mathcal{U}$ is a sufficiently small neighborhood of 0 in $\mathrm{d}^*\Omega^2(M) \oplus \mathrm{H}^1(M, \mathbb{R})$. In this case, the chain (4.2.70) is independent of $B \in \mathcal{S}$ and coincides (up to a sign) with the elliptic de Rham complex

$$0 \longrightarrow \Omega^0(M) \xrightarrow{\mathrm{d}} \Omega^1(M) \xrightarrow{\mathrm{d}} \Omega^2(M) \longrightarrow 0. \tag{4.2.71}$$

The closedness assumptions of the subspaces $\mathfrak{g} \, . \, m$ and $\operatorname{Im} \delta_m J$ were checked in Example 4.2.17 for the model under consideration. In summary, Theorem 4.2.27 implies that $J$ can be brought into a strong MGS normal form. In the present setting, the normal form is particularly simple. Indeed, for every connection $A$, the curvature map restricted to the slice $\mathcal{S}$ at $A$ is already in the MGS normal form

$$\mathcal{J} : \mathcal{S} = A + \mathcal{U} \mapsto \Omega^2(M), \qquad A + \alpha_1 + \alpha_2 \mapsto F_A + \mathrm{d}\alpha_1, \tag{4.2.72}$$

where $\alpha_1 \in \mathrm{d}^*\Omega^2(M)$ and $\alpha_2 \in \mathrm{H}^1(M, \mathbb{R})$. We read off that the linear part of $\mathcal{J}$ is given by the isomorphism $\mathrm{d} \colon \mathrm{d}^*\Omega^2(M) \to \mathrm{d}\Omega^1(M)$ and that the singular part $\mathcal{J}_{\mathrm{sing}}$ vanishes. The latter should not come as a surprise as the $\mathcal{G}au(P)$-action on $\mathcal{C}(P)$ has only one orbit type, namely $\mathrm{U}(1)$. ◊

The linear action of a compact group on a symplectic vector space provides another situation where we can show the existence of an MGS normal form. To see this, let $(X, \omega)$ be a symplectic Fréchet space endowed with a continuous linear symplectic action of the compact finite-dimensional Lie group $G$. Recall from Example 4.2.5, that in this case, the momentum map $J \colon X \to \mathfrak{g}^*$ for the $G$-action is given by

$$\kappa(J(x), \xi) = \frac{1}{2}\omega(x, \xi \, . \, x), \tag{4.2.73}$$

for $\xi \in \mathfrak{g}$. Moreover, $J$ is equivariant with respect to the coadjoint action of $G$ on $\mathfrak{g}^*$.

Theorem 4.2.29 *Let $(X, \omega)$ be a symplectic Fréchet space endowed with a continuous linear symplectic action of the compact Lie group $G$. Then, the equivariant momentum map $J\colon X \to \mathfrak{g}^*$ defined in (4.2.73) can be brought into a strong MGS normal form at every point of $X$.* ◊

*Proof.* Let $x \in X$ and let $\mu = J(x)$. The stabilizer subgroup $G_\mu$ is a compact Lie subgroup of $G$, because $G$ is compact and thus finite-dimensional. Therefore, the $G_\mu$-action on $X$ is proper and, by Theorem A.2.4, it admits a linear slice at $x$. Since $\mathfrak{g}$ is finite-dimensional, the image of $T_x J$ is (weakly) closed and the orbit $\mathfrak{g} \,.\, x$ is symplectically closed due to Lemma 4.1.8. Define $\rho\colon \mathfrak{g}^* \to \mathfrak{g}^*$ by $\rho(\nu) = \nu - \mu$. Clearly, $\rho$ is $\mathrm{Ad}^*_{G_\mu}$-invariant and satisfies $T_\nu \rho = \mathrm{id}_{\mathfrak{g}^*}$ for every $\nu \in \mathfrak{g}^*$. Now, Theorem 3.1.9 and Theorem 4.2.25 imply that $J$ can be brought into an MGS normal form at $x$. Since the symplectic form $\omega$ on $X$ is constant, an argument similar to the one in the proof of Proposition 4.2.19 shows that $J_{\mathrm{sing}}$ satisfies (4.2.54). □

## 4.3 Singular symplectic reduction

In this section, we are concerned with the geometric structure of the symplectic quotient

$$\check{M}_\mu \equiv M /\!\!/_\mu G = J^{-1}(\mu)/G_\mu. \tag{4.3.1}$$

In finite dimensions, the geometry of $\check{M}_\mu$ attracted a lot of attention starting with the work of Meyer [Mey73] and Marsden and Weinstein [MW74] on the regular case and of Arms, Gotay, and Jennings [AGJ90] and Sjamaar and Lerman [SL91] on the singular case. Extending these classical structure theorems to our infinite-dimensional context, we will show that the decomposition of $\check{M}_\mu$ into orbit types is a stratification and that each stratum carries a natural symplectic structure. These results are obtained under the standing assumption that $J$ can be brought into an MGS normal form. Moreover, we will see that the frontier condition of the stratification relies on the strong MGS normal form combined with a certain approximation property.

We continue to work in the setting of the previous section. Let $(M, \omega)$ be a symplectic $G$-manifold with proper $G$-action and with equivariant momentum map $J\colon M \to G^*$, where $\kappa(G^*, G)$ is a dual pair and $G^*$ carries a given coconjugation action. Recall from Appendix A.2 that $M$ decomposes into the orbit type subsets $M_{(H)}$ and that $M_{(H)}$ are submanifolds of $M$ given that the group action admits a slice at every point. The following is the symplectic version of this fact.

Proposition 4.3.1 *Let $(M, \omega)$ be a symplectic $G$-manifold with proper $G$-action and with equivariant momentum map $J\colon M \to G^*$. Let $\mu \in G^*$. If $J$ can be brought into* *an MGS normal form at every point of* $M_\mu := J^{-1}(\mu)$*, then, for every orbit type* $(H)$ *of the* $G_\mu$*-action on* $M_\mu$*, the set*

$$M_{(H),\mu} := M_{(H)} \cap J^{-1}(\mu) \tag{4.3.2}$$

*is a smooth submanifold of* $M$*. Moreover, there exists a unique smooth manifold structure on the quotient*

$$\check{M}_{(H),\mu} := M_{(H),\mu} / G_\mu \tag{4.3.3}$$

*such that the natural projection* $\pi_{(H),\mu} \colon M_{(H),\mu} \to \check{M}_{(H),\mu}$ *is a smooth surjective submersion. Furthermore, for every* $m \in M_{(H),\mu}$*, the projection* $\mathrm{T}_m \pi_{(H),\mu}$ *restricts to an isomorphism*

$$\mathrm{T}_m \pi_{(H),\mu} \colon E_{G_m} \to \mathrm{T}_{[m]} \check{M}_{(H),\mu}, \tag{4.3.4}$$

*where* $E_{G_m}$ *denotes the set of fixed points of* $E = \mathrm{T}_m S \cap \operatorname{Ker} \delta_m J \subseteq \mathrm{T}_m M$ *under the* $G_m$*-action and* $S$ *is the* $G_\mu$*-slice in the MGS normal form at* $m$*.* ◇

*Proof.* Let $m \in M_{(H),\mu}$. By assumption, $J$ can be brought into an MGS normal form $(G_m, X, \mathfrak{g}^*, \hat{J}, J_{\mathrm{sing}})$ using maps $\chi^{\mathrm{T}} \colon G_\mu \times_{G_m} S \to M$ and $\rho \colon G^* \supseteq V' \to \mathfrak{g}^*$. Thus, we have

$$J \circ \chi^{\mathrm{T}} \circ (\mathrm{id}_{G_\mu} \times \iota_S)\big([g, x_1, x_2]\big) = \rho^{\mathrm{T}}\big([g, \hat{J}(x_2) + J_{\mathrm{sing}}(x_1, x_2)]\big) \tag{4.3.5}$$

for $g \in G_\mu$ and $(x_1, x_2) \in U \subseteq X$, where $\iota_S \colon X \supseteq U \to S$ is the slice diffeomorphism. Since $(\rho^{\mathrm{T}})^{-1}(\mu) = G_\mu \times_{G_m} \{0\}$ and since $\hat{J}$ is an isomorphism, we find that $\chi^{\mathrm{T}}$ identifies the level set $M_\mu = J^{-1}(\mu)$ with the set

$$G_\mu \times_{G_m} \big\{x_1 \in U \cap \operatorname{Ker} : \; J_{\mathrm{sing}}(x_1) = 0\big\}. \tag{4.3.6}$$

Moreover, $\chi^{\mathrm{T}}$ is $G_\mu$-invariant and so $M_{(H),\mu}$ is locally identified with

$$G_\mu \times_{G_m} \big(U_{(H)} \cap \operatorname{Ker} \cap J_{\mathrm{sing}}^{-1}(0)\big), \tag{4.3.7}$$

where for $U_{(H)}$ the stabilizer is taken with respect to the linear $G_m$-action but conjugation is understood with respect to $G_\mu$. Since the action is proper, we have $S_{(H)} = S_{G_m}$ and therefore $U_{(H)} = U_{G_m}$. Hence, for every $x_1 \in U_{(H)} \cap \operatorname{Ker}$, $w \in \operatorname{Ker}$ and $\xi \in \mathfrak{g}_m$, we find

$$\kappa\big(\mathrm{T}_{x_1} J_{\mathrm{sing}}(w), \xi\big) = -\bar{\omega}_{x_1}(\xi \,.\, x_1, w) = 0 \tag{4.3.8}$$

and so $\mathrm{T}_{x_1} J_{\mathrm{sing}}(w) = 0$. By possibly shrinking $U$, we may assume that $U_{G_m} \cap \operatorname{Ker}$ is convex. Now the Fundamental Theorem of Calculus [Nee06, Proposition

I.2.3.2] implies

$$J_{\text{sing}}(x_1) = J_{\text{sing}}(0) + \int_0^1 \mathrm{T}_{(tx_1)} J_{\text{sing}}(x_1)\, \mathrm{d}t = 0 \tag{4.3.9}$$

for every $x_1 \in U_{G_m} \cap \operatorname{Ker}$. Thus, in summary, $M_{(H),\mu}$ is locally identified with

$$G_\mu \times_{G_m} \left(U_{(H)} \cap \operatorname{Ker} \cap J_{\text{sing}}^{-1}(0)\right) = G_\mu \times \left(U_{G_m} \cap \operatorname{Ker}\right), \tag{4.3.10}$$

from which we conclude that $M_{(H),\mu}$ is a smooth submanifold of $M$. Moreover, $\check{M}_{(H),\mu}$ is locally identified with $U_{G_m} \cap \operatorname{Ker}$ and so carries a smooth manifold structure modeled on $\operatorname{Ker}_{G_m}$. In these coordinates, the projection $\pi_{(H),\mu} \colon M_{(H),\mu} \to \check{M}_{(H),\mu}$ corresponds to the projection onto the second factor and thus is a smooth submersion. Finally, note that the $G_m$-equivariant isomorphism $\mathrm{T}_0 \iota_S \colon X \to \mathrm{T}_m S$ identifies $\operatorname{Ker}_{G_m}$ with $E_{G_m}$. □

Remark 4.3.2 In contrast to the regular case, the level set $M_\mu$ is in general not a smooth manifold. Indeed, as we have seen in (4.3.6), $M_\mu$ is locally modeled on spaces of the type $\operatorname{Ker} \cap J_{\text{sing}}^{-1}(0)$. If the MGS normal form is strong, then $(J_{\text{sing}})\restriction_{\operatorname{Ker}}$ is a quadratic form and thus $M_\mu$ has conic singularities[1]. ◇

Proposition 4.3.3 *In the setting of Proposition 4.3.1, assume additionally that* $\mathfrak{g} \,.\, m$ *is symplectically closed for every* $m \in M_{(H),\mu}$. *Then, there exists a symplectic form* $\check{\omega}_{(H),\mu}$ *on* $\check{M}_{(H),\mu}$ *uniquely determined by*

$$\pi_{(H),\mu}^* \check{\omega}_{(H),\mu} = \iota_{(H),\mu}^* \omega, \tag{4.3.11}$$

*where* $\iota_{(H),\mu} \colon M_{(H),\mu} \to M$ *is the natural injection.* ◇

*Proof.* To proof that (4.3.11) uniquely defines a 2-form $\check{\omega}_{(H),\mu}$ on $\check{M}_{(H),\mu}$ it suffices to show that $\iota_{(H),\mu}^* \omega$ is $G_\mu$-invariant and horizontal with respect to the smooth submersion $\pi_{(H),\mu}$. Invariance under $G_\mu$ is clear, because $\iota_{(H),\mu}$ is $G_\mu$-equivariant and $\omega$ is $G$-invariant. Furthermore, for every $\xi \in \mathfrak{g}_\mu$, $m \in M_{(H),\mu}$ and $v \in \mathrm{T}_m M_{(H),\mu}$, we find

$$\left(\iota_{(H),\mu}^* \omega\right)_m (\xi \,.\, m, v) = \omega_m(\xi \,.\, m, v) = -\kappa(\delta_m J(v), \xi) = 0, \tag{4.3.12}$$

because $J$ is constant on $M_{(H),\mu}$. In summary, $\iota_{(H),\mu}^* \omega$ is $G_\mu$-invariant and horizontal, and thus descends to a smooth 2-form $\check{\omega}_{(H),\mu}$ on $\check{M}_{(H),\mu}$ which, by definition, satisfies (4.3.11). Moreover, the identity (4.3.11) shows that

[1] The fact that the level set of the momentum map has conic singularities is well-known in finite dimensions (see, e.g., [OR03, Proposition 8.1.2]) and was first observed by Arms, Marsden, and Moncrief [AMM81, Theorem 5].

$\check{\omega}_{(H),\mu}$ is closed. It remains to prove that $\check{\omega}_{(H),\mu}$ is non-degenerate. For this purpose, recall from Proposition 4.3.1 that the projection $\mathrm{T}_m\pi_{(H),\mu}$ restricts to an isomorphism of $\mathrm{T}_{[m]}\check{M}_{(H),\mu}$ with $E_{G_m}$, where $E = \mathrm{T}_mS \cap \operatorname{Ker} \delta_m J$ and $S$ is the $G_\mu$-slice in MGS normal form at $m$. Equation (4.3.11) shows that, under this isomorphism, $(\check{\omega}_{(H),\mu})_{[m]}$ coincides with the restriction of $\omega_m$ to $E_{G_m}$. By Proposition 4.2.15, $E$ is a symplectic subspace of $(\mathrm{T}_mM, \omega)$. Since $G_m$ is compact, Proposition 4.1.15 implies that $E_{G_m}$ is symplectic as well. Thus, we conclude that $\check{\omega}_{(H),\mu}$ is non-degenerate. □

At this point, we know that $\check{M}_\mu$ decomposes into orbit type manifolds, every one of which carries a symplectic structure. We will now see that the pieces fit together in a particularly nice way. The reader might want to recall from Appendix A.2 the notion of a stratified space.

Proposition 4.3.4 *Let $(M, \omega)$ be a symplectic $G$-manifold with proper $G$-action and with equivariant momentum map $J\colon M \to G^*$. Let $\mu \in G^*$. Assume that $J$ can be brought into a strong MGS normal form $(H, X, \mathfrak{g}, \hat{J}, J_{\mathrm{sing}}, \bar{\omega})$ at every point $m \in M_\mu$ such that the intersection*

$$U_{(K)} \cap \operatorname{Ker} \cap J_{\mathrm{sing}}^{-1}(0) \tag{4.3.13}$$

*is non-empty for every orbit type $(K)$ of the $G_\mu$-action on $M_\mu$ satisfying $(K) \leq (G_m)$. Then, the decomposition of $M_\mu$ and $\check{M}_\mu$ into orbit type subsets $M_{(H),\mu}$ and $\check{M}_{(H),\mu}$, respectively, is a stratification.* ◊

A strong MGS normal form $(H, X, \mathfrak{g}, \hat{J}, J_{\mathrm{sing}}, \bar{\omega})$ satisfying the assumption of Proposition 4.3.4 is said to *have the approximation property*.

*Proof.* Let $m \in M_\mu$ and let $(K)$ be an orbit type of the $G_\mu$-action on $M_\mu$ with $(K) \leq (H)$, where $H = G_m$. By assumption, there exists a strong MGS normal form $(H, X, \mathfrak{g}, \hat{J}, J_{\mathrm{sing}}, \bar{\omega})$ at $m$ such that $Y \equiv U_{(K)} \cap \operatorname{Ker} \cap J_{\mathrm{sing}}^{-1}(0)$ is non-empty. A point $x \in U \cap \operatorname{Ker}$ lies in $Y$ if and only if $(H_x) = (K)$ and $J_{\mathrm{sing}}(x) = 0$. Since the $H$-action on $X$ is linear, we have $H_x = H_{\alpha x}$ for all $\alpha \in \mathbb{R}_{>0}$. Moreover, by (4.2.54), we obtain

$$\kappa\big(J_{\mathrm{sing}}(\alpha x), \xi\big) = \frac{1}{2}\bar{\omega}_0\big(\xi \,.\, (\alpha x), \alpha x\big) = \alpha^2\, \kappa\big(J_{\mathrm{sing}}(x), \xi\big) \tag{4.3.14}$$

Thus, for every $x \in Y$ and $\alpha \in \mathbb{R}_{>0}$, the point $\alpha x$ lies in $Y$ as well. By letting $\alpha \to 0$, we conclude that the point $m$ lies in the closure of $Y$ in $X$. Thus, the claim follows from Proposition 3.3.1. □

In summary, we obtain the following result concerning the structure of the symplectic quotient $\check{M}_\mu$.

Theorem 4.3.5 (Singular Reduction Theorem) *Let $(M, \omega)$ be a symplectic $G$-manifold with proper action and equivariant momentum map $J\colon M \to G^*$. Let $\mu \in G^*$. Assume that $G_\mu$ is a Lie subgroup of $G$ and that $\mathfrak{g} \,.\, m$ is symplectically closed for* *every $m \in M_\mu = J^{-1}(\mu)$. If $J$ can be brought into an MGS normal form at every point of $M_\mu$, then the following holds:*

*(i) For every orbit type $(H)$ of the $G_\mu$-action on $M_\mu$, the orbit type subset $M_{(H),\mu}$ is a smooth submanifold of $M$. Moreover, there exists a unique smooth manifold structure on the quotient*

$$\check{M}_{(H),\mu} := M_{(H),\mu}/G_\mu \tag{4.3.15}$$

*such that the natural projection $\pi_{(H),\mu}\colon M_{(H),\mu} \to \check{M}_{(H),\mu}$ is a smooth submersion.*

*(ii) If, additionally, the MGS normal forms can be chosen to be strong and to have the approximation property, then the decomposition of $M_\mu$ and $\check{M}_\mu$ into orbit type subsets $M_{(H),\mu}$ and $\check{M}_{(H),\mu}$, respectively, is a stratification.*

*(iii) For every orbit type $(H)$ of the $G_\mu$-action on $M_\mu$, there exists a symplectic form $\check{\omega}_{(H),\mu}$ on $\check{M}_{(H),\mu}$ uniquely determined by*

$$\pi^*_{(H),\mu}\,\check{\omega}_{(H),\mu} = \iota^*_{(H),\mu}\omega, \tag{4.3.16}$$

*where $\iota_{(H),\mu}\colon M_{(H),\mu} \to M$ is the natural injection.* ◇

Example 4.3.6 Continuing in the setting of Example 4.2.6, let $P \to M$ be a principal U(1)-bundle on the closed surface $M$. As we have seen, the momentum map $\mathcal{J}\colon \mathcal{C}(P) \to \Omega^2(M)$ for the $\mathcal{G}au(P)$-action on $\mathcal{C}(P)$ is given by (minus) the curvature. Thus, the symplectic quotient at 0,

$$\check{\mathcal{C}}_0(P) \equiv \mathcal{J}^{-1}(0)/\mathcal{G}au(P), \tag{4.3.17}$$

coincides with the moduli space of flat connections. Moreover, by the discussion in Example 4.2.28, $\mathcal{J}$ can be brought into a strong MGS normal form at every $A \in \mathcal{C}(P)$. Since the $\mathcal{G}au(P)$-action is free, Theorem 4.3.5 implies that $\check{\mathcal{C}}_0(P)$ is a symplectic manifold.

In the present setting, we can use the method of invariants to get a direct description of $\check{\mathcal{C}}_0(P)$. Let $(\gamma_i)$ be a family of closed piecewise smooth paths in $M$ generating $\pi_1(M)$. The holonomy $\mathrm{Hol}_A(\gamma_i)$ of $\gamma_i$ with respect to a connection $A$ furnishes a map

$$\mathcal{K}\colon \mathcal{C}(P) \to \mathrm{H}^1(M, \mathrm{U}(1)), \quad A \mapsto \big(\gamma_i \mapsto \mathrm{Hol}_A(\gamma_i)\big), \tag{4.3.18}$$

where we used the Universal Coefficient Theorem and the Hurewicz Theorem to identify $\mathrm{H}^1(M, \mathrm{U}(1))$ with $\mathrm{Hom}(\pi_1(M), \mathrm{U}(1))$. Since $\mathcal{K}$ is $\mathcal{G}au(P)$-invariant, it descends to a map $\check{\mathcal{K}}$ from $\check{\mathcal{C}}_0(P)$ to $\mathrm{H}^1(M, \mathrm{U}(1))$. Moreover, by standard arguments, $\check{\mathcal{K}}$ is a diffeomorphism. Under the identification given by $\check{\mathcal{K}}$, the reduced symplectic form on $\check{\mathcal{C}}_0(P)$ coincides with the intersection form on $\mathrm{H}^1(M, \mathbb{R})$. The corresponding story for a principal $G$-bundle $P$ with non-abelian structure group $G$ will be discussed below in Section 4.4. It turns out that, in this case, the moduli space $\check{\mathcal{C}}_0(P)$ of flat connections has singularities and is a stratified symplectic space. ◊

It is quite remarkable that, in the linear setting, the usual finite-dimensional result about symplectic strata directly generalizes to the infinite-dimensional realm without any further assumptions.

Corollary 4.3.7 *Let $(X, \omega)$ be a symplectic Fréchet space endowed with a continuous linear symplectic action of the compact Lie group $G$. Then, the $G$-action has a unique (up to a constant) equivariant momentum map $J\colon X \to \mathfrak{g}^*$ given by*

$$\kappa(J(x), \xi) = \frac{1}{2}\omega(x, \xi \,.\, x) \tag{4.3.19}$$

*for $x \in X$ and $\xi \in \mathfrak{g}$. Moreover, for every orbit type $(H)$, the subset $X_{(H),0} \coloneqq X_{(H)} \cap J^{-1}(0)$ is a smooth submanifold of $X$ and there exists a unique smooth manifold structure on $\check{X}_{(H),0} \coloneqq X_{(H),0}/G$ such that the natural projection $\pi_{(H)}\colon X_{(H),0} \to \check{X}_{(H),0}$ is a smooth submersion. Furthermore, $\check{X}_{(H),0}$ carries a symplectic form $\check{\omega}_{(H)}$ uniquely determined by*

$$\pi^*_{(H)}\check{\omega}_{(H)} = \omega\restriction_{X_{(H),0}}. \tag{4.3.20}$$

◊

*Proof.* Since $G$ is finite-dimensional, (4.3.19) defines a smooth map $J\colon X \to \mathfrak{g}^*$. That $J$ is indeed a momentum map follows from a routine calculation, which we leave to the reader. As $G$ is compact, the $G$-action on $X$ is proper. Moreover, for every $x$, the orbit $\mathfrak{g} \,.\, x$ is a finite-dimensional subspace of $X$ and thus is symplectically closed due to Lemma 4.1.8. According to Theorem 4.2.29, $J$ can be brought into an MGS normal form. Hence, the claims follow from the Singular Reduction Theorem 4.3.5. □

Let us now pass from the kinematic picture presented so far to dynamics.

Proposition 4.3.8 *Let $h$ be a $G$-invariant Hamiltonian on the symplectic $G$-manifold $(M, \omega)$ with equivariant momentum map $J\colon M \to G^*$. Assume*[1] *that the associated Hamiltonian vector field $X_h$ exists and that it has a unique flow $\mathrm{FL}^h_t$. Let $(H)$ be an orbit type and $\mu \in G^*$. Then,*

*(i) the flow $\mathrm{FL}^h_t$ is $G$-equivariant and leaves $M_{(H),\mu}$ invariant and, hence, it projects to a flow $\check{\mathrm{FL}}^h_t$ on $\check{M}_{(H),\mu}$,*

[1] Recall that the Hamiltonian vector field associated to the Hamiltonian $h$ may not exist in infinite dimensions and that vector fields on Fréchet manifolds need not have flows. The latter is more or less equivalent to local in time solutions of the corresponding partial differential equation.

*(ii) the projected flow* $\check{\mathrm{FL}}^h_t$ *is Hamiltonian with respect to the smooth function* $\check{h}_{(H)}$ *on* $\check{M}_{(H),\mu}$ *defined by*

$$\pi^*_{(H),\mu}\check{h}_{(H)} = h_{\restriction M_{(H),\mu}} \,. \tag{4.3.21}$$

◇

*Proof.* Since $h$ is $G$-invariant, the associated Hamiltonian vector field $X_h$ is invariant, too. The calculation

$$\frac{\mathrm{d}}{\mathrm{d}t}\bigg|_t \mathrm{FL}^h_t(g\cdot m) = (X_h)_{g\cdot m} = g\,.\,(X_h)_m = \frac{\mathrm{d}}{\mathrm{d}t}\bigg|_t g\cdot \mathrm{FL}^h_t(m) \tag{4.3.22}$$

shows that the flow $\mathrm{FL}^h_t$ is $G$-equivariant (since, by assumption, it exists and is unique). Moreover, the Noether Proposition 4.2.11 implies that the flow $\mathrm{FL}^h_t$ leaves $M_\mu$ invariant. Hence, $\mathrm{FL}^h_t$ preserves $M_{(H),\mu}$ and so it projects onto a flow $\check{\mathrm{FL}}^h_t$ on $\check{M}_{(H),\mu}$. Denote the induced vector field on $\check{M}_{(H),\mu}$ by $\check{X}_{(H)}$. Since $h$ is $G$-invariant and since $\pi_{(H)}$ a surjective submersion, $h$ descends to a smooth function $\check{h}_{(H)}$ on $\check{M}_{(H),\mu}$. That $\check{X}_{(H)}$ is Hamiltonian with respect to $\check{h}_{(H)}$, indeed, is verified by a routine calculation. □

## 4.4 Application: Yang–Mills equation over a Riemannian surface

In this section, we are concerned with the moduli space of Yang–Mills connections on a principal bundle over a closed surface. This moduli space was extensively studied both from the geometric and algebraic point of view by Atiyah and Bott [AB83]. They described the structure of this moduli space using infinite-dimensional techniques inspired by symplectic reduction. We rework and extend, in the framework of Fréchet manifolds, the approach of Atiyah and Bott. As an application of the Reduction Theorem 4.2.27, we show that (a variant of) the Yang–Mills moduli space is a symplectic stratified space. Moreover, Atiyah and Bott established that Yang–Mills connections are in bijective correspondence with certain conjugacy classes of group homomorphisms. This alternative finite-dimensional model has been further investigated in [Hue93; HJ94; Jef94]. The extended moduli space construction developed in these papers yields a stratified symplectic structure on such spaces of homomorphisms. In the second part of this section, we prove that the above bijection between the Yang–Mills moduli space and these spaces of homomorphisms is an isomorphism of stratified symplectic spaces, that is, it is compatible with the additional topological, smooth and symplectic structure of both spaces under consideration. The results described in this section are based on joint work with J. Huebschmann published in extended form in [DH18].

Let $G$ be a compact connected Lie group and let $P \to M$ be a principal $G$-bundle over a closed Riemannian surface $(M, g)$. Fix an $\mathrm{Ad}_G$-invariant pairing on the Lie algebra $\mathfrak{g}$ of $G$. We are interested in connections $A \in \mathcal{C}(P)$ satisfying the Yang–Mills equation

$$\mathrm{d}_A * F_A = 0, \tag{4.4.1}$$

where $*$ refers to the Hodge star operator associated to the Riemannian metric $g$ on $M$. A special class of Yang–Mills connections is provided by connections $A$ whose curvature is of the form

$$F_A = \xi \cdot \mathrm{vol}_g, \tag{4.4.2}$$

where $\xi$ is an element of the Lie algebra $\mathfrak{z}$ of the center of $G$. We call such a connection a *central Yang–Mills connection* and refer to $\xi$ as the charge of $A$. The importance of central Yang–Mills connections for the study of the solution space of the Yang–Mills equation on a Riemannian surface comes from the following observation.

Proposition 4.4.1 ([AB83, p. 560]) *Every Yang–Mills connection $A$ on $P$ is reducible to a central Yang–Mills connection $A_\xi$ on a subbundle $P_\xi \subseteq P$.* ◇

*Proof.* The statement follows directly from the theory of symmetry reduction with $* F_A \in \Gamma^\infty(\mathrm{Ad}P)$ playing the role of the Higgs field, cf. [RS17, Section 1.6]. For completeness, let provide more details. The section $* F_A$ of $\mathrm{Ad}P$ can be viewed as a $G$-equivariant map $P \to \mathfrak{g}$. For $\xi \in \mathfrak{g}$ in the image of $* F_A$, consider the subbundle

$$P_\xi = \{p \in P : (* F_A)(p) = \xi\} \subseteq P. \tag{4.4.3}$$

It is straightforward to verify that $P_\xi$ is a reduction of $P$ to the stabilizer subgroup $G_\xi$ of $\xi$. Moreover, the Yang–Mills equation entails that $A$ is reducible to a connection $A_\xi$ on $P_\xi$. By construction, the curvature of $A_\xi$ satisfies $* F_{A_\xi} = \xi$, which shows that $A_\xi$ is a central Yang–Mills connection on $P_\xi$. □

Thus, the analysis of the moduli space of Yang–Mills connections is divided into two steps: first, for all $\xi \in \mathfrak{g}$, determine the possible reductions $P_\xi$ of $P$ to the stabilizer subgroup $G_\xi$; second, investigate the moduli space of central Yang–Mills connections on $P_\xi$. As $\xi$ is determined by the topological type of $P_\xi$ (according to Chern–Weil theory), the first step has a topological flavor and has been extensively discussed in [AB83, Section 6]. In the following, we focus on the second step from the point of view of symplectic reduction.

Recall from Example 4.3.6, that in the case $G = \mathrm{U}(1)$ the moduli space of flat connections was realized as a symplectic quotient. In a similar vein, we now describe the moduli space of central Yang–Mills connections on $P$ as a symplectically reduced space. The space $\mathcal{C}(P)$ of connections on $P$ is an affine space modeled on the tame Fréchet space $\Omega^1(M, \mathrm{Ad}P)$ of 1-forms on $M$ with values in the adjoint bundle $\mathrm{Ad}P$. As in Example 4.2.2, the 2-form $\omega$ on $\mathcal{C}(P)$ defined by the integration pairing

$$\omega_A(\alpha, \beta) = \int_M \langle \alpha \wedge \beta \rangle \tag{4.4.4}$$

for $\alpha, \beta \in \Omega^1(M, \mathrm{Ad}P)$ is a symplectic form, where $\langle \cdot \wedge \cdot \rangle$ denotes the wedge product relative to the $\mathrm{Ad}_G$-invariant pairing on $\mathfrak{g}$. The natural action on $\mathcal{C}(P)$ of the group $\mathcal{G}au(P)$ of gauge transformations of $P$ is smooth and preserves the symplectic structure $\omega$. In order to determine the associated momentum map, note that the natural pairing[1]

$$\kappa \colon \Gamma^\infty(\mathrm{Ad}P) \times \Gamma^\infty(\mathrm{Ad}P) \to \mathbb{R}, \qquad (\phi, \varrho) \mapsto -\int_M \langle \phi, \varrho \rangle \operatorname{vol}_g \tag{4.4.5}$$

identifies $\Gamma^\infty(\mathrm{Ad}P)$ as the dual of $\mathfrak{gau}(P)$. A straightforward calculation similar to the one in Example 4.2.6 verifies that the map

$$\mathcal{J} \colon \mathcal{C}(P) \to \Gamma^\infty(\mathrm{Ad}P), \qquad A \mapsto * F_A \tag{4.4.6}$$

is a momentum map for the $\mathcal{G}au(P)$-action on $\mathcal{C}(P)$. Thus, the symplectic quotient

$$\mathcal{J}^{-1}(\xi) / \mathcal{G}au(P) \tag{4.4.7}$$

at $\xi \in \mathfrak{z}$ (viewed as a constant section of $\mathrm{Ad}P$) coincides with the moduli space $\check{\mathcal{C}}_\xi(P)$ of central Yang–Mills connections with charge $\xi$. Recall that the $\mathcal{G}au(P)$-action on $\mathcal{C}(P)$ is in general not free and thus the symplectic quotient is not a smooth symplectic manifold. The following is the next best thing one could hope for.

Theorem 4.4.2 *For every $\xi \in \mathfrak{z}$, the orbit type subsets of the moduli space $\check{\mathcal{C}}_\xi(P)$ are finite-dimensional symplectic manifolds.* ◇

*Proof.* For $A \in \mathcal{J}^{-1}(\xi)$, let us verify the assumptions of Theorem 4.2.27 for the momentum map $\mathcal{J}$:

(i) Since $\xi$ is a central element, its stabilizer coincides with the whole group $\mathcal{G}au(P)$, which is a geometric tame Fréchet Lie group.

(ii) The $\mathcal{G}au(P)$-action on $\mathcal{C}(P)$ is proper and admits a slice $\mathcal{S}$ at $A$ as discussed in Section 3.4.

(iii) Since the pairing $\kappa$ is $\mathrm{Ad}_{\mathcal{G}au(P)}$-invariant, the coadjoint action coincides with the adjoint action and thus is clearly linear.

[1] The sign in front of the integral turns out to be a convenient choice in the sequel.

(iv) The chain (4.2.70) takes the following form here:

$$0 \longrightarrow \Omega^0(M, \mathrm{Ad}P) \xrightarrow{-\mathrm{d}_B} \Omega^1(M, \mathrm{Ad}P) \xrightarrow{*\mathrm{d}_B} \Omega^0(M, \mathrm{Ad}P), \longrightarrow 0 \tag{4.4.8}$$

where the connection $B$ on $P$ is an element of the slice $\mathcal{S}$. This is clearly a chain of differential operators tamely parametrized by $B \in \mathcal{S}$. Moreover, for $B = A$, this chain is an elliptic complex, because we have

$$\mathrm{d}_A \, \mathrm{d}_A \eta = [F_A, \eta] = [\xi, \eta] \, \mathrm{vol}_g = 0 \tag{4.4.9}$$

for every $\eta \in \Gamma^\infty(\mathrm{Ad}P)$.

(v) Arguments similar to those used in Example 4.2.17 show that $\mathfrak{gau}(P) \,.\, A$ is symplectically closed in $\mathrm{T}_A\mathcal{C}(P)$ and that the image of $\delta_A J$ is $\mathrm{L}^2$-closed in $\mathfrak{gau}(P)$.

Thus, by Theorem 4.2.27, the momentum map $\mathcal{J}$ can be bought into a strong MGS normal form at $A$. Now the claim follows from the Singular Reduction Theorem 4.3.5. □

The statement that the top stratum of $\check{\mathcal{C}}_\xi(P)$ is endowed with a natural symplectic reduction was already obtained by Atiyah and Bott [AB83, p. 587]. The symplectic nature of the singular strata has been established in [Hue96, Theorem 1.2] in the Sobolev framework by reducing the problem to a *finite-dimensional* symplectic quotient.

Remark 4.4.3 Let us spell out the local structure of $\check{\mathcal{C}}_\xi(P)$ at a point $[A]$ as given by Theorems 4.2.27 and 4.3.5. According to Proposition 4.2.15, the space Ker in the strong MGS normal form is identified with the middle homology of the chain (4.4.8) at $B = A$, that is,

$$\mathrm{Ker} = \mathrm{Ker}\big(*\,\mathrm{d}_A : \Omega^1(M, \mathrm{Ad}P) \to \Omega^0(M, \mathrm{Ad}P)\big) / \mathrm{Im}\, \mathrm{d}_A \equiv \mathrm{H}^1_A(M, \mathrm{Ad}P). \tag{4.4.10}$$

Note that this (co)homology group is only well-defined because $A$ is a central Yang–Mills connection. Similarly, the Lie algebra of the stabilizer subgroup $\mathcal{G}au_A(P)$ of $A$ is given by

$$\mathfrak{gau}_A(P) = \mathrm{Ker}\big(\mathrm{d}_A : \Omega^0(M, \mathrm{Ad}P) \to \Omega^1(M, \mathrm{Ad}P)\big) \equiv \mathrm{H}^0_A(M, \mathrm{Ad}P). \tag{4.4.11}$$

The symplectic structure $\omega_A$ on $\mathrm{T}_A\mathcal{C}(P)$ restricts to a symplectic structure $\bar{\omega}_A$ on $\mathrm{H}^1_A(M, \mathrm{Ad}P)$. By (4.4.4), we have

$$\bar{\omega}_A\big([\alpha], [\beta]\big) = \int_M \langle \alpha \wedge \beta \rangle. \tag{4.4.12}$$

In other words, $\bar{\omega}_A$ is the non-abelian generalization of the intersection form. Nondegeneracy of $\bar{\omega}_A$ can be verified by a direct calculation or using the Bifurcation Lemma 4.2.14 (ii). It is clear that the action of $\mathcal{G}au(P)$ on $\mathcal{C}(P)$ induces a symplectic action of $\mathcal{G}au_A(P)$ on $\mathrm{H}^1_A(M, \mathrm{Ad}P)$:

$$\lambda \cdot \alpha = \frac{\mathrm{d}}{\mathrm{d}\varepsilon}\bigg|_0 \lambda \cdot (A + \varepsilon\alpha) = \mathrm{Ad}_\lambda \alpha \tag{4.4.13}$$

for $\lambda \in \mathcal{G}au_A(P)$ and $[\alpha] \in \mathrm{H}^1_A(M, \mathrm{Ad}P)$. The infinitesimal action of $\mathfrak{gau}_A(P)$ is given by the Lie bracket so that the associated momentum map is defined by

$$\mathcal{J}_{\mathrm{sing}} \colon \mathrm{H}^1_A(M, \mathrm{Ad}P) \to \mathrm{H}^0_A(M, \mathrm{Ad}P), \qquad [\alpha] \mapsto \frac{1}{2} *[\alpha \wedge \alpha], \tag{4.4.14}$$

where $[\cdot \wedge \cdot]$ denotes the wedge product of $\mathrm{Ad}P$-valued differential forms relative to the Lie bracket. Indeed, by Example 4.2.5, we have

$$\begin{aligned}\kappa\big(\mathcal{J}_{\mathrm{sing}}([\alpha]), \xi\big) &= \frac{1}{2}\bar{\omega}_A\big([\alpha], [\mathrm{ad}_\xi \alpha]\big) \\ &= \frac{1}{2}\int_M \langle \alpha \wedge [\xi, \alpha]\rangle \\ &= -\frac{1}{2}\int_M \langle [\alpha \wedge \alpha], \xi\rangle \\ &= \frac{1}{2}\kappa\big(*[\alpha \wedge \alpha], \xi\big).\end{aligned} \tag{4.4.15}$$

Now, (4.3.6) entails that, locally, the moduli space of central Yang–Mills connections is modeled on

$$\mathcal{J}^{-1}_{\mathrm{sing}}(0) / \mathcal{G}au_A(P). \tag{4.4.16}$$

Such a local description of $\check{\mathcal{C}}_\xi(P)$ has already been obtained by Huebschmann [Hue95, Theorem 2.32] in the framework of Sobolev spaces using techniques similar to those employed in our construction of the MGS normal form. ◇

Recall from Example 4.3.6 that, for $G = \mathrm{U}(1)$, the moduli space of flat connections is identified with the space of group homomorphisms from $\pi_1(M)$ to $\mathrm{U}(1)$. Following [AB83], we discuss now a similar identification of $\check{\mathcal{C}}_\xi(P)$ for general compact groups $G$. For this purpose, fix a point $m_0 \in M$, denote the genus of $M$ by $l$, and choose a canonical system $(u_1, v_1, \ldots, u_l, v_l)$ of closed curves in $M$ based at $m_0$. That is to say, cutting $M$ successively along these curves yields a planar polygon. In this picture, the curve

$$u_1 \cdot v_1 \cdot u_1^{-1} v_1^{-1} \cdots u_l \cdot v_l \cdot u_l^{-1} v_l^{-1} \tag{4.4.17}$$

corresponds to the boundary of the polygon. We will denote the homotopy class of this boundary curve by $r$. As a consequence of the Seifert–van Kampen theorem, the fundamental class $\pi_1(M, m_0) \equiv \pi_1$ of $M$ is generated by the homotopy classes $([u_1], [v_1], \ldots, [u_l], [v_l])$ subject to the condition $r = 1$. Consider the finite group $\Gamma$ generated by $([u_1], [v_1], \ldots, [u_l], [v_l])$ and an element $J$ subject to the relation $r = J$. Note that $J$ is central in $\Gamma$. The natural projection from $\Gamma$ to $\pi_1$ is a group homomorphism with kernel $\mathbb{Z}\langle r\rangle$ (the free group generated by $r$). Relative to the embedding $\mathbb{Z}\langle r\rangle \to \mathbb{R}$ of abelian groups defined by $r \mapsto 1$, let $\Gamma_{\mathbb{R}}$ be the 1-dimensional Lie group characterized by the requirement that

$$\begin{array}{ccccccccc} 1 & \longrightarrow & \mathbb{Z}\langle r\rangle & \longrightarrow & \Gamma & \longrightarrow & \pi_1 & \longrightarrow & 1 \\ & & \downarrow & & \downarrow & & \downarrow & & \\ 0 & \longrightarrow & \mathbb{R} & \longrightarrow & \Gamma_{\mathbb{R}} & \longrightarrow & \pi_1 & \longrightarrow & 1 \end{array} \tag{4.4.18}$$

be a commutative diagram of central extensions.

Remark 4.4.4 By [Mor91, Theorem 2.2], the group $\Gamma_{\mathbb{R}}$ admits the following geometric description. Let $\mathcal{L}^{\infty}_{m_0}(M)$ be the space of closed, piecewise smooth loops in $M$ based at $m_0$. Define the equivalence relation $\sim_{\mathrm{vol}_g}$ on $\mathcal{L}^{\infty}_{m_0}(M)$ by declaring two piecewise smooth loops $\gamma_1$ and $\gamma_2$ to be equivalent if and only if there exists a homotopy $h\colon [0,1] \times [0,1] \to M$ between $\gamma_1$ and $\gamma_2$ such that

$$\int_{[0,1]\times[0,1]} h^* \mathrm{vol}_g = 0. \tag{4.4.19}$$

Then, $\Gamma_{\mathbb{R}}$ is identified with the quotient of $\mathcal{L}^{\infty}_{m_0}(M)$ relative to the equivalence relation $\sim_{\mathrm{vol}_g}$. ◊

We continue to denote the image of the central element $r \in \Gamma$ under the embedding $\Gamma \to \Gamma_{\mathbb{R}}$ by $r$. With this notation, the embedding of $\mathbb{R}$ into $\Gamma_{\mathbb{R}}$ is given by $t \mapsto tr$ for $t \in \mathbb{R}$. Note that this embedding yields a natural isomorphism of $\mathbb{R}/\mathbb{Z}$ with $\Gamma_{\mathbb{R}}/\Gamma$. Therefore, the projections onto $\Gamma_{\mathbb{R}}/\Gamma$ and $\pi_1$ yield an exact sequence

$$1 \longrightarrow \mathbb{Z}\langle r\rangle \longrightarrow \Gamma_{\mathbb{R}} \longrightarrow \mathbb{R}/\mathbb{Z} \times \pi_1 \longrightarrow 1 \tag{4.4.20}$$

In order to return to the geometric picture, assume that $\mathrm{vol}_g$ is normalized in such a way that the total volume is 1. Thus, $\mathrm{vol}_g$ is a closed 2-form that defines an integral class in $\mathrm{H}^2(M, \mathbb{R})$ and we can always find a principal $\mathbb{R}/\mathbb{Z}$-bundle $Q \to M$ and a connection on $Q$ with $\mathrm{vol}_g$ as its curvature. This construction is a basic result of the theory of geometric quantization and, in this context, $Q$ is referred to as the prequantum bundle, see e.g. [Woo97, Proposition 8.3.1]. The group $\Gamma$ is isomorphic to the fundamental group of $Q$. Indeed, since the curvature determines the connection only up to an element of $\mathrm{H}^1(M, \mathbb{R}/\mathbb{Z})$, we can choose the connection on $Q$ in such a way that the holonomies of the curves $u_i$ and $v_i$ are trivial (with respect to a fixed point of $Q$ in the fiber over $m_0$). The horizontal lifts of $u_i$ and $v_i$ are then closed curves and we can take their homotopy classes as generators of $\pi_1(Q) \simeq \Gamma$. Since the universal covering $\tilde{M} \to M$ of $M$ is a flat principal $\pi_1$-bundle, the product connection on the principal $\mathbb{R}/\mathbb{Z} \times \pi_1$-bundle $Q \times_M \tilde{M}$ still has curvature $\mathrm{vol}_g$. Lifting the structure group to $\Gamma_\mathbb{R}$ along (4.4.20) yields a principal $\Gamma_\mathbb{R}$-bundle $Q_M \to M$ with connection $A_M$ and curvature $F_{A_M} = \mathrm{vol}_g$. Comparing with (4.4.2), we see that $A_M$ a central Yang–Mills connection with charge $1 \in \mathbb{R}$. The pair $(Q_M, A_M)$ plays a universal role as we shall now explain.

For every $\xi \in \mathfrak{z}$, let $\mathrm{Hom}_\xi(\Gamma_\mathbb{R}, G)$ denote the space of Lie group homomorphisms $\chi\colon \Gamma_\mathbb{R} \to G$ with the property that

$$\chi(tr) = \exp(t\xi) \tag{4.4.21}$$

for all $t \in \mathbb{R}$. The evaluation of a homomorphism $\chi \in \mathrm{Hom}_\xi(\Gamma_\mathbb{R}, G)$ on $(u_1, v_1, \dots, u_l, v_l)$ yields an embedding of $\mathrm{Hom}_\xi(\Gamma_\mathbb{R}, G)$ into $G^{2l}$ as the subset

$$\left\{(a_1, b_1, \dots, a_l, b_l) : \; a_1 b_1 a_1^{-1} b_1^{-1} \cdots a_l b_l a_l^{-1} b_l^{-1} = \exp(\xi)\right\}. \tag{4.4.22}$$

Throughout, we endow $\mathrm{Hom}_\xi(\Gamma_\mathbb{R}, G)$ with the relative topology and identify it with its image under the embedding into $G^{2l}$. Note that $G$ acts naturally on $\mathrm{Hom}_\xi(\Gamma_\mathbb{R}, G)$ by conjugation. This action corresponds to the diagonal conjugation action under the embedding of $\mathrm{Hom}_\xi(\Gamma_\mathbb{R}, G)$ into $G^{2l}$.

Given $\chi \in \mathrm{Hom}_\xi(\Gamma_\mathbb{R}, G)$, we can form the associated principal $G$-bundle $P_\chi = Q_M \times_\chi G$. It turns out that the map $\chi \mapsto P_\chi$ defines a bijection between the connected components of $\mathrm{Hom}_\xi(\Gamma_\mathbb{R}, G)$ and the topological types of principal $G$-bundles on $M$ having $\xi \in \mathfrak{z}$ as its corresponding characteristic class, see [DH18, Proposition 3.1]. We denote by $\mathrm{Hom}_\xi(\Gamma_\mathbb{R}, G)_P$ the connected component that corresponds to the principal G-bundle $P$. The reference connection $A_M$ on $Q_M$ yields a principal connection $A_\chi$ on $P_\chi$. Since the Lie algebra homomorphism $\mathrm{T}_r\chi\colon \mathbb{R} \to \mathfrak{g}$ associated to $\chi \in \mathrm{Hom}_\xi(\Gamma_\mathbb{R}, G)$ satisfies $\mathrm{T}_r\chi(1) = \xi$, we find

$$F_{A_\chi} = \xi \cdot \mathrm{vol}_g\,. \tag{4.4.23}$$

Thus, $A_\chi$ is a central Yang–Mills connection with charge $\xi$. If $\chi$ is an element of $\mathrm{Hom}_\xi(\Gamma_\mathbb{R}, G)_P$, then $P_\chi$ has the same topological type as $P$ so that there exists a vertical automorphism $P_\chi \to P$ of principal $G$-bundles. The push-forward of $A_\chi$ along this automorphism yields a central Yang–Mills connection on $P$, which we continue to denote by $A_\chi$. Note that the connection $A_\chi$ on $P$ depends on the choice of the vertical automorphism $P_\chi \to P$ and thus only its class modulo gauge transformations is intrinsically well-defined. The crucial insight of Atiyah and Bott [AB83, Theorem 6.7] (see also [Hue94, Theorem 2.1]) was that every central Yang–Mills connection on $P$ can be obtained in this way (up to gauge transformations).

Proposition 4.4.5 *For every $\xi \in \mathfrak{z}$, the assignment $\chi \mapsto A_\chi$ yields a bijection between the space $\mathrm{Hom}_\xi(\Gamma_\mathbb{R}, G)_P / G$ of conjugacy classes of homomorphism $\Gamma_\mathbb{R} \to G$ and the moduli space $\check{\mathcal{C}}_\xi(P)$ of central Yang–Mills connections on $P$ with charge $\xi$.* ◇

Let us recall how, conversely, every central Yang–Mills connection $A$ on $P$ gives rise to a Lie group homomorphism $\chi_A : \Gamma_\mathbb{R} \to G$. Fix a point $q_0 \in Q_M$ in the fiber over $m_0$. For every $x \in \Gamma_\mathbb{R}$, there exists an $A_M$-horizontal path $\hat{\gamma}$ in $Q_M$ from $q_0$ to $q_0 \cdot x$. Project $\hat{\gamma}$ to a necessarily closed curve $\gamma$ in $M$ and define $\chi_A(x) \in G$ to be the holonomy of $\gamma$. By construction, $\chi_A(r)$ coincides with the holonomy of the boundary curve $v_1 \cdot u_1^{-1} v_1^{-1} \cdots u_l \cdot v_l \cdot u_l^{-1} v_l^{-1}$ and thus is equal to

$$\chi_A(r) = \int_M F_A = \int_M \xi \cdot \mathrm{vol}_g = \xi. \tag{4.4.24}$$

Hence, the resulting group homomorphism $\chi_A : \Gamma_\mathbb{R} \to G$ is an element of $\mathrm{Hom}_\xi(\Gamma_\mathbb{R}, G)$. Since $\Gamma_\mathbb{R}$ is abelian, the homomorphism $\chi_A$ does not depend on the choice of $q_0$. However, choosing a different reference point in the fiber of $P$ over $m_0$ for the holonomy $\chi_A(x) \in G$ yields a group homomorphism $\Gamma_\mathbb{R} \to G$ conjugate to $\chi_A$.

Let $\mathrm{Rep}_\xi(\Gamma_\mathbb{R}, G)_P = \mathrm{Hom}_\xi(\Gamma_\mathbb{R}, G)_P / G$. The space $\mathrm{Rep}_\xi(\Gamma_\mathbb{R}, G)_P$ decomposes naturally into orbit types according to the $G$-action on $\mathrm{Hom}_\xi(\Gamma_\mathbb{R}, G)_P$ by conjugation. Moreover, the extended moduli space construction developed in [Hue93; HJ94; Jef94] endows $\mathrm{Rep}_\xi(\Gamma_\mathbb{R}, G)_P$ with a natural stratified symplectic structure. In [AMM98], the extended moduli space construction was reinterpreted as a quotient procedure for so-called quasi-Hamiltonian $G$-space. On the other hand, we have seen that the moduli space $\check{\mathcal{C}}_\xi(P)$ arises as a symplectic quotient of an infinite-dimensional symplectic space. Note that the character of both constructions is very different: reduction of a finite-dimensional quasi-Hamiltonian system vs. infinite-dimensional symplectic reduction. Thus, it is natural to ask whether the bijection of Proposition 4.4.5 respects the manifold and the symplectic structure of both quotients.

The first step towards an answer to this question is the following, which we quote without proof.

Proposition 4.4.6 ([DH18, Theorem 5.2]) *Let $p_0 \in P$ be a point in the fiber over $m_0$. For every closed curve $\gamma$ in $M$ based at $m_0$, the map*

$$\mathrm{Hol}^{p_0}_\gamma : \mathcal{C}(P) \to G, \quad A \mapsto \mathrm{Hol}^{p_0}_\gamma(A) \tag{4.4.25}$$

*that assigns to $A$ the holonomy $\mathrm{Hol}^{p_0}_\gamma(A)$ of $\gamma$ is a smooth map with left logarithmic derivative*

$$\delta_A \mathrm{Hol}^{p_0}_\gamma(\alpha) = \int_{\gamma_A} \alpha, \tag{4.4.26}$$

*where $\gamma_A$ denotes the $A$-horizontal lift of $\gamma$ to $p_0$ and $\alpha \in \mathrm{T}_A\mathcal{C}(P) \simeq \Omega^1(M, \mathrm{Ad}P)$ is viewed as a $\mathfrak{g}$-valued 1-form on $P$.* ◊

Recall the canonical system $(u_1, v_1, \ldots, u_l, v_l)$ of closed curves in $M$ based at $m_0$. The Wilson loop map $\mathcal{W}_{p_0} \colon \mathcal{C}(P) \to G^{2l}$ is defined by

$$\mathcal{W}_{p_0}(A) = \big(\mathrm{Hol}^{p_0}_{u_1}(A), \mathrm{Hol}^{p_0}_{v_1}(A), \ldots, \mathrm{Hol}^{p_0}_{u_l}(A), \mathrm{Hol}^{p_0}_{v_l}(A)\big). \tag{4.4.27}$$

As a consequence of Proposition 4.4.6, $\mathcal{W}_{p_0}$ is a smooth map. Recall from (4.4.22) that $\mathrm{Hom}_\xi(\Gamma_\mathbb{R}, G)$ is identified with a subset of $G^{2l}$. If $A$ is a central Yang–Mills connection with charge $\xi$, then the holonomy of the boundary curve $r$ of the fundamental polygon is given by

$$\mathrm{Hol}^{p_0}_r(A) = \exp\left(\int_M F_A\right) = \exp(\xi). \tag{4.4.28}$$

Thus, $\mathcal{W}_{p_0}(A)$ lies in $\mathrm{Hom}_\xi(\Gamma_\mathbb{R}, G) \subseteq G^{2l}$ for every central Yang–Mills connection $A$. It is not hard to see that the corresponding homomorphism $\Gamma_\mathbb{R} \to G$ coincides with the homomorphism $\chi_A$ that was constructed above.

We maintain the choice of a point $p_0 \in P$ in the fiber over $m_0$. Recall that the evaluation at $p_0$ of a gauge transformation $\lambda \in \mathcal{G}au(P)$, seen as a $G$-equivariant map $P \to G$, yields a Lie group homomorphism

$$\mathrm{ev}_{p_0} \colon \mathcal{G}au(P) \to G. \tag{4.4.29}$$

The group $\mathcal{G}au_{m_0}(P)$ of *pointed gauge transformations* is the kernel of $\mathrm{ev}_{p_0}$. This group is independent of the particular choice of $p_0$ and depends only on $m_0$; this justifies our notation $\mathcal{G}au_{m_0}(P)$. By [Nee06, Proposition IV.3.4], $\mathcal{G}au_{m_0}(P)$ is a normal, locally exponential Lie subgroup of $\mathcal{G}au(P)$. Note that the natural action of $\mathcal{G}au_{m_0}(P)$ on $\mathcal{C}(P)$ is free. Moreover, the slices for the $\mathcal{G}au(P)$-action can be modified to yield slices for the $\mathcal{G}au_{m_0}(P)$-action, see [DH18, Proposition 11.9]. Hence, according to Proposition A.2.7, the orbit space

$$\check{\mathcal{C}}_{m_0}(P) = \mathcal{C}(P)/\mathcal{G}au_{m_0}(P) \tag{4.4.30}$$

is a smooth Fréchet manifold. Using a similar notation, the quotient of the space of central Yang–Mills connections with charge $\xi$ by the group of pointed gauge transformation will be denoted by $\check{\mathcal{C}}_{\xi,m_0}(P)$. Clearly, the Wilson loop map $\mathcal{W}_{p_0}$ is invariant under the action of pointed gauge transformations and thus descends to a smooth map $\check{\mathcal{W}}_{p_0} \colon \check{\mathcal{C}}_{m_0}(P) \to G^{2l}$.

Proposition 4.4.7 *The Wilson loop map $\mathcal{W}_{p_0} \colon \mathcal{C}(P) \to G^{2l}$ yields a smooth map*

$\check{\mathcal{W}}_{p_0}\colon \check{\mathcal{C}}_{m_0}(P) \to G^{2l}$ *which restricts to a G-equivariant homeomorphism*

$$\check{\mathcal{C}}_{\xi,m_0}(P) \to \mathrm{Hom}_{\xi}(\Gamma_{\mathbb{R}}, G)_P. \tag{4.4.31}$$

*Accordingly,* $\check{\mathcal{W}}_{p_0}$ *induces a homeomorphism*

$$\check{\mathcal{C}}_{\xi}(P) \to \mathrm{Rep}_{\xi}(\Gamma_{\mathbb{R}}, G)_P. \tag{4.4.32}$$

*Moreover, with respect to the G-orbit type stratifications on both sides, this homeomorphism is an isomorphism of stratified spaces which, on each stratum, restricts to a diffeomorphism onto the corresponding stratum of the target space.* ◇

*Proof.* We only sketch the main arguments and refer the reader to [DH18, Theorem 7.1] for details.

Arguments similar to the ones discussed in connection with Proposition 4.4.5 show that $\check{\mathcal{W}}_{p_0}$ restricts to a bijection between $\check{\mathcal{C}}_{\xi,m_0}(P)$ and $\mathrm{Hom}_{\xi}(\Gamma_{\mathbb{R}}, G)_P$. Since $\check{\mathcal{W}}_{p_0}$ is a smooth map, its restriction to $\check{\mathcal{C}}_{\xi,m_0}(P)$ is continuous. By the Strong Uhlenbeck Compactness Theorem [Uhl82], for any sequence $(A_j)$ of smooth Yang–Mills connections with constant curvature[1], there exists a subsequence $(A_{j_i})$ and a sequence $(\lambda_i)$ of smooth gauge transformations such that the sequence $(\lambda_i \cdot A_{j_i})$ converges uniformly in the Fréchet topology to a smooth connection $A \in \mathcal{C}(P)$. Therefore, the space $\check{\mathcal{C}}_{\xi,m_0}(P)$ is compact (in the Fréchet topology) and, consequently, the Wilson loop map restricts to a homeomorphism onto the compact subspace $\mathrm{Hom}_{\xi}(\Gamma_{\mathbb{R}}, G)_P$ of the compact manifold $G^{2l}$ as asserted. Moreover, the $G$-equivariant homeomorphism $\check{\mathcal{C}}_{\xi,m_0}(P) \to \mathrm{Hom}_{\xi}(\Gamma_{\mathbb{R}}, G)_P$ induces a homeomorphism $\check{\mathcal{C}}_{\xi}(P) \to \mathrm{Rep}_{\xi}(\Gamma_{\mathbb{R}}, G)_P$ which is compatible with the orbit type stratifications on both sides. Since the Wilson loop map is smooth, it restricts, on each stratum, to a diffeomorphism onto the corresponding stratum of the target space. □

Remark 4.4.8 The usage of the Uhlenbeck Compactness Theorem in the above proof can be avoided using a detailed slice analysis, see [DH18, Lemma 7.5]. ◇

Finally, let us comment on the behavior of the isomorphism $\check{\mathcal{C}}_{\xi}(P) \to \mathrm{Rep}_{\xi}(\Gamma_{\mathbb{R}}, G)_P$ with regard to the symplectic structure on both quotients. For simplicity, we restrict attention to the abelian case $G = \mathrm{U}(1)$ and to $l = 1$. In this case, the extended moduli space construction yields a symplectic form $\sigma$ on $\mathrm{U}(1) \times \mathrm{U}(1)$ given by

$$\sigma_{(a,b)}(a \,.\, \xi_1 + b \,.\, \eta_1, a \,.\, \xi_2 + b \,.\, \eta_2) = \langle \xi_1, \eta_2 \rangle - \langle \xi_2, \eta_1 \rangle \tag{4.4.33}$$

for $a, b \in \mathrm{U}(1)$ and $\xi_1, \xi_2, \eta_1, \eta_2 \in \mathfrak{u}(1)$. That is, $\sigma$ coincides with the (normalized) volume form on the torus $\mathrm{U}(1) \times \mathrm{U}(1)$. Using the formula (4.4.26) for the

[1] In fact, the hypothesis that the curvature be uniformly $L^2$-bounded suffices.

derivative of the Wilson loop map, we obtain for the pull-back of $\sigma$ to $\mathcal{C}(P)$:

$$(\mathcal{W}^*_{p_0}\sigma)_A(\alpha,\beta) = \left\langle \int_u \alpha, \int_v \beta \right\rangle - \left\langle \int_u \beta, \int_v \alpha \right\rangle, \tag{4.4.34}$$

where $A \in \mathcal{C}(P)$, $\alpha, \beta \in \Omega^1(M)$ and $u, v$ are the canonical closed curves on the torus $M$. Moreover, if $\alpha$ and $\beta$ are closed, a simple calculation using Hodge theory yields the identity

$$(\mathcal{W}^*_{p_0}\sigma)_A(\alpha,\beta) = \int_M \alpha \wedge \beta = \omega(\alpha,\beta). \tag{4.4.35}$$

Note that the space of closed 1-forms is the tangent space to $\mathcal{J}^{-1}(\xi)$. We thus find

$$(\mathcal{W}^*_{p_0}\sigma)_{\restriction \mathcal{J}^{-1}(\xi)} = \omega. \tag{4.4.36}$$

Upon projecting with respect to the action of $\mathcal{G}au(P)$, we see that the diffeomorphism $\check{\mathcal{C}}_\xi(P) \longrightarrow \mathrm{Rep}_\xi(\Gamma_\mathbb{R}, \mathrm{U}(1))_P$ intertwines the symplectic structures. A similar reasoning works for higher genus $l \geq 1$ and non-abelian structure group $G$. We thus obtain in summary.

Theorem 4.4.9 *The Wilson loop map $\check{\mathcal{W}}_{p_0}$ induces an isomorphism*

$$\check{\mathcal{C}}_\xi(P) \longrightarrow \mathrm{Rep}_\xi(\Gamma_\mathbb{R}, G)_P \tag{4.4.37}$$

*of stratified symplectic spaces.* ◇

Thus, we conclude that the bijection $\check{\mathcal{C}}_\xi(P) \longrightarrow \mathrm{Rep}_\xi(\Gamma_\mathbb{R}, G)_P$ established by Atiyah and Bott [AB83] is compatible with the additional stratified symplectic structure of both quotients.

# Singular Cotangent Bundle Reduction 5

In most applications in physics, the phase space is a cotangent bundle $\mathrm{T}^*Q$ over the configuration space $Q$ of the system. When performing symplectic reduction for that case, it is of interest to connect the geometry of the symplectically reduced space $\mathrm{T}^*Q /\!\!/_{\mu} G$ to the properties of the quotient $Q/G$. In finite dimensions and for proper free $G$-actions, this connection is well understood: the reduction $\mathrm{T}^*Q /\!\!/_0 G$ at zero is symplectomorphic to the cotangent bundle $\mathrm{T}^*(Q/G)$ (with its canonical symplectic form) of the reduced configuration space $Q/G$, see, for example, [Mar+07, Section 2] for details. The singular case is more complicated due to the occurrence of new phenomena that are absent in the regular case. In finite dimensions, Perlmutter, Rodriguez-Olmos, and Sousa-Dias [PRS07] have shown that the fibered structure of the cotangent bundle yields a refinement of the usual orbit-momentum type strata into so-called seams. The principal seam is symplectomorphic to a cotangent bundle while the singular seams are coisotropic submanifolds of the corresponding symplectic stratum. In this chapter, we extend these results concerning singular symplectic reduction of cotangent bundles to the case of infinite-dimensional Fréchet[1] manifolds. In the finite-dimensional setting, the proof that the seams are manifolds follows from an appropriate symplectic slice theorem, see [Sch07; RT17]. The latter theorem provides a normal form both for the symplectic structure and for the momentum map, and additionally the local symplectomorphism bringing the system to normal form is adapted to the fiber structure of the cotangent bundle. When passing to the Fréchet context, it is impossible to follow the approach of [PRS07; Sch07] because a few essential tools required in the construction are not readily available for Fréchet manifolds, see Remark 5.2.1 for details. Similar to the strategy in Chapter 4, at the root of our approach lies the observation that, for the cotangent bundle reduction procedure, it is not essential to bring the symplectic structure into a normal form. Instead, we exploit the additional structure of the cotangent bundle to directly construct the normal form of the momentum map for the lifted action. Using the harmonic oscillator as an example, we also discuss the influence of the singular seams on the dynamics. Finally, we apply our general theory to the singular cotangent bundle reduction of Yang–Mills–Higgs theory. The content of this chapter will be published in a slightly modified version as part of [DR18b].

[1] Most results of this chapter readily generalize to infinite-dimensional manifolds modeled on Mackey complete, locally convex spaces.

## 5.1 Problems in infinite dimensions. The setting

Let $Q$ be a Fréchet manifold. The tangent bundle $\mathrm{T}Q$ of $Q$ is a smooth manifold in such a way that the projection $\mathrm{T}Q \to Q$ is a smooth locally trivial bundle. However, the topological dual bundle $\mathrm{T}'Q := \bigsqcup_{q \in Q} (\mathrm{T}_q Q)'$ is not a *smooth* fiber bundle for non-Banach manifolds $Q$, cf. [Nee06, Remark I.3.9] and Appendix A.3. As a substitute, we say that a smooth Fréchet bundle $\overset{\star}{\tau} \colon \mathrm{T}^* Q \to Q$ is a cotangent bundle if there exists a fiberwise non-degenerate pairing with $\mathrm{T}Q$, see Appendix A.3 for details. We note that according to this definition the cotangent bundle is no longer canonically associated to $Q$ but requires the choice of a bundle $\mathrm{T}^*Q$ and of a pairing $\mathrm{T}^*Q \times_Q \mathrm{T}Q \to \mathbb{R}$. Similarly to the finite-dimensional case, for a given cotangent bundle $\mathrm{T}^*Q$, the formula

$$\theta_p(v) = \langle p, \mathrm{T}_p \overset{\star}{\tau}(v) \rangle, \qquad v \in \mathrm{T}_p(\mathrm{T}^*Q)\,, \tag{5.1.1}$$

defines a smooth 1-form $\theta$ on $\mathrm{T}^*Q$. Furthermore, $\omega = \mathrm{d}\theta$ is a symplectic form.

Next, assume that a Fréchet Lie group $G$ acts smoothly on $Q$. We continue to use the dot notation to write this action as $(g, q) \mapsto g \cdot q$. Throughout this chapter, we assume that the $G$-action on $Q$ is proper and that it admits a slice at every point. By linearization, we also get a smooth action of $G$ on the tangent bundle $\mathrm{T}Q$, which we write using the lower dot notation as $g \,.\, v \in \mathrm{T}_{g \cdot q} Q$ for $g \in G$ and $v \in \mathrm{T}_q Q$. The action on $\mathrm{T}Q$ induces a $G$-action on $\mathrm{T}^*Q$ by requiring that the pairing be left invariant, that is,

$$\langle g \cdot p, v \rangle = \langle p, g^{-1} \,.\, v \rangle \tag{5.1.2}$$

for $p \in \mathrm{T}^*_{g^{-1} \cdot q} Q$ and $v \in \mathrm{T}_q Q$. In order for this equation to define a smooth action on $\mathrm{T}^*Q$, the action $\mathrm{T}_{g^{-1} \cdot q} Q \to \mathrm{T}_q Q$ needs to be weakly continuous with respect to the pairing $\langle \cdot\,, \cdot \rangle$ for every $g \in G$.

In finite dimensions, there always exists a $G$-equivariant diffeomorphism between $\mathrm{T}Q$ and $\mathrm{T}^*Q$. Therefore, the orbit types of $\mathrm{T}Q$ and $\mathrm{T}^*Q$ coincide. In the infinite-dimensional setting, a vector space may not be isomorphic to its dual and thus the orbit types may differ. Indeed, we now give an example where the action on the dual space has more orbit types.

Example 5.1.1 Consider the space $\ell^1$ of real-valued doubly infinite sequences $(x_n)_{n \in \mathbb{Z}}$ satisfying

$$\|(x_n)\|_1 := \sum_{n \in \mathbb{Z}} |x_n| < \infty. \tag{5.1.3}$$

With respect to the norm $\|\cdot\|_1$, $\ell^1$ is a Banach space. The topological dual is isomorphic to the space $\ell^\infty$ of all bounded sequences $(\alpha_n)_{n \in \mathbb{Z}}$, which is a Banach space with respect to the uniform norm

$$\|(\alpha_n)\|_\infty := \sup_{n \in \mathbb{Z}} |\alpha_n|. \tag{5.1.4}$$

The pairing between $\ell^1$ and $\ell^\infty$ is given by

$$\langle (\alpha_n), (x_n) \rangle := \sum_{n \in \mathbb{Z}} \alpha_n x_n, \qquad (\alpha_n) \in \ell^\infty, (x_n) \in \ell^1. \tag{5.1.5}$$

The group $G = \mathbb{Z}$ of integers acts on $\ell^1$ via the shift operators $T_j \colon \ell^1 \longrightarrow \ell^1$, with $j \in \mathbb{Z}$, defined by

$$T_j(x_n) := (x_{n+j}) . \tag{5.1.6}$$

This action is linear and self-adjoint in the sense that the dual action on $\ell^\infty$ is also given by the shift operators $T_j \colon \ell^\infty \longrightarrow \ell^\infty$. The actions on $\ell^1$ and $\ell^\infty$ are continuous for every topology on $\mathbb{Z}$ (in particular, we may endow $\mathbb{Z}$ with the cofinite topology in which it is compact). Recall that all subgroups of $\mathbb{Z}$ are of the form $k\mathbb{Z}$ for some integer $k \geq 0$. A sequence $(x_n)$ has stabilizer $k\mathbb{Z}$ if and only if it is $k$-periodic. The subgroups $\{0\}$ ($k = 0$) and $\mathbb{Z}$ ($k = 1$) can easily be realized as stabilizer subgroups of some $(x_n) \in \ell^1$. However, none of the other subgroups $k\mathbb{Z}$ with $k > 1$ occurs as the stabilizer subgroup of the action on $\ell^1$. Indeed, if a sequence $(x_n)$ is non-zero and $k$-periodic for $k > 1$, then it is divergent and hence never an element of $\ell^1$. On the other hand, every periodic sequence is bounded and thus all subgroups $k\mathbb{Z}$ for $k \geq 0$ occur as stabilizers of the action on the dual space $\ell^\infty$. ◊

In order to exclude such pathological phenomena that make it impossible to connect the orbit types of $Q$ with the ones for the lifted action on $\mathrm{T}^*Q$, in the following, we will assume that there exists a $G$-equivariant diffeomorphism[1] between $\mathrm{T}Q$ and $\mathrm{T}^*Q$; an assumption that holds in the applications we are interested in. Note that for a Hilbert manifold $Q$, a $G$-invariant scalar product yields such a $G$-equivariant diffeomorphism.

By Proposition A.3.3, the lifted action of $G$ on the cotangent bundle $\mathrm{T}^*Q$ preserves the canonical symplectic form $\omega$. In order to define the momentum map, we need to specify a dual pair $\kappa(\mathfrak{g}^*, \mathfrak{g})$. With respect to this data, the momentum map $J \colon \mathrm{T}^*Q \longrightarrow \mathfrak{g}^*$, if it exists, satisfies

$$\kappa(J(p), \xi) = \langle p, \xi \,.\, q \rangle \tag{5.1.7}$$

for $p \in \mathrm{T}^*_q Q$ and $\xi \in \mathfrak{g}$. Note that the right-hand side, viewed as a functional on $\mathfrak{g}$, may not be representable by an element $J(p) \in \mathfrak{g}^*$. In this case, a $\mathfrak{g}^*$-valued momentum map for the lifted action does not exist. The existence of the momentum map mainly depends on the chosen duality $\kappa$ and not so much on the cotangent bundle. In fact, we now give an example for an action of an infinite-dimensional Lie group acting on a *finite-dimensional* cotangent bundle that does not possess a momentum map.

[1] Without mentioning it in the future, we assume that the diffeomorphism $\mathrm{T}Q \simeq \mathrm{T}^*Q$ is fiber-preserving.

Example 5.1.2 Let $M$ be a compact finite-dimensional manifold endowed with a volume form vol. Consider a finite-dimensional Lie group $G$ that acts on a finite-dimensional manifold $Q$. Fix a point $m_0 \in M$ and let the current group $C^\infty(M, G)$ act on $Q$ via the evaluation group homomorphism $\mathrm{ev}_{m_0}\colon C^\infty(M, G) \to G$, that is,

$$\lambda \cdot q := \lambda(m_0) \cdot_G q \tag{5.1.8}$$

for $\lambda \in C^\infty(M, G)$ and $q \in Q$. The induced action of $C^\infty(M, G)$ on the cotangent bundle $T^*Q$ is symplectic but it does not possess a momentum map with respect to the natural dual pairing

$$\kappa\colon C^\infty(M, \mathfrak{g}^*) \times C^\infty(M, \mathfrak{g}) \ni (\mu, \xi) \mapsto \int_M \langle \mu, \xi \rangle \, \mathrm{vol} \in \mathbb{R}. \tag{5.1.9}$$

If the momentum map $J$ were to exist, then it would necessarily satisfy the relation

$$\int_M \langle J(p), \xi \rangle \, \mathrm{vol} = \kappa(J(p), \xi) = \langle p, \xi(m_0) \,.\, q \rangle = \langle J_G(p), \xi(m_0) \rangle \tag{5.1.10}$$

for all $p \in T^*Q$ and $\xi \in C^\infty(M, \mathfrak{g})$, where $J_G\colon T^*Q \to \mathfrak{g}^*$ denotes the momentum map for the $G$-action. This is only possible if $J(p)$ is a delta distribution localized at the point $m_0$. However, by definition, the dual $C^\infty(M, \mathfrak{g}^*)$ only contains regular distributions.

In other words, the reason for the non-existence of the momentum map is that the evaluation map $\mathrm{ev}_{m_0}\colon C^\infty(M, \mathfrak{g}) \to \mathfrak{g}$ is *not* weakly continuous with respect to the pairings $\kappa$ and $\langle \cdot\,, \cdot \rangle$, and hence does not have an adjoint $\mathfrak{g}^* \to C^\infty(M, \mathfrak{g}^*)$. ◊

In the following, we always assume that the lifted action on $T^*Q$ has a momentum map.

Finally, there is another basic problem typical to the infinite-dimensional setting. As the map $\mathrm{ev}_{m_0}\colon C^\infty(M, \mathfrak{g}) \to \mathfrak{g}$ of Example 5.1.2 shows, the adjoint of a linear mapping between Fréchet spaces may not exist. However, our strategy for the construction of the normal form will involve dualizing. Thus, in all constructions of Chapter 5 that rely on dualizing, we will need to assume that the corresponding adjoint maps exist. Moreover, even if the adjoint exists and the original map is injective, the adjoint is in general not surjective but it only has a dense image, see Proposition B.2.3 (iv). We will also assume that such adjoints are actually surjective. As we will see in Section 5.6, these assumptions are fulfilled for Yang–Mills theory because the maps involved are (elliptic) differential operators.

## 5.2 Lifted slices and normal form

The fibered structure of the cotangent bundle yields a refinement of the usual orbit-momentum type strata into so-called seams. As in the context of standard symplectic reduction, the manifold structure of these seams follows from an appropriate normal form theorem for the momentum map. To arrive at such a normal form result, it would seem natural to start from the construction of the MGS normal form in Section 4.2.3 and adapt it to the geometry of the cotangent bundle case. However, this approach would lead to a theorem with rather strong assumptions in addition to the ones necessary for Theorem 4.2.27. It turns out to be more efficient to incorporate the additional structure of the cotangent bundle from the outset and to pursue a completely different strategy for the construction of the normal form:

(i) By using a slice $S$ at $q \in Q$, reduce the problem to $\mathrm{T}^*(G \times_{G_q} S)$.

(ii) Establish an equivariant diffeomorphism $\mathrm{T}^*(G \times_{G_q} S) \simeq G \times_{G_q} (\mathfrak{m}^* \times \mathrm{T}^* S)$, where $\mathfrak{m}$ is a complement of $\mathfrak{g}_q$ in $\mathfrak{g}$.

(iii) Calculate the momentum map under these identifications.

As a by-product, our theory provides an approach much simpler than the construction of [PRS07; Sch07] in the finite-dimensional setting.

Remark 5.2.1 We also note that it is impossible to directly extend the approach of [PRS07; Sch07] in the finite-dimensional setting to the Fréchet context according to the following serious obstacles: one needs the existence of (orthogonal) complements and an appropriate version of regular symplectic reduction to construct the normal form reference system; moreover, the Inverse Function Theorem is central to the construction of the local symplectomorphism bringing the system into the normal form. All these tools are not readily available for Fréchet manifolds. ◇

Let $p \in \mathrm{T}^* Q$ be a point in the fiber over $q \in Q$. Assume that the $G$-action on $Q$ has a slice $S$ at $q$. According to the strategy outlined above, the first step is to reduce the problem of determining the local structure of the momentum map near $p$ to a problem on $\mathrm{T}^*(G \times_{G_q} S)$. This reduction is accomplished by applying the following standard result to the tube map[1] $\chi^{\mathrm{T}} : G \times_{G_q} S \to Q$, which plainly extends from the finite- to the infinite-dimensional setting.

Proposition 5.2.2 (Lifting point transformations) *Let $C$ and $Q$ be Fréchet manifolds and let $\phi : C \to Q$ be a diffeomorphism. Assume that the lift*

$$\mathrm{T}^*\phi : \mathrm{T}^* Q \to \mathrm{T}^* C, \quad p \mapsto \phi^* p, \tag{5.2.1}$$

[1] Recall from Proposition A.2.3 that $\chi^{\mathrm{T}}$ is a diffeomorphism onto an open neighborhood of $q$ in $Q$. Moreover, it is $G$-equivariant with respect to the action of $G$ on $G \times_{G_q} S$ by left translation on the $G$-factor.

*of $\phi$ is a diffeomorphism. Then,* $\mathrm{T}^*\phi$ *is a symplectomorphism. If, moreover, $\phi$ is $G$-equivariant and the lifted action on* $\mathrm{T}^*Q$ *has a momentum map $J_Q$, then $J_C := J_Q \circ \mathrm{T}^*\phi^{-1} \colon \mathrm{T}^*C \to \mathfrak{g}^*$ is a momentum map for the lifted $G$-action on* $\mathrm{T}^*C$. ◇

According to step two of the general strategy, we now establish a convenient identification of the cotangent bundle of $G \times_{G_q} S$. Generalizing the finite-dimensional case, we say that a Lie subgroup $H \subseteq G$ is *reductive* if its Lie algebra $\mathfrak{h}$ has an $\mathrm{Ad}_H$-invariant complement in $\mathfrak{g}$. Since we consider only proper actions of $G$ on $Q$, the stabilizer $G_q$ is always compact and, hence, the following lemma shows that $G_q$ is a reductive Lie subgroup of $G$.

Lemma 5.2.3 *For every compact Lie subgroup $H \subseteq G$, there exists an* $\mathrm{Ad}_H$*-invariant complement $\mathfrak{m}$ of $\mathfrak{h}$ in $\mathfrak{g}$. Moreover, there exists a weakly complementary decomposition* $\mathfrak{g}^* = \mathfrak{m}^* \oplus \mathfrak{h}^*$. ◇

*Proof.* Every compact Lie group is finite-dimensional, because every locally compact topological vector space is finite-dimensional, see [Köt83, Proposition 8.7.1]. As a finite-dimensional subspace, $\mathfrak{h}$ is automatically closed and has a topological complement by [Köt83, Proposition 15.5.2 and 20.5.5]. The complement can be chosen to be $\mathrm{Ad}_H$-invariant by taking the average over the projection using the invariant Haar measure, see Lemma A.2.5. By [Köt83, Proposition 20.5.1], there exists a weakly complementary decomposition $\mathfrak{g}^* = \mathfrak{g}_q^0 \oplus \mathfrak{m}^0$, where the subscript denotes the annihilator. Choose $\mathfrak{g}_q^0 = \mathfrak{m}^*$ and $\mathfrak{m}^0 = \mathfrak{h}^*$. □

Let $\iota \colon G \times S \to G \times_{G_q} S$ be the natural projection and, for $a \in G$ and $s \in S$, let $\iota_a \colon S \to G \times_{G_q} S$ and $\iota_s \colon G \to G \times_{G_q} S$ denote the induced embeddings, respectively.

Definition 5.2.4 Let $S$ be a slice at $q$ and let $\mathfrak{m}$ be an $\mathrm{Ad}_{G_q}$-invariant complement of $\mathfrak{g}_q$ in $\mathfrak{g}$. For a given choice of $\mathrm{T}^*(G \times_{G_q} S)$, the slice $S$ will be called *compatible with the cotangent bundle structures* if it fulfills the following additional requirements:

(ST1) The lift $\mathrm{T}^*\chi^{\mathrm{T}}$ of the tube map $\chi^{\mathrm{T}} \colon G \times_{G_q} S \to Q$ exists and is a diffeomorphism onto its image.

(ST2) The injective maps

$$\mathrm{T}_a \iota_s(a \,.\, \cdot) \colon \mathfrak{m} \to \mathrm{T}_{[a,s]}(G \times_{G_q} S), \qquad \mathrm{T}_s \iota_a \colon \mathrm{T}_s S \to \mathrm{T}_{[a,s]}(G \times_{G_q} S) \quad (5.2.2)$$

have surjective adjoints. ◇

Whether a slice is compatible with the cotangent bundle structures strongly depends on the dual pairings involved, so that one cannot hope to find a general criterion for the existence of such a slice. Note that, in the finite-dimensional setting, every slice is compatible with the cotangent bundle structures. In the sequel, we will show that there is a natural choice for $\mathrm{T}^*(G \times_{G_q} S)$.

Lemma 5.2.5 *For every choice of an* $\mathrm{Ad}_{G_q}$*-invariant complement* $\mathfrak{m}$ *of* $\mathfrak{g}_q$ *in* $\mathfrak{g}$*, there exists a* $G$*-equivariant diffeomorphism*

$$\mathrm{T}(G \times_{G_q} S) \simeq G \times_{G_q} (\mathfrak{m} \times \mathrm{T}S). \tag{5.2.3}$$

◊

*Proof.* It is well-known that the choice of an Ad-invariant complement $\mathfrak{m}$ yields a homogeneous connection in $G \to G/G_q$. Indeed, the horizontal space at $a \in G$ is, by definition, $a \,.\, \mathfrak{m}$ and the horizontal lift of a vector $[a, \xi] \in G \times_{G_q} \mathfrak{m} \simeq \mathrm{T}(G/G_q)$ to the point $a \in G$ is $a \,.\, \xi$. Accordingly, the tangent bundle to the associated bundle $G \times_{G_q} S$ splits into its horizontal and vertical parts so that the $G$-equivariant map defined by

$$G \times \mathfrak{m} \times \mathrm{T}S \to \mathrm{T}(G \times_{G_q} S), \qquad (a, \varsigma, Y_s) \mapsto \mathrm{T}_a \iota_s(a \,.\, \varsigma) + \mathrm{T}_s \iota_a(Y_s), \tag{5.2.4}$$

is a $G_q$-invariant submersion which descends to a $G$-equivariant diffeomorphism between $G \times_{G_q} (\mathfrak{m} \times \mathrm{T}S)$ and $\mathrm{T}(G \times_{G_q} S)$. □

By dualizing the isomorphism (5.2.3), we get a natural choice for the cotangent bundle $\mathrm{T}^*(G \times_{G_q} S)$. For that purpose, we choose $\mathfrak{m}^*$ as constructed in the proof of Lemma 5.2.3. Moreover, viewing $S$ as a submanifold of $Q$, let $\mathrm{T}^*S$ be the image of $\mathrm{T}S$ under the $G$-equivariant diffeomorphism $\mathrm{T}Q \to \mathrm{T}^*Q$. Now, the cotangent bundle is provided by the $G$-equivariant diffeomorphism

$$\phi\colon G \times_{G_q} (\mathfrak{m}^* \times \mathrm{T}^*S) \to \mathrm{T}^*(G \times_{G_q} S) \tag{5.2.5}$$

defined by

$$\langle \phi([a, (\nu, \alpha_s)]), \mathrm{T}_a \iota_s(a \,.\, \varsigma) + \mathrm{T}_s \iota_a(Y_s) \rangle = \kappa(\nu, \varsigma) + \langle \alpha_s, Y_s \rangle \tag{5.2.6}$$

for $\varsigma \in \mathfrak{m}$ and $Y_s \in \mathrm{T}_s S$. With this choice of $\mathrm{T}^*(G \times_{G_q} S)$ the second condition (ST2) in Definition 5.2.4 is automatically satisfied, because the adjoints are the identical mappings.

Remark 5.2.6 Alternatively, one could define $\mathrm{T}^*(G \times_{G_q} S)$ by dualizing the $G$-equivariant (local) diffeomorphism $\mathrm{T}\chi^{\mathrm{T}}\colon \mathrm{T}(G \times_{G_q} S) \to \mathrm{T}Q$ given by the tube diffeomorphism $\chi^{\mathrm{T}}$. In this case, the first condition (ST1) in Definition 5.2.4 is automatically satisfied. ◊

Remark 5.2.7 In the finite-dimensional context, Schmah [Sch07, Proposition 13] established a similar identification of $\mathrm{T}^*(G \times_{G_q} S)$ using regular cotangent bundle reduction. The starting point is the cotangent bundle $\mathrm{T}^*(G \times S)$. Using left translation, identify $\mathrm{T}^*G$ with $G \times \mathfrak{g}^*$. Denote points in $\mathrm{T}^*(G \times S) \simeq G \times \mathfrak{g}^* \times \mathrm{T}^*S$ by tuples $(a, \mu, \alpha_s)$. A straightforward calculation shows that the lift of the twisted $G_q$-action on $G \times S$ has the momentum map

$$J_{G_q^{\mathrm{T}}}(a, \mu, \alpha_s) = -\mu\restriction_{\mathfrak{g}_q} + J_{G_q}(\alpha_s), \tag{5.2.7}$$

where $J_{G_q}\colon \mathrm{T}^*S \longrightarrow \mathfrak{g}_q^*$ is the momentum map for the lifted $G_q$-action on $\mathrm{T}^*S$. Define the map

$$\varphi\colon G \times \mathfrak{g}^* \times \mathrm{T}^*S \supseteq J^{-1}_{G_q^{\mathrm{T}}}(0) \longrightarrow \mathrm{T}^*(G \times_{G_q} S) \tag{5.2.8}$$

by

$$\langle \varphi(a, \mu, \alpha_s), \mathrm{T}\iota(a\,.\,\xi, Y_s)\rangle = \kappa(\mu, \xi) + \langle \alpha_s, Y_s\rangle, \tag{5.2.9}$$

for $\xi \in \mathfrak{g}$ and $Y_s \in \mathrm{T}_s S$. The regular cotangent bundle reduction theorem [OR03, Theorem 6.6.1] shows that $\varphi$ is a submersion and that it descends to a symplectomorphism $\check{\varphi}$ of $J^{-1}_{G_q^{\mathrm{T}}}(0)/G_q$ with $\mathrm{T}^*(G \times_{G_q} S)$. In order to establish the link with the isomorphism $\phi$ discussed above (taken in the finite-dimensional case), we define

$$\psi\colon G \times (\mathfrak{m}^* \times \mathrm{T}^*S) \longrightarrow G \times (\mathfrak{g}^* \times \mathrm{T}^*S), \tag{5.2.10}$$

$$(a, \nu, \alpha_s) \longmapsto (a, \nu + J_{G_q}(\alpha_s), \alpha_s). \tag{5.2.11}$$

By construction, $\psi$ takes values in $J^{-1}_{G_q^{\mathrm{T}}}(0)$ and on this set it has a smooth inverse,

$$J^{-1}_{G_q^{\mathrm{T}}}(0) \longrightarrow G \times (\mathfrak{m}^* \times \mathrm{T}^*S), \quad (a, \mu, \alpha_s) \longmapsto (a, \mu\restriction_{\mathfrak{m}}, \alpha_s). \tag{5.2.12}$$

Now, it is an easy exercise in chasing identifications to see that the following diagram commutes

$$\begin{array}{ccccc}
J^{-1}_{G_q^{\mathrm{T}}}(0) & & \xleftarrow{\quad\quad\psi\quad\quad} & & G \times (\mathfrak{m}^* \times \mathrm{T}^*S) \\
\downarrow & \searrow^{\varphi} & & \swarrow & \downarrow \\
J^{-1}_{G_q^{\mathrm{T}}}(0)/G_q & \xrightarrow{\check{\varphi}} & \mathrm{T}^*(G \times_{G_q} S) & \xleftarrow{\phi} & G \times_{G_q} (\mathfrak{m}^* \times \mathrm{T}^*S).
\end{array} \tag{5.2.13}$$

Since the regular cotangent bundle reduction theorem does not directly generalize to infinite dimensions[1], the approach of [Sch07] is not available to us in our infinite-dimensional context. One can, however, read (5.2.13) as a proof that the regular reduction of $\mathrm{T}^*(G \times S)$ coincides with $\mathrm{T}^*(G \times_{G_q} S)$, indeed. ◇

Lemma 5.2.8 *Under the diffeomorphism $\phi$ defined in* (5.2.6)*, the momentum map $J\colon \mathrm{T}^*(G \times_{G_q} S) \longrightarrow \mathfrak{g}^*$ for the lifted $G$-action is identified with the map*

$$(J \circ \phi)([a, (\nu, \alpha_s)]) = \mathrm{Ad}^*_a(\nu + J_{G_q}(\alpha_s)), \tag{5.2.14}$$

*where $J_{G_q}\colon \mathrm{T}^*S \longrightarrow \mathfrak{g}_q^*$ is the momentum map for the lifted $G_q$-action on $S$. In particular, the condition $(J \circ \phi)([a, (\nu, \alpha_s)]) = 0$ is equivalent to $J_{G_q}(\alpha_s) = 0$ and $\nu = 0$.* ◇

[1] We do prove a regular reduction theorem below in Theorem 5.3.8. However, this result relies on the existence of the normal form and thus cannot be used to construct it.

*Proof.* First, we note that the lifted $G_q$-action on $\mathrm{T}^*S$ has a momentum map $J_{G_q} \colon \mathrm{T}^*S \to \mathfrak{g}_q^*$, because, by properness of the action, $G_q$ is compact and hence finite-dimensional. The canonical $G$-action $g \cdot [a, s] = [ga, s]$ on $G \times_{G_q} S$ has the fundamental vector field

$$\xi \,.\, [a, s] = \mathrm{T}_a \iota_s(\xi \,.\, a) = \mathrm{T}_a \iota_s \left(a \,.\, (\mathrm{Ad}_a^{-1} \xi)\right), \quad \xi \in \mathfrak{g}. \tag{5.2.15}$$

Let $\beta = \phi([a, (\nu, \alpha_s)]) \in \mathrm{T}^*_{[a,s]}(G \times_{G_q} S)$. Then, for every $\xi \in \mathfrak{g}$, the momentum map satisfies

$$\begin{aligned}
\kappa(J(\beta), \mathrm{Ad}_a \xi) &= \left\langle \beta, (\mathrm{Ad}_a \xi) \,.\, [a, s] \right\rangle \\
&= \left\langle \beta, \mathrm{T}_a \iota_s(a \,.\, \xi) \right\rangle \\
&= \left\langle \beta, \mathrm{T}_a \iota_s \left(a \,.\, \xi_{\mathfrak{g}_q}\right) \right\rangle + \left\langle \beta, \mathrm{T}_a \iota_s \left(a \,.\, \xi_{\mathfrak{m}}\right) \right\rangle \\
&= \left\langle \beta, \mathrm{T}_s \iota_a \left(\xi_{\mathfrak{g}_q} \,.\, s\right) \right\rangle + \left\langle \beta, \mathrm{T}_a \iota_s \left(a \,.\, \xi_{\mathfrak{m}}\right) \right\rangle \\
&= \left\langle \alpha_s, \xi_{\mathfrak{g}_q} \,.\, s \right\rangle + \kappa(\nu, \xi_{\mathfrak{m}}) \\
&= \kappa\left(J_{G_q}(\alpha_s), \xi_{\mathfrak{g}_q}\right) + \kappa(\nu, \xi_{\mathfrak{m}}),
\end{aligned} \tag{5.2.16}$$

where we have decomposed $\xi = \xi_{\mathfrak{g}_q} + \xi_{\mathfrak{m}}$ into $\xi_{\mathfrak{g}_q} \in \mathfrak{g}_q$ and $\xi_{\mathfrak{m}} \in \mathfrak{m}$. In the line before the last line we have used (5.2.6). Hence,

$$\mathrm{Ad}^*_{a^{-1}}(J \circ \phi)([a, (\nu, \alpha_s)]) \equiv \mathrm{Ad}^*_{a^{-1}} J(\beta) = J_{G_q}(\alpha_s) + \nu \tag{5.2.17}$$

and so $J(\beta) = 0$ if and only if $J_{G_q}(\alpha_s) = 0$ and $\nu = 0$. □

Combining the diffeomorphism $\phi$ with the local tube diffeomorphism

$$\mathrm{T}^*\chi^{\mathrm{T}} \colon \mathrm{T}^*Q \to \mathrm{T}^*(G \times_{G_q} S) \tag{5.2.18}$$

yields a convenient normal form for the momentum map of the lifted $G$-action on $\mathrm{T}^*Q$.

Theorem 5.2.9 (Normal form) *Let $Q$ be a Fréchet $G$-manifold with proper $G$-action. Assume that $\mathrm{T}^*Q$ is a Fréchet manifold, which is $G$-equivariantly diffeomorphic to $\mathrm{T}Q$. Let $p \in \mathrm{T}^*_q Q$ and assume that the $G$-action on $Q$ admits a slice at $q$ compatible with the cotangent bundle structures. Then, the map $\Phi \colon G \times_{G_q} (\mathfrak{m}^* \times \mathrm{T}^*S) \to \mathrm{T}^*Q$ defined by the condition*

$$\langle \Phi([a, (\nu, \alpha_s)]), (\mathrm{Ad}_a \rho) \,.\, (a \cdot s) + a \,.\, Y_s \rangle = \kappa(\nu, \rho) + \langle \alpha_s, Y_s \rangle, \tag{5.2.19}$$

*for all $\rho \in \mathfrak{m}$ and $Y_s \in \mathrm{T}_s S$, is a diffeomorphism onto an open neighborhood of $p$ in $\mathrm{T}^*Q$. Moreover, $\Phi$ is $G$-equivariant with respect to left translation on the $G$-factor and the lifted action on $\mathrm{T}^*Q$. Assume, in addition, that the momentum map $J \colon \mathrm{T}^*Q \to \mathfrak{g}^*$* *for the lifted G-action exists. Then, under* $\Phi$*, it is identified with the map*

$$(J \circ \Phi)([a, (\nu, \alpha_s)]) = \mathrm{Ad}^*_a(\nu + J_{G_q}(\alpha_s)), \tag{5.2.20}$$

*where* $J_{G_q} \colon \mathrm{T}^*S \to \mathfrak{g}^*_q$ *is the momentum map for the lifted* $G_q$*-action on* $\mathrm{T}^*S$*.* ◊

*Proof.* Let the map $\Phi = \mathrm{T}^*(\chi^{\mathrm{T}})^{-1} \circ \phi$ be defined as the composition of the diffeomorphism $\phi \colon G \times_{G_q} (\mathfrak{m}^* \times \mathrm{T}^*S) \to \mathrm{T}^*(G \times_{G_q} S)$ with the local diffeomorphism $\mathrm{T}^*(\chi^{\mathrm{T}})^{-1} \colon \mathrm{T}^*(G \times_{G_q} S) \to \mathrm{T}^*Q$. The composition of (5.2.4) with $\mathrm{T}\chi^{\mathrm{T}}$ yields the map

$$G \times_{G_q} (\mathfrak{m} \times \mathrm{T}S) \to \mathrm{T}Q, \quad [a, (\rho, Y_s)] \mapsto \mathrm{T}_{a,s}(\chi^{\mathrm{T}} \circ \iota)(a \,.\, \rho, Y_s). \tag{5.2.21}$$

Note that we have

$$\chi^{\mathrm{T}} \circ \iota(a, s) = \chi^{\mathrm{T}}([a, s]) = a \cdot s. \tag{5.2.22}$$

Thus, the above map reads

$$G \times_{G_q} (\mathfrak{m} \times \mathrm{T}S) \to \mathrm{T}Q, \quad [a, (\rho, Y_s)] \mapsto (\mathrm{Ad}_a \rho) \,.\, (a \cdot s) + a \,.\, Y_s, \tag{5.2.23}$$

which is a diffeomorphism onto its image, the latter being a subbundle of $\mathrm{T}Q$. By dualizing, we obtain the expression (5.2.19) for $\Phi$. The asserted properties follow immediately from the previous discussion. □

We should note, however, that the semi-global diffeomorphism $\Phi$ does not bring the momentum map $J$ in an MGS normal form in the sense of Definition 4.2.18. In fact, it does not even provide a slice for the lifted action, because we have taken the quotient by the 'wrong' stabilizer group, i.e., the model space around $p \in \mathrm{T}^*_q Q$ is of the form $G \times_{G_q} V$ instead of $G \times_{G_p} W$. In finite dimensions, much work in the study of singular cotangent bundle reduction is devoted to constructing a bona fide symplectic slice adapted to the cotangent bundle structure (see, e.g., [Sch07; RT17]). Converting the normal form $G \times_{G_q} (\mathfrak{m}^* \times \mathrm{T}^*S)$ into a symplectic slice requires a detailed analysis of the Witt–Artin decomposition in the cotangent bundle setting and then extending these infinitesimal results to a local statement using the Inverse Function Theorem. It is not clear if and how these steps generalize to the infinite-dimensional context as they lead to delicate issues of analytic nature. It turns out, however, that the simple normal form (5.2.20) of the momentum map we have constructed so far is sufficient for most questions concerning singular cotangent bundle reduction.

Recall that the choice of a complement $\mathfrak{m}$ of $\mathfrak{g}_q$ in $\mathfrak{g}$ yields a canonical identification of $\mathrm{T}_{[e]}(G/G_q)$ with $\mathfrak{m}$.

Proposition 5.2.10 *Under the diffeomorphism* $\Phi$ *constructed in Theorem 5.2.9 and under the identification* $\mathrm{T}_{[a,(\nu,\alpha_s)]}(G \times_{G_q} (\mathfrak{m}^* \times \mathrm{T}^*S)) \simeq \mathfrak{m} \times \mathfrak{m}^* \times \mathrm{T}_{\alpha_s}(\mathrm{T}^*S)$*, the* *canonical symplectic form $\omega$ on $\mathrm{T}^*Q$ takes the following form:*

$$\begin{aligned}&(\Phi^*\omega)_{[a,(\nu,\alpha_s)]}\big((\xi^1,\eta^1,Z^1),(\xi^2,\eta^2,Z^2)\big)\\&\quad=\kappa(\eta^1,\xi^2)-\kappa(\eta^2,\xi^1)-\kappa\Big(\nu+J_{G_q}(\alpha_s),[\xi^1,\xi^2]\Big)+\omega^S_{\alpha_s}(Z^1,Z^2),\end{aligned}\tag{5.2.24}$$

*where $\xi^i\in\mathfrak{m}$, $\eta^i\in\mathfrak{m}^*$ and $Z^i\in\mathrm{T}_{\alpha_s}(\mathrm{T}^*S)$ for $i=1,2$ and $\omega^S$ is the canonical symplectic form on $\mathrm{T}^*S$. In particular, at points $[a,(\nu,\alpha_s)]$ with $J\circ\Phi([a,(\nu,\alpha_s)])=0$, the third term on the right-hand side vanishes and $\Phi^*\omega$ is the direct sum of the canonical symplectic forms on $\mathfrak{m}\times\mathfrak{m}^*$ and $\mathrm{T}^*S$.* ◇

*Proof.* Let $\overset{\star}{\tau}\colon \mathrm{T}^*Q\longrightarrow Q$ be the canonical projection. Then, $\overset{\star}{\tau}\circ\Phi([a,(\nu,\alpha_s)])=a\cdot s$. Using (5.2.19), for the pull-back of the canonical 1-form $\theta$ we find:

$$\begin{aligned}&(\Phi^*\theta)_{[a,(\nu,\alpha_s)]}([a\,.\,\xi,(\eta,Z)])\\&\qquad=\langle\Phi([a,(\nu,\alpha_s)]),\mathrm{T}_{[a,(\nu,\alpha_s)]}(\overset{\star}{\tau}\circ\Phi)([a\,.\,\xi,(\eta,Z)])\rangle,\\&\qquad=\langle\Phi([a,(\nu,\alpha_s)]),(\mathrm{Ad}_a\,\xi)\,.\,(a\cdot s)+a\,.\,\mathrm{T}_{\alpha_s}\overset{\star}{\tau}(Z)\rangle\\&\qquad=\langle\Phi([a,(\nu,\alpha_s)]),(\mathrm{Ad}_a\,\xi_{\mathfrak{m}})\,.\,(a\cdot s)+a\,.\,(\xi_{\mathfrak{g}_q}\,.\,s)+a\,.\,\mathrm{T}_{\alpha_s}\overset{\star}{\tau}(Z)\rangle\\&\qquad=\kappa(\nu,\xi_{\mathfrak{m}})+\langle\alpha_s,\xi_{\mathfrak{g}_q}\,.\,s\rangle+\langle\alpha_s,\mathrm{T}_{\alpha_s}\overset{\star}{\tau}(Z)\rangle\\&\qquad=\kappa(\nu,\xi_{\mathfrak{m}})+\kappa(J_{G_q}(\alpha_s),\xi_{\mathfrak{g}_q})+\theta^S_{\alpha_s}(Z),\end{aligned}\tag{5.2.25}$$

where $\xi=\xi_{\mathfrak{m}}+\xi_{\mathfrak{g}_q}\in\mathfrak{g}$, $\eta\in\mathfrak{m}^*$ and $Z\in\mathrm{T}_{\alpha_s}(\mathrm{T}^*S)$, and $\theta^S$ is the canonical 1-form on $\mathrm{T}^*S$. If we introduce the left Maurer–Cartan form $\vartheta\in\Omega^1(G,\mathfrak{g})$ by $\vartheta_a(a\,.\,\xi)=\xi$, then we get

$$(\Phi^*\theta)_{[a,(\nu,\alpha_s)]}=\kappa\Big(\nu+J_{G_q}(\alpha_s),\vartheta_a(\cdot)\Big)+\theta^S_{\alpha_s}(\cdot).\tag{5.2.26}$$

The Maurer–Cartan equation $\mathrm{d}\vartheta=-\frac{1}{2}[\vartheta\wedge\vartheta]$ yields

$$\begin{aligned}&(\mathrm{d}\Phi^*\theta)_{[a,(\nu,\alpha_s)]}\big([a\,.\,\xi^1,(\eta^1,Z^1)],[a\,.\,\xi^2,(\eta^2,Z^2)]\big)\\&\qquad=\kappa(\eta^1,\xi^2_{\mathfrak{m}})-\kappa(\eta^2,\xi^1_{\mathfrak{m}})+\kappa(\nu+J_{G_q}(\alpha_s),(\mathrm{d}\vartheta)_a(a\,.\,\xi^1,a\,.\,\xi^2))\\&\qquad\quad+\kappa(\mathrm{T}_{\alpha_s}J_{G_q}(Z^1),\xi^2_{\mathfrak{g}_q})-\kappa(\mathrm{T}_{\alpha_s}J_{G_q}(Z^2),\xi^1_{\mathfrak{g}_q})\\&\qquad\quad+(\mathrm{d}\theta^S)_{\alpha_s}(Z^1,Z^2)\\&\qquad=\kappa(\eta^1,\xi^2_{\mathfrak{m}})-\kappa(\eta^2,\xi^1_{\mathfrak{m}})-\kappa\Big(\nu+J_{G_q}(\alpha_s),[\xi^1,\xi^2]\Big)\\&\qquad\quad+\omega^S_{\alpha_s}(Z^1,\xi^2_{\mathfrak{g}_q}\,.\,\alpha_s)-\omega^S_{\alpha_s}(Z^2,\xi^1_{\mathfrak{g}_q}\,.\,\alpha_s)+\omega^S_{\alpha_s}(Z^1,Z^2).\end{aligned}\tag{5.2.27}$$

Now the identification $\mathrm{T}_{[a,(\nu,\alpha_s)]}(G\times_{G_q}(\mathfrak{m}^*\times\mathrm{T}^*S))\simeq\mathfrak{m}\times\mathfrak{m}^*\times\mathrm{T}_{\alpha_s}(\mathrm{T}^*S)$ amounts to setting the $\mathfrak{g}_q$-component of $\xi^i$ to zero and we thus arrive at the claimed formula (5.2.24). □

## 5.3 Symplectic stratification

In Section 4.3, we have seen that the decomposition of the symplectically reduced phase space into orbit type manifolds yields a stratification. For the reduction scheme, we heavily relied on the MGS normal form to prove that the orbit type subsets carry a smooth manifold structure. As we will see now, the Normal Form Theorem 5.2.9 for the momentum map established in the previous section allows to make similar statements concerning the local structure in the context of cotangent bundles.

Theorem 5.3.1 *Let $Q$ be a Fréchet $G$-manifold. Assume that the $G$-action is proper, that it admits at every point a slice compatible with the cotangent bundle structures and that the decomposition of $Q$ into orbit types satisfies the frontier condition. Moreover, assume that $\mathrm{T}^*Q$ is a Fréchet manifold, which is $G$-equivariantly diffeomorphic to $\mathrm{T}Q$, and that the lifted action on $\mathrm{T}^*Q$, endowed with its canonical symplectic form $\omega$, has a momentum map $J$. Then, the following holds:*

*(i) The set of orbit types of $P := J^{-1}(0)$ with respect to the lifted $G$-action coincides with the set of orbit types for the $G$-action on $Q$.*

*(ii) The reduced phase space $\check{P} := J^{-1}(0)/G$ is stratified into orbit type manifolds $\check{P}_{(K)} := (J^{-1}(0))_{(K)}/G$.*

*(iii) Assume, additionally, that the orbit $\mathfrak{g} \,.\, p$ is symplectically closed for all $p \in P$. Then, for every orbit type $(K)$, the manifold $\check{P}_{(K)}$ carries a symplectic form $\check{\omega}_{(K)}$ uniquely determined by*

$$\pi^*_{(K)}\check{\omega}_{(K)} = \omega\restriction_{P_{(K)}}, \tag{5.3.1}$$

*where $\pi_{(K)}\colon P_{(K)} \to \check{P}_{(K)}$ is the natural projection.* ◊

In the sequel, we will prove this theorem by means of a series of lemmas. Let us start with the observation that the orbit types of the lifted action are tightly connected to the ones for the action on the base manifold. The following results extend the determination of orbit types of the lifted action [Rod06; PRS07, Theorem 5] to the infinite-dimensional setting.

Lemma 5.3.2 *Under the assumptions of Theorem 5.3.1, for every orbit type $(K)$ of $\mathrm{T}^*Q$ there exist a triple $(H, \tilde{H}, L)$ such that $(K)$ is represented by $K = \tilde{H} \cap L$, where $\tilde{H} \subseteq H$ are stabilizer subgroups of $Q$ and $L$ is a stabilizer subgroup of the $H$-action on $(\mathfrak{g}/\mathfrak{h})^*$ with $\mathfrak{h}$ being the Lie algebra of $H$.*

*Conversely, for every triple $(H, (\tilde{H}), L)$, where $H$ is a stabilizer subgroup of $Q$, $(\tilde{H})$ is an orbit type of $Q$ fulfilling $(\tilde{H}) \leq (H)$ and $L$ is a stabilizer subgroup of the $H$-action on $(\mathfrak{g}/\mathfrak{h})^*$, there exists a representative $\tilde{H}$ of $(\tilde{H})$ and a stabilizer subgroup $K$ of $\mathrm{T}^*Q$ such that $K = \tilde{H} \cap L$.* ◊

*Proof.* Let $p \in \mathrm{T}^*_q Q$ and denote the stabilizer of $p$ and $q$ by $K$ and $H$, respectively. Choose an $\mathrm{Ad}_H$-invariant complement $\mathfrak{m}$ of $\mathfrak{h} \equiv \mathfrak{g}_q$ in $\mathfrak{g}$, which according to Lemma 5.2.3 is possible because $H$ is compact. We will first show that there exists $\tilde{q} \in Q$ (close to $q$) and $\nu \in \mathfrak{m}^*$ such that

$$K = G_{\tilde{q}} \cap H_\nu. \tag{5.3.2}$$

For this purpose, let $S$ be a slice at $q$. Denote the model space of $S$ by $X$. By Theorem 5.2.9, there exists $\nu \in \mathfrak{m}^*$ and $\alpha \in \mathrm{T}^*_q S$ such that $\Phi([e, (\nu, (q, \alpha))]) = p$. Since $\Phi$ is $G$-equivariant, the common stabilizer of $\nu$ and $\alpha$ under the $H$-action on $\mathfrak{m}^* \times \mathrm{T}^*_q S$ is $K$, that is, $K = H_\alpha \cap H_\nu$. By assumption, there exists a $G$-equivariant diffeomorphism between $\mathrm{T}Q$ and $\mathrm{T}^*Q$. The latter induces an $H$-equivariant diffeomorphism between $\mathrm{T}S$ and $\mathrm{T}^*S$. Thus, there exists $x \in X \simeq \mathrm{T}_q S$ whose stabilizer under the $H$-action is $H_\alpha$. Since the topology of $X$ is generated by absorbent sets and $S$ is (diffeomorphic to) an open neighborhood of 0 in $X$, there exists $r \in \mathbb{R}$ such that $x = r\tilde{q}$ for some $\tilde{q} \in S$. By linearity of the $H$-action, both $x$ and $\tilde{q}$ have the same stabilizer $H_\alpha$. As $\tilde{q} \in S$, (SL2) of Definition A.2.2 shows that $G_{\tilde{q}} \subseteq H$ and so $G_{\tilde{q}} = H_{\tilde{q}}$. Hence, we have found $\tilde{q}$ and $\nu$ satisfying (5.3.2).

In the converse direction, let $H$ be a stabilizer subgroup of $Q$, $(\tilde{H})$ be an orbit type of $Q$ fulfilling $(\tilde{H}) \leq (H)$ and $L$ be a stabilizer subgroup of the $H$-action on $(\mathfrak{g}/\mathfrak{h})^*$. Let $q \in Q$ be such that $G_q = H$. Choose a slice $S$ at $q$. Since $(\tilde{H}) \leq (H)$, the frontier condition and (SL3) of Definition A.2.2 imply that there exists $\tilde{q} \in S$ with $(H_{\tilde{q}}) = (G_{\tilde{q}}) = (\tilde{H})$. Choose an $\mathrm{Ad}_H$-invariant complement $\mathfrak{m}$ of $\mathfrak{h}$ in $\mathfrak{g}$ and $\nu \in \mathfrak{m}^*$ with $H_\nu = L$. Since $\Phi$ is equivariant, the point $p = \Phi([e, (\nu, (\tilde{q}, 0))])$ has stabilizer

$$G_p = H_{\tilde{q}} \cap H_\nu, \tag{5.3.3}$$

which completes the proof. □

Corollary 5.3.3 *Under the assumptions of Theorem 5.3.1, the set of orbit types of $P = J^{-1}(0)$ with respect to the lifted $G$-action coincides with the set of orbit types for the $G$-action on $Q$, that is, a subgroup $K$ is the stabilizer $K$ of some $p \in P$ if and only if there exists a point $\tilde{q} \in Q$ such that $G_{\tilde{q}} = K$.* ◇

*Proof.* One direction is clear, as the zero section of $\mathrm{T}^*Q$ has the same orbit types as $Q$ and is contained in the zero level set of $J$. The conclusion in the converse direction follows by similar arguments as in the proof of Lemma 5.3.2. Indeed, $J(p) = 0$ implies that $\nu = 0$ by (5.2.20) and hence $G_{\tilde{q}} = K$ by (5.3.2). □

Lemma 5.3.4 *Under the assumptions of Theorem 5.3.1 the following holds. For every orbit type $(K)$ of $\mathrm{T}^*Q$, the orbit type stratum $(\mathrm{T}^*Q)_{(K)}$ is a submanifold of $\mathrm{T}^*Q$ and the quotient $(\mathrm{T}^*Q)_{(K)}/G$ is a Fréchet manifold.* ◇

*Proof.* Let $p \in (\mathrm{T}^*Q)_{(K)}$ and denote its base point by $q \in Q$. By Theorem 5.2.9, it is enough to show that

$$\Big(G \times_{G_q} (\mathfrak{m}^* \times \mathrm{T}^*S)\Big)_{(G_p)} = G \times_{G_q} \Big((\mathfrak{m}^* \times \mathrm{T}^*S)_{(G_p)}\Big) \tag{5.3.4}$$

is a submanifold of $G \times_{G_q} (\mathfrak{m}^* \times \mathrm{T}^*S)$. In this expression, $(G_p)$ still denotes the conjugacy class of $G_p$ in $G$ and not in $G_q$. Since the action is linear and $G_q$ is compact, Theorem A.2.4 implies that the $G_q$-action on $\mathfrak{m}^* \times \mathrm{T}^*S$ admits a slice at every point. Using the existence of these slices and Lemma A.2.1, we conclude that $(\mathfrak{m}^* \times \mathrm{T}^*S)_{(G_p)}$ is a submanifold of $\mathfrak{m}^* \times \mathrm{T}^*S$, cf. Proposition A.2.7. Finally, the $G$-quotient is a smooth manifold, because it is locally identified with $(\mathfrak{m}^* \times \mathrm{T}^*S)_{(G_p)}/G_q$. □

The normal form for the momentum map of the lifted action established in Theorem 5.2.9 yields in a similar manner the following.

Lemma 5.3.5 *Under the assumptions of Theorem 5.3.1 the following holds. For every orbit type $(K)$ of the lifted G-action, the set*

$$P_{(K)} := (\mathrm{T}^*Q)_{(K)} \cap J^{-1}(0) = (J^{-1}(0))_{(K)} \tag{5.3.5}$$

*is a smooth submanifold of $(\mathrm{T}^*Q)_{(K)}$. Moreover, there exists a unique smooth manifold structure on*

$$\check{P}_{(K)} := P_{(K)}/G \tag{5.3.6}$$

*such that the natural projection $\pi_{(K)}\colon P_{(K)} \longrightarrow \check{P}_{(K)}$ is a smooth submersion.* ◇

*Proof.* Let $p \in P_{(K)}$ and let $q \in Q$ be its base point. Under the local diffeomorphism $\Phi$ established in Theorem 5.2.9, the set $P_{(K)}$ is identified with

$$\Big(G \times_{G_q} (\mathfrak{m}^* \times \mathrm{T}^*S)\Big)_{(K)} \cap (J \circ \Phi)^{-1}(0) = G \times_{G_q} \Big(\{0\} \times (J_{G_q}^{-1}(0))_{(K)}\Big). \tag{5.3.7}$$

This equality is a direct consequence of (5.3.4) and (5.2.20). By Corollary 4.3.7, the orbit type subset $(J_{G_q}^{-1}(0))_{(K)} = (\mathrm{T}^*S)_{(K)} \cap J_{G_q}^{-1}(0)$ is a smooth submanifold of $\mathrm{T}^*S$. By the same corollary, the quotient space $(J_{G_q}^{-1}(0))_{(K)}/G_q$, which is the model space of $\check{P}_{(K)}$, is a smooth manifold, too. □

Lemma 5.3.6 *Under the assumptions of Theorem 5.3.1, the orbit type decompositions of $\check{Q}$, $P$ and $\check{P}$ satisfy the frontier condition.* ◇

*Proof.* The quotient $\check{Q}$ inherits the frontier condition from $Q$ by [DR18c, Theorem 4.6]. Similarly, $\check{P}$ inherits the frontier condition from $P$. Let $p \in P$. We have to show that $p \in \overline{P_{(K)}}$ if and only if $(G_p) \geq (K)$. Denote the base point of $p$ by $q$, abbreviate $H \equiv G_q$ and let $\Phi$ be the local diffeomorphism constructed in Theorem 5.2.9. As we have seen in the proof of Lemma 5.3.5, $P_{(K)}$ is locally identified with

$$G \times_H \left(\{0\} \times (J_H^{-1}(0))_{(K)}\right). \tag{5.3.8}$$

Moreover, $\overline{P_{(K)}}$ is locally identified with

$$G \times_H \Big(\{0\} \times \overline{(J_H^{-1}(0))_{(K)}}\Big), \tag{5.3.9}$$

because the quotient map $G \times \mathrm{T}^*S \to G \times_H \mathrm{T}^*S$ is open and $f^{-1}(\overline{A}) = \overline{f^{-1}(A)}$ for every open continuous map $f\colon Y \to Z$ and every subset $A \subseteq X$. Write $p = \Phi([e, (0, \alpha)])$, where $\alpha \in \mathrm{T}^*_q S$ satisfies $J_H(\alpha) = 0$. Since $G_p = H_\alpha$, it is enough to show that $\alpha$ lies in the closure of $(J_H^{-1}(0))_{(K)}$ if and only if $(H_\alpha) \geq (K)$.

First, suppose that $\alpha \in \overline{(J_H^{-1}(0))_{(K)}}$. Since $H$ is compact, there exists a slice at $\alpha$ and hence a neighborhood of $\alpha$ in $\mathrm{T}^*S$ such that every point in this neighborhood has a stabilizer subconjugate to $H_\alpha$. However, by assumption, $(J_H^{-1}(0))_{(K)}$ has to intersect this neighborhood and thus $(K) \leq (H_\alpha)$.

For the converse direction, we first need to establish a result about the orbit type decomposition of $S$. Since the orbit type stratification of $Q$ satisfies the frontier condition, we have $\bigcup_{(H) \geq (K)} Q_{(H)} \subseteq \overline{Q_{(K)}}$ for every orbit type $(K)$ of $Q$. We will show now that we get a similar approximation property in the slice. For this purpose, let $(H) > (K)$ be orbit types of the action on $Q$ and let $s \in S_{(H)}$. For every open neighborhood $V$ of $s$ in $S$, the image $\chi^S(U \times V)$ under the local slice diffeomorphism $\chi^S\colon U \times S \to Q$ is an open neighborhood of $s$ in $Q$. Since $Q_{(H)} \subseteq \overline{Q_{(K)}}$, the intersection

$$\chi^S(U \times V) \cap Q_{(K)} = \chi^S(U \times V_{(K)}) \tag{5.3.10}$$

is non-empty. Thus, $V_{(K)}$ is non-empty and we have shown that $S_{(H)} \subseteq \overline{S_{(K)}}$ for all orbit types $(H) > (K)$.

Now, suppose that $(K)$ is an orbit type of $P$ with $(K) \leq (H_\alpha)$. By Corollary 5.3.3, $(K)$ is also an orbit type of the $G$-action on $Q$. We will show that every open neighborhood $W$ of $\alpha$ in $\mathrm{T}^*_q S$ has non-empty intersection with $(J_H^{-1}(0))_{(K)}$. Note that the momentum map $J_H$ vanishes on the whole fiber $\mathrm{T}^*_q S$, because $q$ has stabilizer $H$. Thus, it suffices to show that $\alpha$ lies in the closure of $(\mathrm{T}^*_q S)_{(K)}$. Since $\mathrm{T}^*Q$ is $G$-equivariantly diffeomorphic to $\mathrm{T}Q$, there exists an $H$-equivariant diffeomorphism of $\mathrm{T}^*_q S$ and $\mathrm{T}_q S$. Let $v \in \mathrm{T}_q S$ be the image of $\alpha$ under this diffeomorphism. Since $S$ is (diffeomorphic to) an open subset of $\mathrm{T}_q S$, there exists a non-zero $r \in \mathbb{R}$ such that $rv \in S$. Note that scaling by $r$ is an $H$-equivariant diffeomorphism of $\mathrm{T}_q S$. In particular, the stabilizer of $rv$ coincides with $H_\alpha$. Thus, in summary, we have reduced the problem to showing that $rv$ lies in the closure of $(\mathrm{T}_q S)_{(K)}$. But, as $(H_{rv}) = (H_\alpha) \geq (K)$, we have

$$rv \in \overline{S_{(K)}} \subseteq \overline{(\mathrm{T}_q S)_{(K)}}, \tag{5.3.11}$$

using the approximation property $S_{(H_\alpha)} \subseteq \overline{S_{(K)}}$. □

Remark 5.3.7 In the paper [PRS07] it was silently taken for granted that the decomposition of $Q$ into orbit types always satisfies the frontier condition. However, in examples, there may be an orbit type whose closure contains some but not all fixed points; and hence the frontier condition is violated in these cases (see e.g., [CM18, Example 17 and Remark 13]). ◇

Finally, the last point in Theorem 5.3.1 follows from Proposition 4.3.3. This completes the proof of Theorem 5.3.1.

For the special case, when the $G$-action on $Q$ has only one orbit type[1], we obtain the infinite-dimensional counterpart to the well-known cotangent bundle reduction theorem for one orbit type [ER90, Theorem 1].

Theorem 5.3.8 *In the setting of Theorem 5.3.1, assume additionally that the $G$-action on $Q$ has only one orbit type. Then, $\check{P} = J^{-1}(0)/G$ is symplectomorphic to* $\mathrm{T}^*(Q/G)$ *with its canonical symplectic structure.* ◇

*Proof.* Let $(H)$ denote the orbit type of $Q$. By assumption, we have $Q_{(H)} = Q$ and thus the existence of slices ensures that the quotient $\check{Q} \equiv Q/G$ is a smooth manifold, see Proposition A.2.7. Theorem 5.3.1 entails that every point of $J^{-1}(0)$ has orbit type $(H)$ and that, moreover, the reduced space $\check{P} \equiv J^{-1}(0)/G$ is a smooth manifold and carries a closed 2-form[2] $\check{\omega}$. Let $\pi\colon Q \to \check{Q}$ denote the natural projection. We claim that the map

$$\tau\colon \mathrm{T}^*\check{Q} \to \check{P}, \quad \alpha_{[q]} \mapsto [\pi^*\alpha_{[q]}] \tag{5.3.12}$$

is a well-defined diffeomorphism and that $\tau^*\check{\omega}$ coincides with the canonical symplectic form on $\mathrm{T}^*\check{Q}$. First, for every $q \in Q$ and $\alpha \in \mathrm{T}^*_{[q]}\check{Q}$

$$\langle \pi^*\alpha_{[q]}, \xi \,.\, q\rangle = \langle \alpha_{[q]}, \mathrm{T}_q\pi(\xi \,.\, q)\rangle = 0 \tag{5.3.13}$$

shows that $\pi^*\alpha_{[q]}$ lies in the zero level set of $J$ and thus $\tau$ is well defined as a map with values in $\check{P}$. It is straightforward to see that $\tau$ is bijective as $\pi$ is a submersion with $\operatorname{Ker} \mathrm{T}_q\pi = \mathfrak{g} \,.\, q$. In order to show that $\tau$ is a diffeomorphism, let $q \in Q$ with $G_q = H$ and choose a slice $S$ at $q$. In slice coordinates, we have $Q \simeq G \times_H S$ and so $\check{Q} \simeq S$, because the assumption that $Q$ has only the orbit type $(H)$ implies $S = S_H$. Thus, locally $\mathrm{T}^*\check{Q} \simeq \mathrm{T}^*S$. On the other hand, in the proof of Lemma 5.3.5 we have seen that $\check{P}$ is locally identified with $J_H^{-1}(0)_{(H)}/H$, where $J_H\colon \mathrm{T}^*S \to \mathfrak{h}^*$ is the momentum map for the $H$-action on $\mathrm{T}^*S$. Now, $S = S_H$ implies that $J_H^{-1}(0)_{(H)}/H$ simply coincides with $\mathrm{T}^*S$. Going back to the definitions, it is easy to see that in these coordinates $\tau$ is the identity map on $\mathrm{T}^*S$ and thus a diffeomorphism. Finally, by Proposition 5.2.10, under these identifications the canonical symplectic form on $\mathrm{T}^*\check{Q}$ and the form $\check{\omega}$ both are equal to the canonical symplectic form on $\mathrm{T}^*S$. □

[1] This assumption includes, of course, also the case of a free action.

[2] A priori, $\check{\omega}$ may be degenerate as we have not assumed that the orbits be symplectically closed.

## 5.4 Secondary stratification

Examples of singular cotangent bundle reduction (such as the ones coming from lattice gauge theory [FRS07]) show that the natural projection from the reduced phase space $\check{P} = \mathrm{T}^*Q /\!\!/_0 G$ to the reduced configuration space $\check{Q} = Q/G$ is not a morphism of stratified spaces, i.e., there exist strata in $\check{P}$ that project onto different strata in $\check{Q}$. In the finite-dimensional context, Perlmutter, Rodriguez-Olmos, and Sousa-Dias [PRS07] have refined the symplectic stratification of $\check{P}$ in such a way that the projection $\check{P} \to \check{Q}$ becomes a morphism of stratified spaces. This so-called secondary stratification has also the advantage of identifying certain strata in the reduced phase space as cotangent bundles. In this section, we construct this refined stratification in our infinite-dimensional setting. In particular, we show that the secondary strata are submanifolds (for finite-dimensional manifolds this is shown in [PRS07, Theorem 7], but the proof there does not directly translate to the infinite-dimensional setting). Moreover, we study how the secondary strata behave with respect to the ambient symplectic structure and investigate the frontier condition. This discussion culminates in Theorem 5.4.8. In this section, we continue to work in the setting of Theorem 5.3.1.

For orbit types $(K)$ and $(H)$ of $\mathrm{T}^*Q$ and $Q$, respectively, consider the subset

$$(\mathrm{T}^*Q)^{(K)}_{(H)} \coloneqq \{p \in \mathrm{T}^*_q Q : q \in Q_{(H)}, p \in (\mathrm{T}^*Q)_{(K)}\} \tag{5.4.1}$$

of the orbit type stratum $(\mathrm{T}^*Q)_{(K)}$. Since the projection $\mathrm{T}^*Q \to Q$ is $G$-equivariant, $(\mathrm{T}^*Q)^{(K)}_{(H)}$ is non-empty only if $(K) \leq (H)$. Moreover, the union of $(\mathrm{T}^*Q)^{(K)}_{(H)}$ over all orbit types $(H)$ fulfilling this condition yields the orbit type stratum $(\mathrm{T}^*Q)_{(K)}$. Note that the stabilizer $G_p$ of $p \in \mathrm{T}^*_q Q$ under the lifted $G$-action on $\mathrm{T}^*Q$ coincides with the stabilizer $[G_q]_p$ of $p$ under the dual isotropy action of $G_q$ on the fiber $\mathrm{T}^*_q Q$. Whence, we equivalently have

$$(\mathrm{T}^*Q)^{(K)}_{(H)} \coloneqq \{p \in \mathrm{T}^*_q Q : q \in Q_{(H)}, [G_q]_p \sim K\}. \tag{5.4.2}$$

Lemma 5.4.1 *Assume that the $G$-action on $Q$ is proper and that it admits a slice compatible with the cotangent bundle structures at every point. Moreover, assume that $\mathrm{T}^*Q$ is a Fréchet manifold, which is $G$-equivariantly diffeomorphic to $\mathrm{T}Q$. Let $(K)$ be an orbit type of $\mathrm{T}^*Q$. Then, for every orbit type $(H)$ of $Q$ fulfilling $(H) \geq (K)$, the sets $(\mathrm{T}^*Q)^{(K)}_{(H)}$ and $(\mathrm{T}^*Q)^{(K)}_{(H)}/G$ are submanifolds of $(\mathrm{T}^*Q)_{(K)}$ and $(\mathrm{T}^*Q)_{(K)}/G$, respectively.* ◊

We call the sets $(\mathrm{T}^*Q)^{(K)}_{(H)}$ the *secondary strata* and the decomposition of $\mathrm{T}^*Q$ into these secondary strata is referred to as the *secondary orbit type stratification*. As far as we know, the above result is novel even for the finite-dimensional case.

*Proof.* Let $p \in (\mathrm{T}^*Q)^{(K)}_{(H)}$ and denote its base point by $q \in Q$. Without loss of generality we may assume that $G_p = K$ and $G_q = H$ with $G_p \subseteq G_q$. By Theorem 5.2.9, it is enough to show that the corresponding subset

$$\Big(G \times_{G_q} (\mathfrak{m}^* \times \mathrm{T}^*S)\Big)^{(G_p)}_{(G_q)} \subseteq G \times_{G_q} (\mathfrak{m}^* \times \mathrm{T}^*S) \tag{5.4.3}$$

is a submanifold. By definition, $S$ is diffeomorphic to an open subset of a Fréchet space $X$ and, hence, we may identify $\mathrm{T}^*S$ with $S \times X^*$. By Lemma A.2.1, for every point $s$ in the slice, $G_s$ is conjugate to $G_q$ if and only if $G_s = G_q$. We thus find

$$\Big(G \times_{G_q} (\mathfrak{m}^* \times \mathrm{T}^*S)\Big)^{(G_p)}_{(G_q)} = G \times_{G_q} \Big((\mathfrak{m}^* \times X^*)_{(G_p)} \times S_{G_q}\Big). \tag{5.4.4}$$

In this expression, $(G_p)$ clearly denotes the conjugacy class of $G_p$ in $G$ and not in $G_q$. Now, since $G_q$ is compact, $(\mathfrak{m}^* \times X^*)_{(G_p)}$ is a submanifold of $\mathfrak{m}^* \times X^*$ (the proof follows by the same arguments as in the proof of Lemma 5.3.4). Finally, the $G$-quotient is again a smooth manifold, because it is locally identified with $(\mathfrak{m}^* \times X^*)_{(G_p)}/G_q \times S_{G_q}$. □

For the study of the interaction of the secondary orbit type stratification with the momentum map geometry we need the following basic result about linear cotangent bundle reduction.

Lemma 5.4.2 *Let $\langle X^*, X\rangle$ be a dual pair of Fréchet spaces and let $G$ be a compact Lie group acting linearly on $X$ and, by duality, also on $X^*$. Then, the lifted $G$-action on $\mathrm{T}^*X$ has an equivariant momentum map $J\colon \mathrm{T}^*X \longrightarrow \mathfrak{g}^*$. Moreover, for every orbit type $(K)$,*

$$(\mathrm{T}^*X)^{(K)}_{(G)} \cap J^{-1}(0) \simeq X_G \times X^*_{(K)} \tag{5.4.5}$$

*is a submanifold of* $\mathrm{T}^*X$. ◇

*Proof.* Under the identification $\mathrm{T}^*X \simeq X \times X^*$, the canonical 1-form takes the form

$$\theta_{x,\alpha}(y,\beta) = \langle \alpha, y\rangle, \qquad x, y \in X, \alpha, \beta \in X^*. \tag{5.4.6}$$

Since $G$ is compact and, hence, finite-dimensional, the linear $G$-action has a momentum map $J$ defined by

$$\kappa(J(x,\alpha), \xi) = \langle \alpha, \xi \,.\, x\rangle, \qquad \xi \in \mathfrak{g}. \tag{5.4.7}$$

Note that $J(x,\alpha) = 0$ if $x \in X_G$. Since, by definition of the secondary strata,

$$(\mathrm{T}^*X)^{(K)}_{(G)} \simeq X_{(G)} \times X^*_{(K)} = X_G \times X^*_{(K)} \tag{5.4.8}$$

holds, we obtain

$$(\mathrm{T}^*X)^{(K)}_{(G)} \cap J^{-1}(0) \simeq X_G \times X^*_{(K)}. \tag{5.4.9}$$

Since $G$ is compact, the $G$-action on $X^*$ admits a slice at every point according to Theorem A.2.4 and thus the orbit type manifold $X^*_{(K)}$ is a submanifold of $X^*$, see Proposition A.2.7. Therefore, $(\mathrm{T}^*X)^{(K)}_{(G)} \cap J^{-1}(0)$ is a submanifold of $\mathrm{T}^*X$. □

We now return to the general non-linear setting. Given two orbit types $(K) \leq (H)$, following [PRS07], we call the set

$$P^{(K)}_{(H)} \coloneqq (\mathrm{T}^*Q)^{(K)}_{(H)} \cap J^{-1}(0) \tag{5.4.10}$$

a *preseam* and the quotient $\check{P}^{(K)}_{(H)} \coloneqq P^{(K)}_{(H)}/G$ a *seam*.

Lemma 5.4.3 *For every pair of orbit types $(H) \geq (K)$, the preseam $P^{(K)}_{(H)}$ is a smooth submanifold of $\mathrm{T}^*Q$ and the seam $\check{P}^{(K)}_{(H)}$ is a smooth submanifold of $\check{P}_{(K)}$ and of $(\mathrm{T}^*Q)_{(K)}/G$. Moreover, $\check{P}^{(K)}_{(H)}$ is a smooth fiber bundle over $\check{Q}_{(H)}$.* ◇

*Proof.* Let $p \in (\mathrm{T}^*Q)^{(K)}_{(H)}$ and denote its base point by $q \in Q$. Let $S$ be a slice at $q$ and let $X$ be the model space of $S$. Using the local diffeomorphism $\Phi$ of Theorem 5.2.9 and the isomorphism of equation (5.4.5), in a neighborhood of $p$ we can identify the preseam with the submanifold

$$\begin{aligned} \Big(G \times_{G_q} (\mathfrak{m}^* \times \mathrm{T}^*S)\Big)^{(G_p)}_{(G_q)} \cap J^{-1}(0) &= G \times_{G_q} \Big(\{0\} \times J^{-1}_{G_q}(0)^{(G_p)}_{(G_q)}\Big) \\ &\simeq G \times_{G_q} \Big(\{0\} \times S_{G_q} \times X^*_{(G_p)}\Big) \end{aligned} \tag{5.4.11}$$

of $G \times_{G_q} (\mathfrak{m}^* \times \mathrm{T}^*S)$. Similarly, the seam $\check{P}^{(K)}_{(H)}$ has locally the same structure as the smooth manifold $S_{G_q} \times (X^*_{(G_p)}/G_q)$. Under these identifications, the quotient map $\check{P}^{(K)}_{(H)} \to \check{Q}_{(H)}$ corresponds to the projection onto the first factor and is thus a locally trivial submersion. □

We now come to the interaction of the seams with the symplectic geometry. For this purpose, denote by $\check{\omega}^{(K)}_{(H)}$ the restriction of $\check{\omega}_{(K)}$ to $\check{P}^{(K)}_{(H)} \subseteq \check{P}_{(K)}$. The injection $\mathrm{T}(Q_{(H)}) \to (\mathrm{T}Q)\restriction_{Q_{(H)}}$ induces a surjective map $\mathrm{pr}\colon (\mathrm{T}^*Q)\restriction_{Q_{(H)}} \to \mathrm{T}^*(Q_{(H)})$ and, thereby, a map $\bar{\pi}\colon P^{(K)}_{(H)} \to \mathrm{T}^*(Q_{(H)})$. Let $\bar{\omega}_{(H)}$ denote the canonical symplectic form on $\mathrm{T}^*(Q_{(H)})$. With this notation, we can give a characterization of $\check{\omega}^{(K)}_{(H)}$ similar to the one for the reduced symplectic form.

Lemma 5.4.4 *The restriction* $\check{\omega}^{(K)}_{(H)}$ *of* $\check{\omega}_{(K)}$ *to* $\check{P}^{(K)}_{(H)} \subseteq \check{P}_{(K)}$ *is uniquely characterized by*

$$\left(\pi^{(K)}_{(H)}\right)^{*} \check{\omega}^{(K)}_{(H)} = \bar{\pi}^{*}\, \bar{\omega}_{(H)}, \tag{5.4.12}$$

*where* $\pi^{(K)}_{(H)} \colon P^{(K)}_{(H)} \to \check{P}^{(K)}_{(H)}$ *is the canonical projection.* ◊

*Proof.* We first note that the restriction of $\omega$ to $(\mathrm{T}^{*}Q)\restriction_{Q_{(H)}}$ coincides with the pull-back $\mathrm{pr}^{*}\, \bar{\omega}_{(H)}$. In fact, the commutative diagram

$$\begin{array}{ccccc} \mathrm{T}^{*}(Q_{(H)}) & \xleftarrow{\ \mathrm{pr}\ } & (\mathrm{T}^{*}Q)\restriction_{Q_{(H)}} & \longrightarrow & \mathrm{T}^{*}Q \\ \downarrow & & \downarrow & & \downarrow \\ Q_{(H)} & \xleftarrow{\ \mathrm{id}\ } & Q_{(H)} & \longrightarrow & Q \end{array} \tag{5.4.13}$$

and a straightforward calculation show that the pull-back $\mathrm{pr}^{*}\, \bar{\theta}_{(H)}$ of the canonical 1-form on $\mathrm{T}^{*}(Q_{(H)})$ coincides with the restriction of $\theta$ to $(\mathrm{T}^{*}Q)\restriction_{Q_{(H)}}$. Now the claim follows by chasing along the following commutative diagram:

$$\begin{array}{ccccc} \mathrm{T}^{*}(Q_{(H)}) & \xleftarrow{\ \mathrm{pr}\ } & (\mathrm{T}^{*}Q)\restriction_{Q_{(H)}} & \longrightarrow & \mathrm{T}^{*}Q \\ & \nwarrow^{\bar{\pi}} & \uparrow & & \uparrow \\ & & P^{(K)}_{(H)} & \longrightarrow & P_{(K)} \\ & & \downarrow{\scriptstyle \pi^{(K)}_{(H)}} & & \downarrow \\ & & \check{P}^{(K)}_{(H)} & \longrightarrow & \check{P}_{(K)}. \end{array} \tag{5.4.14}$$

□

The construction above provides additional insight into the structure of the seam $\check{P}^{(H)}_{(H)}$. To see this, note that $\bar{\pi}$ takes values in the zero level set of the momentum map

$$\bar{J}_{(H)} \colon \mathrm{T}^{*}(Q_{(H)}) \to \mathfrak{g}^{*}\,. \tag{5.4.15}$$

Moreover, $\bar{\pi}$ is $G$-equivariant and thus descends to a map

$$\check{\bar{\pi}} \colon \check{P}^{(K)}_{(H)} \to \bar{J}^{-1}_{(H)}(0)/G. \tag{5.4.16}$$

By Theorem 5.3.8, the target space $\left(\bar{J}^{-1}_{(H)}(0)\right)_{(H)}/G$ is symplectomorphic to $\mathrm{T}^{*}(\check{Q}_{(H)})$.

Proposition 5.4.5 *For every orbit type* $(H)$*, the restriction of* $\check{\omega}_{(H)}$ *to* $\check{P}^{(H)}_{(H)}$ *is symplectic and* $\check{\bar{\pi}}$ *is a symplectomorphism between* $\check{P}^{(H)}_{(H)}$ *and* $\mathrm{T}^{*}(\check{Q}_{(H)})$ *with its* *canonical symplectic structure*[1]. ◊

We call $\check{P}^{(H)}_{(H)}$ the *principal seam.*

*Proof.* As we have seen in the proof of Lemma 5.4.3, the seam $\check{P}^{(H)}_{(H)}$ is locally identified with $S_H \times (X^*_{(H)}/H) = S_H \times X^*_H \simeq \mathrm{T}^*(S_H)$. On the other hand, $\mathrm{T}^*(\check{Q}_{(H)})$ is locally equivalent to $\mathrm{T}^*(S_H)$, see the proof of Theorem 5.3.8. It is straightforward to see that, in these coordinates, $\check{\bar{\pi}}$ is the identity map on $\mathrm{T}^*(S_H)$ and hence a diffeomorphism. Finally, by Lemma 5.4.4, $\check{\bar{\pi}}$ intertwines the closed 2-form $\check{\omega}^{(H)}_{(H)}$ with the canonical symplectic form on $\mathrm{T}^*(\check{Q}_{(H)})$. Since $\check{\bar{\pi}}$ is a diffeomorphism, $\check{\omega}^{(H)}_{(H)}$ is symplectic. □

In finite dimensions, one can show that the seams $\check{P}^{(K)}_{(H)}$ are coisotropic with respect to the reduced symplectic form $\check{\omega}_{(K)}$ on $\check{P}_{(K)}$, see [PRS07, Corollary 9]. The proof, however, relies on counting dimensions and thus does not generalize to the infinite-dimensional setting. A different idea to show that the seams are coisotropic is to use a Witt–Artin decomposition adapted to the cotangent bundle case. In the finite-dimensional context, such a decomposition was established in [PRS08], but its extension to infinite dimensions is not immediate and will be left to future work.

Let us now discuss the frontier condition. For this purpose, we endow the set of pairs of orbit types $((K),(H))$ satisfying $(H) \geq (K)$ with the partial order

$$((K),(H)) \leq ((K'),(H')) \quad \text{if and only if} \quad (K) \leq (K') \text{ and } (H) \leq (H').$$

Lemma 5.4.6 *Assume that the orbit type decomposition of $Q$ satisfies the frontier condition. For every pair of orbit types $(K)$ and $(H)$ fulfilling $(H) \geq (K)$, we have*

$$\overline{\check{P}^{(K)}_{(H)}} = \bigcup_{((K'),(H')) \geq ((K),(H))} \check{P}^{(K')}_{(H')}, \tag{5.4.17}$$

*so that the decomposition of $\check{P}$ into seams also satisfies the frontier condition.* ◊

*Proof.* Since the orbit type decomposition of $Q$ satisfies the frontier condition, the orbit type decomposition of $P$ and $\check{P}$ share this property by Lemma 5.3.6. Let $\check{P}^{(K')}_{(H')}$ be a seam that has a non-empty intersection with the closure of $\check{P}^{(K)}_{(H)}$. In particular, $\check{P}^{(K')}_{(H')}$ intersects the closure of $\check{P}_{(K)}$ and thus $(K') \geq (K)$ as the orbit type decomposition satisfies the frontier condition. Since the canonical projection $\check{P} \to \check{Q}$ is continuous, a similar argument shows that $(H') \geq (H)$.

[1] One might expect that $(\mathrm{T}^*Q)^{(H)}_{(H)}$ is symplectomorphic to $\mathrm{T}^*(Q_{(H)})$. However, simple examples like $G = \mathrm{SO}(n)$ acting on $\mathbb{R}^n$ show that this is not the case.

For the converse direction, let $\check{P}^{(K)}_{(H)}$ and $\check{P}^{(K')}_{(H')}$ be seams with $(K') \geq (K)$ and $(H') \geq (H)$. We have to show that $\check{P}^{(K')}_{(H')}$ lies in the closure of $\check{P}^{(K)}_{(H)}$. Let $[p] \in \check{P}^{(K')}_{(H')}$ and choose $p \in P^{(K')}_{(H')}$. Denote the base point of $p$ by $q$. We will show that every neighborhood of $p$ in $P$ has a non-trivial intersection with the preseam $P^{(K)}_{(H)}$. Since this is a local question, it is enough to consider it in a tubular neighborhood of the form $G \times_{G_q} (\mathfrak{m}^* \times \mathrm{T}^*S)$. That is, it is enough to show that every neighborhood of $p \restriction_{\mathrm{T}_q S} \in \mathrm{T}^*_q S$ in $J^{-1}_{G_q}(0)$ contains a point $(s, \alpha) \in S \times X^* \simeq \mathrm{T}^*S$ with $(G_s) = (H)$ and $(G_\alpha) = (K)$. The existence of such a point follows from the fact that $P$ and $Q$ satisfy the frontier condition. Indeed, since $(H) \leq (H') = (G_q)$, the frontier condition implies that every neighborhood of $q$ in $Q$ contains a point $q'$ such that $(G_{q'}) = (H)$. Without loss of generality, we may assume that we work in slice coordinates and that there are, thus, $g \in G$ and $s \in S$ with $q' = g \cdot s$. By equivariance of the stabilizer, we have $(G_s) = (G_{q'}) = (H)$. The construction of $\alpha \in J^{-1}_{G_q}(0)$ with $(G_\alpha) = (K)$ follows from similar arguments using the frontier condition for $\check{P}$. □

Corollary 5.4.7 *Assume that the orbit type decomposition of $Q$ satisfies the frontier condition. Then, for every pair of orbit types $(K)$ and $(H)$ with $(H) \geq (K)$, the closure of $\check{P}^{(K)}_{(H)}$ in $\check{P}_{(K)}$ is the union of $\check{P}^{(K)}_{(H')}$ over all orbit types $(H') \geq (H)$. In particular, the decomposition of $\check{P}_{(K)}$ into seams satisfies the frontier condition.* ◊

Let us summarize.

Theorem 5.4.8 (Singular cotangent bundle reduction at zero) *Let $Q$ be a Fréchet $G$-manifold. Assume that the $G$-action is proper, that it admits at every point a slice compatible with the cotangent bundle structures and that the decomposition of $Q$ into orbit types satisfies the frontier condition. Moreover, assume that $\mathrm{T}^*Q$ is a Fréchet manifold, which is $G$-equivariantly diffeomorphic to $\mathrm{T}Q$, and that the lifted action on $\mathrm{T}^*Q$, endowed with its canonical symplectic form $\omega$, has a momentum map $J$. Assume, additionally, that the orbit $\mathfrak{g} \,.\, p$ is symplectically closed for all $p \in J^{-1}(0)$. Then, the following holds:*

*(i) The set of orbit types of $J^{-1}(0)$ with respect to the lifted $G$-action coincides with the set of orbit types for the $G$-action on $Q$.*

*(ii) The reduced phase space $\check{P} = J^{-1}(0)/G$ is stratified into orbit type subsets $\check{P}_{(K)} = (J^{-1}(0))_{(K)}/G$. For every orbit type $(K)$, the set $\check{P}_{(K)}$ is a smooth manifold and carries a symplectic form $\check{\omega}_{(K)}$.*

*(iii) Every symplectic stratum* $\check{P}_{(K)}$ *is further stratified as*

$$\check{P}_{(K)} = \bigsqcup_{(H)\geq(K)} \check{P}^{(K)}_{(H)}, \tag{5.4.18}$$

*where each seam* $\check{P}^{(K)}_{(H)}$ *is a smooth fiber bundle over* $\check{Q}_{(H)}$.

*(iv) For every orbit type* $(H)$*, the principal seam* $\check{P}^{(H)}_{(H)}$ *endowed with the restriction of the symplectic form* $\check{\omega}_{(H)}$ *is symplectomorphic to* $\mathrm{T}^*(\check{Q}_{(H)})$ *endowed with its canonical symplectic structure.*

*(v) The decomposition*

$$\check{P} = \bigsqcup_{(H)\geq(K)} \check{P}^{(K)}_{(H)} \tag{5.4.19}$$

*is a stratification of* $\check{P}$ *called the* secondary stratification. *Moreover, the projection* $\mathrm{T}^*Q \longrightarrow Q$ *induces a stratified surjective submersion* $\check{P} \longrightarrow \check{Q}$ *with respect to the secondary stratification of* $\check{P}$ *and the orbit type stratification of* $\check{Q}$. ◊

## 5.5 Dynamics

Let us now pass from the kinematic picture presented so far to dynamics. According to Proposition 4.3.8, the Hamiltonian flow on the original phase space $\mathrm{T}^*Q$ generated by a $G$-invariant Hamiltonian $h\colon \mathrm{T}^*Q \longrightarrow \mathbb{R}$ descends to a Hamiltonian flow on each symplectic stratum of the reduced phase space $\check{P}$. In this sense, the symplectic stratification of $\check{P}$ represents an additional conversed quantity. The interaction of dynamics with the secondary stratification is more complicated. The seams are in general not preserved by the Hamiltonian flow. The following example suggests that, for each orbit type $(K)$, the singular seams $\check{P}^{(K)}_{(H)}$ with $(H) > (K)$ stitch together the dynamics in the cotangent bundle $\check{P}^{(K)}_{(K)} \simeq \mathrm{T}^*(\check{Q}_{(K)})$.

Example 5.5.1 Consider a two-dimensional isotropic harmonic oscillator, whose coordinates are $q = (q_1, q_2)$ and the corresponding momenta are $p = (p_1, p_2)$. Consequently, the phase space is $\mathrm{T}^*\mathbb{R}^2$ and the Hamiltonian of the system is given by

$$H(q,p) = \frac{1}{2}\|p\|^2 + \frac{1}{2}\|q\|^2. \tag{5.5.1}$$

Note that U(1) acting by rotation in the $(q_1, q_2)$-plane is a symmetry of $H$. The angular momentum

$$J(q,p) = q_1 p_2 - q_2 p_1 \tag{5.5.2}$$

is the momentum map for the lift of this action to $T^*\mathbb{R}^2$. Hence, $J(q,p) = 0$ if and only if $q$ and $p$ are parallel. The U(1)-action on $Q$ is free except at the origin, which has stabilizer U(1). Consequently, the secondary stratification of $P = J^{-1}(0)$ is

$$P = \underbrace{\{(0,0)\}}_{P^{\mathrm{U}(1)}_{\mathrm{U}(1)}} \cup \underbrace{\{(0, p \neq 0)\}}_{P^{\{e\}}_{\mathrm{U}(1)}} \cup \underbrace{\{(q \neq 0, p) : \; q \parallel p\}}_{P^{\{e\}}_{\{e\}}} . \tag{5.5.3}$$

In order to identify the reduced phase space, consider the map

$$K \colon T^*\mathbb{R}^2 \to \mathbb{R}^3, \quad (q,p) \mapsto \begin{pmatrix} E_+ \\ E_- \\ H \end{pmatrix} := \begin{pmatrix} \frac{1}{2}\|p\|^2 - \frac{1}{2}\|q\|^2 \\ q \cdot p \\ \frac{1}{2}\|p\|^2 + \frac{1}{2}\|q\|^2 \end{pmatrix}. \tag{5.5.4}$$

The reason for this notation will become clear in a moment. On the way, we note that the combination of $K$ and $J$ yields the momentum map for the U(2)-symmetry[1], see [CB97, I.3.3]. The Hamiltonian $H$ is clearly non-negative and a direct calculation shows that $H^2 - J^2 = E_+^2 + E_-^2$. Hence, the image of $P = J^{-1}(0)$ under $K$ is the upper cone (with origin included), see Figure 5.1. Moreover, $K$ is U(1)-invariant and descends to a homeomorphism

$$\check{K} \colon \check{P} \to C \tag{5.5.5}$$

of the reduced phase space $\check{P} = J^{-1}(0)/\mathrm{U}(1)$ with the upper cone $C \subseteq \mathbb{R}^3$. The image of $\check{P}^{\mathrm{U}(1)}_{\mathrm{U}(1)}$ under $\check{K}$ is the origin and the seam $\check{P}^{\{e\}}_{\mathrm{U}(1)}$ gets mapped onto the line $L$ determined by $E_- = 0$, $H = E_+$ and $H > 0$. The remaining part of the cone corresponds to $\check{P}^{\{e\}}_{\{e\}}$.

We now pass to the symplectic structure. The Poisson brackets of the components $E_\pm$ and $H$ of $K$ (viewed as real-valued functions on $T^*\mathbb{R}^2$) are given by

$$\{H, E_\pm\} = \mp 2E_\mp \qquad \{E_+, E_-\} = 2H. \tag{5.5.6}$$

These relations are identical with the commutation relations of $\mathfrak{sl}(2,\mathbb{R})$. Hence, $K$ is a Poisson map from $T^*\mathbb{R}^2$ to $\mathbb{R}^3 \simeq \mathfrak{sl}(2,\mathbb{R})$, where the latter space carries the usual Lie–Poisson structure, i.e., the one given by the bivector field

$$\Pi = -2E_- \,\partial_H \wedge \partial_{E_+} + 2E_+ \,\partial_H \wedge \partial_{E_-} + 2H \,\partial_{E_+} \wedge \partial_{E_-}. \tag{5.5.7}$$

The symplectic top stratum $\check{P}_{\{e\}}$, i.e. the cone without the origin, is a coadjoint orbit of $\mathrm{SL}(2,\mathbb{R})$ and thus carries the Kostant–Kirillov–Souriau symplectic form.

---

[1] To be more precise, to arrive at the momentum map of Cushman and Bates [CB97, I.3.3] one has to exchange the coordinates $q_2$ and $p_1$. The U(2) symmetry is a good starting point for the qualitative discussion of the dynamics using the energy-momentum map; a topic which will not be further developed here.

As the singular stratum $\check{P}_{\mathrm{U}(1)}$ is zero-dimensional, its symplectic form vanishes.

Recall from Proposition 5.4.5 the construction of the symplectomorphism $\check{\bar{\pi}}$ between $\check{P}^{\{e\}}_{\{e\}}$ and $\mathrm{T}^*(\check{Q}_{\{e\}})$. Moreover, the map $[q] \mapsto \frac{1}{2}\|q\|^2$ identifies $\check{Q}_{\{e\}}$ with $\mathbb{R}_{>0}$ and thus yields a symplectomorphism of $\mathrm{T}^*(\check{Q}_{\{e\}})$ with $\mathrm{T}^*\mathbb{R}_{>0} \simeq \mathbb{R}_{>0} \times \mathbb{R}$. A straightforward calculation shows that the combined symplectomorphism

$$\psi\colon \check{P}^{\{e\}}_{\{e\}} \to \mathrm{T}^*(\check{Q}_{\{e\}}) \to \mathrm{T}^*\mathbb{R}_{>0} \tag{5.5.8}$$

is given by

$$\psi([q,p]) = \left(\frac{1}{2}\|q\|^2, \frac{q \cdot p}{\|q\|^2}\right) = \left(\frac{1}{2}(H - E_+), \frac{E_-}{H - E_+}\right). \tag{5.5.9}$$

Consider the map

$$I\colon \mathrm{T}^*\mathbb{R}_{>0} \to \mathbb{R}^3, \quad (\bar{q}, \bar{p}) \mapsto \begin{pmatrix} \bar{q}(\bar{p}^2 - 1) \\ 2\bar{q}\bar{p} \\ \bar{q}(\bar{p}^2 + 1) \end{pmatrix}. \tag{5.5.10}$$

The image of $I$ is the upper cone $C$ without the line $L$ and the following diagram commutes

$$\begin{array}{ccc} \check{P}^{\{e\}}_{\{e\}} & \xrightarrow{\check{K}} & C \setminus L \\ \Big\downarrow \check{\bar{\pi}} & \overset{\psi}{\searrow} & \Big\uparrow I \\ \mathrm{T}^*(\check{Q}_{\{e\}}) & \longrightarrow & \mathrm{T}^*\mathbb{R}_{>0}\,. \end{array} \tag{5.5.11}$$

Moreover, a direct calculation shows that the components of $I$ again satisfy the commutation relations (5.5.6) and hence $I$ is a Poisson immersion of $\mathrm{T}^*\mathbb{R}_{>0}$ into $(\mathfrak{sl}(2,\mathbb{R}), \Pi)$. To summarize the kinematic picture, we have decomposed the reduced phase space into the symplectic strata $\check{P}_{\mathrm{U}(1)}$ and $\check{P}_{\{e\}}$. The symplectic stratum $\check{P}_{\{e\}}$ further decomposes into the cotangent bundle $\mathrm{T}^*\mathbb{R}_{>0}$ and the line $L$. This decomposition is in accordance with Theorem 5.4.8.

Let us now discuss the dynamics. Using (5.5.7), the Hamiltonian vector field $X_H = \mathrm{d}H \lrcorner \Pi$ on $\mathbb{R}^3$ generated by $H$ is

$$X_H = -2E_- \,\partial_{E_+} + 2E_+ \,\partial_{E_-}. \tag{5.5.12}$$

Hence, the time evolution is given by rotation in the $(E_+, E_-)$-plane with $H = \text{const}$. In particular, the flow periodically hits the line $L$, i.e., the seam $\check{P}^{\{e\}}_{\mathrm{U}(1)}$. It is interesting to compare this behavior to the Hamiltonian flow on $\mathrm{T}^*\mathbb{R}_{>0}$ generated by

$$\bar{H} \coloneqq I^* H = \bar{q}(\bar{p}^2 + 1). \tag{5.5.13}$$

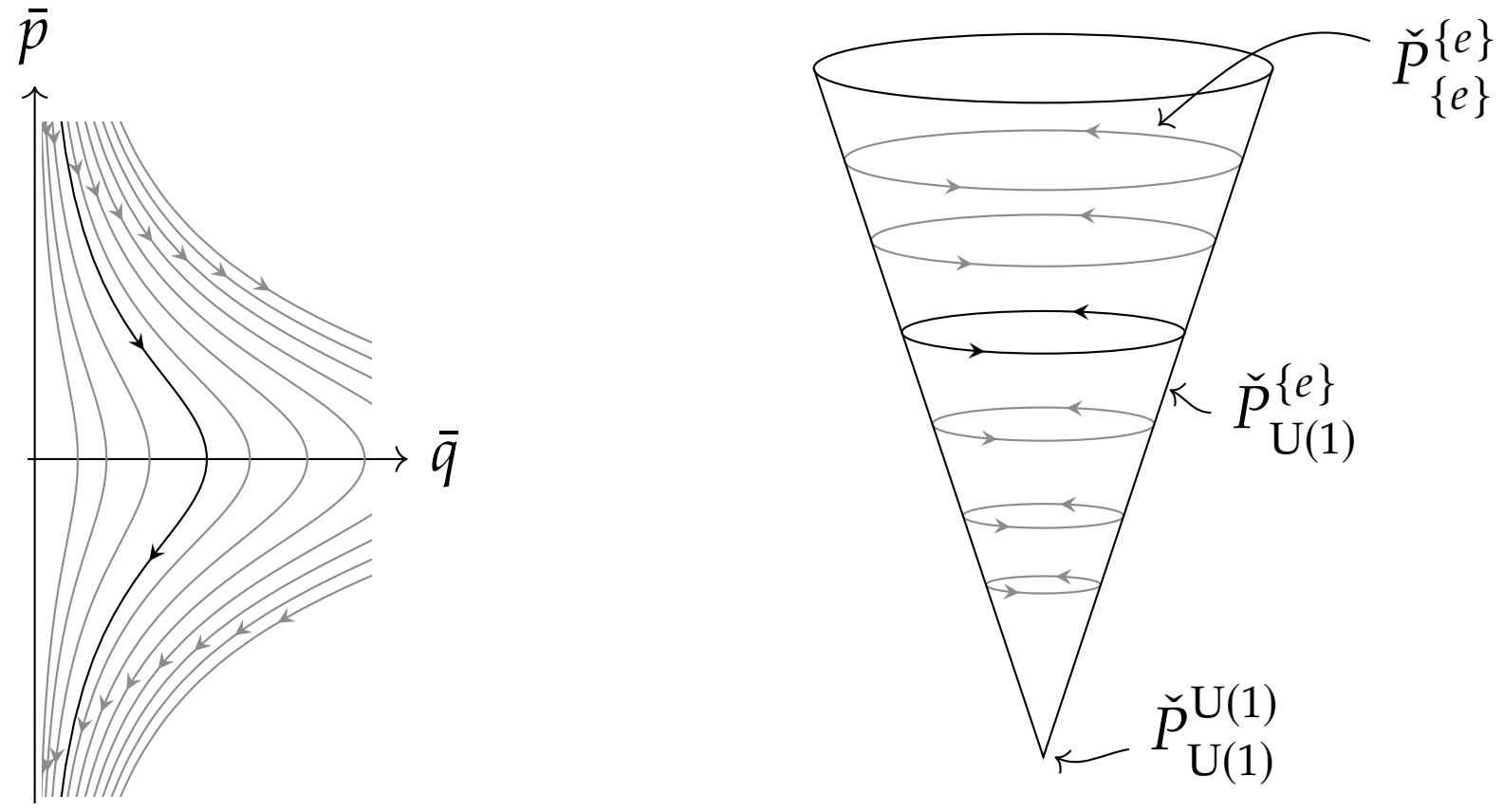

(a) Integral curves of $Y_{\bar{H}}$ in $\mathrm{T}^*\mathbb{R}_{>0}$. (b) Integral curves of $X_H$ in $\check{P}$.

Figure 5.1: Comparison of the Hamiltonian flows in $\mathrm{T}^*\mathbb{R}_{>0}$ and $\check{P}$.

The associated Hamiltonian vector field has the form

$$Y_{\bar{H}} = 2\bar{q}\bar{p}\ \partial_{\bar{q}} - (\bar{p}^2 + 1)\ \partial_{\bar{p}} \tag{5.5.14}$$

and, hence, the integral curves are given by

$$\bar{q}(t) = \bar{H}_0 \cos^2(t + t_0), \qquad \bar{p}(t) = -\tan(t + t_0), \tag{5.5.15}$$

where $\bar{H}_0$ and $t_0$ are determined by the initial conditions. Under the map $I$, they read

$$t \mapsto \bar{H}_0 \begin{pmatrix} -2\cos(2(t + t_0)) \\ -2\sin(2(t + t_0)) \\ 1 \end{pmatrix}. \tag{5.5.16}$$

Note that in $\mathrm{T}^*\mathbb{R}_{>0}$ the flow is not defined at times $t_c = \frac{\pi}{2} + k\pi - t_0$ with $k \in \mathbb{N}$. At these times, $\bar{p}$ explodes, i.e., the trajectory periodically tries to quickly leave the configuration space $\mathbb{R}_{>0}$. On the other hand, the flow under $I$ continuously extends to $t \in \mathbb{R}$. In other words, the map $I$ plays the role of regularizing the dynamics in $\mathrm{T}^*\mathbb{R}_{>0}$. In this sense, the singular seam $\check{P}^{\{e\}}_{\mathrm{U}(1)}$ stitches together the singular solution in $\check{P}^{\{e\}}_{\{e\}} \simeq \mathrm{T}^*\mathbb{R}_{>0}$ to a nice periodic flow. See Figure 5.1 for a visual comparison of the flows in $\mathrm{T}^*\mathbb{R}_{>0}$ and $\check{P}$. ◊

## 5.6 Application: Yang–Mills–Higgs theory

In the sequel, we will investigate the stratified structure of the reduced phase space of the Yang–Mills–Higgs theory. As a starting point, we use the Hamiltonian picture for this model as developed in [DR18a] based on the (3 + 1)-splitting of the configuration space and on a geometric constraint analysis, cf. also [Śni99].

First, we check that the model meets the assumptions made in the general theory developed in the previous sections. Thus, the singular cotangent bundle reduction theorem holds here and implies that the reduced phase space of the theory is a stratified symplectic space. Next, we analyze the normal form and the stratification in some detail. In particular, we find that including the singular strata leads to a refinement of what is called the resolution of the Gauß constraint in the physics literature. Our analysis of the normal form improves upon earlier work [Arm81; AMM81] on the singular geometry of the momentum map level set. We also describe the orbit types for the model after symmetry breaking, which leads to a more transparent picture of the stratification. Finally, we further analyze the secondary stratification in the concrete example of the Higgs sector of the Glashow–Weinberg–Salam model.

Let $(M, g)$ be a 3-dimensional oriented manifold with time-dependent Riemannian metric, which plays the role of a Cauchy surface in the $(3+1)$-splitting. Denote the induced volume form by $\mathrm{vol}_g$. The geometry underlying Yang–Mills-Higgs theory is that of a principal $G$-bundle $P \to M$, where $G$ is a connected compact Lie group. As usual, the gauge potential is a connection $A$ on $P$ and a bosonic matter field is a section $\varphi$ of the associated vector bundle $F = P \times_G \underline{F}$, where the typical fiber $\underline{F}$ carries a unitary $G$-representation. Thus, the space of configurations $\mathcal{Q}$ of Yang–Mills–Higgs theory consists of pairs $(A, \varphi)$. It obviously is the product of the infinite-dimensional affine Fréchet space $\mathcal{C}$ of connections[1] on $P$ and the Fréchet space $\mathcal{F}$ of sections of $F$. Let $\underline{V} \colon \underline{F} \to \mathbb{R}$ be a $G$-invariant function and denote the induced function on $F$ by $V$ (the Higgs potential).

In order to underline that the Hodge dual is defined in terms of a linear functional on the space of differential forms, we use the convention that the Hodge dual of a vector-valued differential form $\alpha \in \Omega^k(M, F)$ is the *dual-valued* differential form $* \alpha \in \Omega^{3-k}(M, F^*)$. Moreover, we use the diamond product[2]

$$\diamond \colon \Omega^k(M, F) \times \Omega^{3-r-k}(M, F^*) \to \Omega^{3-r}(M, \mathrm{Ad}^* P), \tag{5.6.1}$$

which is defined by

$$\langle \xi \wedge (\alpha \diamond \beta) \rangle = \langle (\xi \mathbin{\dot\wedge} \alpha) \wedge \beta \rangle \in \Omega^{\dim M}(M) \tag{5.6.2}$$

for all $\xi \in \Omega^r(M, \mathrm{Ad} P)$, where $\dot\wedge \colon \Omega^r(M, \mathrm{Ad} P) \times \Omega^k(M, F) \to \Omega^{r+k}(M, F)$ is the natural operation obtained by combining the Lie algebra action $\mathfrak{g} \times \underline{F} \to \underline{F}$, $(\xi, f) \mapsto \xi \,.\, f$ with the wedge product operation. Moreover, for $\alpha \in \Omega^k(M, F)$ and $\beta \in \Omega^{3-k}(M, F^*)$, we have denoted by $\langle \alpha \wedge \beta \rangle$ the real-valued top-form which arises from combining the wedge product with the natural pairing $\langle \cdot, \cdot \rangle : F \times F^* \to \mathbb{R}$.

[1] If there is no risk of confusion, we write $\mathcal{C}$ instead of $\mathcal{C}(P)$ for the space of connection on $P$. Similarly, the group of gauge transformations of $P$ will be denoted by $\mathcal{G}au(P)$ or simply by $\mathcal{G}au$.

[2] This is the natural extension to differential forms of the diamond product $\diamond \colon F \times F^* \to \mathfrak{g}^*$ that plays an important role in the study of Lie–Poisson systems.

Since $\mathcal{Q}$ is an affine space, its tangent bundle is trivial with fiber $\Omega^1(M, \mathrm{Ad}P) \times \Gamma^\infty(F)$. We will denote points in $\mathrm{T}\mathcal{Q}$ by tuples $(A, \alpha, \varphi, \zeta)$ with $\alpha \in \Omega^1(M, \mathrm{Ad}P)$ and $\zeta \in \Gamma^\infty(F)$. A natural choice for the cotangent bundle $\mathrm{T}^*\mathcal{Q}$ is the trivial bundle over $\mathcal{Q}$ with fiber $\Omega^2(M, \mathrm{Ad}^*P) \times \Omega^3(M, F^*)$. We denote elements of this fiber by pairs $(D, \Pi)$. Then, the natural pairing with $\mathrm{T}\mathcal{Q}$ is given by integration over $M$,

$$\langle (D, \Pi), (\alpha, \zeta) \rangle = \int_M \langle D \wedge \alpha \rangle + \int_M \langle \Pi \wedge \zeta \rangle. \tag{5.6.3}$$

The equations of motion are derived from their 4-dimensional covariant counterpart in the temporal gauge, see [DR18a]. In this $(3+1)$ formulation, $M$ is a Cauchy surface and the Lorentzian metric on the spacetime $\mathbb{R} \times M$ is of the form $-\ell(t)^2\, \mathrm{d}t^2 + g(t)$, where $\ell$ is the lapse function, that is, $\ell(t) \in C^\infty(M)$. With this notation at hand, the evolutionary form of the Yang–Mills–Higgs equations is given by:

$$\partial_t D = -\,\mathrm{d}_A(\ell * F_A) - \ell\varphi \diamond *(\mathrm{d}_A\varphi), \tag{5.6.4a}$$
$$\partial_t A = \ell * D, \tag{5.6.4b}$$
$$\partial_t \Pi = \mathrm{d}_A(\ell * \mathrm{d}_A\varphi) - \ell\, V'(\varphi)\, \mathrm{vol}_g, \tag{5.6.4c}$$
$$\partial_t \varphi = \ell * \Pi, \tag{5.6.4d}$$
$$\mathrm{d}_A D + \varphi \diamond \Pi = 0\,. \tag{5.6.4e}$$

By [DR18a, Theorem 3.4], the evolution equations (5.6.4a) to (5.6.4d) are Hamiltonian with respect to the Hamiltonian

$$\begin{aligned}\mathcal{H}(A, D, \varphi, \Pi) = \int_M \frac{\ell}{2}\Big(&\langle D \wedge *D \rangle + \langle F_A \wedge *F_A \rangle + \langle \Pi \wedge *\Pi \rangle \\ &+ \langle \mathrm{d}_A\varphi \wedge *\mathrm{d}_A\varphi \rangle + 2\,V(\varphi)\,\mathrm{vol}_g\Big).\end{aligned} \tag{5.6.5}$$

Moreover, equation (5.6.4e) is the Gauß constraint. In terms of the cotangent bundle geometry, it has the following interpretation, cf. [AMM81; Śni99; DR18a]. On $\mathcal{Q} = \mathcal{C} \times \mathcal{F}$ we have a left action of the group $\mathcal{G}au = \Gamma^\infty(P \times_G G)$ of local gauge transformations,

$$A \mapsto \mathrm{Ad}_\lambda A + \lambda\, \mathrm{d}\lambda^{-1}, \quad \varphi \mapsto \lambda \cdot \varphi, \tag{5.6.6}$$

for $\lambda \in \mathcal{G}au$. The Hamiltonian $\mathcal{H}$ is invariant under the lift of this action to $\mathrm{T}^*\mathcal{Q}$. A straightforward calculation shows that

$$\mathcal{J}(A, D, \varphi, \Pi) = \mathrm{d}_A D + \varphi \diamond \Pi \tag{5.6.7}$$

is the momentum map for the lifted action with respect to the natural choice $\mathfrak{gau}^* = \Omega^3(M, \mathrm{Ad}^* P)$, see [DR18a, Equation 3.10]. Hence, the Gauß constraint (5.6.4e) is equivalent to the momentum map constraint $\mathcal{J} = 0$.

Remark 5.6.1 In [DR18a], we have accomplished a unification of the Hamiltonian evolution equations (5.6.4a) to (5.6.4d) with the constraint (5.6.4e) by developing a novel variational principle (which we called the Clebsch–Lagrange principle). Besides the variation of configuration variables, the latter includes also the variation of the symmetry generators of the system. Here, these generators coincide with the time-component of the gauge potential of the 4-dimensional theory. In this language, the choice of the temporal gauge has the interpretation of being the first step in symplectic reduction by stages. In the sequel, we discuss the reduction of the remaining symmetry of the Cauchy problem. A version of the reduction by stages theorem thus shows that the reduced phase space we obtain coincides with the reduced phase space of the 4-dimensional theory. ◊

Note that the action of the group of gauge transformations on $\mathcal{Q}$ is usually not free. Hence, the model under consideration fits into the general setting of infinite-dimensional *singular* cotangent bundle reduction as discussed in Chapter 5. We now show that all assumptions made in the general discussion are met for the Yang–Mills–Higgs system:

(i) $\mathcal{Q}$ is a Fréchet manifold, because it is an affine space modeled on the Fréchet vector space $\Omega^1(M, \mathrm{Ad}P) \times \Gamma^\infty(F)$.

(ii) $\mathcal{G}au$ is a Fréchet Lie group, because it is realized as the space of sections of the group bundle $P \times_G G$.

(iii) The cotangent bundle $\mathrm{T}^*\mathcal{Q} = \mathcal{Q} \times \Omega^2(M, \mathrm{Ad}^* P) \times \Omega^3(M, F^*)$ is clearly a Fréchet manifold. The Hodge operator yields a fiber-preserving $\mathcal{G}au$-equivariant diffeomorphism between $\mathrm{T}\mathcal{Q}$ and $\mathrm{T}^*\mathcal{Q}$.

(iv) The $\mathcal{G}au$-action on $\mathcal{Q}$ is affine and thus smooth. Moreover, it is proper, see Section 3.4.

(v) The $\mathcal{G}au$-action on $\mathcal{Q}$ admits a slice at every point. First, recall from Section 3.4 that the action of $\mathcal{G}au$ on $\mathcal{C}$ has a slice $\mathcal{S}_{A_0}$ at every[1] $A_0 \in \mathcal{C}$ given in terms of the Coulomb gauge condition:

$$\mathcal{S}_{A_0} := \{A \in \mathcal{U} : \mathrm{d}^*_{A_0}(A - A_0) = 0\}, \tag{5.6.8}$$

where $\mathcal{U}$ is a suitable open neighborhood of $A_0$ in $\mathcal{C}$. Note that this slice fixes the gauge transformations up to elements of the stabilizer $\mathcal{G}au_{A_0}$ of

[1] Here and in the following, the expression $A_0$ should not be confused with the time component of a connection.

$A_0$. Thus, we are left with the $\mathcal{G}au_{A_0}$-action on $\mathcal{F}$. This is a linear action of a finite-dimensional compact group on a Fréchet space and hence has a slice $\mathcal{S}_{\varphi_0}$ at every point $\varphi_0 \in \mathcal{F}$, see Theorem A.2.4. By Proposition A.2.6, the product $\mathcal{S}_{A_0} \times \mathcal{S}_{\varphi_0}$ is a slice at $(A_0, \varphi_0)$ for the $\mathcal{G}au$-action on $\mathcal{Q}$. This slice is compatible with the cotangent bundle structures in the sense of Definition 5.2.4 as will be discussed in Section 5.6.1.

(vi) That every infinitesimal orbit of $\mathcal{G}au$ is symplectically closed will be shown in Lemma 5.6.9 below. Thus, in particular, the strong version of the Bifurcation Lemma 4.2.14 holds.

As a consequence, the Normal Form Theorem 5.2.9 holds. For the Reduction Theorems 5.3.1 and 5.4.8 to hold we assume, additionally, that the frontier condition for the decomposition of $\mathcal{Q}$ into gauge orbit types is satisfied. As discussed in Section 3.4, the orbit type decomposition of $\mathcal{C}$ satisfies the frontier condition. However, including matter fields is a rather delicate issue. As we will see, the stabilizer of the Higgs field is given in terms of a series of intersections of stabilizer groups of the $G$-action on $\underline{F}$ and this intersection is hard to analyze in full generality. For the Glashow–Weinberg–Salam model the frontier condition can be verified by direct inspection, see Proposition 5.6.12. Subsequently, we analyze the content of these theorems for the case under consideration.

### 5.6.1 Normal form

As for the classical Hodge–de Rham complex, elliptic theory[1] gives topological isomorphisms (cf. [RS17, Theorem 6.1.9]):

$$\Omega^0(M, \mathrm{Ad}P) = \operatorname{Im} \mathrm{d}_A^* \oplus \operatorname{Ker} \mathrm{d}_A, \tag{5.6.9a}$$

$$\Omega^1(M, \mathrm{Ad}P) = \operatorname{Im} \mathrm{d}_A \oplus \operatorname{Ker} \mathrm{d}_A^*, \tag{5.6.9b}$$

where $\mathrm{d}_A^*$ is the codifferential. By applying the Hodge star operator, we obtain the dual decompositions:

$$\Omega^3(M, \mathrm{Ad}^*P) = \operatorname{Im} \mathrm{d}_A \oplus \operatorname{Ker} \mathrm{d}_A^*, \tag{5.6.9c}$$

$$\Omega^2(M, \mathrm{Ad}^*P) = \operatorname{Im} \mathrm{d}_A^* \oplus \operatorname{Ker} \mathrm{d}_A. \tag{5.6.9d}$$

For clarity of presentation, we first derive the normal form of Theorem 5.2.9 for $\mathrm{T}^*\mathcal{C}$ and thereby ignore the Higgs part for a moment (see Remark 5.6.3). By the above decompositions, we have a splitting of $\mathrm{T}\mathcal{C}$,

$$\mathrm{T}_A\mathcal{C} \simeq \Omega^1(M, \mathrm{Ad}P) = \operatorname{Im} \mathrm{d}_A \oplus \operatorname{Ker} \mathrm{d}_A^*, \tag{5.6.10}$$

[1] These Hodge-type decompositions are usually derived in a Sobolev context but elliptic regularity implies that the decompositions below hold in the $\mathrm{C}^\infty$-setting.

into the canonical vertical distribution $\operatorname{Im} \mathrm{d}_A$ and the $L^2$-orthogonal complement $\operatorname{Ker} \mathrm{d}_A^*$. This decomposition is basic for the study of the geometry of the stratification, see [RS17, Sections 8.3 and 8.4]. As one expects from Hodge theory, these decompositions satisfy the annihilation relations

$$\left(\operatorname{Im} \mathrm{d}_A^*\right)^\circ = \operatorname{Ker} \mathrm{d}_A^*, \qquad (\operatorname{Ker} \mathrm{d}_A)^\circ = \operatorname{Im} \mathrm{d}_A, \tag{5.6.11}$$

$$(\operatorname{Im} \mathrm{d}_A)^\circ = \operatorname{Ker} \mathrm{d}_A, \qquad \left(\operatorname{Ker} \mathrm{d}_A^*\right)^\circ = \operatorname{Im} \mathrm{d}_A^*, \tag{5.6.12}$$

with respect to the natural pairing between $\Omega^k(M, \mathrm{Ad}P)$ and $\Omega^{3-k}(M, \mathrm{Ad}^*P)$. Thus, for the spaces involved in Theorem 5.2.9 we obtain:

$$\mathfrak{gau}_{A_0} = \operatorname{Ker} \mathrm{d}_{A_0} \colon \Omega^0(M, \mathrm{Ad}P) \to \Omega^1(M, \mathrm{Ad}P), \tag{5.6.13}$$

$$\mathfrak{gau}_{A_0}^* = \operatorname{Ker} \mathrm{d}_{A_0}^* \colon \Omega^3(M, \mathrm{Ad}^*P) \to \Omega^2(M, \mathrm{Ad}^*P), \tag{5.6.14}$$

$$\mathfrak{m} = \operatorname{Im} \mathrm{d}_{A_0}^* \colon \Omega^1(M, \mathrm{Ad}P) \to \Omega^0(M, \mathrm{Ad}P), \tag{5.6.15}$$

$$\mathfrak{m}^* = \operatorname{Im} \mathrm{d}_{A_0} \colon \Omega^2(M, \mathrm{Ad}^*P) \to \Omega^3(M, \mathrm{Ad}^*P), \tag{5.6.16}$$

$$\mathrm{T}_A \mathcal{S}_{A_0} = \operatorname{Ker} \mathrm{d}_{A_0}^* \colon \Omega^1(M, \mathrm{Ad}P) \to \Omega^0(M, \mathrm{Ad}P), \tag{5.6.17}$$

$$\mathrm{T}_A^* \mathcal{S}_{A_0} = \operatorname{Ker} \mathrm{d}_{A_0} \colon \Omega^2(M, \mathrm{Ad}^*P) \to \Omega^3(M, \mathrm{Ad}^*P). \tag{5.6.18}$$

In the sequel, we use $\|$ or $\perp$ to denote objects that are parallel or orthogonal to the gauge orbits, respectively. Accordingly, we can write every $E \in \mathrm{T}_A \mathcal{C}$ as $E = \mathrm{Ad}_\lambda(E^\perp + E^\|)$ with $E^\perp \in \mathrm{T}_A \mathcal{S}_{A_0}$ and $E^\| = -\mathrm{d}_A \chi$ for some $\chi \in \mathfrak{m}$. As a consequence, the local diffeomorphism $\mathcal{G}au \times_{\mathcal{G}au_{A_0}} (\mathfrak{m} \times \mathrm{T}\mathcal{S}_{A_0}) \to \mathrm{T}\mathcal{C}$ defined in (5.2.23) here takes the form

$$[\lambda, (\chi, A, E^\perp)] \mapsto \left(\lambda \cdot A, \mathrm{Ad}_\lambda(E^\perp - \mathrm{d}_A \chi)\right). \tag{5.6.19}$$

In order to determine the dual map $\Phi$, we introduce the Faddeev–Popov operator $\triangle_{AA_0}$ for every $A \in \mathcal{S}_{A_0}$ as the composition

$$\operatorname{Im} \mathrm{d}_{A_0} \hookrightarrow \Omega^3(M, \mathrm{Ad}^*P) \xrightarrow{\mathrm{d}_A \mathrm{d}_{A_0}^*} \Omega^3(M, \mathrm{Ad}^*P) \xrightarrow{\mathrm{pr}} \operatorname{Im} \mathrm{d}_{A_0}. \tag{5.6.20}$$

By possibly shrinking the slice $\mathcal{S}_{A_0}$, we may assume that $\triangle_{AA_0}$ is an invertible operator on $\operatorname{Im} \mathrm{d}_{A_0}$ for all $A \in \mathcal{S}_{A_0}$, see [RS17, p. 655] or [DH18, p. 22]. In particular, we have $\mathrm{d}_A \mathrm{d}_{A_0}^* \circ \triangle_{AA_0}^{-1} = \mathrm{id}_{\operatorname{Im} \mathrm{d}_{A_0}}$. Thus, the operator

$$T_A \colon \mathrm{T}_A^* \mathcal{S}_{A_0} \times \mathfrak{m}^* \to \Omega^2(M, \mathrm{Ad}^*P), \quad (D^\perp, \nu) \mapsto D^\perp + \mathrm{d}_{A_0}^* \triangle_{AA_0}^{-1} \nu \tag{5.6.21}$$

satisfies

$$\mathrm{pr}_{\operatorname{Ker} \mathrm{d}_A} \circ T_A(D^\perp, \nu) = D^\perp, \qquad \mathrm{pr}_{\operatorname{Im} \mathrm{d}_A} \circ \mathrm{d}_A T_A(D^\perp, \nu) = \nu. \tag{5.6.22}$$

Thus, $T_A$ is an isomorphism of Fréchet spaces for every $A \in \mathcal{S}_{A_0}$. Moreover, we find

$$\begin{aligned}\langle E^\perp - \mathrm{d}_A \chi, T_A(D^\perp, \nu)\rangle &= \langle E^\perp, T_A(D^\perp, \nu)\rangle + \langle \chi, \mathrm{d}_A T_A(D^\perp, \nu)\rangle \\ &= \langle E^\perp, D^\perp\rangle + \langle \chi, \nu\rangle.\end{aligned} \tag{5.6.23}$$

In summary, here, the local diffeomorphism $\Phi\colon \mathcal{G}au \times_{\mathcal{G}au_{A_0}} (\mathfrak{m}^* \times \mathrm{T}^*\mathcal{S}_{A_0}) \to \mathrm{T}^*\mathcal{C}$ defined in (5.2.19) has the form

$$\Phi([\lambda, (\nu, A, D^\perp)]) = \Big(\lambda \cdot A, \mathrm{Ad}^*_\lambda (D^\perp + \mathrm{d}^*_{A_0} \vartriangle^{-1}_{AA_0} \nu)\Big)\,. \tag{5.6.24}$$

Thus, every $D \in \mathrm{T}^*_A\mathcal{C}$ decomposes as $D = \mathrm{Ad}^*_\lambda(D^\perp + D^\parallel)$ with $D^\perp \in \mathrm{T}^*_A\mathcal{S}_{A_0}$ and $D^\parallel = \mathrm{d}^*_{A_0} \vartriangle^{-1}_{AA_0} \nu$ for some $\nu \in \mathfrak{m}^*$. That is, $\nu = \mathrm{d}_{A_0} D^\parallel$.

Remark 5.6.2 As we have seen, for every $A \in \mathcal{S}_{A_0}$, the map

$$\mathfrak{m} \times \mathrm{T}_A\mathcal{S}_{A_0} \to \mathrm{T}_A\mathcal{C}, \quad (\chi, E^\perp) \mapsto E^\perp - \mathrm{d}_A\chi \tag{5.6.25}$$

has the adjoint $T_A$. Since $T_A$ is an isomorphism, we see that the slice $\mathcal{S}$ is compatible with the cotangent bundle structures. ◊

In these local coordinates, the momentum map $\mathcal{J}$ of (5.6.7) is expressed as follows, cf. (5.2.20):

$$\mathcal{J}(\Phi([\lambda, (\nu, A, D^\perp)]), \varphi, \Pi) = \mathrm{Ad}^*_\lambda(\mathrm{d}_A D^\perp + \nu) + \varphi \diamond \Pi, \tag{5.6.26}$$

where $A \in \mathcal{S}_{A_0}$ and $D^\perp \in \operatorname{Ker} \mathrm{d}_{A_0}$. We denote the matter charge density $\varphi \diamond \Pi$ by $\rho$ and decompose it with respect to (5.6.9c),

$$\rho = \mathrm{Ad}^*_\lambda(\rho^\parallel + \rho^\perp), \tag{5.6.27}$$

where $\rho^\parallel \in \operatorname{Im} \mathrm{d}_{A_0}$ and $\rho^\perp \in \operatorname{Ker} \mathrm{d}^*_{A_0}$. Note that $\mathrm{d}_A D^\perp \in \operatorname{Ker} \mathrm{d}^*_{A_0}$, because by upper semi-continuity of the kernel of a semi-Fredholm operator [Hör07, Corollary 19.1.6] we have

$$\operatorname{Ker} \mathrm{d}^*_{A_0} \mathrm{d}_A \subseteq \operatorname{Ker} \mathrm{d}^*_{A_0} \mathrm{d}_{A_0} = \operatorname{Ker} \mathrm{d}_{A_0}. \tag{5.6.28}$$

Thus, with respect to the decomposition $\mathfrak{gau}^* = \mathfrak{m}^* \oplus \mathfrak{gau}^*_{A_0}$, the Gauß constraint $\mathcal{J} \circ \Phi = 0$ is equivalent to the following equations:

$$\nu + \rho^\parallel = 0\,, \tag{5.6.29a}$$

$$\mathrm{d}_A D^\perp + \rho^\perp = 0\,. \tag{5.6.29b}$$

In summary, we have decomposed the non-linear Gauß constraint into a linear equation and finitely many non-linear equations according to the fact that the stabilizer $\mathfrak{gau}_{A_0}$ is finite-dimensional. In spirit, this splitting is similar to the Kuranishi method or the Lyapunov–Schmidt construction used in Section 2.2, where we also found convenient coordinates to reduce a non-linear equation to a non-linear equation in finite dimensions. However, in contrast to the construction in Section 2.2, we have found the coordinates here by exploiting the cotangent bundle geometry of the problem and have not directly used the Inverse Function Theorem. Note that the construction of the local diffeomorphism $\Phi$ involves the non-local operator $\triangle^{-1}_{AA_0}$. Hence, the reconstruction of the solution $(A, D)$ of the Gauß constraint from a solution $(\nu, A, D^{\perp})$ of (5.6.29) is a non-local and non-linear operation. In the physics literature, one usually only considers the case of trivial stabilizer $\mathfrak{gau}_{A_0}$ (i.e., irreducible connections). In this case, the above construction reduces the Gauß constraint to the linear equation (5.6.29a). However, if one wants to include non-generic configurations, that is, connections with a non-trivial stabilizer, (5.6.29b) must be taken into account as well.

On the quantum level, our observations suggest that the standard quantization methods, which only take the linear constraint (5.6.29a) into account, fail in the neighborhood of reducible connections and need to be supplemented by the non-linear constraint (5.6.29b). For a quantization program for *lattice* gauge theories, where non-generic gauge orbit strata are included, we refer to [RS17, Chapter 9] and to [HRS09] for a case study.

Remark 5.6.3 In the above construction of the normal form, we have considered only $\mathrm{T}^*\mathcal{C}$ and, thereby, we have ignored the Higgs part. By passing to the normal form of the full cotangent bundle $\mathrm{T}^*\mathcal{Q}$ according to Theorem 5.2.9, we note that the stabilizer $\mathfrak{gau}_{A_0}$ further decomposes into the stabilizer $\mathfrak{gau}_{(A_0,\varphi_0)} = \mathfrak{gau}_{A_0} \cap \mathfrak{gau}_{\varphi_0}$ of $(A_0, \varphi_0) \in \mathcal{Q}$ and some complement $\mathfrak{r}$. Accordingly, the non-linear part (5.6.29b) of the Gauß law further decomposes into a linear equation in $\mathfrak{r}^*$ and a non-linear equation in $\mathfrak{gau}^*_{(A_0,\varphi_0)}$. ◇

### 5.6.2 Orbit types

Next, let us consider Theorem 5.3.1. By point (i) of this theorem, the set of orbit types of $\mathcal{P} = \mathcal{J}^{-1}(0)$ with respect to the lifted $\mathcal{G}au$-action coincides with the set of orbit types for the $\mathcal{G}au$-action on $\mathcal{Q}$. We will now determine these orbit types. First, recall from Section 3.4 that the stabilizer $\mathcal{G}au_A$ of $A \in \mathcal{C}$ under the $\mathcal{G}au$-action is isomorphic to the centralizer $\mathrm{C}_G(\mathrm{Hol}_A)$ of the holonomy group of $A$ (based at some point $p_0 \in P$), because every gauge transformation $\lambda \in \mathcal{G}au_A$ is constant on the holonomy bundle $P_A$ of $A$. It is straightforward to include the matter fields $\varphi$ and to pass, thereby, from $\mathcal{C}$ to $\mathcal{Q}$. Indeed, a gauge transformation $\lambda$ leaves $\varphi$ invariant if and only if $\lambda(p) \in G_{\varphi(p)}$ for all $p \in P$, where $G_{\varphi(p)}$ is the stabilizer of $\varphi(p) \in \underline{F}$ under the $G$-action. The equivariance properties $\lambda(p \cdot g) = g^{-1}\lambda(p)g$ and

$$G_{\varphi(p\cdot g)} = G_{g^{-1}\cdot\varphi(p)} = g^{-1}G_{\varphi(p)}g \tag{5.6.30}$$

show that it is actually enough to test $\lambda(p) \in G_{\varphi(p)}$ only for one point per fiber. In particular, it suffices to let $p$ range over points in the holonomy bundle $P_A$. Since a gauge transformation $\lambda$ in the stabilizer of $A$ is necessarily constant on $P_A$, the evaluation map $\mathrm{ev}_{p_0}\colon \mathcal{G}au \to G$ defined in (3.4.4) restricts to an isomorphism of Lie groups

$$\mathcal{G}au_{A,\varphi} = \mathcal{G}au_A \cap \mathcal{G}au_\varphi \simeq \mathrm{C}_G(\mathrm{Hol}_A) \cap \bigcap_{p \in P_A} G_{\varphi(p)} \,. \tag{5.6.31}$$

To summarize, by point (i) of Theorem 5.3.1, we have completely determined the stabilizer subgroups of points in $\mathcal{P} = \mathcal{J}^{-1}(0) \subseteq \mathrm{T}^*\mathcal{Q}$ with respect to the lifted $\mathcal{G}au$-action.

Finally, let us also determine the orbit types of $\mathrm{T}^*\mathcal{Q}$. An arbitrary point $(A, D) \in \mathcal{C} \times \Omega^2(M, \mathrm{Ad}^* P)$ in the cotangent bundle $\mathrm{T}^*\mathcal{C}$ has stabilizer

$$\mathcal{G}au_{A,D} = \mathcal{G}au_A \cap \mathcal{G}au_D, \tag{5.6.32}$$

where $\mathcal{G}au_A$ and $\mathcal{G}au_D$ denote the stabilizers of $A$ and $D$ under the action of the gauge group, respectively.

To analyze $\mathcal{G}au_{A,D}$, for a moment, consider the stabilizer of a general $k$-form $\alpha$ with values in $F = P \times_G \underline{F}$, where $\underline{F}$ carries a $G$-representation. For $\alpha \in \Omega^k(M, F)$ and $p \in P$, let $V_\alpha(p)$ be the subspace of $\underline{F}$ spanned by all elements of the form $\alpha_p(v_1, \dots, v_k)$, where $v_i \in \mathrm{T}_p P$. Similar arguments as above for $\alpha = \varphi$ yield the following identification

$$\mathcal{G}au_{A,\alpha} \simeq \mathrm{C}_G(\mathrm{Hol}_A) \cap \bigcap_{\substack{p \in P_A \\ f \in V_\alpha(p)}} G_f \,. \tag{5.6.33}$$

Now let us return to the stabilizer of $D$. Denote $K \equiv \mathrm{Hol}_A$. Since the adjoint action of $\mathrm{C}_G(K)$ on the Lie algebra $\mathfrak{k}$ of $K$ is trivial, the subspace $\mathfrak{k} \subseteq \mathfrak{g}$ is $\mathrm{Ad}_{\mathrm{C}_G(K)}$-invariant and thus has a $\mathrm{C}_G(K)$-invariant complement $\mathfrak{p}$ in $\mathfrak{g}$. For $\nu \in \mathfrak{k}^*$, we have

$$\langle \mathrm{Ad}^*_g \nu, \xi \rangle = \langle \nu, \mathrm{Ad}^{-1}_g \xi \rangle = \langle \nu, \xi \rangle \tag{5.6.34}$$

for every $g \in \mathrm{C}_G(K)$ and $\xi \in \mathfrak{k}$. Hence, $G_\nu \subseteq \mathrm{C}_G(K)$ for every $\nu \in \mathfrak{k}^*$. In summary, we have shown:

Proposition 5.6.4 *Let $A \in \mathcal{C}$ and $D \in \mathrm{T}^*_A\mathcal{C}$. Then,*

$$\mathcal{G}au_{A,D} \simeq \mathrm{C}_G(\mathrm{Hol}_A) \cap \bigcap_{\substack{p \in P_A \\ \nu \in V_D(p) \cap \mathfrak{p}^*}} G_\nu \,, \tag{5.6.35}$$

*where $V_D(p)$ is the subspace of $\mathfrak{g}^*$ spanned by all elements of the form $D_p(v_1, v_2)$ for $v_i \in \mathrm{T}_p P$.* ◇

### 5.6.3 Description of orbit types after symmetry breaking

Now, let us pass to a description of the model after symmetry breaking. As usual in the physics literature, let us assume that $\varphi$ takes values in one fixed orbit type $\underline{F}_{(K)}$ of the $G$-action for some stabilizer subgroup $K$. Assume, moreover, that the bundle $\underline{F}_{(K)} \to \check{\underline{F}}_{(K)} = \underline{F}_{(K)}/G$ is trivial and choose a smooth section $f_0 \colon \check{\underline{F}}_{(K)} \to \underline{F}_{(K)}$ which takes values in the subset $\underline{F}_K$ of isotropy type $K$. With respect to this choice, we can write the Higgs field $\varphi$, viewed as a $G$-equivariant map $P \to \underline{F}$, as

$$\varphi(p) = \phi(p) \cdot f_0(\eta(p)), \quad p \in P, \tag{5.6.36}$$

where $\phi \colon P \to G$ and $\eta \colon P \to \check{\underline{F}}_{(K)}$ are smooth maps. Since $f_0$ takes values in $\underline{F}_K$, the map $\phi$ is uniquely defined when viewed as a map with values in $G/K$. Moreover, by $G$-equivariance of $\varphi$, we identify $\phi$ as a section of $P \times_G G/K$ and $\eta$ as a smooth map $M \to \check{\underline{F}}_{(K)}$. The map $\eta$ (or rather $\eta$ shifted by the Higgs vacuum) is the surviving Higgs field. It is straightforward to see that the decomposition (5.6.36) depends smoothly on $\varphi$ and hence establishes a diffeomorphism

$$\mathcal{F} \to \Gamma^\infty(P \times_G G/K) \times C^\infty\big(M, \check{\underline{F}}_{(K)}\big), \qquad \varphi \mapsto (\phi, \eta) \tag{5.6.37}$$

of Fréchet manifolds. In geometric terms, $\phi$ yields a reduction of $P$ to the principal $K$-bundle

$$\hat{P} := \{p \in P : \ \phi(p) = [e]\}. \tag{5.6.38}$$

Next, recall the following geometric version of the Higgs mechanism (cf. [RS17, Proposition 7.3.4]). Since $K$ is compact, there exists an $\mathrm{Ad}_K$-invariant decomposition $\mathfrak{g} = \mathfrak{k} \oplus \mathfrak{p}$. Accordingly, the restriction of $A$ to the $K$-bundle $\hat{P}$ splits into

$$A_{\restriction \hat{P}} = \hat{A} + \tau, \tag{5.6.39}$$

where $\hat{A}$ and $\tau$ take values in $\mathfrak{k}$ and $\mathfrak{p}$, respectively. It is an easy exercise to verify that $\hat{A}$ is a principal $K$-connection in $\hat{P}$ and $\tau$ a horizontal 1-form of type $\mathrm{Ad}_K \, \mathfrak{p}$ on $\hat{P}$, i.e., $\tau \in \Omega^1(M, \hat{P} \times_K \mathfrak{p})$. In the physics language, $\hat{A}$ is the reduced gauge field and $\tau$ is the intermediate vector boson. Let us combine this decomposition of $A$ with the diffeomorphism of (5.6.37). For this purpose, consider the smooth Fréchet bundle $\mathcal{E} \to \Gamma^\infty(P \times_G G/K)$, whose fiber over $\phi$ is $\mathcal{C}(\hat{P}) \times \Omega^1(M, \hat{P} \times_K \mathfrak{p})$. The bundle $\mathcal{E}$ carries a natural action[1] of $\mathcal{G}au(P)$ as $\hat{P}$ is a subbundle of $P$. To summarize we obtain the following.

Proposition 5.6.5 *The map $(A, \varphi) \mapsto (\phi, \hat{A}, \tau, \eta)$ defines a $\mathcal{G}au(P)$-equivariant diffeomorphism of $\mathcal{Q}$ with $\mathcal{E} \times C^\infty(M, \check{\underline{F}}_{(K)})$. In particular, we get an isomorphism of stratified spaces between the gauge orbit space $\check{\mathcal{Q}}$ and $\big(\mathcal{E}/\mathcal{G}au(P)\big) \times C^\infty\big(M, \check{\underline{F}}_{(K)}\big)$.* ◊

[1] Note that the action of $\mathcal{G}au$ on $\mathcal{E}$ mixes the variables $\hat{A}$ and $\tau$, because the decomposition $\mathfrak{g} = \mathfrak{k} \oplus \mathfrak{p}$ is not $\mathrm{Ad}_G$-invariant.

When the bundle $\hat{P}$ is non-trivial, the field $\phi$ representing $\hat{P}$ carries topological data, which may be encoded in various ways, e.g. in terms of Chern classes.

Using this proposition, the gauge orbit types of $\mathcal{Q}$ may be characterized in the following more explicit way. First, note that $\eta$ does not contribute to the orbit type structure. Next, recall that $\lambda \in \mathcal{G}au(P)$ preserves $\varphi$ if and only if $\lambda(p) \cdot \varphi(p) = \varphi(p)$ for all $p \in P$. Thus, every $\lambda \in \mathcal{G}au_\varphi(P)$ restricts to a $K$-gauge transformation on $\hat{P}$. In fact, a moment's reflection shows that every $K$-gauge transformation on $\hat{P}$ can be obtained in that way (use $P \times_G K \simeq \hat{P} \times_K K$). This yields the following.

Lemma 5.6.6 *For every $\varphi$ with associated K-bundle $\hat{P}$, we have $\mathcal{G}au_\varphi(P) = \mathcal{G}au(\hat{P})$.* ◊

Let $(A, \varphi) \in \mathcal{Q}$. By equivariance and (5.6.31), we find[1]

$$\mathcal{G}au_{A,\varphi}(P) \simeq \mathrm{C}_G(\mathrm{Hol}_A) \cap \bigcap_{p \in P_A} \phi(p) K \phi(p)^{-1}. \tag{5.6.40}$$

It is interesting to compare the gauge orbit types $\mathcal{G}au_{A,\varphi}(P)$ to the orbit types of the theory after symmetry reduction. Using Proposition 5.6.5, we have

$$\mathcal{G}au_{A,\varphi}(P) = \mathcal{G}au_{\hat{A},\tau}(\hat{P}) = \mathcal{G}au_{\hat{A}}(\hat{P}) \cap \mathcal{G}au_\tau(\hat{P}), \tag{5.6.41}$$

where $\mathcal{G}au_{\hat{A},\tau}(\hat{P})$ denotes the stabilizer of $(\hat{A}, \tau)$ under the natural $\mathcal{G}au(\hat{P})$-action. Indeed, by Lemma 5.6.6, $\mathcal{G}au_\varphi(P)$ is isomorphic to $\mathcal{G}au(\hat{P})$ and it is straightforward to see that a gauge transformation $\lambda \in \mathcal{G}au_\varphi(P)$ leaves $A$ invariant if and only if $\lambda_{\restriction \hat{P}} \in \mathcal{G}au(\hat{P})$ leaves $A_{\restriction \hat{P}} = \hat{A} + \tau$ invariant. Furthermore, using (5.6.33), we obtain a more explicit description of $\mathcal{G}au_\tau(\hat{P})$, which, in summary, yields the following.

Proposition 5.6.7 *The stabilizer of $(A, \varphi)$ under the $\mathcal{G}au(P)$-action on $\mathcal{Q}$ is given by*

$$\mathcal{G}au_{A,\varphi}(P) \simeq \mathrm{C}_K(\mathrm{Hol}_{\hat{A}}) \cap \bigcap_{\substack{p \in \hat{P}_{\hat{A}} \\ \xi \in V_\tau(p)}} G_\xi\,. \tag{5.6.42}$$

◊

Finally, let us consider a special case which is important for instance in the theory of magnetic monopoles [RS17, Chapter 7]. It is defined by the additional assumption that $\phi$ be covariantly constant with respect to $A$. Then, $A$ reduces to a connection $\hat{A}$ on $\hat{P}$ and so $\tau = 0$. Thus, in this case, (5.6.42) simplifies to

$$\mathcal{G}au_{A,\varphi}(P) \simeq \mathcal{G}au_{\hat{A}}(\hat{P}) \simeq \mathrm{C}_K(\mathrm{Hol}_{\hat{A}}). \tag{5.6.43}$$

[1] In the sequel, we assume that the base point $p_0$ defining $P_A$ lies in $\hat{P}$, which can be always assured by translation with a constant $g \in G$.

As an immediate consequence, we obtain the following analogue of Proposition 3.4.1.

Proposition 5.6.8 *The orbit types of the action of $\mathcal{G}au(P)$ on $\mathcal{Q}$ at points $(A, \varphi)$ with $\mathrm{d}_A\phi = 0$ are in one-to-one correspondence with isomorphism classes of holonomy-induced Howe subbundles of K-reductions $\hat{P}$ of P via the map*

$$[\mathcal{G}au_{A,\varphi}(P)] \mapsto [\hat{P}_{\hat{A}} \cdot \mathrm{C}_K^2(\mathrm{Hol}_{\hat{A}})]. \tag{5.6.44}$$

◇

*Proof.* By (5.6.43), the stabilizer $\mathcal{G}au_{A,\varphi}(P)$ of a point $(A, \varphi) \in \mathcal{Q}$ satisfying $\mathrm{d}_A\phi = 0$ is isomorphic to the stabilizer of $\hat{A}$ under $\mathcal{G}au(\hat{P})$. Now recall from Proposition 3.4.1 that orbit types of connections on $\hat{P}$ are in bijective correspondence with isomorphism classes of holonomy-induced Howe subbundles of $\hat{P}$. □

### 5.6.4 Reduced symplectic structure

Finally, let us come to point (iii) of Theorem 5.3.1. The symplectic structure $\Omega$ on $\mathrm{T}^*\mathcal{Q}$ is defined by

$$\begin{aligned}\Omega_{A,\varphi,D,\Pi}&\big((\delta A_1, \delta\varphi_1, \delta D_1, \delta\Pi_1), (\delta A_2, \delta\varphi_2, \delta D_2, \delta\Pi_2)\big)\\ &= \int_M \langle \delta D_1 \wedge \delta A_2\rangle - \langle \delta D_2 \wedge \delta A_1\rangle + \langle \delta\Pi_1 \wedge \delta\varphi_2\rangle - \langle \delta\Pi_2 \wedge \delta\varphi_1\rangle,\end{aligned} \tag{5.6.45}$$

where $(\delta A_j, \delta D_j, \delta\varphi_j, \delta\Pi_j) \in \mathrm{T}_{(A,D,\varphi,\Pi)}(\mathrm{T}^*\mathcal{Q})$ for $j = 1, 2$.

Lemma 5.6.9 *The orbit $\mathfrak{gau} \,.\, (A, D, \varphi, \Pi)$ is symplectically closed in $\mathrm{T}^*\mathcal{Q}$ for every $(A, D, \varphi, \Pi) \in \mathrm{T}^*\mathcal{Q}$, that is,*

$$\big(\mathfrak{gau} \,.\, (A, D, \varphi, \Pi)\big)^{\Omega\Omega} = \mathfrak{gau} \,.\, (A, D, \varphi, \Pi). \tag{5.6.46}$$

◇

*Proof.* First, note that

$$j(\delta A, \delta D, \delta\varphi, \delta\Pi) := (-*\delta D, *\delta A, -*\delta\Pi, *\delta\varphi) \tag{5.6.47}$$

defines an almost complex structure $j$ on $\mathrm{T}^*\mathcal{Q}$, which intertwines the symplectic structure $\Omega$ with the $\mathrm{L}^2$-scalar product. Clearly, for every vector subspace $W \subseteq \mathrm{T}_{(A,D,\varphi,\Pi)}(\mathrm{T}^*\mathcal{Q})$ we have $W^\Omega = jW^\perp$, where $W^\perp$ is the $\mathrm{L}^2$-orthogonal complement. Hence,

$$W^{\Omega\Omega} = j(jW^\perp)^\perp = W^{\perp\perp} \tag{5.6.48}$$

and it remains to show that $W^{\perp\perp} = W$ for $W = \mathfrak{gau} \,.\, (A, D, \varphi, \Pi)$. By the Bipolar Theorem B.1.8, this holds if and only if the orbit is closed with respect to the weak $\mathrm{L}^2$-topology. Note that the infinitesimal action has the character of a multiplication operator in the $D$- and $\varphi$-direction and hence in these components the orbit is not closed (it may actually be dense). Nonetheless, the 'diagonal' orbit $\mathfrak{gau} \,.\, (A, D, \varphi, \Pi)$ is closed as we will show now.

Let $(\delta A, \delta D, \delta\varphi, \delta\Pi)$ be an element of $\mathrm{T}_{(A,D,\varphi,\Pi)}(\mathrm{T}^* Q)$ that does not lie on the $\mathfrak{gau}$-orbit. It will be convenient to write $\alpha \equiv \delta A$ and abbreviate $(\delta D, \delta\varphi, \delta\Pi)$ by $\beta$. We have to construct an $\mathrm{L}^2$-open neighborhood $\mathcal{U}$ of $(\alpha, \beta)$ in $\mathrm{T}_{(A,D,\varphi,\Pi)}(\mathrm{T}^* \mathcal{Q})$ that is disjoint from the $\mathfrak{gau}$-orbit. First, suppose that $\alpha$ cannot be written as $\mathrm{d}_A \xi$ for some $\xi \in \mathfrak{gau}$. Since the decomposition (5.6.9b) is orthogonal with respect to the $\mathrm{L}^2$-scalar product, the orbit $\mathfrak{gau} \,.\, A$ is closed in $\mathrm{T}_A \mathcal{C}$ with respect to the weak $\mathrm{L}^2$-topology. Thus, there exists an $\mathrm{L}^2$-open neighborhood $\mathcal{U}_\alpha$ of $\alpha$ in $\Omega^1(M, \mathrm{Ad} P)$ which is disjoint from $\mathfrak{gau} \,.\, A$. Then,

$$\mathcal{U} \coloneqq \mathcal{U}_\alpha \times \Omega^2(M, \mathrm{Ad}^* P) \times \Gamma^\infty(F) \times \Omega^3(M, F^*) \tag{5.6.49}$$

is an $\mathrm{L}^2$-open neighborhood of $(\alpha, \beta)$ in $\mathrm{T}_{(A,D,\varphi,\Pi)}(\mathrm{T}^* Q)$ which does not intersect the $\mathfrak{gau}$-orbit. Now suppose that $\alpha = \mathrm{d}_A \xi$ for some $\xi \in \mathfrak{gau}$. Then, $\xi$ is uniquely determined up to an element of the finite-dimensional stabilizer $\mathfrak{gau}_A$. Since the orbit $\mathfrak{gau}_A \,.\, (D, \varphi, \Pi)$ is finite-dimensional, it is automatically closed. Hence, there exists an $\mathrm{L}^2$-open neighborhood $\mathcal{U}_\beta$ of $\beta$ and an $\mathrm{L}^2$-open subset $\mathcal{V}$ containing the orbit $\mathfrak{gau}_A \,.\, (D, \varphi, \Pi)$ such that $\mathcal{U}_\beta$ and $\mathcal{V}$ have empty intersection. Let $\mathcal{W} \subseteq \mathfrak{gau}$ denote the inverse image of $\mathcal{V}$ under the action $\xi \to \xi \,.\, (D, \varphi, \Pi)$. Note that $\mathfrak{gau}_A \subseteq \mathcal{W}$. Now,

$$\mathcal{U}_\alpha \coloneqq \mathrm{d}_A \mathcal{W} \oplus \operatorname{Ker} \mathrm{d}_A^* \tag{5.6.50}$$

is an $\mathrm{L}^2$-open neighborhood of $\alpha$ in $\Omega^1(M, \mathrm{Ad} P)$. Suppose $\alpha' \in \mathcal{U}_\alpha$ and $\beta' \in \mathcal{U}_\beta$ are of the form $\alpha' = \mathrm{d}_A \xi'$ and $\beta' = \xi' \,.\, \beta$ for some $\xi' \in \mathfrak{gau}$. Then, $\xi' \in \mathcal{W}$ and thus $\beta' = \xi' \,.\, \beta \in \mathcal{V}$. However, by assumption, $\beta'$ is an element of $\mathcal{U}_\beta$. Since the latter is disjoint from $\mathcal{V}$, we have constructed a contradiction. In summary, $\mathcal{U} \coloneqq \mathcal{U}_\alpha \times \mathcal{U}_\beta$ is an $\mathrm{L}^2$-open neighborhood of $(\alpha, \beta)$ in $\mathrm{T}_{(A,D,\varphi,\Pi)}(\mathrm{T}^* Q)$, which has empty intersection with the $\mathfrak{gau}$-orbit. □

By this lemma, the strata of the reduced phase space inherit a symplectic form from $(\mathrm{T}^* \mathcal{Q}, \Omega)$, according to point (iii) of Theorem 5.3.1.

Moreover, we have shown that Theorem 5.4.8 holds for Yang–Mills–Higgs theory and thus, in particular, every symplectic stratum further decomposes into seams and a cotangent bundle. This secondary stratification will be further analyzed below in the concrete example of the Glashow–Weinberg–Salam model.

### 5.6.5 Example: Glashow–Weinberg–Salam model

We now specialize the discussion to the Higgs sector of the standard model of electroweak interactions. As in the general setting, the configurations of this model are pairs $(A, \varphi)$ consisting of a connection in a principal bundle $P$ and a section of an associated vector bundle $P \times_G \underline{F}$. For this model, the principal bundle $P$ is an SU(2) × U(1)-bundle over $M$ and the associated bundle $F$ has typical fiber $\underline{F} = \mathbb{C}^2$ carrying the following representation of SU(2) × U(1):

$$\rho_{a,\vartheta}(z_1, z_2) = a \cdot e^{\frac{i}{2}\vartheta} \cdot (z_1, z_2), \qquad a \in \mathrm{SU}(2), \vartheta \in [0, 4\pi). \tag{5.6.51}$$

The Higgs potential $\underline{V}\colon \mathbb{C}^2 \to \mathbb{R}$ has the form

$$\underline{V}(f) = \lambda \left( \|f\|^2 - \frac{\nu^2}{2} \right)^2 \tag{5.6.52}$$

for given $\lambda > 0$ and non-zero $\nu \in \mathbb{R}$.

It is straightforward to see that under the representation $\rho$ the origin is a fixed point and that all other points have a stabilizer conjugate to

$$K := \left\{ \left( \begin{pmatrix} e^{\frac{i}{2}\vartheta} & 0 \\ 0 & e^{-\frac{i}{2}\vartheta} \end{pmatrix}, e^{i\vartheta} \right) \; : \; \vartheta \in [0, 4\pi) \right\}, \tag{5.6.53}$$

which is isomorphic to U(1) and plays the role of the electromagnetic gauge group after symmetry breaking. As common in the physics literature, we assume that $\varphi$ vanishes nowhere, that is, we ignore the singular orbit type in $\underline{F}$. The generic orbits in $\underline{F}$ are three-spheres centered at the origin and hence the quotient $\check{\underline{F}}_{(K)} = \underline{F}_{(K)}/G$ is diffeomorphic to $\mathbb{R}_{>0}$. All the points $\frac{r}{\sqrt{2}}(0, \nu)$ for $r \in \mathbb{R}_{>0}$ have stabilizer $K$. Hence, the map $f_0\colon r \mapsto \frac{r}{\sqrt{2}}(0, \nu)$ is a smooth section of $\underline{F}_{(K)} \to \check{\underline{F}}_{(K)}$ taking values in $\underline{F}_K$. Accordingly, the decomposition (5.6.36) of $\varphi$ simplifies here to

$$\varphi = \frac{\eta}{\sqrt{2}} \phi \cdot \begin{pmatrix} 0 \\ \nu \end{pmatrix}, \tag{5.6.54}$$

where $\phi \in \Gamma^\infty(P \times_G G/K)$ and $\eta \in C^\infty(M, \mathbb{R}_{>0})$. Note that, in this presentation, $V(\varphi) = \frac{\lambda \nu^2}{2}(\eta^2 - 1)^2$.

### 5.6.5.1 *Orbit types of $Q$*

According to the general program, we have to determine the Howe subgroups of SU(2) × U(1). For that purpose, we use the following elementary result, whose proof is a simple calculation that we leave to the reader.

Lemma 5.6.10 *Let $G$ be a group and $L$ be an abelian group. A subgroup $H$ of $G \times L$ is Howe if and only if there exists a subgroup $H'_G$ of $G$ such that*

$$H = C_G(H'_G) \times L. \tag{5.6.55}$$

◊

| $\mathrm{Hol}_A$ | $\mathcal{G}au_A$ | $H$ |
|---|---|---|
| $\{e\}, \mathbb{Z}_2$ | SU(2) | $\mathbb{Z}_2$ |
| U(1) | U(1) | U(1) |
| SU(2) | $\mathbb{Z}_2$ | SU(2) |

Table 5.1: List of all possible holonomy groups for SU(2) up to conjugacy with the corresponding stabilizer $\mathcal{G}au_A = \mathrm{C}_G(\mathrm{Hol}_A)$ and the Howe subgroup $H = \mathrm{C}^2_G(\mathrm{Hol}_A)$.

According to this lemma, we first have to determine the Howe subgroups of SU(2). Clearly, each Howe subgroup of SU(2) is conjugate to the centralizer $\mathbb{Z}_2$, to SU(2) itself or to U(1) (seen as a subgroup via the embedding $\alpha \mapsto \left(\begin{smallmatrix} \alpha & 0 \\ 0 & \bar{\alpha} \end{smallmatrix}\right)$). This corresponds to the trivial group or $\mathbb{Z}_2$, SU(2) and U(1) as holonomy groups, respectively, see Table 5.1. Correspondingly, the Howe subgroups of $\mathrm{SU}(2) \times \mathrm{U}(1)$ are conjugate to $\mathrm{SU}(2) \times \mathrm{U}(1)$, $\mathrm{U}(1) \times \mathrm{U}(1)$ or $\mathbb{Z}_2 \times \mathrm{U}(1)$.

Even in the case when (5.6.55) uniquely determines $H'_G$, there is still room for different subgroups $H'$ of $G \times L$ with $H = \mathrm{C}_{G\times L}(H')$. Indeed, recall that by Goursat's lemma subgroups of $G \times L$ are in bijection with quintuples $(G_1, G_2, L_1, L_2, \theta)$, where $G_2 \trianglelefteq G_1 \subseteq G$, $L_2 \trianglelefteq L_1 \subseteq L$ and $\theta\colon G_1/G_2 \to L_1/L_2$ is an isomorphism. Such a tuple yields the subgroup

$$H' = \{(g, l) \in G_1 \times L_1 : \ \theta(g\, G_2) = l\, L_2\} \subseteq G \times L. \tag{5.6.56}$$

Note that the projection of $H'$ to the $G$-factor coincides with $G_1$. Hence the knowledge of $H'_G$ just determines $G_1 = H'_G$. We now determine the possible choices for the other elements $(G_2, L_1, L_2, \theta)$ that generate the Howe subgroups of $\mathrm{SU}(2) \times \mathrm{U}(1)$, as summarized in Table 5.2.

- The Howe subgroup SU(2) is generated by $G_1 = \{e\}$ or $G_1 = \mathbb{Z}_2$. In the first case, we hence have $G_2 = \{e\}$ and so $L_1 = L_2$ with $\theta$ being trivial. In particular, $L_1 = L_2$ has to be either $\{e\}$, U(1) or the cyclic group $\mathbb{Z}_p$ for some $p \in \mathbb{N}$, since these are the only subgroups of U(1). For $G_1 = \mathbb{Z}_2$ there are two cases: First we may choose $G_2 = \mathbb{Z}_2$, which then again requires $L_1 = L_2$. Secondly, also $G_2 = \{e\}$ is possible. Then $L_1/L_2$ has to be isomorphic to $\mathbb{Z}_2$, which is only possible[1] if $L_1 = \mathbb{Z}_{2p}$ and $L_2 = \mathbb{Z}_p$ for some $p \in \mathbb{N}$.

- The Howe subgroup U(1) is generated by $G_1 = \mathrm{U}(1)$. There are two non-trivial choices for $G_2$. First, $G_2 = \mathrm{U}(1)$ leads again to $L_1 = L_2$. The second choice $G_2 = \mathbb{Z}_q$ for some $q \in \mathbb{N}$ enforces $L_1 = \mathrm{U}(1)$ and $L_2 = \mathbb{Z}_p$. Since the map $z \mapsto z^q$ induces an isomorphism of $\mathrm{U}(1)/\mathbb{Z}_q$ with U(1) and the only automorphisms of U(1) are of the form $z \mapsto z^k$ for some $k \in \mathbb{N}$, isomorphisms $\mathrm{U}(1)/\mathbb{Z}_q \to \mathrm{U}(1)/\mathbb{Z}_p$ are necessarily induced by maps of the form $z \mapsto z^{\frac{kq}{p}}$ for some $k \in \mathbb{N}$.

[1] Note that $\mathbb{Z}_p/\mathbb{Z}_q$ is isomorphic to $\mathbb{Z}_{p/q}$.

| | $H'$ | | | | |
|---|---|---|---|---|---|
| $H$ | $G_1$ | $G_2$ | $L_1$ | $L_2$ | $\theta$ |
| $\mathrm{SU}(2) \times \mathrm{U}(1)$ | $\{e\}$ | $\{e\}$ | $\mathrm{U}(1)$ | $\mathrm{U}(1)$ | trivial |
| | $\{e\}$ | $\{e\}$ | $\mathbb{Z}_p$ | $\mathbb{Z}_p$ | trivial |
| | $\{e\}$ | $\{e\}$ | $\{e\}$ | $\{e\}$ | trivial |
| | $\mathbb{Z}_2$ | $\mathbb{Z}_2$ | $\mathrm{U}(1)$ | $\mathrm{U}(1)$ | trivial |
| | $\mathbb{Z}_2$ | $\mathbb{Z}_2$ | $\mathbb{Z}_p$ | $\mathbb{Z}_p$ | trivial |
| | $\mathbb{Z}_2$ | $\mathbb{Z}_2$ | $\{e\}$ | $\{e\}$ | trivial |
| | $\mathbb{Z}_2$ | $\{e\}$ | $\mathbb{Z}_{2p}$ | $\mathbb{Z}_p$ | $\mathrm{id}_{\mathbb{Z}_2}$ |
| $\mathrm{U}(1) \times \mathrm{U}(1)$ | $\mathrm{U}(1)$ | $\mathrm{U}(1)$ | $\mathrm{U}(1)$ | $\mathrm{U}(1)$ | trivial |
| | $\mathrm{U}(1)$ | $\mathrm{U}(1)$ | $\mathbb{Z}_p$ | $\mathbb{Z}_p$ | trivial |
| | $\mathrm{U}(1)$ | $\mathrm{U}(1)$ | $\{e\}$ | $\{e\}$ | trivial |
| | $\mathrm{U}(1)$ | $\mathbb{Z}_q$ | $\mathrm{U}(1)$ | $\mathbb{Z}_p$ | $z \mapsto z^{\frac{kq}{p}}$ for some $k \in \mathbb{N}$ |
| | $\mathrm{U}(1)$ | $\{e\}$ | $\mathrm{U}(1)$ | $\{e\}$ | $z \mapsto z^k$ for some $k \in \mathbb{N}$ |
| $\mathbb{Z}_2 \times \mathrm{U}(1)$ | many choices | | | | |

Table 5.2: List of all Howe subgroups $H$ of $\mathrm{SU}(2) \times \mathrm{U}(1)$ up to conjugacy with the corresponding generator $H'$ satisfying $\mathrm{C}_{\mathrm{SU}(2)\times\mathrm{U}(1)}(H') = H$.

- The Howe subgroup $\mathbb{Z}_2$ is generated by $G_1 = \mathrm{SU}(2)$ giving rise to a plethora of possible choices for $G_2$, $L_1$ and $L_2$. Fortunately, this case will be of no further interest below, so we do not need to dive into details.

According to Proposition 3.4.1, we have accomplished the algebraic part of the classification of stabilizer subgroups $\mathcal{G}au_A$. Depending on the topologies of $M$ and $P$, some of the Howe subgroups $H$ may possibly not occur in the final classification of orbit types.

According to (5.6.40), in order to calculate the stabilizer $\mathcal{G}au_{A,\varphi}$, we need to determine which conjugates of $K$ intersect $K$ again. Writing $a \in \mathrm{SU}(2)$ as $a = \begin{pmatrix} \alpha & -\bar{\beta} \\ \beta & \bar{\alpha} \end{pmatrix}$ with $\alpha, \beta \in \mathbb{C}$ such that $|\alpha|^2 + |\beta|^2 = 1$, we find

$$a \begin{pmatrix} e^{\frac{i\vartheta}{2}} & 0 \\ 0 & e^{-\frac{i\vartheta}{2}} \end{pmatrix} a^{-1} = \begin{pmatrix} |\alpha|^2 e^{\frac{i\vartheta}{2}} + |\beta|^2 e^{-\frac{i\vartheta}{2}} & -\alpha\bar{\beta}\left(e^{-\frac{i\vartheta}{2}} - e^{\frac{i\vartheta}{2}}\right) \\ \bar{\alpha}\beta\left(e^{\frac{i\vartheta}{2}} - e^{-\frac{i\vartheta}{2}}\right) & |\alpha|^2 e^{-\frac{i\vartheta}{2}} + |\beta|^2 e^{\frac{i\vartheta}{2}} \end{pmatrix}. \tag{5.6.57}$$

Hence $\left(a \begin{pmatrix} e^{\frac{i\vartheta}{2}} & 0 \\ 0 & e^{-\frac{i\vartheta}{2}} \end{pmatrix} a^{-1}, e^{i\vartheta}\right)$ is an element of $K$ if and only if either $\beta = 0$ or

$e^{i\vartheta} = 1$. Thus,

$$\bigcap_{p \in P} \phi(p) K \phi(p)^{-1} = \begin{cases} K & \text{if } \beta(p) = 0 \text{ for all } p \in P_A, \\ \mathbb{Z}_2 & \text{otherwise,} \end{cases} \tag{5.6.58}$$

where $\mathbb{Z}_2$ is viewed as the subgroup

$$\mathbb{Z}_2 \simeq \left\{ \left( \begin{pmatrix} 1 & 0 \\ 0 & 1 \end{pmatrix}, 1 \right), \left( \begin{pmatrix} -1 & 0 \\ 0 & -1 \end{pmatrix}, 1 \right) \right\} \tag{5.6.59}$$

of $K$. Thus, we obtain the following.

Proposition 5.6.11 *The common stabilizer $\mathcal{G}au_{A,\varphi}$ of $(A, \varphi)$ is either $K$ if $\beta(p) = 0$ for all $p \in P_A$ and the stabilizer of $A$ is $\mathrm{SU}(2) \times \mathrm{U}(1)$ or $\mathrm{U}(1) \times \mathrm{U}(1)$; or otherwise it is $\mathbb{Z}_2$.* ◇

We now describe the structure of the orbit types in terms of the fields $(\hat{A}, \tau)$ after symmetry breaking, see (5.6.39). For this purpose, choose the following basis of $\mathfrak{g} = \mathfrak{su}(2) \times \mathfrak{u}(1)$:

$$\left\{ t_a = \frac{i}{2} \sigma_a, i \right\}, \tag{5.6.60}$$

where $\sigma_a$ are the Pauli matrices. In terms of these generators, the induced representation of $\mathfrak{g}$, which will also be denoted by $\rho$, is determined by

$$\rho_{t_a} = t_a, \qquad \rho_i = \frac{i}{2} \mathbb{1}. \tag{5.6.61}$$

Let $t_\pm = t_3 \pm i$. Then, the Lie algebra $\mathfrak{k}$ of $K$ is spanned by $t_+$ and the complement $\mathfrak{p}$ is spanned by $\{t_1, t_2, t_-\}$. With respect to the basis $\{t_a, i\}$, we expand the connection as

$$A = g W^a t_a + i g' B, \tag{5.6.62}$$

where we have introduced the coupling constants $g$ and $g'$. Thus, passing to the basis $\{t_1, t_2, t_\pm\}$ yields the decomposition $A = \hat{A} + \tau$, where

$$\hat{A} = \frac{1}{2}(g W^3 + g' B) t_+, \tag{5.6.63a}$$

$$\tau = g W^1 t_1 + g W^2 t_2 + \frac{1}{2}(g W^3 - g' B) t_-. \tag{5.6.63b}$$

As common in the physics literature, it is convenient to introduce the fields

$$W_\pm = \frac{1}{\sqrt{2}}(W^1 \mp i W^2), \tag{5.6.64a}$$

$$\begin{pmatrix} Z \\ A_\gamma \end{pmatrix} = \frac{1}{\sqrt{g^2 + g'^2}} \begin{pmatrix} g & -g' \\ g' & g \end{pmatrix} \begin{pmatrix} W^3 \\ B \end{pmatrix}. \tag{5.6.64b}$$

Moreover, we denote the elementary charge $e$ by $e = \frac{gg'}{\sqrt{g^2+g'^2}}$ and the Weinberg angle $\theta_W$ by $\tan\theta_W = \frac{g'}{g}$. Then,

$$\hat{A} = \left( eA_\gamma + \frac{g\cos\theta_W - g'\sin\theta_W}{2} Z \right) t_+, \tag{5.6.65a}$$

$$\tau = gW_+ t + gW_- \bar{t} + \frac{gg'}{2e} Z t_- \tag{5.6.65b}$$

with $t := \frac{1}{\sqrt{2}}(t_1 + \mathrm{i}t_2)$ and $\bar{t} := \frac{1}{\sqrt{2}}(t_1 - \mathrm{i}t_2)$. For later use, let us record the commutation relations of the new basis:

$$[t, \bar{t}\,] = \mathrm{i}t_3, \qquad [t, t_\pm] = -\mathrm{i}t, \tag{5.6.66a}$$

$$[t_+, t_-] = 0, \qquad [\bar{t}, t_\pm] = \mathrm{i}\bar{t}. \tag{5.6.66b}$$

By (5.6.41), we have $\mathcal{G}au_{A,\varphi} = \mathcal{G}au_{\hat{A},\tau}(\hat{P})$. Since $K$ is abelian, the stabilizer $\mathcal{G}au_{\hat{A}}(\hat{P})$ is isomorphic to $K$. Thus it remains to determine which gauge transformations leave $\tau$ invariant. For this purpose, let $k = \left( \begin{pmatrix} e^{\frac{i}{2}\vartheta} & 0 \\ 0 & e^{-\frac{i}{2}\vartheta} \end{pmatrix}, e^{i\vartheta} \right) \in K$. A straightforward calculation shows that $k$ acts on the basis $\{t_a, \mathrm{i}\}$ as follows:

$$\mathrm{Ad}_k\, t_1 = \cos\vartheta\, t_1 - \sin\vartheta\, t_2, \tag{5.6.67a}$$

$$\mathrm{Ad}_k\, t_2 = \sin\vartheta\, t_1 + \cos\vartheta\, t_2, \tag{5.6.67b}$$

$$\mathrm{Ad}_k\, t_3 = t_3, \tag{5.6.67c}$$

$$\mathrm{Ad}_k\, \mathrm{i} = \mathrm{i}. \tag{5.6.67d}$$

That is, $k$ acts as a rotation in the $(t_1, t_2)$-plane and acts trivially on $t_\pm$. Thus, by Proposition 5.6.7,

$$\mathcal{G}au_{A,\varphi}(P) = \mathcal{G}au_{\hat{A}}(\hat{P}) \cap \mathcal{G}au_\tau(\hat{P}) \simeq \begin{cases} K & \text{if } \tau \text{ is proportional to } t_-, \\ \mathbb{Z}_2 & \text{otherwise.} \end{cases} \tag{5.6.68}$$

In other words, $(\hat{A}, \tau)$ has non-trivial stabilizer if and only if $W^1 = 0 = W^2$, i.e., if the $Z$-boson is the only non-trivial intermediate vector boson. Hence, on the non-generic stratum only one intermediate boson survives.

Finally, we show that the decomposition of $\mathcal{Q}$ defines a stratification.

Proposition 5.6.12 *The decomposition of $\mathcal{Q}$ into orbit types satisfies the frontier condition.* ◊

*Proof.* As we have seen above, the $\mathcal{G}au$-action on $\mathcal{Q}$ has only the orbit types $(K)$ and $(\mathbb{Z}_2)$. The frontier condition thus requires that every pair $(A, \varphi) \in \mathcal{Q}$ with orbit type $(K)$ can be approximated by a sequence $(A_i, \varphi_i)$ with orbit type $(\mathbb{Z}_2)$. By Proposition 5.6.11, a pair $(A, \varphi)$ has a stabilizer conjugate to $K$ only when $A$ has a stabilizer conjugate to $\mathrm{SU}(2) \times \mathrm{U}(1)$ or to $\mathrm{U}(1) \times \mathrm{U}(1)$. In both cases, the approximation theorem [KR86, Theorem 4.3.5] shows that there is a converging sequence $A_i \longrightarrow A$ of connections $A_i$ with stabilizer conjugate to $\mathbb{Z}_2 \times \mathrm{U}(1)$. By Lemma 5.6.6, the pair $(A_i, \varphi)$ has stabilizer conjugate to $(\mathbb{Z}_2 \times \mathrm{U}(1)) \cap K \simeq \mathbb{Z}_2$ and converges to $(A, \varphi)$ by construction. □

### 5.6.5.2 *Orbit types of* $\mathrm{T}^*\mathcal{Q}$

Next, we determine the secondary stratification of the cotangent bundle. For this purpose, we endow $\mathfrak{g}$ with the $\mathrm{Ad}_G$-invariant scalar product given as the product of (minus) the Killing form on $\mathfrak{su}(2)$ and the usual scalar product on $\mathfrak{u}(1)$. The normalization is such that the generators $\{t_a, \mathrm{i}\}$ form an orthonormal basis. In the sequel, we will always use this scalar product to identify $\mathfrak{g}^*$ with $\mathfrak{g}$. Now, let $D \in \Omega^2(M, \mathrm{Ad}^* P)$. We decompose $D$ according to $\mathfrak{g} = \mathfrak{k} \oplus \mathfrak{p}$ into $D_{\restriction \hat{P}} = D_{\mathfrak{k}} + D_{\mathfrak{p}}$, with $D_{\mathfrak{k}} \in \Omega^2(M, \mathrm{Ad}^* \hat{P}) \simeq \Omega^2(M, \mathfrak{k}^*)$ and $D_{\mathfrak{p}} \in \Omega^2(M, \hat{P} \times_K \mathfrak{p}^*)$. Recall from the discussion above that $\mathcal{G}au_{\hat{A}}(\hat{P})$ is isomorphic to $K$, viewed as constant gauge transformations in $\hat{P}$. Then, similarly to the above reasoning, using (5.6.67) and the fact that $K$ is abelian, we find

$$\mathcal{G}au_D(\hat{P}) \cap \mathcal{G}au_{\hat{A}}(\hat{P}) = \mathcal{G}au_{D_{\mathfrak{p}}}(\hat{P}) \cap \mathcal{G}au_{\hat{A}}(\hat{P}) = \begin{cases} K & \text{if } D_{\mathfrak{p}} \text{ is proportional to } t_-, \\ \mathbb{Z}_2 & \text{otherwise.} \end{cases} \tag{5.6.69}$$

We now turn to the stabilizer of $\Pi \in \mathcal{F}^*$, which we view as $\Pi \in \Gamma^\infty(P \times_G \mathbb{C}^2)$ using the volume form. As we have seen above, a point $f \in \mathbb{C}^2$ has stabilizer $K$ if and only if it is of the form

$$f = \frac{r}{\sqrt{2}} \begin{pmatrix} \mathrm{e}^{\mathrm{i}\alpha} & 0 \\ 0 & \mathrm{e}^{-\mathrm{i}\alpha} \end{pmatrix} \cdot \begin{pmatrix} 0 \\ \nu \end{pmatrix} = \frac{r}{\sqrt{2}} \begin{pmatrix} 0 \\ \mathrm{e}^{-\mathrm{i}\alpha} \nu \end{pmatrix} \tag{5.6.70}$$

for some $r \in \mathbb{R}_{>0}$ and $\mathrm{e}^{\mathrm{i}\alpha} \in \mathrm{U}(1)$, that is, if its first component $f_1$ vanishes. Thus, we have

$$\mathcal{G}au_\Pi(\hat{P}) \cap \mathcal{G}au_{\hat{A}}(\hat{P}) = \begin{cases} K & \text{if } \Pi_1 = 0 \text{ on } \hat{P}, \\ \mathbb{Z}_2 & \text{otherwise.} \end{cases} \tag{5.6.71}$$

Hence, in summary, we find

Proposition 5.6.13 *The stabilizer of $(A, \varphi, D, \Pi) \in \mathrm{T}^*\mathcal{Q}$ is conjugate to $K$ if the following conditions are met on $\hat{P}$ (otherwise it is conjugate to $\mathbb{Z}_2$):*

*(i) $\tau$ is proportional to $t_-$, i.e. $W_\pm = 0$,*

*(ii) $D_{\mathfrak{p}}$ is proportional to $t_-$ and*

*(iii) the first component of $\Pi$ vanishes.* ◇

We note that $\mathrm{T}^*\mathcal{Q}$ thus has the same orbit types as $\mathcal{Q}$, cf. Proposition 5.6.11.

Remark 5.6.14 More generally, instead of (5.6.51) we could consider the representation

$$\rho^Y_{a,\vartheta}(z_1,z_2) = a \cdot \mathrm{e}^{\mathrm{i}Y\vartheta} \cdot (z_1,z_2), \qquad a \in \mathrm{SU}(2), \vartheta \in [0,4\pi), \tag{5.6.72}$$

which changes the weak hypercharge of the Higgs field from $\frac{1}{2}$ to $Y \in \mathbb{Q}$. In this case, the stabilizer group $K$ is replaced by

$$K^Y := \left\{ \left( \begin{pmatrix} \mathrm{e}^{\mathrm{i}Y\vartheta} & 0 \\ 0 & \mathrm{e}^{-\mathrm{i}Y\vartheta} \end{pmatrix}, \mathrm{e}^{\mathrm{i}\vartheta} \right) \; : \; \vartheta \in [0,4\pi) \right\}. \tag{5.6.73}$$

Moreover, the generic orbit type is no longer $\mathbb{Z}_2$ but the subgroup of $\mathbb{Z}_2 \times \mathrm{U}(1)$ generated by the elements

$$\left( \begin{pmatrix} \mathrm{e}^{\mathrm{i}\pi n} & 0 \\ 0 & \mathrm{e}^{-\mathrm{i}\pi n} \end{pmatrix}, \mathrm{e}^{\mathrm{i}\frac{\pi n}{Y}} \right), \qquad n \in \mathbb{Z}. \tag{5.6.74}$$

We see, in particular, that the orbit type stratification of the configuration space $\mathcal{Q}$ depends on the weak hypercharge of the Higgs field. ◇

Note that there might be topological obstructions related to the conditions in Proposition 5.6.13. Thus, the complete classification of gauge orbit types depends on the topology of $P$ and $M$.

Remark 5.6.15 Let us consider the special case $M = S^3$. By the standard principal fiber bundle classification theorem, all $G$-bundles over $S^3$ are trivial, because $\pi_2(G) = 0$. That is, $P$ is trivial in that case. The same applies to any subbundle of $P$ and hence to any holonomy-induced Howe subbundle. As a consequence, the classification problem for that base manifold boils down to the algebraic problem we just solved. ◇

### 5.6.5.3 *Momentum map and reduced phase space*

Let us determine the expression of the momentum map given by (5.6.7),

$$\mathcal{J} \colon \mathrm{T}^*\mathcal{Q} \to \mathfrak{gau}^*, \qquad (A, D, \varphi, \Pi) \mapsto \mathrm{d}_A D + \varphi \diamond \Pi \,\mathrm{vol}_g \,. \tag{5.6.75a}$$

First, we calculate the second summand. For this purpose, we consider the momentum map $\underline{J} \colon \mathbb{C}^2 \times \mathbb{C}^2 \to \mathfrak{g}^*$ for the lifted $G$-action on $\mathrm{T}^*\mathbb{C}^2$. The latter is determined by the equations

$$\langle \underline{J}(z,v), t_a \rangle = \langle v, \rho_{t_a} z \rangle = \mathrm{Re}(v^* t_a z) \,, \tag{5.6.76}$$

$$\langle \underline{J}(z,v), \mathrm{i} \rangle = \langle v, \rho_{\mathrm{i}} z \rangle = -\frac{1}{2} \mathrm{Im}(v^* z) \,, \tag{5.6.77}$$

and hence, in vector form, it is given by

$$\underline{J}(z, v) = -\frac{1}{2}\operatorname{Im}\begin{pmatrix} v_1^* z_2 + v_2^* z_1 \\ \mathrm{i}v_2^* z_1 - \mathrm{i}v_1^* z_2 \\ v_1^* z_1 - v_2^* z_2 \\ v_1^* z_1 + v_2^* z_2 \end{pmatrix}. \tag{5.6.78}$$

Thus, we have

$$\begin{aligned}(\varphi \diamond \Pi)_{\restriction \hat{P}} &= \underline{J}\left(\frac{\eta}{\sqrt{2}}(0, \nu), \Pi\right) = -\frac{\eta\nu}{2\sqrt{2}}\operatorname{Im}\begin{pmatrix} \Pi_1^* \\ -\mathrm{i}\Pi_1^* \\ -\Pi_2^* \\ \Pi_2^* \end{pmatrix} \\ &= -\frac{\eta\nu}{2}\left(\frac{\mathrm{i}\Pi_1}{2}t - \frac{\mathrm{i}\Pi_1^*}{2}\bar{t} + \frac{\operatorname{Im}\Pi_2}{\sqrt{2}}t_-\right).\end{aligned} \tag{5.6.79}$$

If $\Pi$ lies in the singular orbit type $(K)$, then its first component $\Pi_1$ vanishes and thus the current $\varphi \diamond \Pi$ is proportional to $t_-$ in this case.

Next, let us expand $D_{\restriction \hat{P}} \in \Omega^2(M, \hat{P} \times_K \mathfrak{g}^*)$ in the basis $\{t, \bar{t}, t_-, t_+\}$:

$$D_{\restriction \hat{P}} = \underbrace{\frac{D_-}{g}t + \frac{D_+}{g}\bar{t} + \left(\frac{e}{gg'}D_Z - \frac{g\cos\theta_{\mathrm{W}} - g'\sin\theta_{\mathrm{W}}}{2gg'}D_\gamma\right)t_-}_{D_{\mathfrak{p}}} + \underbrace{\frac{D_\gamma}{2e}t_+}_{D_{\mathfrak{k}}}, \tag{5.6.80}$$

where the normalization was chosen in such a way that the symplectic form stays in Darboux form in the new coordinates $(W_\pm, D_\pm, Z, D_Z, A_\gamma, D_\gamma)$ (with respect to the scalar product in which $\{t_a, i\}$ is an orthonormal basis). In this notation, using (5.6.65) and (5.6.66), we find

$$\begin{aligned}(\mathrm{d}_A D)_{\restriction \hat{P}} &= \mathrm{d}D_{\restriction \hat{P}} + [\hat{A}, D_{\restriction \hat{P}}] + [\tau, D_{\restriction \hat{P}}] \\ &= \frac{1}{g}\,\mathrm{d}D_- t + \frac{1}{g}\,\mathrm{d}D_+\bar{t} + \frac{e}{gg'}\,\mathrm{d}D_Z t_- \\ &\quad - \frac{g\cos\theta_{\mathrm{W}} - g'\sin\theta_{\mathrm{W}}}{2gg'}\,\mathrm{d}D_\gamma t_- \\ &\quad + \frac{1}{2e}\,\mathrm{d}D_\gamma t_+ + \mathrm{i}\left(\sin\theta_{\mathrm{W}}A_\gamma + \cos\theta_{\mathrm{W}}Z\right) \wedge (D_- t - D_+\bar{t}) \\ &\quad + \mathrm{i}W_+ \wedge \left(D_+ t_3 - (\sin\theta_{\mathrm{W}}D_\gamma + \cos\theta_{\mathrm{W}}D_Z)t\right) \\ &\quad - \mathrm{i}W_- \wedge \left(D_- t_3 - (\sin\theta_{\mathrm{W}}D_\gamma + \cos\theta_{\mathrm{W}}D_Z)\bar{t}\right).\end{aligned} \tag{5.6.81}$$

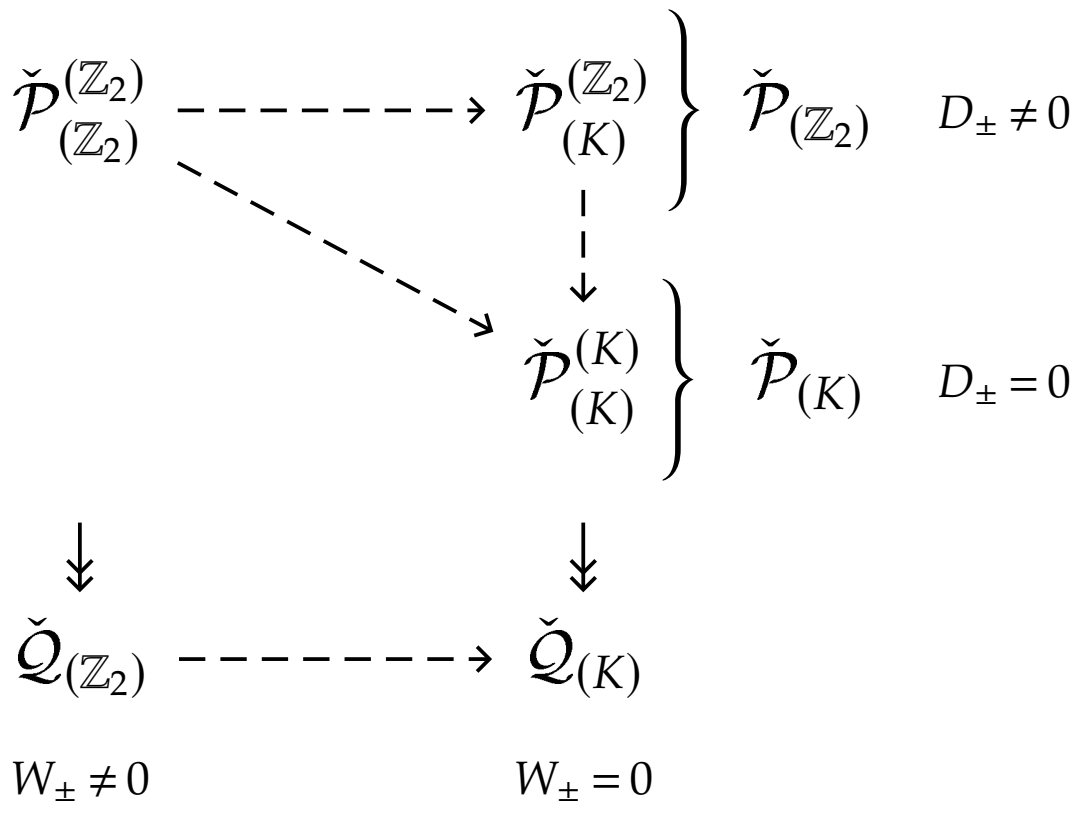

Figure 5.2: Schematic illustration of the secondary stratification of the reduced phase $\check{\mathcal{P}}$ and its relation to the orbit type stratification of the reduced configuration space $\check{\mathcal{Q}}$. Dotted arrows mean that the target lies in the closure of the source.

Hence, the Gauß constraint (5.6.4e) takes the following form:

$$\begin{aligned} &\mathrm{d}D_+ - \mathrm{i}g(\sin\theta_\mathrm{W} A_\gamma + \cos\theta_\mathrm{W} Z) \wedge D_+ \\ &\qquad = -\mathrm{i}g W_- \wedge (\sin\theta_\mathrm{W} D_\gamma + \cos\theta_\mathrm{W} D_Z) - \mathrm{i}\frac{\eta\nu g}{4}\Pi_1^* \,\mathrm{vol}_g, \end{aligned} \tag{5.6.82a}$$

$$\begin{aligned} &\mathrm{d}D_- + \mathrm{i}g(\sin\theta_\mathrm{W} A_\gamma + \cos\theta_\mathrm{W} Z) \wedge D_- \\ &\qquad = \mathrm{i}g W_+ \wedge (\sin\theta_\mathrm{W} D_\gamma + \cos\theta_\mathrm{W} D_Z) + \mathrm{i}\frac{\eta\nu g}{4}\Pi_1 \,\mathrm{vol}_g, \end{aligned} \tag{5.6.82b}$$

$$\mathrm{d}D_Z = \mathrm{i}g\cos\theta_\mathrm{W}(W_- \wedge D_- - W_+ \wedge D_+) + \frac{\eta\nu g g'}{2\sqrt{2}e}\operatorname{Im}\Pi_2 \,\mathrm{vol}_g, \tag{5.6.82c}$$

$$\mathrm{d}D_\gamma = \mathrm{i}e(W_- \wedge D_- - W_+ \wedge D_+). \tag{5.6.82d}$$

Thus, considered on the singular stratum where we have $W_\pm = 0 = D_\pm$ according to Proposition 5.6.13, the Gauß constraint is equivalent to the two equations

$$\mathrm{d}D_Z = \frac{\eta\nu g g'}{2\sqrt{2}e}\operatorname{Im}\Pi_2 \,\mathrm{vol}_g, \qquad \mathrm{d}D_\gamma = 0. \tag{5.6.83}$$

Since these equations are decoupled, the Gauß constraint cuts out a smooth subbundle of $(\mathrm{T}^*\mathcal{Q})^{(K)}_{(K)}$, whose fiber is parametrized by the fields $D_\gamma \in \Omega^2_{\mathrm{cl}}(M)$, $D_Z \in \Omega^2(M)$ and $\frac{\nu}{\sqrt{2}}\operatorname{Re}\Pi_2 \in C^\infty(M)$.

According to Theorem 5.4.8, the reduced phase space decomposes into

$$\check{\mathcal{P}} = \underbrace{\check{\mathcal{P}}^{(K)}_{(K)}}_{\check{\mathcal{P}}_{(K)}} \cup \underbrace{\check{\mathcal{P}}^{(\mathbb{Z}_2)}_{(K)} \cup \check{\mathcal{P}}^{(\mathbb{Z}_2)}_{(\mathbb{Z}_2)}}_{\check{\mathcal{P}}_{(\mathbb{Z}_2)}} \tag{5.6.84}$$

and the strata $\check{\mathcal{P}}_{(K)}$ and $\check{\mathcal{P}}_{(\mathbb{Z}_2)}$ are symplectic. Moreover, $\check{\mathcal{P}}^{(K)}_{(K)}$ is symplectomorphic to the cotangent bundle of $\check{\mathcal{Q}}_{(K)}$. As we have seen, this singular stratum is the (reduced) phase space of the theory of a photon and the $Z$-boson without any other intermediate bosons. In contrast, on the generic stratum $\check{\mathcal{P}}^{(\mathbb{Z}_2)}_{(\mathbb{Z}_2)} \simeq \mathrm{T}^*(\check{\mathcal{Q}}_{(\mathbb{Z}_2)})$ all intermediate vector bosons are present. This cotangent bundle is stitched together by the seam $\check{\mathcal{P}}^{(\mathbb{Z}_2)}_{(K)}$, where no $W$-bosons are present but their conjugate momenta are non-zero. The secondary stratification is schematically illustrated in Figure 5.2. In the next section, we study the structure of the reduced phase space in detail and show that it is similar to that of the harmonic oscillator discussed in Example 5.5.1.

5.6.5.4 *Description of the reduced phase space*

We now give an explicit parameterization of the reduced phase space. The main idea is to identify a part of the configuration space on which the group of gauge transformations acts transitively and thereby to reduce the symmetry to a subgroup. Next, we repeat this process until only a finite-dimensional symmetry remains. To simplify the discussion, we limit our attention to the case $M = S^3$. As noted in Remark 5.6.15, for $M = S^3$ all bundles occurring in the construction are trivial and we can hence represent all geometric objects on the bundle as objects living on $S^3$.

1. Reduction of the gauge group from $\mathrm{C}^\infty(M, G)$ to $\mathrm{C}^\infty(M, K)$: Recall from Proposition 5.6.5 that the Higgs mechanism yields a parametrization of $(A, \varphi) \in \mathcal{Q}$ in terms of the variables $(\phi, \hat{A}, \tau, \eta)$, which are viewed as elements of a bundle over $\Gamma^\infty(P \times_G G/K)$. In the present case of a trivial bundle $P$, we can strengthen this result. For this purpose, recall that $G/K$ is diffeomorphic to $S^3$ via the $G$-action on $\mathbb{C}^2$ defined in (5.6.51). Hence, in particular, the $K$-bundle $G \to G/K$ is trivial. Therefore, every smooth map $M \to G/K$ lifts to a smooth map $M \to G$ and so the action of $\mathcal{G}au(P) = \mathrm{C}^\infty(M, G)$ on $\Gamma^\infty(P \times_G G/K) \simeq \mathrm{C}^\infty(M, G/K)$ is transitive. Moreover, the stabilizer of the constant function taking values in the identity coset is the subgroup $\mathcal{G}au(\hat{P}) = \mathrm{C}^\infty(M, K)$ (this is in accordance with Lemma 5.6.6). Therefore, $\mathrm{C}^\infty(M, G/K)$ is diffeomorphic to the homogeneous space $\mathcal{G}au(P)/\mathcal{G}au(\hat{P})$. In other words, the decomposition (5.6.54) takes the form

$$\varphi = \frac{\eta}{\sqrt{2}} \lambda \cdot \begin{pmatrix} 0 \\ \nu \end{pmatrix}, \tag{5.6.85}$$

where $\eta \in \mathrm{C}^\infty(M, \mathbb{R}_{>0})$ and $\lambda \in \mathcal{G}au(P)$ is unique up to an element of $\mathcal{G}au(\hat{P})$. By implementing the unitary gauge, i.e., by gauging away $\lambda$, we obtain the following refinement of Proposition 5.6.5. Recall the decomposition $\mathfrak{g} = \mathfrak{p} \oplus \mathfrak{k}$ with $\mathfrak{p}$ spanned by $\{t, \bar{t}, t_-\}$ and $\mathfrak{k}$ spanned by $t_+$.

Proposition 5.6.16 *The assignment*[1]

$$(A, \varphi) \mapsto ([\lambda, \hat{A}, \tau], \eta), \tag{5.6.86}$$

*where* $\lambda$ *and* $\eta$ *are determined by* (5.6.85) *and* $(\lambda^{-1} \cdot A) = \hat{A} + \tau$, *defines a* $\mathcal{G}au(P)$*-equivariant diffeomorphism*

$$\mathcal{Q} \to \mathcal{G}au(P) \times_{\mathcal{G}au(\hat{P})} \big(\Omega^1(M, \mathfrak{k}) \times \Omega^1(M, \mathfrak{p})\big) \times C^\infty(M, \mathbb{R}_{>0}). \tag{5.6.87}$$

*Here, on the right hand side,* $\mathcal{G}au(P)$ *acts by left translation on the first factor and* $\mathcal{G}au(\hat{P})$ *acts by gauge transformations on the space of connections* $\Omega^1(M, \mathfrak{k})$ *and via the adjoint action on* $\Omega^1(M, \mathfrak{p})$. *In particular, the gauge orbit space* $\check{\mathcal{Q}}$ *is isomorphic to* $\big(\Omega^1(M, \mathfrak{k}) \times \Omega^1(M, \mathfrak{p})\big) / \mathcal{G}au(\hat{P}) \times C^\infty(M, \mathbb{R}_{>0})$ *in the sense of stratified spaces.* ◊

Recall from the discussion in Section 5.2, that there is a natural description of the cotangent bundle of an associated bundle such that the momentum map is brought into a convenient normal form, cf. Theorem 5.2.9 (in Section 5.2 the focus lies on certain associated bundles where the slice is the typical fiber — however, the discussion there does not really use the properties of a slice). Using the same strategy, we identify

$$\mathcal{G}au(P) \times_{\mathcal{G}au(\hat{P})} \Big(C^\infty(M, \mathfrak{p}) \times \Omega^1(M, \mathfrak{k})^2 \times \Omega^1(M, \mathfrak{p})^2\Big) \times C^\infty(M, \mathbb{R}_{>0}) \times C^\infty(M, \mathbb{R}) \tag{5.6.88}$$

with $\mathrm{T}\mathcal{Q}$ via the map

$$\big([\lambda, (\varsigma, \hat{A}, \delta\hat{A}, \tau, \delta\tau)], \eta, \delta\eta\big) \mapsto \begin{pmatrix} A = \lambda \cdot (\hat{A} + \tau) \\ \delta A = \mathrm{Ad}_\lambda(\delta\hat{A} + \delta\tau) - \mathrm{d}_A \varsigma \\ \varphi = \frac{\eta}{\sqrt{2}} \lambda \cdot \begin{pmatrix} 0 \\ \nu \end{pmatrix} \\ \delta\varphi = \frac{\delta\eta}{\eta}\varphi + \varsigma \,.\, \varphi \end{pmatrix}. \tag{5.6.89}$$

Here, we have denoted $\Omega^1(M, \cdot) \times \Omega^1(M, \cdot) \equiv \Omega^1(M, \cdot)^2$. A straightforward calculation shows that the dual map

$$\begin{aligned} \mathrm{T}^*\mathcal{Q} \to \mathcal{G}au(P) \times_{\mathcal{G}au(\hat{P})} &\Big(\Omega^3(M, \mathfrak{p}^*) \times \Omega^1(M, \mathfrak{k}) \\ &\times \Omega^2(M, \mathfrak{k}^*) \times \Omega^1(M, \mathfrak{p}) \times \Omega^2(M, \mathfrak{p}^*)\Big) \\ &\times C^\infty(M, \mathbb{R}_{>0}) \times \Omega^3(M, \mathbb{R}), \\ (A, D, \varphi, \Pi) \mapsto &\big([\lambda, \nu, \hat{A}, \hat{D}, \tau, D_\tau], \eta, \Pi_\eta\big) \end{aligned} \tag{5.6.90}$$

[1] With a slight abuse of notation, we continue to use the notation $\hat{A}$ and $\tau$. Note, however, that these fields differ from the ones introduced in (5.6.62) by a gauge transformation.

is given by

$$\hat{D} + D_\tau = \operatorname{Ad}^*_{\lambda^{-1}} D, \tag{5.6.91}$$

$$\Pi_\eta = \frac{1}{\eta} \operatorname{Re}(\Pi^* \varphi), \tag{5.6.92}$$

$$\nu = \operatorname{pr}_{\mathfrak{p}^*}(\mathrm{d}_A D + \varphi \diamond \Pi). \tag{5.6.93}$$

Moreover, let us parametrize $(\hat{A}, \tau)$ in terms of the fields $(W_\pm, Z, A_\gamma)$ as in (5.6.65) and $(\hat{D}, D_\tau)$ in terms of the fields $(D_\pm, D_Z, D_\gamma)$ as in (5.6.80). Then, by (5.6.82), the $\mathfrak{k}^*$-projection of the Gauß constraint is given by

$$0 = \operatorname{pr}_{\mathfrak{k}^*}(\mathrm{d}_A D + \varphi \diamond \Pi) = \frac{1}{2e} \mathrm{d} D_\gamma + \operatorname{Im}\big(D_- \wedge W_-\big). \tag{5.6.94}$$

Denote $\mathcal{Q}_{\text{red}} = \Omega^1(M, \mathfrak{k}) \times \Omega^1(M, \mathbb{C})$. Its elements are $(A_\gamma, W_-)$. Moreover, elements of $\mathrm{T}^*_{(A_\gamma, W_-)} \mathcal{Q}_{\text{red}} = \Omega^2(M, \mathfrak{k}^*) \times \Omega^2(M, \mathbb{C})$ are denoted by $(D_\gamma, D_-)$. The right hand side of (5.6.94) is the momentum map for the induced $\mathcal{G}au(\hat{P})$-action[1] on $\mathrm{T}^* \mathcal{Q}_{\text{red}}$, because the momentum map for the lift of the $K$-action (5.6.67) to $\mathrm{T}^*\mathbb{C}$ is given by $\underline{J}_K(z, w) = -\operatorname{Im}(w^* z)$. The upshot of this first symmetry reduction is the following description of the reduced phase space.

Proposition 5.6.17 *The diffeomorphism* (5.6.87) *induces an isomorphism*

$$\check{\mathcal{P}} = \mathrm{T}^* \mathcal{Q} /\!\!/ \mathcal{G}au(P) \simeq \big(\mathrm{T}^* \mathcal{Q}_{red} /\!\!/ \mathcal{G}au(\hat{P})\big) \times \mathrm{T}^*\big(\Omega^1(M, \mathbb{R}\, t_-) \times \mathrm{C}^\infty(M, \mathbb{R}_{>0})\big) \tag{5.6.95}$$

*of symplectic stratified spaces.* ◇

*Proof.* First note that the Gauß constraint (5.6.4e) is equivalent to $\nu = 0$ and the condition (5.6.94), which is the momentum map constraint for the $\mathcal{G}au(\hat{P})$-action as we have discussed above. Now the assertion follows from the decomposition $\mathfrak{p} = \mathbb{C} \oplus \mathbb{R}$ with $\mathbb{C}$ spanned by $\{t, \bar{t}\}$ and $\mathbb{R}$ spanned by $t_-$ and from the fact that $K$ acts trivially on $t_-$, cf. (5.6.67). □

Hence, the singular structure of the phase space is completely encoded in the symplectic reduction of $\mathrm{T}^* \mathcal{Q}_{\text{red}}$ by the $\mathcal{G}au(\hat{P})$-action.

2. Reduction from $\mathcal{G}au(\hat{P})$ to $K$: Since $S^3$ is simply connected, the Hodge decomposition theorem yields

$$\Omega^1(M, \mathfrak{k}) = \mathrm{d}\, \mathrm{C}^\infty(M, \mathfrak{k}) \oplus \mathrm{d}^* \Omega^2(M, \mathfrak{k}). \tag{5.6.96}$$

[1] Note that $\lambda \in \mathcal{G}au(\hat{P})$ acts on $A_\gamma \in \Omega^1(M, \mathfrak{k})$ by $\lambda \cdot A_\gamma = A_\gamma - \frac{1}{e} \lambda^{-1} \mathrm{d}\lambda$.

Note that every map $f\colon M \longrightarrow \mathfrak{k}$ induces a map $\hat{f} = \exp_K \circ f\colon M \longrightarrow K$ such that $\hat{f}^{-1}\,\mathrm{d}\hat{f} = \mathrm{d}f$. Accordingly, every $K$-connection $A_\gamma$ can be written as

$$A_\gamma = \psi^{-1}\,\mathrm{d}\psi + \beta, \tag{5.6.97}$$

where $\psi \in \mathrm{C}^\infty(M, K)$ and $\beta \in \mathrm{d}^*\Omega^2(M, \mathfrak{k})$ is uniquely determined by the curvature $F_{A_\gamma}$ of $A_\gamma$. The action of $\mathcal{G}au(\hat{P})$ on $\mathrm{d}\Omega^0(M, \mathfrak{k})$, viewed as a subspace of the space of $K$-connections, is transitive with kernel consisting of the constant functions. We identify this kernel with $K$. Note that $K$ acts trivially on $\mathrm{d}^*\Omega^2(M, \mathfrak{k})$ and by rotation on $\Omega^1(M, \mathbb{C})$, cf. (5.6.67).

Lemma 5.6.18 *The map*

$$\begin{aligned} \mathcal{Q}_{red} &\longrightarrow \big(\mathcal{G}au(\hat{P}) \times_K \Omega^1(M, \mathbb{C})\big) \times \mathrm{d}\Omega^1(M, \mathfrak{k}) \\ (A_\gamma, W_-) &\longmapsto ([\psi, v], F_{A_\gamma}), \end{aligned} \tag{5.6.98}$$

*with $v = \psi^{-1} W_-$ is a diffeomorphism.* ◇

*Proof.* The inverse map $([\psi, v], \beta) \longmapsto (A_\gamma, W_-)$ is given by

$$A_\gamma = \psi^{-1}\,\mathrm{d}\psi + \mathrm{d}^* \vartriangle^{-1}\beta, \qquad W_- = \psi \,.\, v, \tag{5.6.99}$$

with $\psi \in \mathcal{G}au(\hat{P})$, $\beta \in \mathrm{d}\Omega^1(M, \mathfrak{k})$ and $v \in \Omega^1(M, \mathbb{C})$, because $\mathrm{d}\,\mathrm{d}^* \vartriangle^{-1}$ is the identity operator on $\mathrm{d}\Omega^1(M, \mathfrak{k})$. □

To get a convenient description of the cotangent bundle $\mathrm{T}^*\mathcal{Q}_{\mathrm{red}}$, we follow the same strategy as above. Let $\mathrm{C}^\infty(M, \mathfrak{k})_0$ denote the functions whose average over $M$ vanishes. The surjective linear map

$$\mathrm{C}^\infty(M, \mathfrak{k}) \longrightarrow \mathrm{C}^\infty(M, \mathfrak{k})_0, \qquad \psi \longmapsto \psi - \frac{1}{\mathrm{vol}_M} \int_M \psi\, \mathrm{vol}_g \tag{5.6.100}$$

has kernel $\mathfrak{k}$ and thus yields the identification $\mathrm{C}^\infty(M, \mathfrak{k})/\mathfrak{k} \simeq \mathrm{C}^\infty(M, \mathfrak{k})_0$. Dually, the integral map $\int_M\colon \Omega^3(M, \mathfrak{k}^*) \longrightarrow \mathfrak{k}^*$ has as kernel the space $\Omega^3(M, \mathfrak{k}^*)_0$, which is the dual space to $\mathrm{C}^\infty(M, \mathfrak{k})_0$. Linearizing the reconstruction equations (5.6.99) yields

$$\begin{aligned} \delta A_\gamma &= \mathrm{d}\,\delta\psi + \mathrm{d}^* \vartriangle^{-1} \delta\beta, \\ \delta W_- &= \delta\psi \,.\, (\psi \cdot v) + \psi \,.\, \delta v, \end{aligned} \tag{5.6.101}$$

where $\delta\psi \in \mathrm{C}^\infty(M, \mathfrak{k})_0$, $\delta\beta \in \mathrm{d}\Omega^1(M, \mathfrak{k})$ and $\delta v \in \Omega^1(M, \mathbb{C})$. By dualizing these equations, we get a diffeomorphism

$$\begin{aligned}\mathrm{T}^*\mathcal{Q}_{\mathrm{red}} &\to \mathcal{G}au(\hat{P}) \times_K \big(\Omega^3(M,\mathfrak{k}^*)_0 \times \Omega^1(M,\mathbb{C})^2\big) \\ &\qquad \times \mathrm{d}\Omega^1(M,\mathfrak{k}) \times \mathrm{d}^*\Omega^2(M,\mathfrak{k}^*) \\ (A_\gamma, D_\gamma, W_-, D_-) &\mapsto \big([\psi, \Pi_0, v, D_v], F_{A_\gamma}, D_{F_{A_\gamma}}\big),\end{aligned} \tag{5.6.102}$$

with

$$\begin{aligned}\Pi_0 &= \mathrm{pr}_{\Omega^3(M,\mathfrak{k}^*)_0}\left(\frac{1}{2e}\,\mathrm{d}D_\gamma + \mathrm{Im}(D_- \wedge W_-)\right), \\ D_v &= \psi^{-1} \,.\, D_-, \\ D_{F_{A_\gamma}} &= \mathrm{d}^* \vartriangle^{-1} D_\gamma.\end{aligned} \tag{5.6.103}$$

Moreover, the reduced Gauß constraint (5.6.94) is equivalent to $\Pi_0 = 0$ and

$$0 = \int_M \left(\frac{1}{2e}\,\mathrm{d}D_\gamma + \mathrm{Im}\big(D_- \wedge W_-\big)\right) = \int_M \mathrm{Im}\big(D_- \wedge W_-\big). \tag{5.6.104}$$

The right hand side of this identity is the momentum map for the lifted $K$-action on $\mathrm{T}^*\big(\Omega^1(M,\mathbb{C})\big)$. Thus, the second symmetry reduction yields the following.

Proposition 5.6.19 *The diffeomorphism* (5.6.98) *induces an isomorphism*

$$\mathrm{T}^*\mathcal{Q}_{red} /\!\!/ \mathcal{G}au(\hat{P}) \simeq \Big(\mathrm{T}^*\big(\Omega^1(M,\mathbb{C})\big) /\!\!/ K\Big) \times \mathrm{d}\Omega^1(M,\mathfrak{k}) \times \mathrm{d}^*\Omega^2(M,\mathfrak{k}^*) \tag{5.6.105}$$

*of symplectic stratified spaces.* ◇

Hence, in combination with the first reduction, see Proposition 5.6.17, we find that the singular structure of the reduced phase space $\check{\mathcal{P}}$ is completely determined by the singular cotangent bundle reduction of $\mathrm{T}^*\big(\Omega^1(M,\mathbb{C})\big)$ with respect to the action of the *finite-dimensional* (compact) Lie group $K$. Note that $\mathrm{T}^*\big(\Omega^1(M,\mathbb{C})\big)$ is pointwise the phase space of three (coupled) harmonic oscillators and the $K$ action corresponds to the diagonal U(1)-symmetry. This shows that the singularity structure of the reduced phase space is essentially determined by a finite-dimensional Lie group action. The reduced phase space $\mathrm{T}^*\big(\Omega^1(M,\mathbb{C})\big) /\!\!/ K$ may thus be studied using classical geometric invariant theory for the finite-dimensional reference system. This will be done elsewhere.

#### 5.6.5.5 *Dynamics on the singular stratum*

It is a challenge for further research to study the dynamics of this model on its full stratified phase space. As a first step, we analyze the dynamics on the singular stratum. First, a word of warning concerning conventions is in order. Since we followed the usual physics convention and introduced the coupling constants when writing the connection $A$ in terms of $W$ and $B$ fields, see (5.6.62), we need to use a different scalar product on $\mathfrak{g}$ in the Hamiltonian to counterbalance the coupling constants[1]. Let $\kappa$ be the scalar product on $\mathfrak{g}$ in which $\{t_a, i\}$ is orthogonal with norm $\kappa(t_a, t_a) = \frac{1}{g^2}$ and $\kappa(\mathrm{i}, \mathrm{i}) = \frac{1}{g'^2}$ (and hence $\kappa(t_+, t_\pm) = \frac{g'^2 \pm g^2}{g^2 g'^2} = \kappa(t_-, t_\mp)$). Moreover, let $\kappa^{-1}$ be the inverse of $\kappa$, i.e., we have $\kappa^{-1}(t_a, t_a) = g^2$ and $\kappa^{-1}(\mathrm{i}, \mathrm{i}) = g'^2$. In terms of these scalar products, the Hamiltonian $\mathcal{H}$ defined by (5.6.5) reads as follows:

$$\mathcal{H}(A, D, \varphi, \Pi) = \int_M \frac{\ell}{2} \left( \|D\|^2_{\kappa^{-1}} + \|F_A\|^2_\kappa + \|\Pi\|^2_{\mathbb{C}} + \|\mathrm{d}_A \varphi\|^2_{\mathbb{C}} + 2\, V(\varphi) \right) \mathrm{vol}_g \, . \tag{5.6.106}$$

Proposition 5.6.20 *On the singular stratum* $(\mathrm{T}^* \mathcal{Q})^{(K)}_{(K)}$, *the Hamiltonian* (5.6.106) *has the form*

$$\begin{aligned} &\mathcal{H}(A_\gamma, Z, \eta, D_\gamma, D_Z, \Pi_2) \\ &\quad = \int_M \frac{\ell}{2} \Big( \|D_\gamma\|^2 + \|D_Z\|^2 + \|\mathrm{d}A_\gamma\|^2 + \|\mathrm{d}Z\|^2 + \|\Pi_2\|^2_{\mathbb{C}} \\ &\qquad\qquad + \frac{\nu^2}{2} \|\mathrm{d}\eta\|^2 + \frac{\eta^2 \nu^2 (g^2 + g'^2)}{8} \|Z\|^2 + \lambda \nu^2 (\eta^2 - 1)^2 \Big) \mathrm{vol}_g \, . \end{aligned} \tag{5.6.107}$$

◊

By examining the Hamiltonian (5.6.107), we can read off the particle content over the singular stratum. We obtain a non-interacting system consisting of electrodynamics described by the photon $A_\gamma$, the theory of a massive vector boson described by the $Z$-boson with mass $m_Z^2 = \frac{\eta^2 \nu^2 (g^2 + g'^2)}{4}$ and the theory of a self-interacting real scalar field described by the Higgs boson $\eta$ with mass $m_\eta^2 = -4\lambda\nu^2$.

---

[1] The reader might find it instructive to compare the situation at hand to that of classical mechanics. There, the kinetic part of the Lagrangian is usually written in the form $L = \frac{m}{2} v^2$. But, the mass can be absorbed in the metric $g$ on the configuration space by setting $\|v\|^2_g := m v^2$. Then, the Hamiltonian is given by $H = \frac{1}{2} \|p\|^2_{g^{-1}}$, where the norm is taken with respect to the inverse (or dual) metric $g^{-1}$. In our field theoretic setting, the coupling constants play the role of inverse masses.

*Proof.* Over the singular stratum, we have $D_\pm = 0$ and thus

$$\begin{aligned} \|D\|^2_{\kappa^{-1}} &= \left\| D_Z - \frac{g^2 - g'^2}{2gg'} D_\gamma \right\|^2 + \frac{(g^2 + g'^2)^2}{4g^2 g'^2} \|D_\gamma\|^2 \\ &\quad + \frac{g^2 - g'^2}{gg'} \left\langle D_Z - \frac{g^2 - g'^2}{2gg'} D_\gamma, D_\gamma \right\rangle \\ &= \|D_Z\|^2 + \|D_\gamma\|^2. \end{aligned} \tag{5.6.108}$$

The curvature of $A$ reads in terms of the fields after symmetry breaking as follows:

$$\begin{aligned} (F_A)_{\restriction \hat{P}} &= F_{\hat{A}} + \mathrm{d}_{\hat{A}} \tau + \frac{1}{2} [\tau \wedge \tau] \\ &= F_{\hat{A}} + \mathrm{d}\tau + \mathrm{i} g^2 \left( \sin\theta_{\mathrm{W}} A_\gamma + \cos\theta_{\mathrm{W}} Z \right) \wedge (W_+ t - W_- \bar{t}) \\ &\quad + \mathrm{i} g^2 \, W_+ \wedge W_- \, t_3. \end{aligned} \tag{5.6.109}$$

Hence, on the singular stratum we simply have

$$(F_A)_{\restriction \hat{P}} = \left( e \, \mathrm{d}A_\gamma + \frac{g \cos\theta_{\mathrm{W}} - g' \sin\theta_{\mathrm{W}}}{2} \mathrm{d}Z \right) t_+ + \frac{gg'}{2e} \mathrm{d}Z \, t_-. \tag{5.6.110}$$

The norm of $F_A$ with respect to $\kappa$, is thus, on the singular stratum, given by

$$\|F_A\|^2_\kappa = \|\mathrm{d}A_\gamma\|^2 + \|\mathrm{d}Z\|^2. \tag{5.6.111}$$

According to (5.6.61), we get

$$\mathrm{d}_A \varphi = \left( \mathrm{d} + g W^a t_a + \frac{\mathrm{i} g'}{2} B \mathbb{1} \right) \varphi. \tag{5.6.112}$$

Using the representation (5.6.54) and the definition (5.6.64) of $W_\pm$ and $Z$, we find in terms of the fields after symmetry breaking

$$\begin{aligned} (\mathrm{d}_A \varphi)_{\restriction \hat{P}} &= \frac{\mathrm{d}\eta}{\sqrt{2}} \begin{pmatrix} 0 \\ \nu \end{pmatrix} + \frac{g\eta}{\sqrt{2}} W^a t_a \begin{pmatrix} 0 \\ \nu \end{pmatrix} + \frac{\mathrm{i} g' \eta}{2\sqrt{2}} B \begin{pmatrix} 0 \\ \nu \end{pmatrix} \\ &= \frac{\mathrm{i} g \eta \nu}{2} W_+ \begin{pmatrix} 1 \\ 0 \end{pmatrix} + \left( \frac{\nu}{\sqrt{2}} \mathrm{d}\eta - \frac{\mathrm{i} \eta \nu \sqrt{g^2 + g'^2}}{2\sqrt{2}} Z \right) \begin{pmatrix} 0 \\ 1 \end{pmatrix}. \end{aligned} \tag{5.6.113}$$

Thus, on the singular stratum,

$$\|\mathrm{d}_A \varphi\|^2_{\mathbb{C}} = \frac{\nu^2}{2} \|\mathrm{d}\eta\|^2 + \frac{\eta^2 \nu^2 (g^2 + g'^2)}{8} \|Z\|^2. \tag{5.6.114}$$

Plugging these identities into (5.6.106) yields (5.6.107). □

# 6 Outlook

In this thesis, we have established a new general framework to study moduli spaces and singular symplectic quotients in infinite dimensions. The techniques developed in the previous chapters open many exciting paths of further investigation. We list some relevant open problems and fundamental issues:

- It would be very interesting to extend the discussion of the normal form of a smooth map in Section 2.2 to higher orders. One would expect that the knowledge of higher order terms of the Taylor expansion of a smooth map $f$ yields further control over the behavior of its singular part $f_{\text{sing}}$. This way, one would gain deeper insight into the singular structure of the level sets of $f$. As an application, the singularities of the set of volume-preserving embeddings seem particularly tractable since they are determined by a map between finite-dimensional spaces. In the context of hydrodynamics, the set of volume-preserving embeddings serves as the configuration space of a free boundary fluid flow and it would be particularly important to see which physical phenomena are linked to the singularities of the configuration space.

- The methods developed in Chapter 4 yield a suitable normal form for the kinematic part of a Hamiltonian system with symmetries. It would be very interesting to continue this research and investigate Hamiltonian dynamics in that context. In particular, our theory should be suitable to study equivariant bifurcation phenomena. For this, it might turn out especially helpful that our proof of the MGS Normal Form Theorem is constructive and does not rely on the relative Darboux theorem, cf. Remark 4.2.26.

- Our results concerning infinite-dimensional singular symplectic reduction clarify the structure of the reduced phase space on the classical level. The next step is to pass to the quantum theory by developing a theory of geometric quantization for infinite-dimensional stratified symplectic spaces. First elements of such a theory are presented in [Die+], where prequantum bundles over symplectic section spaces are constructed.

- In the context of singular cotangent bundle reduction, we have seen that each singular stratum of the reduced phase space further decomposes into seams. Our discussion of the harmonic oscillator in Section 5.5 shows that the singular seams have a regularizing effect on the dynamics on the principal seam. It would be interesting to continue the study of the dynamics on the reduced phase space and, in particular, to clarify the meaning of the seams for the interplay between reduction and regularization.

- In Section 5.6.5, we have found an explicit parameterization of the reduced phase space of the Glashow–Weinberg–Salam model. Moreover, we have shown that the singular structure of the reduced phase space is encoded in the action of the finite-dimensional Lie group U(1). To gain a deeper understanding of how these singularities influence the properties of the corresponding quantum theory, it is of interest to extend the discussion to Cauchy surfaces of arbitrary topological type and to further analyze the structure of the singular reduced phase space. The latter can be achieved, for example, by using methods from geometric invariant theory.

# Calculus on Infinite-Dimensional Manifolds

# A

Our references for terminology and notation in the framework of infinite-dimensional differential geometry are [Ham82] for the tame Fréchet category and [Nee06] for the general locally convex setting.

By a manifold we understand a possibly infinite-dimensional smooth manifold $M$ without boundary modeled on locally convex vector spaces. More precisely, different connected components of $M$ are allowed to be modeled on non-isomorphic vector spaces in the same sense as e.g. in [Lan99, Section II.1; MRA02, Definition 3.1.1]. A subset $S$ of $M$ is called a submanifold if, at each point $s \in S$, there exist a chart $\kappa\colon M \supseteq U \to X$ and a closed subspace $Y \subseteq X$ such that $\kappa(U \cap S) = \kappa(U) \cap Y$. The submanifold $S$ is said to be split if, additionally, each subspace $Y$ is topologically complemented in $X$. A Lie group is a group $G$ endowed with a smooth manifold structure such that multiplication and inversion are smooth maps. The Lie algebra associated to a Lie group is written as a lowercase Fraktur letter corresponding to the uppercase Latin letter denoting the group, i.e. $\mathfrak{g}$ is the Lie algebra of the Lie group $G$. A Lie subgroup $H \subseteq G$ is called *principal* if the natural fibration $G \to G/H$ is a locally trivial principal bundle, see [GN, Section 7.1.4].

The derivative of a smooth map $f\colon M \to N$ at $m \in M$ is denoted by $\mathrm{T}_m f\colon \mathrm{T}_m M \to \mathrm{T}_{f(m)} N$. When the target or domain is a Lie group $G$, then it is often convenient to exploit the left trivialization of the tangent bundle of $G$ and introduce the following variations of the derivative of a map. The *left derivative* of a smooth map $\psi\colon G \to N$ at $g \in G$ is the map $\mathrm{T}^{\mathrm{L}}_g \psi\colon \mathfrak{g} \to \mathrm{T}_{\psi(g)} N$ defined by

$$\mathrm{T}^{\mathrm{L}}_g \psi(\xi) = \mathrm{T}_g \psi(g \,.\, \xi), \tag{A.0.1}$$

where $g \,.\, \xi$ denotes the left translation by $g \in G$ of $\xi \in \mathfrak{g}$ to an element of $\mathrm{T}_g G$. Given a smooth map $\phi\colon M \to G$, the *(left) logarithmic derivative* $\delta_m \phi\colon \mathrm{T}_m M \to \mathfrak{g}$ at $m \in M$ is defined by

$$\delta_m \phi(v) = \phi(m)^{-1} \,.\, \mathrm{T}_m(v) \tag{A.0.2}$$

for $v \in \mathrm{T}_m M$. If $f\colon P \to M$ is a smooth map, then the chain rule clearly yields

$$\delta_p(\phi \circ f) = \delta_{f(p)} \phi \circ \mathrm{T}_p f. \tag{A.0.3}$$

Moreover, a straightforward calculation shows that the left derivative and the logarithmic derivative are inverse to each other in the sense that, for every

(local) diffeomorphism $\psi\colon G \to N$, we have

$$\delta_{\psi(g)}\psi^{-1} \circ \mathrm{T}^{\mathrm{L}}_{g}\psi = \mathrm{id}_{\mathfrak{g}}. \tag{A.0.4}$$

## A.1 Inverse Function Theorems

In this section, we give a brief overview of different generalizations of the classical Inverse Function Theorem to the infinite-dimensional setting. The primary focus is on Glöckner's Inverse Function Theorem for maps between Banach spaces with parameters in a locally convex space and on the Nash–Moser theorem in the tame Fréchet category.

As a reference point, let us recall the classical version of the Inverse Function Theorem in the Banach setting.

Theorem A.1.1 (Banach version, [Lang99, Theorem I.5.2]) *Let $X, Y$ be Banach spaces and let $f\colon X \supseteq U \to Y$ be a smooth map defined on an open neighborhood $U$ of $0$ in $X$. If $\mathrm{T}_0 f\colon X \to Y$ is an isomorphism of Banach spaces, then $f$ is a local diffeomorphism at $0$.* ◇

Glöckner [Glö05; Glö06] has established the following generalization of the Banach Inverse Function Theorem to smooth maps depending on parameters in a locally convex space. Similar results have been obtained in [Hil99] (using a slightly stronger notion of differentiability) and in [Tei01] (using the convenient calculus).

Theorem A.1.2 (Banach version with parameters, [Glö06, Theorem 2.3]) *Let $P \subseteq E$ be an open neighborhood of $0$ in the locally convex vector space $E$. Let $X, Y$ be Banach spaces, let $U$ be an open neighborhood of $0$ in $X$ and let $f\colon E \times X \supseteq P \times U \to Y$ be a smooth map. If the partial derivative $\mathrm{T}^2_{(0,0)} f\colon X \to Y$ of $f$ at $(0, 0)$ with respect to the second variable is an isomorphism of Banach spaces, then the map*

$$E \times X \supseteq P \times U \to E \times Y, \quad (p, x) \mapsto \big(p, f(p, x)\big) \tag{A.1.1}$$

*is a local diffeomorphism at $(0, 0)$.* ◇

We now recall the main notions of the tame Fréchet category and the Nash–Moser Inverse Function Theorem, cf. [Ham82]. A Fréchet space $X$ is called graded if it carries a distinguished increasing fundamental system of seminorms $\|\cdot\|_k$. A graded Fréchet space is called *tame* if the seminorms satisfy an additional interpolation property, which formalizes the idea that $X$ admits smoothing operators, see [Ham82, Definition II.1.3.2] for the exact statement. Let $X$ and $Y$ be tame Fréchet spaces. A continuous (possibly non-linear) map $f\colon X \supseteq U \to Y$ defined on an open subset $U \subseteq X$ is *r-tame smooth* if it satisfies a local estimate of the form

$$\|f(x)\|_k \leq C(1 + \|x\|_{k+r}). \tag{A.1.2}$$

Roughly speaking, this means that $f$ has a maximal loss of $r$ derivatives. Moreover, a smooth map $f$ is called *r-tame smooth* if $f$ and all its derivatives $\mathrm{d}^{(j)} f : U \times X^j \to Y$ are $r$-tame.

Theorem A.1.3 (Nash–Moser Inverse Function Theorem, [Ham82, Section III.1]) *Let $X$ and $Y$ be tame Fréchet spaces, let $U$ be an open neighborhood of $0$ in $X$ and let $f : X \supseteq U \to Y$ be a tame smooth map. Assume that the derivative $\mathrm{T}f$ has a tame smooth family $\psi^f$ of inverses, that is, $\psi^f : U \times Y \to X$ is a tame smooth map and the family $\psi^f_x : Y \to X$ is inverse to $\mathrm{T}_x f$ for every $x \in U$. Then, the map $f$ is a tame local diffeomorphism at $0$.* ◊

The important point is that the derivative of $f$ has to be invertible in a *neighborhood* of 0 and that one requires tame estimates for the inverses.

Finally, let us comment on other generalizations of the Inverse Function Theorem that have been obtained in various analytical setups.

Remarks A.1.4

(i) Inspired by analytic problems in symplectic field theory, Hofer–Wysocki–Zehnder have introduced the scale calculus, see e.g. [HWZ17a; HWZ17b] and references therein. In this approach, one works with sequences of Banach spaces $E_i$ connected by compact inclusions $E_{i+1} \to E_i$. With an appropriate notion of scale-differentiability and scale-Fredholm maps, one can show that a version of the Inverse Function Theorem holds for smooth scale-Fredholm maps by reducing the issue to the ordinary Banach Inverse Function Theorem.

    Since $E_\infty = \bigcap_i E_i$ is a Fréchet space, one can reinterpret the scale Inverse Function Theorem as a theorem for maps between the corresponding Fréchet spaces. Scale calculus focuses on the Banach spaces $E_i$ while, in contrast, the Nash–Moser approach is mainly concerned with $E_\infty$ and only the shadow of the inclusions $E_{i+1} \to E_i$ is noticed in form of the tame estimates, see [Ger16] for a detailed comparison of the scale calculus and the tame category. Scale calculus is tailored to the elliptic setting one encounters in symplectic field theorem. We emphasize that the Nash–Moser theorem covers these cases as well, but additionally allows for applications that go well beyond the elliptic setting.

(ii) A local surjectivity theorem for maps between Fréchet spaces was presented by Ekeland [Eke11], which is based on his variational principle in place of the usual Newton iteration procedure.

(iii) For a rather rigid class of smooth maps satisfying a Lipschitz bound, a strong Inverse Function Theorem in the Fréchet setting is given in [Mül08, Theorem 4.10]. ◊

## A.2 Slices and orbit type stratification

Let $M$ be a (locally convex) manifold. Assume a (locally convex) Lie group $G$ acts smoothly on $M$, that is, assume that the action map $G \times M \to M$ is smooth. We refer to this setting by saying that $M$ is a $G$-manifold. The action is often written, using the dot notation, as $(g, m) \mapsto g \cdot m$. Similarly, the induced action of the Lie algebra $\mathfrak{g}$ of $G$ is denoted by $\xi \,.\, m \in \mathrm{T}_m M$ for $\xi \in \mathfrak{g}$ and $m \in M$. Clearly, $m \mapsto \xi \,.\, m$ is the Killing vector field generated by $\xi$. Furthermore, $G \cdot m = \{g \cdot m : g \in G\} \subseteq M$ is the orbit through $m \in M$. The $G$-action is called proper if inverse images of compact subsets under the map

$$G \times M \to M \times M, \qquad (g, m) \mapsto (g \cdot m, m) \tag{A.2.1}$$

are compact.

The subgroup $G_m := \{g \in G : g \cdot m = m\}$ is called the *stabilizer subgroup* of $m \in M$. It is not known, even for Banach Lie group actions, whether $G_m$ is always a *Lie* subgroup, see [Nee06, Problem IX.3.b]. However, for proper actions this is the case, see [DR18c, Lemma 2.11]. The $G$-action is called free if all stabilizer subgroups are trivial. Two subgroups $H$ and $K$ of $G$ are said to be conjugate if there exists $a \in G$ such that $aHa^{-1} = K$; we write $H \sim K$ is this case. In view of the equivariance relation $G_{g \cdot m} = gG_m g^{-1}$, for every $m \in M$ and $g \in G$, we can assign to every orbit $G \cdot m$ the conjugacy class $(G_m)$, which is called the *orbit type* of $m$. We put a preorder on the set of orbit types by declaring $(H) \leq (K)$ for two orbit types, represented by the stabilizer subgroups $H$ and $K$, if there exists $a \in G$ such that $aHa^{-1} \subseteq K$. If the action is proper, this preorder is actually a partial order. This follows from the following helpful property of compact subgroups.

Lemma A.2.1 ([DR18c, Lemma 2.12]) *Let $G$ be a Lie group. Let $H$ and $K$ be two compact Lie subgroups of $G$. If $K$ is conjugate to $H$ and $K \subseteq H$, then $K = H$.* ◇

For every closed subgroup $H \subseteq G$, define the following subsets of $M$:

$$\begin{aligned} M_H &= \{m \in M : G_m = H\}, \\ M_{(H)} &= \{m \in M : (G_m) = (H)\}. \end{aligned}$$

The subset $M_H$ is called the *isotropy type subset* and $M_{(H)}$ is the *subset of orbit type* $(H)$. Analogous definitions hold for every subset $N \subseteq M$, so, for example, $N_H = N \cap M_H$.

As in finite dimensions, the local structure of these sets is studied with the help of slices, cf. Proposition A.2.7 below. Since slices play a fundamental role in the construction of the normal form of an equivariant map, for the convenience of the reader, we now recall the definition of a slice in infinite dimensions, cf. [DR18c].

Definition A.2.2 Let $M$ be a $G$-manifold. A *slice* at $m \in M$ is a submanifold $S \subseteq M$ containing $m$ with the following properties:

(SL1) The submanifold $S$ is invariant under the induced action of the stabilizer subgroup $G_m$, that is $G_m \cdot S \subseteq S$.

(SL2) Any $g \in G$ with $(g \cdot S) \cap S \neq \emptyset$ is necessarily an element of $G_m$.

(SL3) The stabilizer $G_m$ is a principal Lie subgroup of $G$ and the principal bundle $G \to G/G_m$ admits a local section $\chi\colon G/G_m \supseteq U \to G$ defined on an open neighborhood $U$ of the identity coset $[e]$ in such a way that the map

$$\chi^S\colon U \times S \to M, \qquad ([g], s) \mapsto \chi([g]) \cdot s \tag{A.2.2}$$

is a diffeomorphism onto an open neighborhood $V \subseteq M$ of $m$. We call $V$ a *slice neighborhood* of $m$.

(SL4) The partial slice $S_{(G_m)} = \{s \in S : G_s \text{ is conjugate to } G_m\}$ is a closed submanifold of $S$.

(SL5) There exist a continuous representation of $G_m$ on a locally convex vector space $X$ and a $G_m$-equivariant diffeomorphism $\iota_S$ from a $G_m$-invariant open neighborhood of 0 in $X$ onto $S$ such that $\iota_S(0) = m$. ◊

The notion of a slice is closely related to the concept of a tubular neighborhood.

Proposition A.2.3 ([DR18c, Proposition 2.6.2]) *Let $M$ be a $G$-manifold. For every slice $S$ at $m \in M$, the tube map*

$$\chi^{\mathrm{T}}\colon G \times_{G_m} S \to M, \qquad [g, s] \mapsto g \cdot s \tag{A.2.3}$$

*is a $G$-equivariant diffeomorphism onto an open, $G$-invariant neighborhood $W$ of $G \cdot m$ in $M$.* ◊

In the finite-dimensional context, the existence of slices for proper actions is ensured by Palais' slice theorem [Pal61]. Passing to the infinite-dimensional case, this may no longer be true and additional hypotheses have to be made. We refer the reader to [Sub86; DR18c] for general slice theorems in infinite dimensions and to [Ebi70; ACM89; CMM91] for constructions of slices for concrete examples. One of the problems one faces in the infinite-dimensional setting is the failure of the Inverse Function Theorem and one usually is forced to use hard alternatives such as the Nash–Moser theorem. However, for linear actions of compact groups, the situation is better and we have the following existence result.

Theorem A.2.4 ([DR18c, Theorem 3.15]) *Let $X$ be a Fréchet space and let $G$ be a compact Lie group that acts linearly and continuously on $X$. Then, there exists a slice at every point of $X$.* ◊

The proof of the linear slice theorem uses, among other things, the following result concerning the existence of invariant complements.

Lemma A.2.5 ([DR18c, Lemma 3.13]) *Let $G$ be a compact Lie group which acts linearly and continuously on a Mackey complete locally convex vector space $X$. Let $H$ be a Lie subgroup of $G$. Then, every closed $H$-invariant topologically complemented subspace $E \subseteq X$ admits an $H$-invariant complement.* ◊

In Section 5.6, we need to construct a slice for a group action on a product manifold. This situation is covered by the next result.

Proposition A.2.6 ([DR18c, Proposition 3.29]) *Let $G$ be a Lie group that acts smoothly on the manifolds $M$ and $N$. Assume that the $G$-action admits a slice $S_m \subseteq M$ at a given point $m \in M$ and that the $G_m$-action on $N$ admits a slice $S_n$ at the point $n \in N$. Then, $S_m \times S_n$ is a slice at $(m, n)$ for the diagonal action of $G$ on the product $M \times N$.* ◊

As in the finite-dimensional case, the existence of slices implies many nice properties of the orbit space. For example, we have the following.

Proposition A.2.7 ([DR18c, Propositions 4.1 and 4.5]) *Let $M$ be a $G$-manifold with proper $G$-action. If the $G$-action admits a slice at every point of $M$, then $M_{(H)}$ is a submanifold of $M$. Moreover, $\check{M}_{(H)} = M_{(H)}/G$ carries a smooth manifold structure such that the natural projection $\pi_{(H)}\colon M_{(H)} \to \check{M}_{(H)}$ is a smooth submersion.* ◊

If, in addition, a certain approximation property is satisfied, then the orbit type manifolds fit together nicely and the orbit space is a stratified space, see [DR18c, Theorem 4.2]. More generally, we have the following stratification result for subsets of $M$.

Proposition A.2.8 ([DR18c, Proposition 4.7]) *Let $M$ be a $G$-manifold with proper $G$-action and let $P$ be a closed $G$-invariant subset of $M$. Assume that the $G$-action on $M$ admits a slice $S$ at every point $m \in P$ such that the following holds:*

*(i) $P \cap S_{(G_m)}$ is a closed submanifold of $S_{(G_m)}$,*

*(ii) for every orbit type $(K) \leq (G_m)$, the point $m$ lies in the closure of $S_{(K)} \cap P$ in $S$.*

*Then, the induced partition of $P$ into the orbit type subsets $P_{(H)} = P \cap M_{(H)}$ is a stratification. Moreover, under these assumptions, the decomposition of $\check{P} = P/G$ into $\check{P}_{(H)} = P_{(H)}/G$ is a stratification, too.* ◊

For completeness, we include our definition of a stratification here and refer the reader to [DR18c] for further details and comparison with other notions of stratified spaces in the literature.

Definition A.2.9 Let $X$ be Hausdorff topological space. A partition $\mathcal{Z}$ of $X$ into subsets $X_\sigma$ indexed by $\sigma \in \Sigma$ is called a *stratification* of $X$ if the following conditions are satisfied:

(DS1) Every piece $X_\sigma$ is a locally closed, smooth manifold (whose manifold topology coincides with the relative topology). We will call $X_\sigma$ a *stratum* of $X$.

(DS2) (frontier condition) Every pair of disjoint strata $X_\varsigma$ and $X_\sigma$ with $X_\varsigma \cap \overline{X_\sigma} \neq \emptyset$ satisfies:

a) $X_\varsigma$ is contained in the frontier $\overline{X_\sigma} \setminus X_\sigma$ of $X_\sigma$,

b) $X_\sigma$ does not intersect $\overline{X_\varsigma}$.

In this case, we write $X_\varsigma < X_\sigma$ or $\varsigma < \sigma$. ◇

## A.3 Cotangent bundles in infinite dimensions

The tangent bundle $\mathrm{T}Q$ of a manifold $Q$ is itself a smooth manifold again in such a way that the projection $\mathrm{T}Q \to Q$ is a smooth locally trivial bundle. The dual bundle $\mathrm{T}'Q := \bigsqcup_{q \in Q} (\mathrm{T}_q Q)'$ is more problematic, cf. [Nee06, Remark I.3.9]. In order to endow $\mathrm{T}'Q$ with a smooth fiber bundle structure we need a locally convex topology on the dual $X'$ of the model space $X$ of $Q$ such that for every local diffeomorphism $\phi\colon X \to X$ the map

$$\tau_\phi\colon X \times X' \to X', \qquad (x, \alpha) \mapsto \alpha \circ \mathrm{T}_x \phi \tag{A.3.1}$$

is smooth, because otherwise the notion of smoothness of $\mathrm{T}'Q$ is chart-dependent. It is straightforward to construct a map $\phi$ such that $\tau_\phi$ involves the natural pairing $X \times X' \to \mathbb{R}$. However, this pairing is discontinuous in 0 for any vector topology on $X'$ except when $X$ is a Banach space, see [Mai63, Satz 1]. Thus, in summary, we cannot endow $\mathrm{T}'Q$ with a natural smooth bundle structure for non-Banach manifolds $Q$.

Hence, the most important class of examples of symplectic manifolds is not available a priori in infinite dimensions. We now present a substitute, which will play the role of the cotangent bundle.

Definition A.3.1 A *dual pairing between two vector bundles* $E \to Q$ and $F \to Q$ is a smooth map $h\colon E \times_Q F \to \mathbb{R}$ which is fiberwise non-degenerate, i.e., $h_q\colon E_q \times F_q \to \mathbb{R}$ is a non-degenerate bilinear form for every $q \in Q$. If $E$ is a vector bundle in duality to the tangent bundle $F = \mathrm{T}Q$, then we will write $\mathrm{T}^*Q \equiv E$ and call it a *cotangent bundle* of $Q$. Correspondingly, we denote the bundle projection $\mathrm{T}^*Q \to Q$ by $\overset{*}{\tau}$. ◇

We will often denote the duality by brackets and write the dual pair as $\langle \mathrm{T}^*Q, \mathrm{T}Q \rangle$. If $Q$ carries a Riemannian metric $g$, then the pairing $g\colon \mathrm{T}Q \times_Q \mathrm{T}Q \to \mathbb{R}$ identifies the cotangent bundle of $Q$ with $\mathrm{T}Q$. In the following, we assume that a dual pair $\langle \mathrm{T}^*Q, \mathrm{T}Q \rangle$ is fixed and we will simply refer to $\mathrm{T}^*Q$ as *the* cotangent bundle of $Q$. The reader should however keep in mind that there is no smooth *canonical* cotangent bundle and a choice of a dual pair is always required.

Recall that, for a dual pair $\langle X_2, X_1\rangle$ of vector spaces [Köt83], one has a natural embedding of $X_2$ into the topological dual $X_1'$ of $X_1$. Similarly, for a dual pair of vector bundles, we obtain a natural vector bundle injection of $E$ into the dual bundle $F'$, whose fiber over $m$ is the topological dual $F_m'$ of $F_m$. In particular, every cotangent bundle $\mathrm{T}^*Q$ comes with a natural injection into the topological cotangent bundle $\mathrm{T}'Q$.

Proposition A.3.2 *For a cotangent bundle* $\mathrm{T}^*Q$*, the formula*

$$\theta_p(v) = \langle p, \mathrm{T}_p \overset{\star}{\tau}(v)\rangle_q, \qquad p \in \mathrm{T}^*_q Q, v \in \mathrm{T}_p(\mathrm{T}^*Q) \tag{A.3.2}$$

*defines a smooth* 1*-form* $\theta$ *on* $\mathrm{T}^*Q$*. Furthermore,* $\omega = \mathrm{d}\theta$ *is a symplectic form. We say that* $\omega$ *is the* canonical symplectic form *on* $\mathrm{T}^*Q$*.* ◊

*Proof.* Let $U \subseteq Q$ be an open subset over which $\overset{\star}{\tau}\colon \mathrm{T}^*Q \to Q$ as well as $\tau\colon \mathrm{T}Q \to Q$ trivialize. Denote the fiber model space of $\mathrm{T}Q$ and of $\mathrm{T}^*Q$ by $X$ and $X^*$, respectively. Using a chart on $Q$, we identify $U$ as a subspace of $X$. Then, the local chart representation of the pairing is a smooth map $U \times X^* \times X \to \mathbb{R}$. In this chart, the canonical form $\theta\colon (U \times X^*) \times (X \times X^*) \to \mathbb{R}$ is given by

$$\theta_{q,\alpha}(u, \beta) = \langle \alpha, u\rangle_q \tag{A.3.3}$$

and, hence, $\theta$ is a smooth 1-form on $\mathrm{T}^*Q$. For the symplectic structure

$$\omega\colon (U \times X^*) \times (X \times X^*)^2 \to \mathbb{R}, \tag{A.3.4}$$

we find[1]

$$\begin{aligned}\omega_{q,p}(u_1, \beta_1, u_2, \beta_2) = \partial_q \left(\langle p, u_2\rangle_q\right)(u_1) - \partial_q \left(\langle p, u_1\rangle_q\right)(u_2) \\ + \langle \beta_1, u_2\rangle_q - \langle \beta_2, u_1\rangle_q.\end{aligned} \tag{A.3.5}$$

In finite dimensions, this corresponds to the classical formula

$$\omega = g^i{}_j \,\mathrm{d}p_i \wedge \mathrm{d}q^j + \partial_k \left(g^i{}_j\, p_i\right) \mathrm{d}q^k \wedge \mathrm{d}q^j, \tag{A.3.6}$$

where $g^i{}_j(q)$ denote the components of the dual pair $\langle\,\cdot\,,\,\cdot\,\rangle_q$. In either case, we conclude that $\omega$ is non-degenerate because the dual pair $\langle\,\cdot\,,\,\cdot\,\rangle_q$ possesses this property for every $q \in Q$. □

Let $G$ be a Lie group acting smoothly on $Q$. By linearization, we get a smooth action of $G$ on the tangent bundle, which we write using the lower dot notation as $g \,.\, v \in \mathrm{T}_{g \cdot q} Q$ for $g \in G$ and $v \in \mathrm{T}_q Q$. The action on $\mathrm{T}Q$ induces a $G$-action on $\mathrm{T}^* Q$ by requiring that the pairing be left invariant, that is,

$$\langle g \cdot p, v \rangle_q = \langle p, g^{-1} \,.\, v \rangle_{g^{-1} \cdot q}, \qquad p \in \mathrm{T}^*_{g^{-1} \cdot q} Q, v \in \mathrm{T}_q Q. \tag{A.3.7}$$

[1] This local expression for $\omega$ is well-known for the case that $\mathrm{T}Q$ is put in duality with itself using a Riemannian metric on $Q$, see for example [Mar72, p. 591].

In order for this equation to define a smooth action on $\mathrm{T}^* Q$, the action $\mathrm{T}_{g^{-1} \cdot q} Q \to \mathrm{T}_q Q$ needs to be weakly continuous with respect to the pairing $\langle \cdot , \cdot \rangle$ for every $g \in G$. With respect to this action, the cotangent bundle projection $\overset{\star}{\tau}$ is $G$-equivariant and so the action on $\mathrm{T}^* Q$ is a lift of the action on $Q$.

Let $\mathfrak{g}$ be the Lie algebra of $G$. Similarly to the above strategy for the cotangent bundle, the choice of a Fréchet space $\mathfrak{g}^*$ and of a separately continuous non-degenerate bilinear form $\kappa \colon \mathfrak{g}^* \times \mathfrak{g} \to \mathbb{R}$ yields a dual pair, with $\mathfrak{g}^*$ serving as the dual space of $\mathfrak{g}$.

Proposition A.3.3 *Let $G$ be a Lie group that acts smoothly on $Q$. Then, the lifted action of $G$ on the cotangent bundle* $\mathrm{T}^* Q$ *preserves the canonical symplectic form $\omega$. Let, moreover, $\kappa(\mathfrak{g}^*, \mathfrak{g})$ be a dual pair. If, for every $p \in \mathrm{T}^* Q$, the functional $\xi \mapsto \theta_p(\xi \,.\, p)$ on $\mathfrak{g}$ is represented by an element $J(p) \in \mathfrak{g}^*$ with respect to $\kappa$, then the resulting map $J \colon \mathrm{T}^* Q \to \mathfrak{g}^*$ is a smooth $G$-equivariant $\mathfrak{g}^*$-valued momentum map for the lifted $G$-action on* $\mathrm{T}^* Q$. ◇

*Proof.* The canonical form $\theta$ is $G$-invariant because we have

$$\theta_{g \cdot p}(g \,.\, v) = \langle g \cdot p, \mathrm{T}_{g \cdot p} \overset{\star}{\tau}(g \,.\, v) \rangle = \langle g \cdot p, g \,.\, \mathrm{T}_p \overset{\star}{\tau}(v) \rangle = \langle p, \mathrm{T}_p \overset{\star}{\tau}(v) \rangle \tag{A.3.8}$$

for every $g \in G$ and $v \in \mathrm{T}_p(\mathrm{T}^* Q)$. Thus, the action leaves also the symplectic form $\omega = \mathrm{d}\theta$ invariant. By assumption, the functional $\xi \mapsto \theta_p(\xi \,.\, p)$ on $\mathfrak{g}$ is represented by an element $J(p) \in \mathfrak{g}^*$ with respect to the given dual pair $\kappa(\mathfrak{g}^*, \mathfrak{g})$. This is to say,

$$\kappa(J(p), \xi) = \theta_p(\xi \,.\, p) = \langle p, \xi \,.\, \overset{\star}{\tau}(p) \rangle. \tag{A.3.9}$$

On the other hand, we have $\mathfrak{L}_{\xi^*} \theta = 0$, because $\theta$ is $G$-invariant and so

$$0 = \mathfrak{L}_{\xi^*} \theta = \xi^* \lrcorner \, \mathrm{d}\theta + \mathrm{d}(\xi^* \lrcorner \, \theta) = \xi^* \lrcorner \, \omega + \kappa(\mathrm{d}J, \xi), \tag{A.3.10}$$

which shows that $J$ is a $\mathfrak{g}^*$-valued momentum map. The $G$-equivariance property of $J$ is a direct consequence of the $G$-invariance of the pairing $\langle \mathrm{T}^* Q, \mathrm{T}Q \rangle$ and (A.3.9). □

In contrast to the finite-dimensional case, lifted actions of infinite-dimensional Lie groups do not necessarily possess a momentum map, see Example 5.1.2.

# Dual Pairs B

In this appendix, we summarize the relevant material on the theory of dual pairs without proofs. The exposition is based on standard references [RR64; Sch71; Jar81; Köt83; NB10]. The basic idea underlying duality theory is that one can translate a given problem living on a locally convex space $X$ into a question which primarily involves the dual space $X'$. This dual problem often turns out to be easier to solve than the original one. Conversely, sometimes it can be fruitful to convert problems concerning the dual space into questions involving only the original space.

Definition B.0.1 ([Köt83, Section 10.3]) A *dual pair* is a pair of vector spaces $X_1$ and $X_2$ together with a bilinear form $\langle \cdot\,, \cdot \rangle : X_2 \times X_1 \to \mathbb{R}$ which is non-degenerate in the sense that the following conditions are satisfied:

(i) If $\langle \alpha, x \rangle = 0$ for all $x \in X_1$, then $\alpha = 0$.

(ii) If $\langle \alpha, x \rangle = 0$ for all $\alpha \in X_2$, then $x = 0$.

We will often use the shorthand notation $\langle X_2, X_1 \rangle$ to denote the dual pair. ◊

Although the spaces $X_1$ and $X_2$ are treated on an equal footing in Definition B.0.1, we think of $X_2$ as the dual to $X_1$ and this thinking is reflected in our notation by using Greek letters to denote elements of $X_2$. Non-degeneracy of $\langle \cdot\,, \cdot \rangle$ is equivalent to injectivity of the partial maps

$$X_2 \to X_1^*, \qquad \alpha \mapsto \langle \alpha, \cdot \rangle, \tag{B.0.1}$$

$$X_1 \to X_2^*, \qquad x \mapsto \langle \cdot, x \rangle, \tag{B.0.2}$$

where $X^*$ denotes the algebraic dual of the vector space $X$. We call these maps the *natural embeddings*. If the pairing is separately continuous with respect to given locally convex topologies on $X_1$ and $X_2$, then the natural embeddings take values in the topological dual rather than merely in the algebraic one.

Examples B.0.2

(i) For a locally convex space $X$, the natural pairing of $X$ and its topological dual yields a dual pair $(X', X)$ by the Hahn–Banach Theorem.

(ii) Let $(X, \omega)$ be a (weakly) symplectic vector space in the sense of Definition 4.1.1. Since the symplectic form $\omega$ is non-degenerate, it defines a dual pair $\omega(X, X)$. As usual, we use the musical notation $\omega^\flat : X \to X'$ for the natural embedding. ◊

## B.1 Compatible and polar topologies

The natural embedding $X_2 \to X_1^*$ is almost never surjective onto the full algebraic dual. We can use a topology on $X_1$ to restrict the class of admissible functionals and thereby obtain a surjective map in the following sense.

Definition B.1.1 Let $\langle X_2, X_1 \rangle$ be a dual pair. A locally convex topology $\tau$ on $X_1$ is called *compatible* with $\langle \cdot\,, \cdot \rangle$ if $(X_1, \tau)' = X_2$ under the natural embedding $X_2 \to X_1^*$, i.e., every linear $\tau$-continuous functional on $X_1$ is of the form $\langle \alpha, \cdot \rangle$ for some (unique) $\alpha \in X_2$. ◇

Example B.1.2 For a locally convex space $X$, the original topology on $X$ is compatible with the natural pairing $(X', X)$. The topology on a symplectic vector space $(X, \omega)$ is compatible with $\omega$ if and only if $\omega^\flat$ is surjective, i.e., if $\omega$ is strongly non-degenerate, see Proposition 4.1.5. ◇

Let $\tau_s$ be a topology on $X_1$ for which a given dual pair $\langle \cdot\,, \cdot \rangle : X_2 \times X_1 \to \mathbb{R}$ is separately continuous and let $\tau_c$ be a compatible topology on $X_1$. Since every $\tau_c$-continuous linear functional on $X_1$ is also $\tau_s$-continuous, the topology $\tau_c$ is coarser than $\tau_s$. The difference between $\tau_s$ and $\tau_c$ provides a measure of how much the natural embedding fails to be surjective onto $(X_1, \tau_s)'$. Although this observation is almost tautological, it has far-reaching consequences since it converts the algebraic question of surjectivity of the natural embedding into a topological problem.

In order to systematically study topologies compatible with a dual pair $\langle X_2, X_1 \rangle$, we need the concept of a polar topology. Following [Jar81, Section II.8.2], the *polar* $A^\circ \subseteq X_1$ of a subset $A \subseteq X_2$ is defined by

$$A^\circ := \{x \in X_1 : |\langle \alpha, x \rangle| \leq 1 \text{ for all } \alpha \in A\}. \tag{B.1.1}$$

The polar $B^\circ \subseteq X_2$ of a subset $B \subseteq X_1$ is defined in a similar manner. For a *linear* subspace $Y \subseteq X_2$, the polar $Y^\circ$ coincides with the orthogonal

$$Y^\perp = \{x \in X_1 : \langle \alpha, x \rangle = 0 \text{ for all } \alpha \in Y\}. \tag{B.1.2}$$

Let us record some basic facts about polars.

Proposition B.1.3 ([Jar81, Proposition 8.2.1]) *For every dual pair $\langle X_2, X_1 \rangle$, the following holds:*

*(i)* $\{0\}^\circ = X_1$

*(ii) (Anti-monotony) For subsets $A_1 \subseteq A_2$ of $X_2$, we have $A_1^\circ \supseteq A_2^\circ$.*

*(iii) For all non-zero $\lambda \in \mathbb{R}$ and a subset $A \subseteq X_2$, the identity $(\lambda A)^\circ = \frac{1}{\lambda} A^\circ$ holds.*

*(iv) For linear subspaces $Y_1$ and $Y_2$ of $X_2$, we have $(Y_1 + Y_2)^\perp = Y_1^\perp \cap Y_2^\perp$.*

◊

A polar topology on $X_1$ is a topology whose 0-neighborhood base is given in terms of polars of a suitable family of subsets of $X_2$. The following proposition clarifies the minimal requirements on the subsets $A \subseteq X_2$ such that their polars define a locally convex topology on $X_1$.

Proposition B.1.4 ([Köt83, Proposition 21.1.1; NB10, Theorem 8.3.5]) *Let $\langle X_2, X_1\rangle$ be a dual pair and let $A \subseteq X_2$ be a subset. Then, the following are equivalent:*

*(i) The Minkowski functional*

$$\|x\|_A := \inf\{\kappa \in \mathbb{R} : x \in \kappa \cdot A^\circ\} = \sup_{\alpha \in A} |\langle \alpha, x\rangle| \tag{B.1.3}$$

*of the polar $A^\circ$ is a seminorm on $X_1$.*

*(ii) $A^\circ$ is an absorbing subset of $X_1$.*

*(iii) For all $x \in X_1$, the set of values $\{\langle \alpha, x\rangle : \alpha \in A\}$ is bounded in $\mathbb{R}$.* ◊

A subset $A \subseteq X_2$ satisfying these equivalent conditions is called *weakly bounded*. Accordingly, the seminorms $\|\cdot\|_A$ associated to a family $\mathcal{A}$ of weakly bounded subsets $A$ of $X_2$ determine a topology on $X_1$, which we call the *polar topology* generated by $\mathcal{A}$. The subsets $\varepsilon A^\circ$ for $\varepsilon > 0$ and $A \in \mathcal{A}$ form a 0-neighborhood subbase of this polar topology. Equivalently, a net $(x_\gamma)$ in $X_1$ converges to $x \in X_1$ with respect to the polar topology if and only if $\sup_{\alpha \in A} |\langle \alpha, x - x_\gamma\rangle|$ tends to 0 for all $A \in \mathcal{A}$. Hence, the polar topology is the topology of uniform convergence on the sets of $\mathcal{A}$. The polar topology is Hausdorff and thus locally convex if and only if 0 is the only element orthogonal to the union $\cup_{\mathcal{A}} A$, see [Köt83, Proposition 21.2.2].

Examples B.1.5

(i) (Weak topology) The polar topology on $X_1$ generated by the collection $\mathcal{A}$ of finite subsets of $X_2$ is called the *weak topology* and is denoted by $\sigma(X_1, X_2)$. The weak topology is the coarsest topology on $X_1$ for which the functionals $\langle \alpha, \cdot\rangle$ with $\alpha \in X_2$ are continuous. In other words, the weak topology is the weakest topology on $X_1$ compatible with the dual pair $\langle X_2, X_1\rangle$. The weak topology is equivalently characterized by the following convergence condition: a net $(x_\lambda)$ in $X_1$ converges in the weak topology to $x \in X_1$ if and only if $\langle \alpha, x_\lambda\rangle$ converges to $\langle \alpha, x\rangle$ for all $\alpha \in X_2$.

(ii) (Mackey topology) The *Mackey topology* $\tau(X_1, X_2)$ is the polar topology generated by all convex, circled[1], $\sigma(X_2, X_1)$-compact subsets of $X_2$. The Mackey topology is compatible with the dual pair.

[1] A subset $A \subseteq X$ of a vector space $X$ is called circled if $\{\lambda \in \mathbb{R} : |\lambda| \leq 1\} \cdot A = A$.

(iii) (Strong topology) The family of all weakly bounded subsets of $X_2$ yields the strongest polar topology on $X_1$, which hence is called the *strong topology*. We will denote it by $\beta(X_1, X_2)$. The strong topology is in general not compatible with the dual pair. ◊

Intuitively it is clear that the spectrum of compatible topologies is bounded on both sides: there exists a coarsest compatible topology, because the functionals $\langle \alpha, \cdot \rangle$ have to stay continuous; on the other hand, there has to be a finest compatible topology, because otherwise there are too many continuous functionals to be of the form $\langle \alpha, \cdot \rangle$. Polar topologies are a suitable tool to characterize these topologies that are compatible with a given dual pair.

Theorem B.1.6 (Mackey–Arens, [Sch71, Theorem 3.3]) *Let $\langle X_2, X_1 \rangle$ be a dual pair. A topology $\tau$ on $X_1$ is compatible with $\langle X_2, X_1 \rangle$ if and only if it is a polar topology satisfying*

$$\sigma(X_1, X_2) \leq \tau \leq \tau(X_1, X_2), \tag{B.1.4}$$

*where $\sigma(X_1, X_2)$ and $\tau(X_1, X_2)$ denote the weak and the Mackey topology, respectively.* ◊

As a consequence of the definition, all compatible topologies on $X_1$ have the same set of linear continuous functionals. Hence, it is not surprising that certain other topological properties do not depend on the chosen compatible topology, too.

Proposition B.1.7 ([Wil78, Corollary 8.3.6 and Theorem 8.4.1]) *Let $\langle X_2, X_1 \rangle$ be a dual pair. All compatible topologies on $X_1$ have the same closed convex subsets, the same closed linear subspaces and the same bounded subsets.* ◊

Let $B \subseteq X_1$ be a subset. The polar $B^\circ \subseteq X_2$ of $B$ is defined in the same manner as in (B.1.1):

$$B^\circ := \{\alpha \in X_2 : \ |\langle \alpha, x \rangle| \leq 1 \text{ for all } x \in B\}. \tag{B.1.5}$$

Taking the polar twice yields the *bipolar* $B^{\circ\circ} \subseteq X_1$. While the definition of the bipolar is completely algebraic, the following theorem establishes a deep relation to compatible topologies.

Theorem B.1.8 (Bipolar Theorem, [Köt83, Proposition 20.3.2]) *Let $\langle X_2, X_1 \rangle$ be a dual pair. For every subset $B \subseteq X_1$, the bipolar $B^{\circ\circ}$ coincides with the absolutely convex closure of $B$ with respect to a compatible topology on $X_1$. In particular, for a linear subspace $Y \subseteq X$ we have $Y^{\perp\perp} = \overline{Y}$, where the closure is taken with respect to a compatible topology.* ◊

## B.2 Adjoint maps

Let $\langle X_2, X_1 \rangle$ and $(Y_2, Y_1)$ be dual pairs. The *adjoint* of a linear map $T\colon X_1 \longrightarrow Y_1$ is a linear map $T^*\colon Y_2 \longrightarrow X_2$ satisfying

$$(\beta, Tx) = \langle T^*\beta, x \rangle \tag{B.2.1}$$

for all $x \in X_1$ and $\beta \in Y_2$. The adjoint map $T^*$ is uniquely determined by this relation.

Lemma B.2.1 ([Sch71, Proposition IV.2.1]) *The adjoint of $T\colon X_1 \longrightarrow Y_1$ exists if and only if $T$ is continuous with respect to the weak topologies on $X_1$ and $Y_1$.* ◇

*Proof.* If $T$ is weakly continuous, then, for every $\beta \in Y_2$, the functional $X_1 \ni x \longmapsto (\beta, Tx) \in \mathbb{R}$ is weakly continuous. Since the weak topology is compatible with $\langle X_2, X_1 \rangle$, this functional is represented by an element of $X_2$, say $T^*\beta$. This prescription defines a linear map $T^*\colon Y_2 \longrightarrow X_2$ satisfying $(\beta, Tx) = \langle T^*\beta, x \rangle$. Conversely, if the adjoint $T^*$ exists, then the functional

$$x \longmapsto (\beta, Tx) = \langle T^*\beta, x \rangle \tag{B.2.2}$$

on $X_1$ is weakly continuous for every $\beta \in Y_2$. By the definition of weak topologies, this implies that $T$ is weakly continuous. □

Example B.2.2 Let $\kappa(X, Y)$ be a separately continuous dual pair. Then, the natural embedding takes values in the topological dual space $Y'$ of $Y$, i.e.,

$$\kappa^\flat\colon X \longrightarrow Y', \qquad x \longmapsto \kappa(x, \cdot). \tag{B.2.3}$$

Obviously, the identity

$$(\kappa^\flat x, y) = \kappa(x, y) \tag{B.2.4}$$

holds with respect to the natural dual pair $(Y', Y)$. Hence, $\kappa^\flat$ has the identity on $Y$ as its adjoint relative to the dual pairs $\kappa(X, Y)$ and $(Y', Y)$. According to Lemma B.2.1, the map $\kappa^\flat$ is continuous with respect to the weak-$\kappa$ topology on $X$ and the weak topology on $Y'$. ◇

Proposition B.2.3 ([Sch71, Proposition IV.2.3]) *Let $\langle X_2, X_1 \rangle$ and $(Y_2, Y_1)$ be dual pairs. For every weakly continuous linear map $T\colon X_1 \longrightarrow Y_1$ with adjoint $T^*\colon Y_2 \longrightarrow X_2$, the following holds:*

*(i)* $\big(T(A)\big)^\circ = (T^*)^{-1}(A^\circ)$ *for every subset $A \subseteq X_1$.*

*(ii)* $\operatorname{Ker} T = (\operatorname{Im} T^*)^\perp$.

*(iii)* $\overline{\operatorname{Im} T} = (\operatorname{Ker} T^*)^\perp$, *where the closure is taken with respect to a compatible topology on $Y_1$.*

(iv) *The adjoint $T^*$ is injective if and only if* $\operatorname{Im} T$ *is dense with respect to a compatible topology on* $Y_1$. ◇

We will sometimes write the adjoint of a linear map $T\colon X_1 \to Y_1$ with respect to dual pairs $\langle X_2, X_1\rangle$ and $(Y_2, Y_1)$ in the following diagrammatic form:

$$\begin{array}{ccc} X_1 & \xrightarrow{\;T\;} & Y_1 \\ \times_{\langle\cdot,\cdot\rangle} & & \times_{(\cdot,\cdot)} \\ X_2 & \xleftarrow[T^*]{} & Y_2. \end{array} \tag{B.2.5}$$

Let $\langle Z_1, Z_2\rangle$ be another dual pair and let $S\colon Y_1 \to Z_1$ be a weakly continuous map with adjoint $S^*\colon Z_2 \to Y_2$. Since the adjoint of the composite map $S \circ T$ is $T^* \circ S^*$, see [Sch71, Proposition IV.2.2], we obtain a formalism that allows rudimentary "diagram chasing":

$$\begin{array}{ccccc} X_1 & \xrightarrow{\;T\;} & Y_1 & \xrightarrow{\;S\;} & Z_1 \\ \times_{\langle\cdot,\cdot\rangle} & & \times_{(\cdot,\cdot)} & & \times_{\langle\cdot,\cdot\rangle} \\ X_2 & \xleftarrow[T^*]{} & Y_2 & \xleftarrow[S^*]{} & Z_2. \end{array} \tag{B.2.6}$$

# Own Publications

# List of Symbols

| | |
|---|---|
| $\mathrm{Ad}^* P$ | Coadjoint bundle $\mathrm{Ad}^* = P \times_G \mathfrak{g}^*$ associated to the principal $G$-bundle $P$. |
| $\mathrm{Ad}_g$ | Adjoint action of $g \in G$ on the Lie algebra $\mathfrak{g}$ of $G$. |
| $\mathrm{ad}_A$ | Adjoint action of $A \in \mathfrak{g}$ on the Lie algebra $\mathfrak{g}$. |
| $\mathrm{Ad} P$ | Adjoint bundle $\mathrm{Ad} P = P \times_G \mathfrak{g}$ associated to the principal $G$-bundle $P$. |
| $\mathrm{C}(M)$ | Space of functions on the manifold $M$. |
| $\mathrm{C}^0(X)$ | Space of continuous functions on the topological space $X$. |
| $\mathrm{C}^\infty(M)$ | Space of smooth functions on the manifold $M$. |
| $\mathrm{C}_G(A)$ | Centralizer of the subset $A$ of $G$. |
| $\mathrm{c}_k(E)$ | The $k$-th Chern class of the vector bundle $E \longrightarrow M$. |
| $\mathrm{Ad}^*_g$ | Coadjoint action of $g \in G$ on the dual $\mathfrak{g}^*$ of the Lie algebra $\mathfrak{g}$ of $G$ defined by $\langle \mathrm{Ad}^*_g \, \mu, \xi \rangle = \langle \mu, \mathrm{Ad}_{g^{-1}} \, \xi \rangle$. |
| $\mathrm{ad}^*_A$ | Coadjoint action of $A \in \mathfrak{g}$ on the dual $\mathfrak{g}^*$ of $\mathfrak{g}$. |
| $\mathrm{Coim}\, T$ | Coimage of the linear map $T$. |
| $\mathrm{Coker}\, T$ | Cokernel of the linear map $T$. |
| const | Constant value. |
| $X \lrcorner \alpha$ | Interior product of the vector field $X$ with the differential form $\alpha$. |
| $\hat{\mathrm{H}}^k(M, \mathrm{U}(1))$ | $k$-th Cheeger–Simons differential cohomology group of the manifold $M$. |
| $\mathrm{H}^k(M, R)$ | $k$-th de Rham cohomology group of the manifold $M$ with values in $R$. |
| $\mathcal{D}iff(M)$ | Group of diffeomorphisms of the manifold $M$. |
| $\mathrm{ev}_m \phi$ | Evaluation of the map $\phi$ at the point $m$. |
| $\mathrm{FL}^X$ | Flow of the vector field $X$. |
| $\Gamma^\infty(E)$ | Space of smooth sections of the fiber bundle $E$. |
| $\mathcal{G}au(P)$ | Group of gauge transformation of the principal bundle $P$. |
| $\mathfrak{gau}(P)$ | Lie algebra of infinitesimal gauge transformation of the principal bundle $P$. |
| $\mathrm{H}E$ | Horizontal bundle of the fiber bundle $E$ induced by a connection in $E$. |
| $\mathrm{Hol}^\star_A(\gamma)$ | Holonomy of a closed path $\gamma$ with respect to the connection $A$ and the point $\star$ in the bundle. |
| $\mathrm{Hol}^\star(A)$ | Holonomy group of the connection $A$ relative to the point $\star$ in the bundle. |
| $\mathrm{Hom}(G, H)$ | Space of homomorphisms from $G$ to $H$. |
| $\mathrm{id}_X$ | Identity map on the space $X$. |
| $\operatorname{Im} f$ | Image of the map $f$. |
| $\operatorname{ind} T$ | Index of the Fredholm operator $T$. |
| $\mathrm{j}^r_m \sigma$ | $r$-th jet of the map $\sigma$ at the point $m$. |
| $\mathrm{J}^r E$ | $r$-th jet bundle of the fiber bundle $E$. |
| $\operatorname{Ker} T$ | Kernel of the linear map $T$. |
| $\mathrm{L}^p(X)$ | Space of functions on the measure space $X$ for which the $p$-th power of the absolute value is Lebesgue integrable. |
| $\mathrm{L}^2(X)$ | Space of square-integrable functions on the measure space $X$. |
| $\mathcal{L}^{\propto}_{\star}(M)$ | Space of all closed, piecewise smooth paths in $M$ starting and ending at $\star \in M$. |
| $\mathrm{L}(X, Y)$ | Space of (continuous) linear bundle maps between (topological) vector bundles $X$ and $Y$. |
| $\mathrm{L}(X, Y)$ | Space of (continuous) linear maps between (topological) vector spaces $X$ and $Y$. |
| $\Omega^k(M; E)$ | Space of differential $k$-forms on the manifold $M$ with values in the vector bundle $E \to M$. |
| $\Omega^k_{\mathrm{cl}}(M; \mathbb{R})$ | Space of closed real-valued differential $k$-forms on the manifold $M$. |
| $\Omega^k_{\mathrm{cl},\mathbb{Z}}(M; \mathbb{R})$ | Space of closed real-valued differential $k$-forms on the manifold $M$ whose periods are integral. |
| $\mathrm{per}_\alpha$ | Period homomorphism of $\alpha \in \Omega^1(M, \mathfrak{g})$. |
| $\mathrm{p}_k(E)$ | $k$-th Pontryagin class of the vector bundle $E \to M$. |
| $\mathrm{pr}$ | Projection, for example the natural projection $\mathrm{pr}_M \colon M \times N \to M$. |
| $\operatorname{rk} T$ | Rank of the linear map $T$. |
| $\operatorname{Rep}(G)$ | Space of representations of the group $G$. |
| $\mathrm{H}^k(M, R)$ | $k$-th singular cohomology group of the manifold $M$ with values in the ring $R$. |
| $\ell^p$ | Space of sequences $(x_n)$ satisfying $\sum_n \lvert x_n \rvert^p < \infty$. |
| $\mathrm{H}_k(M, R)$ | $k$-th singular homology group of the manifold $M$ with values in the ring $R$. |
| $\sigma_P$ | Principal symbol of the differential operator $P$. |
| $\mathfrak{sl}(n)$ | Special linear Lie algebra of degree $n$. |
| $\mathrm{SL}(n)$ | Special linear group of degree $n$. |
| $\mathrm{SO}(n)$ | Special orthogonal group of degree $n$. |
| $\mathrm{Sp}(n)$ | Symplectic group of degree $n$. |
| $\mathfrak{su}(n)$ | Special unitary Lie algebra of degree $n$. |
| $\mathrm{SU}(n)$ | Special unitary group of degree $n$. |
| $\mathrm{T}^* M$ | Cotangent bundle of the manifold $M$. |
| $\overset{\star}{\tau}$ | The canonical projection $\overset{\star}{\tau} \colon \mathrm{T}^* M \to M$ of the cotangent bundle. |
| $\mathrm{T}M$ | Tangent bundle of the manifold $M$. |
| $\mathfrak{u}(n)$ | Unitary Lie algebra of degree $n$. |
| $\mathrm{U}(n)$ | Unitary group of degree $n$. |
| $\mathrm{V}E$ | Vertical bundle of the fiber bundle $E$. |
| $\mathrm{vol}_g$ | Volume element induced by the Riemannian metric $g$. |
| $\mathrm{span}\{v_1, \ldots, v_n\}$ | Linear span of the vectors $v_i$. |
| $\mathfrak{X}(M)$ | Space of vector fields on the manifold $M$. |
| $\mathfrak{z}(\mathfrak{g})$ | Center of the Lie algebra $\mathfrak{g}$. |